\documentclass[a4paper,12pt]{amsart}
\usepackage{amsmath,amssymb,amsfonts,amsthm}
\usepackage[]{amsrefs}
\usepackage{bm}
\usepackage[top=25mm,bottom=25mm,left=25mm,right=25mm]{geometry}
\usepackage{array}
\usepackage{comment}
\usepackage{ascmac}
\usepackage[all]{xy}
\usepackage{amsthm}
\usepackage{relsize}
\usepackage{tikz}
\usetikzlibrary{trees}
\usetikzlibrary{cd,arrows,matrix,backgrounds,shapes,positioning,calc,decorations.markings,decorations.pathmorphing,decorations.pathreplacing}
\usepackage{multirow}
\usetikzlibrary{intersections,arrows,calc,decorations.markings}
\newtheorem{theorem}{Theorem}[section]
\usepackage{mathtools,leftidx}

\newtheorem{proposition}[theorem]{Proposition}

\newtheorem{conjecture}[theorem]{Conjecture}
\newtheorem{corollary}[theorem]{Corollary}
\newtheorem{lemma}[theorem]{Lemma}

\theoremstyle{definition}
\newtheorem{remark}[theorem]{Remark}
\newtheorem{example}[theorem]{Example}
\newtheorem{definition}[theorem]{Definition}

\newtheorem{question}[theorem]{Question}

\makeatletter
 
 \@addtoreset{equation}{section}
\makeatother

\def\ZZ{\mathbb Z}
\def\QQ{\mathbb{Q}}

\allowdisplaybreaks[4]
\title[Generalized discrete Markov spectra]{Generalized discrete Markov spectra}
\author{Yasuaki Gyoda}
\address[Yasuaki Gyoda]{Institute for Advanced Research, Nagoya University Furo-cho, Chikusa-ku, Nagoya-shi, 464-8601, Japan}
\email{ygyoda@math.nagoya-u.ac.jp}
\makeatletter
\@ifundefined{subjclassname@2020}{%
  \@namedef{subjclassname@2020}{\textup{2020} Mathematics Subject Classification}}{}
\makeatother
\begin{document}
\begin{abstract}
We develop a generalized Markov theory for the Markov--Lagrange and Markov spectra.  The classical discrete Markov spectrum is governed by Markov numbers, the positive integers occurring in solutions of the Markov equation.  We show that this relation admits a cluster-combinatorial extension governed by generalized Markov numbers.  Replacing the Christoffel-word formalism by snake graphs, we construct generalized discrete Markov spectra attached to the generalized Markov equations
\[
x^2+y^2+z^2+k_1yz+k_2zx+k_3xy=(3+k_1+k_2+k_3)xyz.
\]
Every element of these spectra is realized simultaneously as a Lagrange constant of a quadratic irrational and as a Markov constant of a real indefinite binary quadratic form.  We also prove structural results for these spectra, determine their contribution in the transition interval below Freiman's constant, and identify the boundary value obtained from irrational-slope lines that avoid all lattice points and the midpoints of all horizontal, vertical, and slope-$-1$ unit edges, again realizing it both as a Lagrange constant and as a Markov constant.
\end{abstract}
\maketitle
\tableofcontents
\section{Introduction and Main Results}\label{section:intro}
\subsection{Background}
The starting point of this paper is the classical link between Diophantine spectra and Markov numbers.  We recall it first in a form that will later be generalized by a cluster-combinatorial construction.  The \emph{Markov spectrum} is the set
\[
\mathcal{M}:= \left\{ \sup_{(x,y)\in \mathbb Z^2 \setminus \{0\}} \frac{\sqrt{D}}{|Q(x,y)|} \;\middle|\; Q(x,y) = ax^2 + bxy + cy^2, \; a,b,c\in \mathbb R,\; D = b^2 - 4ac > 0 \right\},
\]
which was studied by Markov around 1880. The Markov constant
\[
\mathcal M(Q):=\sup_{(x,y)\in \mathbb Z^2 \setminus \{0\}} \frac{\sqrt{D}}{|Q(x,y)|}
\]
is understood to be $+\infty$ if $Q(x,y)=0$ for some $(x,y)\in\mathbb Z^2\setminus\{0\}$.  It measures, in this setting, how closely the lines defined by the zeros of the quadratic form $Q(x,y)$ approach lattice points.  Markov identified the part of the spectrum below $3$ in terms of the integer solutions of a single cubic equation.  More precisely, a \emph{Markov number} is a positive integer that appears in a solution $(x,y,z) \in \mathbb Z_{>0}^3$ of
\[
x^2 + y^2 + z^2 = 3xyz.
\]
Markov proved the following theorem.
\begin{theorem}[Markov's theorem \cites{mar1,mar2}]
Let 
\begin{align*}\mathcal M_d& := \left\{ \frac{\sqrt{9n^2 - 4}}{n} \;\middle|\; n \text{ is a Markov number} \right\}\\
\mathcal M_{<3}& :=\{\mathcal M(Q) \mid Q: \text{real indefinite quadratic form  and } \mathcal M(Q) < 3\}\subset \mathcal M.
\end{align*}
Then $\mathcal M_d=\mathcal M_{< 3}$.
\end{theorem}
In this paper, we call the set $\mathcal{M}_d$ the \emph{discrete Markov spectrum}.  The same set also occurs on the approximation-theoretic side.  The \emph{Markov--Lagrange spectrum} is
\[
\mathcal{L} := \left\{ \limsup_{q \to \infty} \frac{1}{q \| q \alpha \|} \;\middle|\; \alpha \in \mathbb{R}\setminus\mathbb Q \right\},
\]
where $\| q \alpha \|$ denotes the distance from $q\alpha$ to the nearest integer.  The Lagrange constant
\[
\mathcal L(\alpha):=\limsup_{q \to \infty} \frac{1}{q \| q \alpha \|}
\]
is the supremum of those $L$ for which there exist infinitely many rationals $\frac{p}{q}$ satisfying
\[
\left| \alpha - \frac{p}{q} \right| <\frac{1}{L q^2}.
\]
Thus the Markov--Lagrange spectrum records the best rational approximations to irrational numbers on the $q^{-2}$ error scale.  Hurwitz's approximation theorem \cite{hurwitz} gives the universal lower bound $\mathcal L(\alpha)\geq\sqrt5$.  The full description below $3$ is part of Markov's theory: if
\[
\mathcal L_{<3}:=\{\mathcal L(\alpha)\mid \alpha\in\mathbb R\setminus\mathbb Q
\text{ and }\mathcal L(\alpha)<3\},
\]
then
\[
\mathcal L_{<3}=\mathcal M_{<3}=\mathcal M_d.
\]
The equality above combines Markov's descriptions on the approximation-theoretic and quadratic-form sides with the standard inclusion $\mathcal L\subset\mathcal M$; see, for example, \cites{aig,cusick-flahive}.
This coincidence is the classical bridge between Markov numbers, binary quadratic forms, and Diophantine approximation.

The regions of $\mathcal L$ and $\mathcal M$ above $3$ also remain central subjects of research.  Both spectra decompose into three regions: the discrete region below $3$, governed by Markov numbers; the transition interval $[3,c_F)$; and the continuous region above \emph{Freiman's constant} $c_F$, called \emph{Hall's ray}, where
\[
c_F:=\dfrac{2221564096 + 283748 \sqrt{462}}{491993569} \approx 4.5278295661\dots.
\]
Hall first proved that the spectra contain a half-line, and Freiman later determined the precise beginning of Hall's ray, namely the above value $c_F$ \cites{hall1947,freiman-hall-ray,freiman-markoff-spectrum}.  Outside the transition interval the Markov--Lagrange spectrum and the Markov spectrum coincide completely.  Freiman also proved that the two spectra are not equal by constructing Markov constants that are not Lagrange constants \cite{freiman-noncoincidence}.  Since the spectra coincide outside the transition interval, this gives the strict inclusion
\[
\mathcal L\cap[3,c_F)\subsetneq\mathcal M\cap [3,c_F).
\]
For further historical background, see \cites{aig, cusick-flahive, reutenauer}.
\subsection{Main Results}
The aim of this paper is to show that the classical bridge just recalled is not an isolated phenomenon.  We develop a generalized Markov theory in which generalized Markov numbers produce explicit Lagrange constants and Markov constants.  The construction is governed by the $(k_1,k_2,k_3)$-\emph{generalized Markov (GM) equation}
\[
x^2 + y^2 + z^2 + k_1 yz + k_2 zx + k_3 xy = (3 + k_1 + k_2 + k_3) xyz,
\]
where $k_1, k_2, k_3$ are nonnegative integers. A positive integer that appears in a solution to this equation is called a $(k_1,k_2,k_3)$-\emph{generalized Markov number}, or \emph{$(k_1,k_2,k_3)$-GM number} for short. This equation was introduced by \cites{gyo21,Bana} for the case $k_1=k_2=k_3=1$ and by \cite{gyomatsu} for the general case in the context of generalized cluster algebras. It has since been studied from several perspectives as a Markov-type Diophantine equation (for example, \cites{banaian, banaian-sen,bana-gyo,chen-li,courcy-litman-mizuno,gyoda-maruyama-sato,chen-li2,bao-li,chen-jia}).

Our first main result is a spectral realization theorem.  In the classical case, every element of $\mathcal M_d=\mathcal L_{<3}=\mathcal M_{<3}$ is realized on both the Lagrange and Markov sides.  The following theorem extends this simultaneous realization to GM numbers and is the basic mechanism from which the generalized discrete Markov spectra are defined.

\begin{theorem}[Theorem \ref{thm:markov-value-gen}, Corollaries \ref{cor:markov-lagrange}, \ref{cor:inL}]\label{thm:main1}
Fix $(k_1,k_2,k_3) \in \mathbb Z_{\geq 0}^3$. Let $n$ be a $(k_1,k_2,k_3)$-GM number, and let $i\in \{1,2,3\}$ be the position of $n$ in a positive integer solution to the $(k_1,k_2,k_3)$-GM equation. Set
\[\Delta(n,i):=((3+k_1+k_2+k_3)n-k_i)^2-4.\]
Then there exists a quadratic irrational $\alpha\in \mathbb Q(\sqrt{\Delta(n,i)})$ and the corresponding quadratic form $Q=(x-\alpha y)(x-\alpha' y)$, where $\alpha'$ is the quadratic conjugate of $\alpha$, such that 
\[\mathcal L(\alpha)=\mathcal M(Q)=\dfrac{\sqrt{\Delta(n,i)}}{n}.\]
In particular, for any $(k_1,k_2,k_3)\in\mathbb Z_{\geq 0}^3$, 
\[\mathcal{M}_{k_1,k_2,k_3} = \left\{ \frac{\sqrt{\Delta(n,i)}}{n} \;\middle|\; \begin{aligned}
&\text{the number $n$ is a $(k_1,k_2,k_3)$-GM number, and} \\
& i \text{ is the position of $n$ in a positive integer solution to}\\
& \text{the $(k_1,k_2,k_3)$-GM equation}
\end{aligned}\right\}
\] is contained in both $\mathcal L$ and $\mathcal M$.
\end{theorem}
For $(k_1,k_2,k_3)=(0,0,0)$ this recovers $\mathcal{M}_{d}$, that is, $\mathcal{M}_{d}=\mathcal M_{0,0,0}$.  We therefore call $\mathcal{M}_{k_1,k_2,k_3}$ the \emph{$(k_1,k_2,k_3)$-generalized discrete Markov spectrum}.

The proof is constructive.  We attach sign rules to rational-slope lines in the plane; the resulting triangle-crossing and edge-crossing data produce a periodic sequence of positive integers, and the corresponding periodic continued fraction gives the quadratic irrational $\alpha$ in Theorem \ref{thm:main1}.  The generalized Markov numbers are realized, at the same time, as numbers of perfect matchings of planar graphs assembled by joining square tiles in a chain, called \emph{snake graphs}.  This planar model is the central combinatorial object of the paper and is also the bridge to the irrational-slope boundary theorem below.

Classical treatments of Markov's theorem often use \emph{Christoffel words} (also called \emph{Cohn words}); see, for instance, \cite{aig}, \cite{bombieri}, and \cite{reutenauer}.  In that setting, Cohn matrices and strongly admissible sequences can both be read from Christoffel words, giving an efficient one-dimensional coding of the classical theory.

In the generalized setting one must additionally retain the types of the crossed edges and the parameters attached to them.  The planar line model records these data, and the associated snake graph converts them into the continued-fraction and perfect-matching quantities used below.  Thus the same geometric construction treats the generalized equations and, after setting all parameters equal to zero, specializes to an alternative construction of the classical sequences.

Markov's theorem characterizes the classical discrete Markov spectrum as the part of the Markov spectrum below $3$.  The generalized spectra constructed here therefore lead naturally to the following problem.

\begin{question}[Question \ref{ques:characterization}]
Is there a characterization of $\mathcal M_{k_1,k_2,k_3}$, or of their union, that generalizes Markov’s theorem?    
\end{question}

We next present two structural results about the generalized discrete Markov spectra that are independent of Question \ref{ques:characterization}.  The first relates the classical discrete Markov spectrum to the $(2,2,2)$-generalized discrete Markov spectrum.  It was shown in \cite{gyomatsu} that classical Markov numbers and $(2,2,2)$-GM numbers are connected through a relation involving squaring and taking square roots.  We prove the corresponding relation at the spectral level.
\begin{theorem}[Theorem \ref{thm:0-2-relation}]\label{thm:main2}
If $r\in \mathcal M_{0,0,0}$, then $3r\in \mathcal M_{2,2,2}$. Conversely, if $R\in \mathcal M_{2,2,2}$, then $\frac R3\in \mathcal M_{0,0,0}$.
\end{theorem}

We also determine exactly which generalized discrete Markov constants lie in the classical transition interval below Freiman's constant.

\begin{theorem}[Theorem \ref{thm:freiman-interval}]\label{thm:main3}
Set
\[\mathcal M':=\bigcup_{(k_1,k_2,k_3)\in \mathbb Z_{\geq 0}^3}\mathcal M_{k_1,k_2,k_3}.\]
Then we have
\[\mathcal M'\cap [3,c_F)=(\mathcal M_{0,0,1}\setminus\{\sqrt 5\})\cup\{2\sqrt 5\}.\]
\end{theorem}
The simultaneous realization of every number in $(\mathcal M_{0,0,1}\setminus\{\sqrt 5\})\cup\{2\sqrt 5\}$ as both a Lagrange constant and a Markov constant appears to be new.  At the same time, other Lagrange constants occur in the transition interval; in particular,
\[
\mathcal M'\cap [3,c_F)\subsetneq \mathcal L\cap [3,c_F).
\]
See Proposition \ref{rem:non-contain}.

The final main result passes from rational slopes to irrational slopes.  It identifies the boundary value approached by the rational-slope construction and realizes this value simultaneously on the Lagrange and Markov sides.  For a bi-infinite positive integer sequence $\mathbf a$, let $L(\mathbf a)$ denote the two-sided continued-fraction value defined in Section~\ref{section:Classical Markov Spectrum}.

\begin{theorem}[Theorem \ref{thm:irrational-slope-boundary-value}, Lemma \ref{lem:rational-values-approach-boundary}, Corollaries \ref{cor:irrational-slope-lagrange-value}, \ref{cor:irrational-slope-markov-value}]\label{thm:main4}
Fix $(k_1,k_2,k_3)\in\mathbb Z_{\geq 0}^3$ and $\sigma\in\mathfrak S_3$, and put
\[
K=3+k_1+k_2+k_3.
\]
The parameters assigned to horizontal, diagonal, and vertical edges are $k_{\sigma(1)}$, $k_{\sigma(2)}$, and $k_{\sigma(3)}$, respectively; see Section~\ref{section:Generalized Cohn matrices and their decompositions} for the sign rules. Let $l$ be a line of positive irrational slope that avoids every lattice point and the midpoint of every horizontal, vertical, and slope-$-1$ unit edge. Applying the triangle-crossing and edge-crossing sign rules along $l$, and then recording the lengths of the maximal consecutive blocks of equal signs, we obtain a bi-infinite positive integer sequence
\[
\mathbf b(l)=(b_n)_{n\in\mathbb Z},
\]
where the origin of the index is chosen arbitrarily. For every $r\in\mathbb Z$, set
\[
\alpha_r=[b_r;b_{r+1},b_{r+2},\ldots],
\qquad
\beta_r=-[0;b_{r-1},b_{r-2},\ldots],
\]
and define the real indefinite quadratic form
\[
Q_r(x,y)=(x-\alpha_r y)(x-\beta_r y).
\]
Then for every $r\in\mathbb Z$, we have
\[
L(\mathbf b(l))=\mathcal L(\alpha_r)=\mathcal M(Q_r)=K.
\]
In particular, $K\in\mathcal L\cap\mathcal M$. 
\end{theorem}

We end by formulating a uniqueness problem suggested by the spectral construction.  The classical starting point is Frobenius's uniqueness conjecture.

\begin{conjecture}[Conjecture \ref{uniqueness-conjecture}]
 For any Markov number $c$, there exists a unique Markov triple $(a,b,c)$ such that $a\leq b\leq c$.    
\end{conjecture}

A direct extension to GM numbers leads to the following question.

\begin{question}[Question \ref{uniqueness-question}]
\emergencystretch=2em
 For any $(k_1,k_2,k_3)$-GM number $c$, is there a unique $(k_1,k_2,k_3)$-GM triple $(a,b,c)$ such that $a\leq b\leq c$?   
\end{question}
However, this direct question has already been answered in the negative by \cite{gyomatsu}. The classical formulations in terms of Lagrange constants and injectivity of the Markov tree, both equivalent to the uniqueness conjecture, become separate questions after generalization. We answer both questions in the negative in general and, in particular, construct an infinite family of unequal-parameter counterexamples to injectivity (Questions \ref{uniqueness-question2} and \ref{uniqueness-question-gen}, and Proposition \ref{prop:injectivity-counterexamples}). By contrast, the injectivity conjecture restricted to $k_1=k_2=k_3$ remains open (Conjecture \ref{conj:equal-parameter-injectivity}).

\subsection{Organization}
Section \ref{section:Classical Markov Spectrum} reviews the classical discrete Markov spectrum and isolates the role played by admissible sequences and Cohn matrices.  It also compares the usual Christoffel-word construction with the planar construction used in this paper.  Section \ref{section:Generalized Markov equation and GM numbers} introduces the generalized Markov equation and generalized Markov numbers.  Section \ref{section:Generalized Cohn matrices and their decompositions} constructs generalized Cohn matrices and their matrix factorizations.  Section \ref{section:Generalized Markov length and generalized Markov distance} develops generalized Markov length and generalized Markov distance, the main tools in the proof of Theorem \ref{thm:main1}.  Section \ref{section:Generalized discrete Markov spectra} proves the spectral realization theorem, establishes the structural results above, and treats the irrational-slope boundary value and its realizations as Lagrange constants and Markov constants. Section \ref{section:Generalized Uniqueness Conjecture} discusses generalized uniqueness questions, constructs an infinite family of counterexamples to injectivity, and formulates the surviving equal-parameter injectivity conjecture.

\subsection*{Acknowledgements}
The author is deeply grateful to Esther Banaian for communicating Theorem \ref{thm:distance-theorem} and many valuable comments. The author would like to thank the participants of the Markov’s Theorem Seminar at the University of Tsukuba from January 27 to 31, 2025: Kota Saito, Hajime Kaneko, Tadahisa Hamada, Katsuki Ito, Satoru Oshima, and Takafumi Tsurumaki. 
Special thanks are due to Toshiki Matsusaka and Shigeki Akiyama for their valuable comments. 
The author is deeply grateful to Shoma Nakabayashi for providing Proposition \ref{prop:injectivity-counterexamples}, which gives an infinite family of counterexamples to Question \ref{uniqueness-question-gen}.
This work was supported by JSPS KAKENHI Grant Number JP25K17224.

\section{Classical Discrete Markov Spectrum}\label{section:Classical Markov Spectrum}
This section recalls several well-known facts about the classical Markov spectrum and reformulates them in a way that is compatible with the generalized theory developed later. 

\subsection{Markov Numbers and Markov Tree}
We consider the \emph{Markov equation}
\[x^2+y^2+z^2=3xyz.\]
A positive integer solution of the Markov equation is called a \emph{Markov triple}, and a positive integer that appears in a Markov triple is a \emph{Markov number.}  
We consider the following binary tree $\mathrm{M}\mathbb T$:
\begin{itemize}\setlength{\leftskip}{-10pt}
\item[(1)] The root vertex is $(1,2,1)$.
\item[(2)] For a vertex $(a,b,c)$, its two children are defined as
\[\begin{xy}(0,0)*+{(a,b,c)}="1",(30,-15)*+{\left(b,\dfrac{b^2+c^2}{a},c\right).}="2",(-30,-15)*+{\left(a,\dfrac{a^2+b^2}{c},b\right)}="3", \ar@{-}"1";"2"\ar@{-}"1";"3"
\end{xy}\]
\end{itemize}
This tree is called the \emph{Markov tree}. The first few vertices in $\mathrm{M}\mathbb{T}$ are as follows: 
\begin{align*}
\begin{xy}(10,0)*+{(1,2,1)}="1",(25,16)*+{(2,5,1)}="2",(25,-16)*+{(1,5,2)}="3", (50,24)*+{(5,13,1)}="4",(50,8)*+{(2,29,5)}="5",(50,-8)*+{(5,29,2)}="6",(50,-24)*+{(1,13,5)}="7", (85,28)*+{(13,34,1)\cdots}="8",(85,20)*+{(5,194,13)\cdots}="9",(85,12)*+{(29,433,5)\cdots}="10",(85,4)*+{(2,169,29)\cdots}="11",(85,-4)*+{(29,169,2)\cdots}="12",(85,-12)*+{(5,433,29)\cdots}="13",(85,-20)*+{(13,194,5)\cdots}="14",(85,-28)*+{(1,34,13)\cdots}="15",\ar@{-}"1";"2"\ar@{-}"1";"3"\ar@{-}"2";"4"\ar@{-}"2";"5"\ar@{-}"3";"6"\ar@{-}"3";"7"\ar@{-}"4";"8"\ar@{-}"4";"9"\ar@{-}"5";"10"\ar@{-}"5";"11"\ar@{-}"6";"12"\ar@{-}"6";"13"\ar@{-}"7";"14"\ar@{-}"7";"15"
\end{xy}
\end{align*}
The Markov tree has the following standard properties.
\begin{proposition}[\cite{gyoda-maruyama-sato}*{Proposition 3.2}]\label{prop:all-markov}
The following statements hold:
\begin{itemize}\setlength{\leftskip}{-10pt}
    \item[(1)] Each vertex $(a,b,c)$ in $\mathrm{M}\mathbb T$ is a Markov triple with $b > \max\{a,c\}$.
    \item[(2)] Every Markov triple $(a,b,c)$ with $b > \max\{a,c\}$ appears exactly once in $\mathrm{M}\mathbb T$.
\end{itemize}
\end{proposition}
Markov triples also satisfy the following coprimality property.
\begin{proposition}[\cite{aig}*{Corollary 3.4}]\label{relatively-prime}
For any Markov triple $(a,b,c)$, every pair among $a,b,c$ is relatively prime.
\end{proposition}

We next introduce the fraction labeling of Markov numbers. For this purpose, we recall Farey triples and the Farey tree, beginning with irreducible fractions in $\mathbb Q_{\geq 0}\cup \{\infty\}$.
\begin{definition}\label{irreducible-fraction}
Let $q\in\QQ_{\geq 0}\cup\{\infty\}$ and $n$ and $d\in\ZZ_{\geq 0}$. The symbol $\frac{n}{d}$ is called the \emph{reduced expression} of $q$ if $n$ and $d$ are relatively prime and $q = \frac{n}{d}$, where $\frac{n}{d}$ is regarded as $\infty$ when $d = 0$ and $n > 0$.
A fraction $\frac{n}{d}$ is called \emph{irreducible} if it is the reduced expression of some $q\in\QQ_{\geq 0}\cup\{\infty\}$.
\end{definition}
\begin{definition}
For $\frac{a}{b}$ and $\frac{c}{d}$, we denote $ad-bc$ by $\det\!\left(\frac{a}{b},\frac{c}{d}\right)$.  
A triple $\left(\frac{a}{b},\frac{c}{d},\frac{e}{f}\right)$ is called a \emph{Farey triple} if
\begin{itemize}\setlength{\leftskip}{-10pt}
    \item[(1)] $\frac{a}{b},\frac{c}{d},\frac{e}{f}$ are irreducible fractions, and
    \item[(2)] $\left|\det\!\left(\frac{a}{b},\frac{c}{d}\right)\right|=\left|\det\!\left(\frac{c}{d},\frac{e}{f}\right)\right|=\left|\det\!\left(\frac{e}{f},\frac{a}{b}\right)\right|=1$.
\end{itemize}
\end{definition}

We define the \emph{Farey tree} $\mathrm{F}\mathbb T$ as follows:
\begin{itemize}\setlength{\leftskip}{-10pt}
\item[(1)] The root vertex is $\left(\frac{0}{1},\frac{1}{1},\frac{1}{0}\right)$.
\item[(2)] Each vertex $\left(\frac{a}{b},\frac{c}{d},\frac{e}{f}\right)$ has the following two children:
\[\begin{xy}(0,0)*+{\left(\dfrac{a}{b},\dfrac{c}{d},\dfrac{e}{f}\right)}="1",(-30,-15)*+{\left(\dfrac{a}{b},\dfrac{a}{b}\oplus\dfrac{c}{d},\dfrac{c}{d}\right)}="2",(30,-15)*+{\left(\dfrac{c}{d},\dfrac{c}{d}\oplus\dfrac{e}{f},\dfrac{e}{f}\right).}="3", \ar@{-}"1";"2"\ar@{-}"1";"3"
\end{xy}\]
where $\frac{a}{b}\oplus\frac{c}{d}=\frac{a+c}{b+d}$.
\end{itemize}

The first few vertices of $\mathrm{F}\mathbb T$ are as follows:
\relsize{+1}
\begin{align*}
\begin{xy}(0,0)*+{\left(\frac{0}{1},\frac{1}{1},\frac{1}{0}\right)}="1",(20,-14)*+{\left(\frac{0}{1},\frac{1}{2},\frac{1}{1}\right)}="2",(20,14)*+{\left(\frac{1}{1},\frac{2}{1},\frac{1}{0}\right)}="3", 
(50,-24)*+{\left(\frac{0}{1},\frac{1}{3},\frac{1}{2}\right)}="4",(50,-8)*+{\left(\frac{1}{2},\frac{2}{3},\frac{1}{1}\right)}="5",(50,8)*+{\left(\frac{1}{1},\frac{3}{2},\frac{2}{1}\right)}="6",(50,24)*+{\left(\frac{2}{1},\frac{3}{1},\frac{1}{0}\right)}="7",(85,-28)*+{\left(\frac{0}{1},\frac{1}{4},\frac{1}{3}\right)\cdots}="8",(85,-20)*+{\left(\frac{1}{3},\frac{2}{5},\frac{1}{2}\right)\cdots}="9",(85,-12)*+{\left(\frac{1}{2},\frac{3}{5},\frac{2}{3}\right)\cdots}="10",(85,-4)*+{\left(\frac{2}{3},\frac{3}{4},\frac{1}{1}\right)\cdots}="11",(85,4)*+{\left(\frac{1}{1},\frac{4}{3},\frac{3}{2}\right)\cdots}="12",(85,12)*+{\left(\frac{3}{2},\frac{5}{3},\frac{2}{1}\right)\cdots}="13",(85,20)*+{\left(\frac{2}{1},\frac{5}{2},\frac{3}{1}\right)\cdots}="14",(85,28)*+{\left(\frac{3}{1},\frac{4}{1},\frac{1}{0}\right)\cdots}="15",\ar@{-}"1";"2"\ar@{-}"1";"3"\ar@{-}"2";"4"\ar@{-}"2";"5"\ar@{-}"3";"6"\ar@{-}"3";"7"\ar@{-}"4";"8"\ar@{-}"4";"9"\ar@{-}"5";"10"\ar@{-}"5";"11"\ar@{-}"6";"12"\ar@{-}"6";"13"\ar@{-}"7";"14"\ar@{-}"7";"15"
\end{xy}
\end{align*}
\relsize{-1}

\begin{proposition}[\cite{aig}*{Section 3.2}]\label{prop:property-farey}
The following hold:\begin{itemize}\setlength{\leftskip}{-10pt}
    \item[(1)] If $\left(r,t,s\right)$ is a Farey triple, then so are $\left(r,r\oplus t,t\right)$ and $\left(t,t\oplus s,s\right)$. In particular, every vertex in $\mathrm{F}\mathbb T$ is a Farey triple.
    \item[(2)] For every irreducible fraction $t \in (0,\infty)$, there exists a unique Farey triple $F$ in $\mathrm{F}\mathbb T$ such that $t$ is the second entry of $F$.
    \item[(3)] For $\left(r,t,s\right)$ in $\mathrm{F}\mathbb T$, the inequalities $r<t<s$ hold. 
\end{itemize}
\end{proposition}

We define $n_t$ to be the Markov number in $\mathrm{M}\mathbb T$ that corresponds to the irreducible fraction $t \in [0,\infty]$ at the corresponding position in $\mathrm{F}\mathbb T$.

For any Markov triple $(n_r,n_t,n_s)$ in $\mathrm{M}\mathbb T$, we consider an integer $u_t$ satisfying the following conditions:
\[\begin{cases}
n_ru_t\equiv n_s \pmod{n_t},\\
0<u_t<n_t.
\end{cases}\]
The integer $u_t$ is uniquely determined because $n_r$, $n_t$, and $n_s$ are pairwise relatively prime (Proposition \ref{relatively-prime}). Moreover, $u_t$ depends only on $t$, since $(r,t,s)$ is uniquely determined by $t$ according to Proposition \ref{prop:property-farey} (2).  
In this way, $u_t$ is defined for $t\in (0,\infty)$.  
We set $u_{\frac{0}{1}}=0$ and $u_{\frac{1}{0}}=1$.  
We call $u_t$ the \emph{characteristic number} of $t$. 
\begin{lemma}\label{lem:t-1/t-relation}
For any irreducible fraction $t \in [0,\infty]$,
\[
 n_{\frac1t}=n_t,
 \qquad
u_{\frac1t} = n_t - u_t.
\]
\end{lemma}
Lemma \ref{lem:t-1/t-relation} is essentially proved in \cite{gyoda-maruyama-sato}*{Proposition 7.22 (2)} (the case $k=0$). Indeed, if $(r,t,s)\in \mathrm{F}\mathbb{T}$, then $(\frac{1}{s},\frac{1}{t},\frac{1}{r})$ appears in $\mathrm{F}\mathbb{T}$ at the position symmetric to $(r,t,s)$ with respect to the central axis of the tree. The Markov tree has the same reflection symmetry, which gives $(n_\frac{1}{s}, n_\frac1t, n_\frac1r) = (n_s,n_t,n_r)$. Under this identification, the element $v_t^+$ in \cite{gyoda-maruyama-sato} coincides with $u_{\frac1t}$.
\subsection{Cohn Words, Cohn Matrices, and Cohn Tree}
Let $\mathfrak M_3$ be the free (non-commutative) monoid of rank $3$ generated by $\{p,q,r\}$.  
We consider the following tree, whose vertices are $\mathfrak M_3^3$:
\begin{itemize}\setlength{\leftskip}{-10pt}
\item[(1)] The root vertex is $(p,q,r)\in \mathfrak M_3^3$.
\item[(2)] For a vertex $(a,b,c)$, its two children are defined as
\[\begin{xy}(0,0)*+{(a,b,c)}="1",(-20,-10)*+{(a,ab,b)}="2",(20,-10)*+{(b,bc,c).}="3", \ar@{-}"1";"2"\ar@{-}"1";"3"
\end{xy}\]
\end{itemize}
This tree is called the \emph{Cohn word tree}, and we denote it by $\mathrm{CoW}\mathbb T$.  
We call the word of an element in $\mathfrak M_3$ appearing in $\mathrm{CoW}\mathbb T$ a \emph{Cohn word}.  

The first few Cohn words in $\mathrm{CoW}\mathbb T$ are as follows:
\begin{align*}
\begin{xy}(0,0)*+{(p,q,r)}="0",(20,15)*+{\left(q,qr,r\right)}="1",(20,-15)*+{\left(p,pq,q\right)}="1'",(45,22.5)*+{\left(qr,qr^2,r\right)}="20",(45,7.5)*+{\left(q,q^2r,qr\right)}="21",(45,-7.5)*+{\left(pq,pq^2,q\right)}="22",(45,-22.5)*+{\left(p,p^2q,pq\right)}="23",(80,26.25)*+{\left(qr^2,qr^3,r\right)}="30",(80,18.75)*+{\left(qr,qrqr^2,qr^2\right)}="31",(80,11.25)*+{\left(q^2r,q^2rqr,qr\right)}="32",(80,3.75)*+{\left(q,q^3r,q^2r\right)}="33",(80,-3.75)*+{\left(pq^2,pq^3,q\right)}="34",(80,-11.25)*+{\left(pq,pqpq^2,pq^2\right)}="35",(80,-18.75)*+{\left(p^2q,p^2qpq,pq\right)}="36",(80,-26.25)*+{\left(p,p^3q,p^2q\right)}="37",\ar@{-}"0";"1"\ar@{-}"0";"1'"\ar@{-}"1";"20"\ar@{-}"1";"21"\ar@{-}"1'";"22"\ar@{-}"1'";"23"\ar@{-}"20";"30"\ar@{-}"20";"31"\ar@{-}"21";"32"\ar@{-}"21";"33"\ar@{-}"22";"34"\ar@{-}"22";"35"\ar@{-}"23";"36"\ar@{-}"23";"37"
\end{xy}.
\end{align*}

We next give a geometric realization of Cohn words.  
For each irreducible fraction $t\in [0,\infty]$, let $c(t)$ denote the element of $\mathfrak M_3$ in $\mathrm{CoW}\mathbb T$ that corresponds to the position of $t$ in $\mathrm{F}\mathbb T$.  
We first construct a word $c_t$ geometrically and then identify it with $c(t)$.

We define a word $c_t$ for every irreducible fraction $t\in[0,\infty]$.  Set
\[
c_{\frac01}:=p,\qquad c_{\frac10}:=r.
\]
For $t=\frac{a}{b}\in(0,\infty)$, consider the line segment in $\mathbb{R}^2$ from $(0,0)$ to $(b,a)$, and denote it by $L_t$. We give $L_t$ the orientation from $(0,0)$ to $(b,a)$. Each time $L_t$ intersects a line of the integer lattice, we take the lattice point lying immediately to the right of the intersection point (with respect to the orientation of $L_t$), and record these in order as $p_1, p_2, \dots, p_n$. Here we set $p_1 = (0,0)$ and $p_n = (b,a)$. If the same lattice point appears more than once, we discard duplicates so that all $p_i$ are distinct. Connecting the points $p_1, p_2, \dots, p_n$ in order by line segments, we assign to each segment a letter according to its slope:  
- if the slope is $0$, we assign $p$;  
- if the slope is $1$, we assign $q$;  
- if the slope is $\infty$, we assign $r$.  

The sequence of these letters, taken in order, is defined to be $c_t$.

\begin{example}
The word $c_{\frac{2}{5}}$ is $p^2qpq$, and $c_{\frac{5}{2}}$ is $qrqr^2$. See Figure \ref{fig:christoffel}.
\begin{figure}[ht]
    \centering
    \begin{tikzpicture}
    % Draw grid
    \draw[thick] (0,0) grid (5,2);
    % Draw line with slope 2/5
    \draw[red,thick] (0,0) -- (5,2);
    \draw[line width=2pt] (0,0) -- (2,0)-- (3,1) -- (4,1) -- (5,2);
    \node at (0.5,-0.5) {$p$};
    \node at (1.5,-0.5) {$p$};
    \node at (2.5,-0.5) {$q$};
    \node at (3.5,-0.5) {$p$};
    \node at (4.5,-0.5) {$q$};
    \node[red] at (1.5,0.6) {\rotatebox{20}{$>$}};
    \node[red] at (3.5,1.4) {\rotatebox{20}{$>$}};
\end{tikzpicture}
\hspace{1cm}
 \begin{tikzpicture}
    % Draw grid
    \draw[thick] (0,0) grid (2,5);
    % Draw line with slope 2/5
    \draw[red,thick] (0,0) -- (2,5);
    \draw[line width=2pt] (0,0) -- (1,1)-- (1,2) -- (2,3) -- (2,5);
    \node at (2.5,0.5) {$q$};
    \node at (2.5,1.5) {$r$};
    \node at (2.5,2.5) {$q$};
    \node at (2.5,3.5) {$r$};
    \node at (2.5,4.5) {$r$};
    \node[red] at (0.6,1.5) {\rotatebox{65}{$>$}};
    \node[red] at (1.4,3.5) {\rotatebox{65}{$>$}};
\end{tikzpicture}
    \caption{The word $c_{\frac{2}{5}}$ and $c_{\frac{5}{2}}$}
    \label{fig:christoffel}
\end{figure}
\end{example}

\begin{remark}
The word $c_t$ is also known as a \emph{Christoffel word}.
\end{remark}

Let $\iota$ be the involution of the alphabet $\{p,q,r\}$ defined by
\[
\iota(p)=r,\qquad \iota(q)=q,\qquad \iota(r)=p.
\]
For a word $w=w_1\cdots w_m$, set
\[
w^\dagger:=\iota(w_m)\cdots\iota(w_1).
\]
Thus $\dagger$ is an involutive anti-automorphism of the free monoid.

\begin{lemma}\label{lem:cohn-word-reciprocity}
For every irreducible fraction $t\in[0,\infty]$,
\[
c_{1/t}=c_t^\dagger
\qquad\text{and}\qquad
c(1/t)=c(t)^\dagger.
\]
\end{lemma}

\begin{proof}
Reflection in the line $y=x$ interchanges horizontal and vertical steps.  Since it reverses the side from which the polygonal lattice path defining $c_t$ is read, it also reverses the order of its steps.  Applied to the geometric definition of $c_t$, this gives $c_{1/t}=c_t^\dagger$.

For the tree-defined words, reflect the Farey tree about its central axis.  A vertex $(r,t,s)$ is then replaced by $(1/s,1/t,1/r)$.  Simultaneously replace a Cohn-word vertex $(a,b,c)$ by $(c^\dagger,b^\dagger,a^\dagger)$.  This fixes the root $(p,q,r)$, and the identity $(uv)^\dagger=v^\dagger u^\dagger$ shows that it interchanges the two child operations.  Induction on the level of the tree therefore gives $c(1/t)=c(t)^\dagger$.
\end{proof}

We use the following standard description of Cohn words.

\begin{theorem}[\cite{aig}*{Theorem 7.6}]
For any irreducible fraction $t\in [0,\infty]$, $c_t=c(t)$ holds.    
\end{theorem}

\begin{proof}
The cited theorem proves the assertion for $t\in[0,1]$.  If $t\in(1,\infty]$, then $1/t\in[0,1)$.  Applying the cited result to $1/t$ and then Lemma \ref{lem:cohn-word-reciprocity} to both constructions gives
\[
c_t=c_{1/t}^\dagger=c(1/t)^\dagger=c(t).
\]
\end{proof}

Let  
\[
P=\begin{bmatrix}3&-1\\1&0\end{bmatrix},\quad
Q=\begin{bmatrix}5&2\\2&1\end{bmatrix},\quad
R=\begin{bmatrix}2&1\\1&1\end{bmatrix}.
\]
For any irreducible fraction $t\in [0,\infty]$, we define a matrix $C_t$ by
\[
C_t:=c_t\big|_{p\mapsto P,\ q\mapsto Q,\ r\mapsto R},
\]
and call it a \emph{Cohn matrix}.  
We define the \emph{Cohn tree} $\mathrm{Co}\mathbb{T}$ by
\[
\mathrm{Co}\mathbb{T}:=\mathrm{CoW}\mathbb{T}\big|_{(c_r,c_t,c_s)\mapsto(C_r,C_t,C_s)}.
\]

\begin{remark}\label{rem:cohn-mannar}
The Cohn matrix defined here differs from that in \cite{aig} in the arrangement of its entries. More precisely, the present definition is obtained from the one in \cite{aig} by swapping the diagonal entries.
\end{remark}

The first few vertices in $\mathrm{Co}\mathbb{T}$ are as follows:
\begin{align*}
\begin{xy}(0,0)*+{\left(\begin{bmatrix}3&-1\\1&0\end{bmatrix},\begin{bmatrix}5&{2}\\2&1\end{bmatrix},\begin{bmatrix}2&{1}\\1&1\end{bmatrix}\right)}="1",(40,-16)*+{\left(\begin{bmatrix}3&-1\\1&0\end{bmatrix},\begin{bmatrix}13&{5}\\5&2\end{bmatrix},\begin{bmatrix}5&{2}\\2&1\end{bmatrix}\right)}="2",(40,16)*+{\left(\begin{bmatrix}5&{2}\\2&1\end{bmatrix},\begin{bmatrix}12&{7}\\5&3\end{bmatrix},\begin{bmatrix}2&{1}\\1&1\end{bmatrix}\right)}="3", 
(100,-8)*+{\left(\begin{bmatrix}13&{5}\\5&2\end{bmatrix},\begin{bmatrix}75&31\\29&12\end{bmatrix},\begin{bmatrix}5&{2}\\2&1\end{bmatrix}\right)}="4",(100,-24)*+{\left(\begin{bmatrix}3&{-1}\\1&0\end{bmatrix},\begin{bmatrix}34&{13}\\13&5\end{bmatrix},\begin{bmatrix}13&{5}\\5&2\end{bmatrix}\right)}="5",(100,24)*+{\left(\begin{bmatrix}12&{7}\\5&3\end{bmatrix},\begin{bmatrix}31&{19}\\13&8\end{bmatrix},\begin{bmatrix}2&{1}\\1&1\end{bmatrix}\right)}="6",(100,8)*+{\left(\begin{bmatrix}5&{2}\\2&1\end{bmatrix},\begin{bmatrix}70&{41}\\29&17\end{bmatrix},\begin{bmatrix}12&{7}\\5&3\end{bmatrix}\right)}="7",\ar@{-}"1";"2"\ar@{-}"1";"3"\ar@{-}"2";"4"\ar@{-}"2";"5"\ar@{-}"3";"6"\ar@{-}"3";"7"
\end{xy}
\end{align*}

\begin{theorem}\label{thm:Ct-description}
For any irreducible fraction $t\in[0,\infty]$, we have
\[
C_t=\begin{bmatrix}
    3n_t-u_t & \dfrac{3u_tn_t-u_t^2-1}{n_t} \\[6pt]
    n_t & u_t
\end{bmatrix}.
\]
In particular, we have $\det (C_t)=1$ and $\mathrm{tr}(C_t)=3n_t$.
\end{theorem}

\begin{proof}
For $t\in[0,1]$, this is \cite{aig}*{Theorem 4.13}, after applying the change of matrix convention described in Remark \ref{rem:cohn-mannar}.  It remains to justify the extension to $t>1$.

Set
\[
H=\begin{bmatrix}2&1\\1&0\end{bmatrix}.
\]
A direct calculation gives
\[
HP^{\mathsf T}H^{-1}=R,\qquad
HQ^{\mathsf T}H^{-1}=Q,\qquad
HR^{\mathsf T}H^{-1}=P.
\]
Consequently, Lemma \ref{lem:cohn-word-reciprocity} implies
\[
C_t=H C_{1/t}^{\mathsf T}H^{-1}.
\]
By Lemma \ref{lem:t-1/t-relation}, $n_{1/t}=n_t$ and $u_{1/t}=n_t-u_t$.  Substituting these identities and the formula already proved for $C_{1/t}$ into the preceding matrix identity gives
\[
C_t=
\begin{bmatrix}
3n_t-u_t & \dfrac{3u_tn_t-u_t^2-1}{n_t}\\[6pt]
n_t&u_t
\end{bmatrix}.
\]
The determinant and trace assertions follow immediately.
\end{proof}

%--------------------------------------
\subsection{Continued Fraction Matrix Decompositions of Cohn Matrices}

In this subsection, we use geometric information to decompose a Cohn matrix into a product of matrices of the form $\begin{bmatrix}a&1\\1&0\end{bmatrix}$.  
We consider the square lattice in $\mathbb{R}^2$ (including its grid lines), subdivided by lines of slope $-1$ that pass through lattice points. We denote this modified lattice by $\widetilde{\mathbb{R}^2}$.  

For a given irreducible fraction $t=\frac{a}{b}\in (0,\infty)$, consider again the line segment in $\widetilde{\mathbb{R}^2}$ from $(0,0)$ to $(b,a)$, denoted $L_t$.  

\begin{figure}[ht]
    \centering
    \includegraphics[scale=0.3]{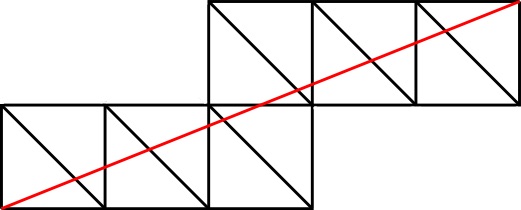}
    \caption{The line segment $L_t$ for $t=\frac{2}{5}$}
    \label{fig:ex-presnakegraph}
\end{figure}

We now define $\overline{L_t}$ as a sufficiently small left push-off of the oriented segment $L_t$, relative to the triangulation.  More precisely, choose a piecewise $C^1$ curve in a small tubular neighborhood of $L_t$ whose interior contains no lattice point, whose intersections with the edges are transverse and avoid all edge midpoints, and whose two endpoints are treated as follows: at the initial lattice point the curve enters through the lower-left edge, while at the terminal lattice point it ends in the interior of the upper-right edge.  At every interior lattice point of $L_t$, and whenever $L_t$ passes through an edge midpoint, the push-off is required to pass locally on the left-hand side of the oriented segment.  Because the segment is compact and has only finitely many intersections with the locally finite triangulation, all sufficiently small generic push-offs satisfying these requirements meet the same ordered list of triangle sides and lie in the same components of the edges with their midpoints removed.  Hence all sign data below are independent of the particular push-off.  We denote any such representative by $\overline{L_t}$.

\begin{figure}[ht]
    \centering
    \includegraphics[scale=0.079]{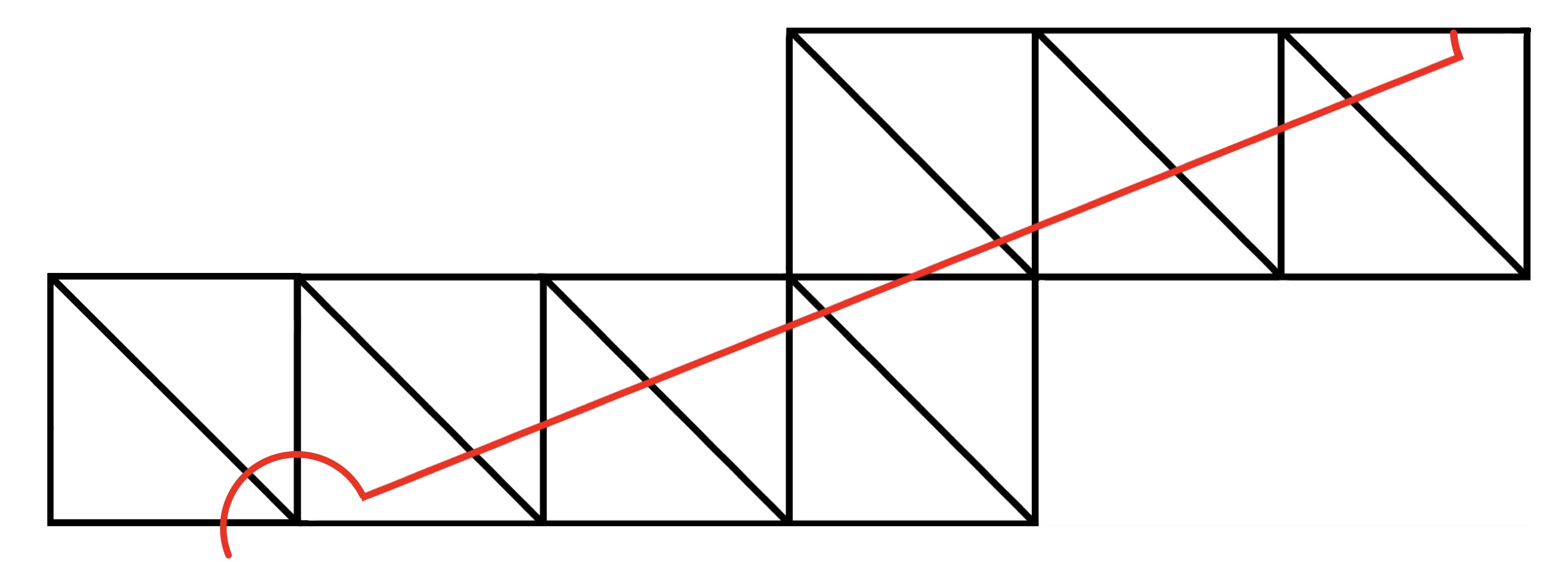}
    \caption{The line segment $\overline{L_t}$ associated with $t=\frac{2}{5}$}
    \label{fig:ex-extendedpresnakegraph}
\end{figure}

In Figure~\ref{fig:ex-extendedpresnakegraph}, this push-off is drawn with a larger deformation for visual clarity. The same convention is used in all subsequent figures.

\begin{remark}\label{rem:one-punctured-torus}
One can also obtain $\overline{L_t}$ using a once-punctured torus:  project $L_t$ onto an arc on the triangulated once-punctured torus, slightly push the arc upward to obtain a loop, and view $\overline{L_t}$ as the embedding of that loop back into $\widetilde{\mathbb R^2}$. Figure \ref{fig:interpretation-of-Lt} shows the case where $t=\frac{1}{1}$.
\begin{figure}[ht]
\centering
\includegraphics[scale=0.13]{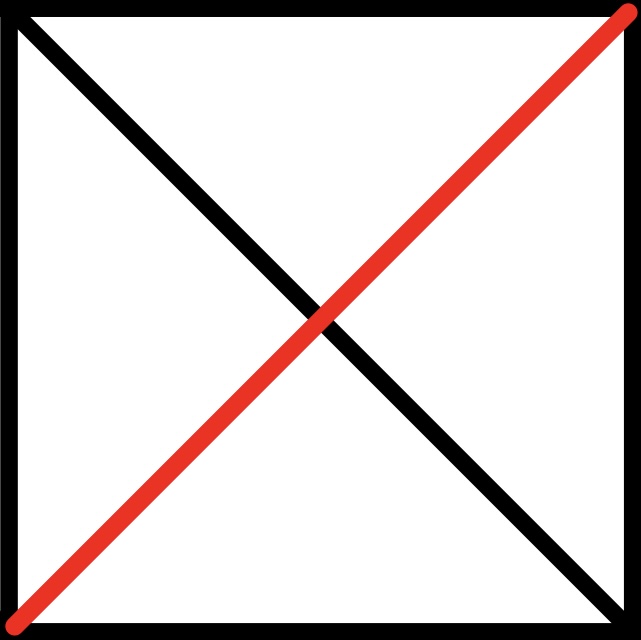}
\hspace{2mm}\raisebox{4\height}{$\mapsto$}
\scalebox{0.65}{\begin{tikzpicture}[baseline=-12mm]
\coordinate (u) at (0,1.5);
\coordinate (d) at (0,-1.5);
\coordinate (m1) at (-3,0);
\coordinate (m2) at (-1.5,0.4);
\coordinate (m3) at (1.5,0.4);
\coordinate (m4) at (3,0);
\coordinate(m5)at (-1.1,0.15);
\coordinate(m6)at (1.1,0.15);
\coordinate (p) at (-0.4,-0.8);
\coordinate (m7) at (0,-0.15);
 \fill(p) circle(3pt);
 \draw [thick](u) to [out=180,in=90] (m1);
 \draw [thick](m1) to [out=-90,in=180] (d);
 \draw [thick](d) to [out=0,in=-90] (m4);
 \draw [thick](m4) to [out=90,in=0] (u);
 \draw [thick](m2) to [out=-40,in=-140] (m3); 
 \draw [thick](m5) to [out=20,in=160] (m6); 
 \draw[thick] (d) to [out=170,in=-170] (m7);
 \draw[thick,dashed] (d) to [out=10,in=-10] (m7);
 \draw[thick] (p).. controls (-3,-0.5)and (-3, -0.5) .. (-3,0);
 \draw[thick] (p).. controls (3,-0.5)and (3, -0.5) .. (3,0);
 \draw[thick][dashed] (-3,0) .. controls (-3,0)and (-3,0.8) .. (0,1) to (0,1).. controls (3,0.8)and (3,0) .. (3,0);
\draw [thick](p) .. controls (1.6,-1.2)and (1.6, -1.2) ..  (2.1,-1.2); 
\draw [thick,dashed](2.1,-1.2)  .. controls (2.8,-0.8)and (2.8, -0.8) .. (1.3,0.3); 
\draw [thick](1.3,0.3)  .. controls (-2.5,2)and (-3.5, 0) .. (p); 
\draw [ultra thick,red](p) .. controls (-1.6,-1.2)and (-1.6, -1.2) ..  (-2.1,-1.2); 
\draw [ultra thick,dashed,red](-2.1,-1.2)  .. controls (-2.8,-0.8)and (-2.8, -0.8) .. (-1.3,0.3); 
\draw [ultra thick,red](-1.3,0.3)  .. controls (3,2)and (3.5, 0) .. (p); 
\end{tikzpicture}}
\raisebox{4\height}{$\mapsto$}
\scalebox{0.65}{\begin{tikzpicture}[baseline=-12mm]
\coordinate (u) at (0,1.5);
\coordinate (d) at (0,-1.5);
\coordinate (m1) at (-3,0);
\coordinate (m2) at (-1.5,0.4);
\coordinate (m3) at (1.5,0.4);
\coordinate (m4) at (3,0);
\coordinate(m5)at (-1.1,0.15);
\coordinate(m6)at (1.1,0.15);
\coordinate (p) at (-0.4,-0.8);
\coordinate (m7) at (0,-0.15);
 \fill(p) circle(3pt);
 \draw [thick](u) to [out=180,in=90] (m1);
 \draw [thick](m1) to [out=-90,in=180] (d);
 \draw [thick](d) to [out=0,in=-90] (m4);
 \draw [thick](m4) to [out=90,in=0] (u);
 \draw [thick](m2) to [out=-40,in=-140] (m3); 
 \draw [thick](m5) to [out=20,in=160] (m6); 
 \draw[thick] (d) to [out=170,in=-170] (m7);
 \draw[thick,dashed] (d) to [out=10,in=-10] (m7);
 \draw[thick] (p).. controls (-3,-0.5)and (-3, -0.5) .. (-3,0);
 \draw[thick] (p).. controls (3,-0.5)and (3, -0.5) .. (3,0);
 \draw[thick][dashed] (-3,0) .. controls (-3,0)and (-3,0.8) .. (0,1) to (0,1).. controls (3,0.8)and (3,0) .. (3,0);
\draw [thick](p) .. controls (1.6,-1.2)and (1.6, -1.2) ..  (2.1,-1.2); 
\draw [thick,dashed](2.1,-1.2)  .. controls (2.8,-0.8)and (2.8, -0.8) .. (1.3,0.3); 
\draw [thick](1.3,0.3)  .. controls (-2.5,2)and (-3.5, 0) .. (p); 
\draw [ultra thick,red](-0.4,-0.6) .. controls (-1.6,-1.2)and (-1.6, -1.2) ..  (-2.1,-1.2); 
\draw [ultra thick,dashed,red](-2.1,-1.2)  .. controls (-2.8,-0.8)and (-2.8, -0.8) .. (-1.3,0.3); 
\draw [ultra thick,red](-1.3,0.3)  .. controls (3,2)and (3.5, 0) .. (-0.4,-0.6); 
\end{tikzpicture}}
\raisebox{4\height}{$\mapsto$}\hspace{2mm}
\includegraphics[scale=0.13]{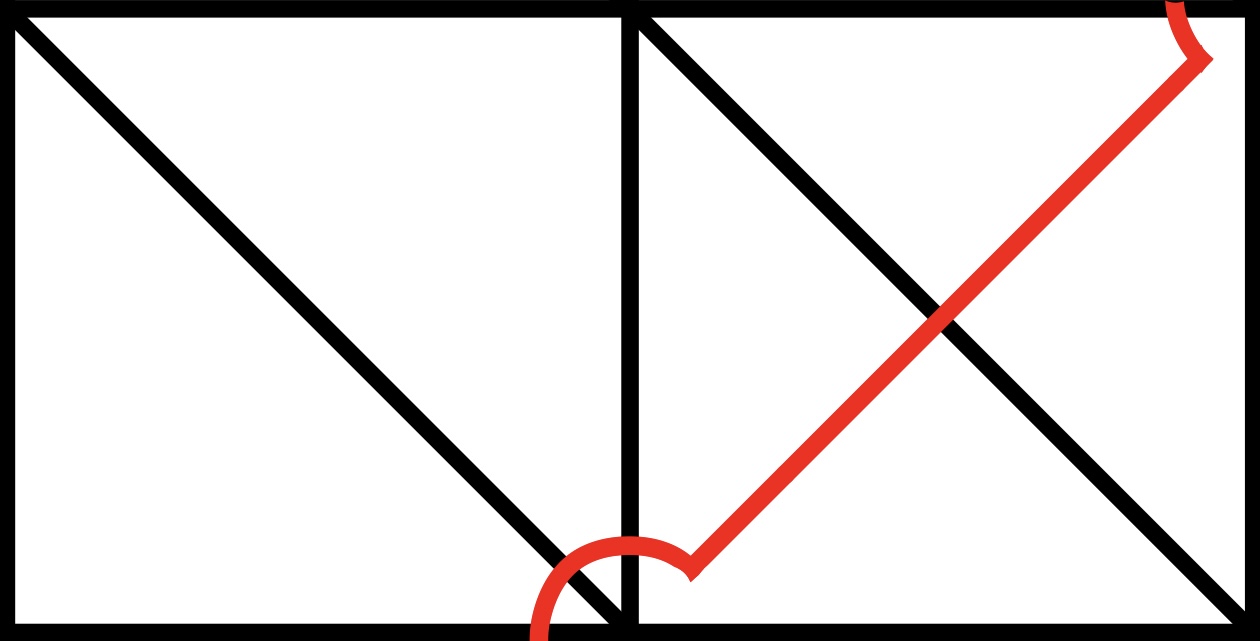}
\caption{Interpretation of $\overline{L_t}$ as a loop on the torus for $t=\frac 11$}\label{fig:interpretation-of-Lt}
\end{figure}
\end{remark}
We use $\overline{L_t}$ to construct an integer sequence for any irreducible fraction $t\in (0,\infty)$. For this purpose, we introduce the sign rule of triangles in $\widetilde{\mathbb R^2}$ associated with a curve segment. 

Let $\gamma\colon[0,1]\to\widetilde{\mathbb R^2}$ be a continuous, piecewise $C^1$ curve segment. Unless stated otherwise, a curve segment is assumed to satisfy the following conditions:
\begin{itemize}\setlength{\leftskip}{-10pt}
    \item $\gamma((0,1))$ does not contain a lattice point;
    \item whenever $\gamma((0,1))$ meets an edge of $\widetilde{\mathbb R^2}$, it intersects the edge transversely and locally passes from one adjacent triangle to the other;
    \item for every elementary triangle $\Delta$ and every connected component $I$ of $\gamma^{-1}(\operatorname{int}\Delta)$, the local branch $\gamma(\overline I)\cap\Delta$ is an embedded proper arc in $\Delta$: its interior lies in $\operatorname{int}\Delta$ and its two distinct endpoints lie on $\partial\Delta$;
    \item as the triangles and edges range over the triangulation, the sets $\gamma^{-1}(\operatorname{int}\Delta)$ have only finitely many connected components in total and $\gamma$ meets the edges at only finitely many parameter values;
    \item $\gamma$ has only finitely many self-intersections.
\end{itemize}

Here, for an edge $e$, $\operatorname{relint}(e)$ denotes the relative interior of $e$, namely, the edge with its two endpoints removed.

For an elementary triangle $\Delta$, a \emph{triangle-passage occurrence} of $\gamma$ is a pair $(\Delta,I)$, where $I$ is a connected component of $\gamma^{-1}(\operatorname{int}\Delta)$. For an edge $e$ and a parameter $u\in(0,1)$ such that $\gamma(u)\in\operatorname{relint}(e)$, the pair $(e,u)$ is called an \emph{edge-crossing occurrence}. These occurrences are ordered by the parameter of $\gamma$. Distinct connected components $I$, or distinct parameter values $u$, define distinct occurrences even if their underlying triangle or edge is the same.

\begin{definition}
Let $\gamma$ be an oriented curve segment on $\widetilde{\mathbb R^2}$. For every triangle-passage occurrence $(\Delta,I)$, the corresponding local branch $\gamma(\overline I)\cap\Delta$ cuts $\Delta$ into two pieces by the standing embeddedness assumption. Apply the following rule and assign a sign in $\{+,-\}$:
\begin{itemize}\setlength{\leftskip}{-10pt}
    \item[(i)] Assign a sign $(-)$ if the left-hand piece of $\Delta$, obtained by cutting it with this local branch, is a quadrilateral (see Figure \ref{fig:minus-righttriangles}):
    \begin{figure}[ht]
        \centering
        \begin{tikzpicture}[baseline=0mm]
        \draw (0,0) -- (1,0) -- (0,1) -- cycle;
        \draw[red,thick] (0.2,0.01) -- (0.5,0.5);
        \node[rotate=238,red] at (0.33,0.22) {$\mathbf <$};
        \end{tikzpicture}\hspace{0cm}
        \begin{tikzpicture}[baseline=0mm]
        \draw (1,0) -- (1,1) -- (0,1) -- cycle;
        \draw[red,thick] (1,0.8) -- (0.6,0.4);
        \node[rotate=230,red] at (0.8,0.6) {$\mathbf <$};
        \end{tikzpicture}\hspace{0.3cm}
        \rotatebox{180}{
        \begin{tikzpicture}[baseline=10mm]
        \draw (1,0) -- (1,1) -- (0,1) -- cycle;
        \draw[red,thick] (1,0.8) -- (0.6,0.4);
        \node[rotate=50,red] at (0.8,0.6) {$\mathbf >$};
        \end{tikzpicture}}
        \hspace{-0.3cm}
        \rotatebox{180}{\begin{tikzpicture}[baseline=10mm]
        \draw (0,0) -- (1,0) -- (0,1) -- cycle;
        \draw[red,thick] (0.2,0.01) -- (0.5,0.5);
        \node[rotate=58,red] at (0.33,0.22) {$\mathbf >$};
        \end{tikzpicture}}\hspace{0.3cm}
        \begin{tikzpicture}[baseline=0mm]
        \draw (0,0) -- (1,0) -- (0,1) -- cycle;
        \draw[red,thick] (0,0.5) -- (0.5,0);
        \node[rotate=315,red] at (0.25,0.25) {$\mathbf >$};
        \end{tikzpicture}\hspace{0cm}
        \begin{tikzpicture}[baseline=0mm]
        \draw (1,0) -- (1,1) -- (0,1) -- cycle;
        \draw[red,thick] (1,0.5) -- (0.5,1);
        \node[rotate=315,red] at (0.75,0.75) {$\mathbf <$};
        \end{tikzpicture}
        \caption{Right-angled triangles with $-$}
        \label{fig:minus-righttriangles}
    \end{figure}
    
    \item[(ii)] Assign a sign $(+)$ to every other triangle-passage occurrence (see Figure \ref{fig:plus-righttriangles}):
    \begin{figure}[ht]
        \centering
        \begin{tikzpicture}[baseline=0mm]
        \draw (0,0) -- (1,0) -- (0,1) -- cycle;
        \draw[red,thick] (0.2,0.01) -- (0.5,0.5);
        \node[rotate=238,red] at (0.33,0.22) {$\mathbf >$};
        \end{tikzpicture}\hspace{0cm}
        \begin{tikzpicture}[baseline=0mm]
        \draw (1,0) -- (1,1) -- (0,1) -- cycle;
        \draw[red,thick] (1,0.8) -- (0.6,0.4);
        \node[rotate=230,red] at (0.8,0.6) {$\mathbf >$};
        \end{tikzpicture}\hspace{0.3cm}    
        \rotatebox{180}{
        \begin{tikzpicture}[baseline=10mm]
        \draw (1,0) -- (1,1) -- (0,1) -- cycle;
        \draw[red,thick] (1,0.8) -- (0.6,0.4);
        \node[rotate=50,red] at (0.8,0.6) {$\mathbf <$};
        \end{tikzpicture}}
        \hspace{-0.3cm}
        \rotatebox{180}{\begin{tikzpicture}[baseline=10mm]
        \draw (0,0) -- (1,0) -- (0,1) -- cycle;
        \draw[red,thick] (0.2,0.01) -- (0.5,0.5);
        \node[rotate=58,red] at (0.33,0.22) {$\mathbf <$};
        \end{tikzpicture}}
        \hspace{0.3cm}
        \begin{tikzpicture}[baseline=0mm]
        \draw (0,0) -- (1,0) -- (0,1) -- cycle;
        \draw[red,thick] (0,0.5) -- (0.5,0);
        \node[rotate=315,red] at (0.25,0.25) {$\mathbf <$};
        \end{tikzpicture}\hspace{0cm}
        \begin{tikzpicture}[baseline=0mm]
        \draw (1,0) -- (1,1) -- (0,1) -- cycle;
        \draw[red,thick] (1,0.5) -- (0.5,1);
        \node[rotate=315,red] at (0.75,0.75) {$\mathbf >$};
        \end{tikzpicture}
        \caption{Right-angled triangles with $+$}
        \label{fig:plus-righttriangles} 
    \end{figure}
\end{itemize}
This rule is called the \emph{triangle-crossing rule} of $\gamma$.
\end{definition}

Next, we define the strongly admissible sequence $s(t)$\footnote{In \cite{aig}, the term “strongly admissible sequence” is used for the bi-infinite sequence ${}^\infty s(t)^\infty$ in the sense of this paper, whereas here we use this name for its single block $s(t)$.
}. For $t=\frac{0}{1}$ or $t=\frac{1}{0}$, set $s(\frac{0}{1})=s(\frac{1}{0})=(1,1)$. For $t\in (0,\infty)$, we define $s(t)$ as follows: 
\begin{itemize}\setlength{\leftskip}{-10pt}
    \item[(1)] Set the orientation of $\overline{L_t}$ from left to right, and arrange the signs assigned to its triangle-passage occurrences by the triangle-crossing rule in their order along $\overline{L_t}$.
    \item[(2)] A \emph{run} in the sign word from (1) is a maximal consecutive block of identical signs. If the runs in that word have lengths $a_1,\dots,a_\ell$, then $(a_1,\dots,a_\ell)$ is called the \emph{strongly admissible sequence} associated with $t$ and is denoted by $s(t)$.
\end{itemize}
The stability statement in the definition of $\overline{L_t}$ shows that $s(t)$ does not depend on the chosen sufficiently small generic push-off.
\begin{example}
Let $t=\frac{2}{5}$. From Figure \ref{fig:signextendedpresnake1}, the sign sequence is
\[
(-,-,+,-,+,-,+,+,-,-,+,-,+,+),
\]
and therefore
\[
s(\tfrac{2}{5})=(2,1,1,1,1,2,2,1,1,2).
\]

\begin{figure}[ht]
    \centering
    \includegraphics[scale=0.081]{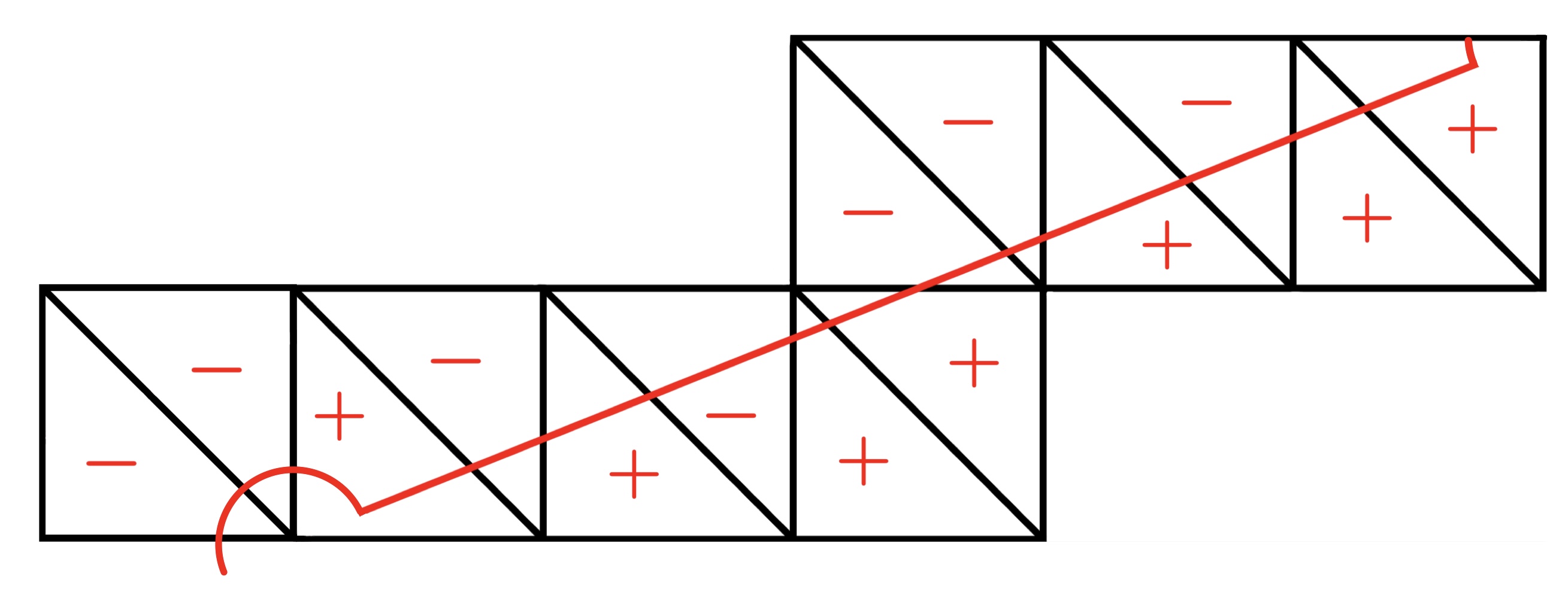}
    \caption{Signs of triangles intersecting $\overline{L_t}$ for $t=\frac{2}{5}$}
    \label{fig:signextendedpresnake1}
\end{figure}

Similarly, for $t=\frac{3}{2}$, Figure \ref{fig:signextendedpresnake2} shows the sign sequence
\[
(-,-,+,+,-,-,+,+,-,+),
\]
so that
\[
s(\tfrac32)=(2,2,2,2,1,1).
\]
\begin{figure}[ht]
    \centering
    \includegraphics[scale=0.079]{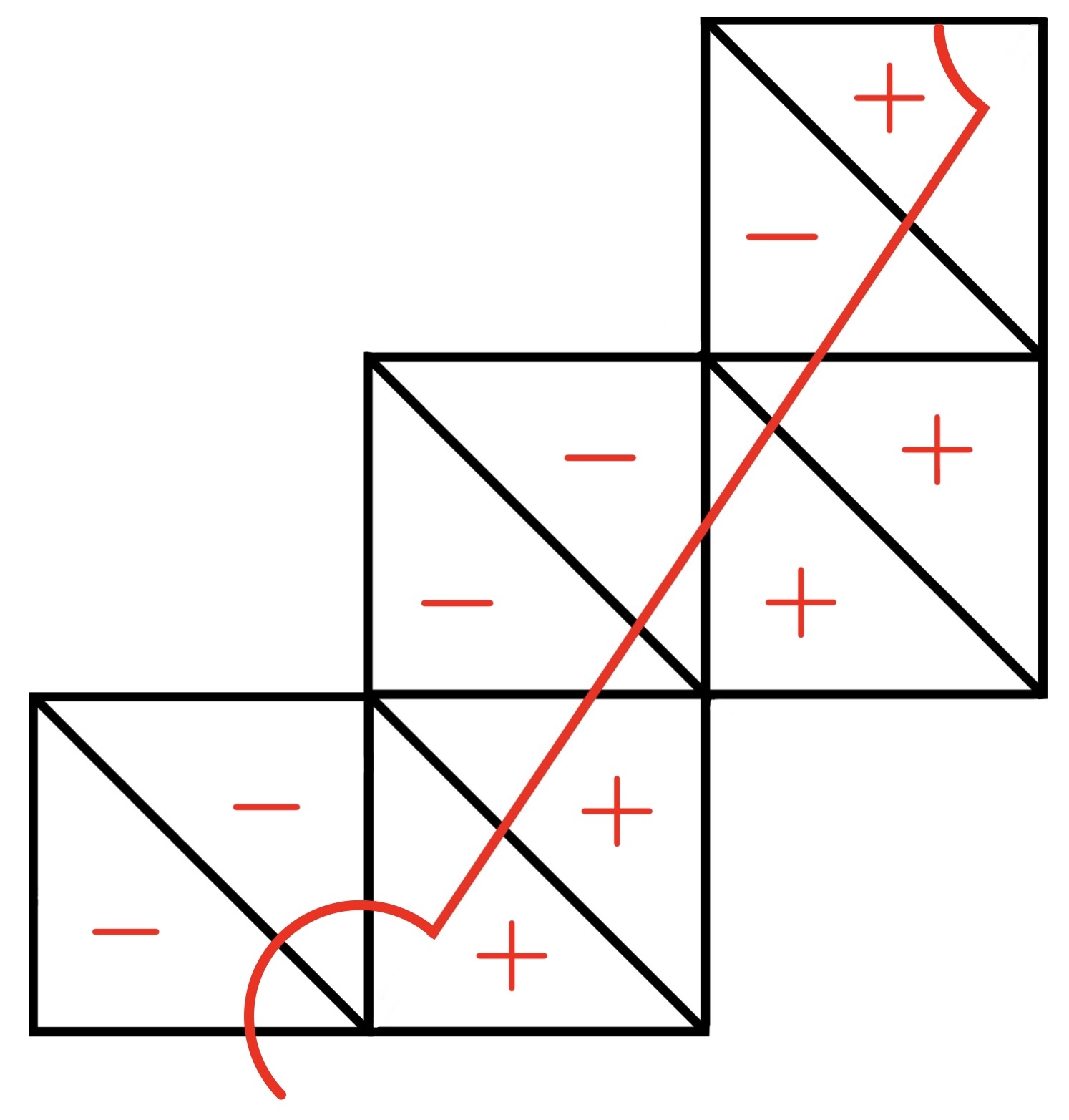}
    \caption{Signs of triangles intersecting $\overline{L_t}$ for $t=\frac{3}{2}$}
    \label{fig:signextendedpresnake2}
\end{figure}
\end{example}
\begin{proposition}\label{rem:t>1remark}
The following properties hold.
\begin{itemize}\setlength{\leftskip}{-10pt}
    \item[(1)] If $t\in[1,\infty]$, then $c_t$ consists of letters $q$ and $r$, and $s(t)$ is obtained from $c_t$ by replacing each $q$ with ``$2,2$" and each $r$ with ``$1,1$".
    \item[(2)] If $t\in (1,\infty)$, $s(t)$ has the form $(2,2,a_3,\dots,a_{n-2},1,1)$. Moreover, in this case,
    \[
    s\left(\tfrac{1}{t}\right)=(2,1,1,a_3,\dots,a_{n-2},2).
    \]
    Here the block $(a_3,\dots,a_{n-2})$ is understood to be empty when $n=4$.
\end{itemize}
\end{proposition}

In \cites{aig,bombieri}, $s(t)$ is constructed from $c_t$ using the replacement in (1).

\begin{proof}
For $t\geq1$, the polygonal lattice path defining $c_t$ has only diagonal and vertical steps, hence only the letters $q$ and $r$.  Inspecting one diagonal step and one vertical step in the triangulated square shows, respectively, the sign blocks
\[
q\longmapsto(2,2),\qquad r\longmapsto(1,1),
\]
which proves (1) by concatenation.

Now let $t>1$.  The first step of this path is diagonal and the last is vertical, so (1) gives
\[
s(t)=(2,2,a_3,\ldots,a_{n-2},1,1).
\]
Reflecting the push-off across $y=x$ produces the push-off for $1/t$ with the endpoint convention interchanged.  The reflection preserves the order of the interior local passages, and a direct check of the two reflected triangle types shows that it preserves their decomposition into runs.  Restoring the prescribed endpoint convention replaces the initial block $(2,2)$ by $(2,1,1)$ and the terminal block $(1,1)$ by $(2)$.  The interior block remains $(a_3,\ldots,a_{n-2})$, giving
\[
s(1/t)=(2,1,1,a_3,\ldots,a_{n-2},2).
\]
The statement also covers $n=4$, when the interior block is empty.
\end{proof}

For $a\in \mathbb Z_{\geq 1}$, the matrix
 $\begin{bmatrix}a & 1\\1 & 0\end{bmatrix}$
is called a \emph{continued fraction matrix}.  
For any finite positive integer sequence $S$, define $CF_S$ as the matrix product obtained by substituting each $a \in S$ with the corresponding continued fraction matrix.  
The following theorem shows that, for every irreducible fraction $t\in(0,\infty]$, the Cohn matrix $C_t$ is a product of continued fraction matrices.
\begin{theorem}\label{continued-fraction-theorem}
For any irreducible fraction $t \in (0,\infty] $, we have
\[
C_t = CF_{s(t)}.
\]
\end{theorem}

\begin{proof}[Proof of Theorem \ref{continued-fraction-theorem} for {$t \in [1,\infty]$}]
Since
\[
Q = \begin{bmatrix}5 & 2\\2 & 1\end{bmatrix} = \begin{bmatrix}2 & 1\\1 & 0\end{bmatrix}^2, \quad
R = \begin{bmatrix}2 & 1\\1 & 1\end{bmatrix} = \begin{bmatrix}1 & 1\\1 & 0\end{bmatrix}^2,
\]
and using Proposition \ref{rem:t>1remark} (1), the theorem follows.
\end{proof}

We now consider the case $t\in (0,1)$.

\begin{lemma}\label{lem:CF(1/t)-description}
For any irreducible fraction $t\in (1,\infty)$,
\[
CF_{s(\frac 1t)} = \begin{bmatrix}2 n_t + u_t & \dfrac{2 n_t^2 - n_t u_t - u_t^2 - 1}{n_t} \\ n_t & n_t - u_t \end{bmatrix}.
\]
\end{lemma}

\begin{proof}
By Proposition \ref{rem:t>1remark} (2), $s(\frac1t)$ is obtained from $s(t)$ by replacing the initial $2,2$ with $2,1,1$ and the final $1,1$ with $2$. The lemma follows from direct matrix computation (see also \cite{bombieri2}*{Corollary 24}).
\end{proof}

\begin{proof}[Proof of Theorem \ref{continued-fraction-theorem} for $t\in (0,1)$]   
By Theorem \ref{thm:Ct-description} and Lemmas \ref{lem:CF(1/t)-description} and \ref{lem:t-1/t-relation}, the \((2,1)\)- and \((2,2)\)-entries of \(C_t\) and \(CF_{s(t)}\) coincide. Since \(\det(C_t) = \det(CF_{s(t)}) = 1\) and \(\mathrm{tr}(C_t) = \mathrm{tr}(CF_{s(t)}) = 3n_t\), the claim follows.
\end{proof}
\begin{remark}
The proof of Theorem \ref{continued-fraction-theorem} uses the correspondence between the Cohn (Christoffel) word $c_t$ and the strongly admissible sequence $s(t)$. In Section~\ref{section:Generalized Cohn matrices and their decompositions}, we give another proof in a more general framework. 
\end{remark}

\subsection{Markov--Lagrange Spectrum and Markov Spectrum}
This section recalls the parts of the Markov--Lagrange spectrum and the Markov spectrum needed later. We first recall the Lagrange constant of an irrational number.
\begin{definition}
For $\alpha\in \mathbb R\setminus\mathbb Q$, define the \emph{Lagrange constant} of $\alpha$ by
\[
\mathcal L(\alpha):=\limsup_{q\to\infty}\frac{1}{q\|q\alpha\|}\in \mathbb R\cup\{+\infty\},
\]
where $\|q\alpha\|$ denotes the distance from $q\alpha$ to the nearest integer. Equivalently, $\mathcal L(\alpha)$ is the supremum of the real numbers $L$ for which
\[
\left| \alpha - \frac{p}{q} \right| < \frac{1}{L q^2}
\]
holds for infinitely many rational numbers $p/q$. We define
\[
\mathcal L := \{\mathcal L(\alpha) \mid \alpha \,\in\, \mathbb R \setminus \mathbb Q\},
\]
and call it the \emph{Markov--Lagrange spectrum}.
\end{definition}

We recall the Perron identity in a form suited to two-sided continued-fraction sequences. Let
$\mathbf a=(a_n)_{n\in\mathbb Z}$ be a bi-infinite positive integer sequence. For each $r\in\mathbb Z$, put
\[
\ell_r(\mathbf a)
:=[a_r;a_{r+1},a_{r+2},\ldots]+[0;a_{r-1},a_{r-2},\ldots].
\]
We then set
\[
L(\mathbf a):=\sup_{r\in\mathbb Z}\ell_r(\mathbf a).
\]

It is often convenient to specify the cut without referring to an index. If
$P=(p_1,p_2,\ldots)$ and $Q=(q_0,q_1,q_2,\ldots)$ are right-infinite positive integer sequences, we write $P^\ast$ for the sequence $P$ placed to the left of a cut, in the order
\[
P^\ast=(\ldots,p_2,p_1),
\]
and define
\[
\ell(P^\ast|Q):=[q_0;q_1,q_2,\ldots]+[0;p_1,p_2,\ldots].
\]
Thus $\ell(P^\ast|Q)$ is the index-free notation for $\ell_r(\mathbf a)$: when
$P=(a_{r-1},a_{r-2},\ldots)$ and $Q=(a_r,a_{r+1},\ldots)$, we have
\[
\ell_r(\mathbf a)=\ell(P^\ast|Q).
\]

Finally, let $T=(t_1,\ldots,t_m)$ be a finite positive integer sequence. We write
\[
T^\ast=(t_m,\ldots,t_1)
\]
for the reversed sequence, and we use
\[
{}^\infty T^\infty=\cdots TTT\cdots, \qquad
{}^\infty T=\cdots TT, \qquad
T^\infty=TT\cdots
\]
for the periodic extensions. The value $L({}^\infty T^\infty)$ is the maximum of $\ell(P^\ast|Q)$ over all cyclic cuts $P^\ast|Q$ of the bi-infinite periodic sequence ${}^\infty T^\infty$.

The following theorem is the Perron identity for the Markov--Lagrange spectrum.
\begin{theorem}[Perron identity \cites{perron1,perron2}]\label{thm:lagrange-quadratic}
Let $\alpha\in\mathbb R\setminus\mathbb Q$, and write its regular continued fraction expansion as
\[
\alpha=[a_0;a_1,a_2,a_3,\dots],
\]
where $a_0\in\mathbb Z$ and $a_i\in\mathbb Z_{>0}$ for $i\geq1$. Then
\[
\mathcal L(\alpha)=\limsup_{n\to\infty}
\left(
[a_{n+1};a_{n+2},a_{n+3},\dots]+[0;a_n,a_{n-1},\dots,a_1]
\right).
\]
In particular, for any finite positive integer sequence $T$, if $\alpha:=[T^\infty]$, then
\[
\mathcal L(\alpha)=L({}^\infty T^\infty).
\]
\end{theorem}
Thus, for a periodic continued fraction, the Lagrange constant is obtained by taking the largest cut value of the corresponding bi-infinite periodic sequence.

Next, we introduce the Markov spectrum.
\begin{definition}
Let $Q(x,y)=ax^2+bxy+cy^2$ be a real indefinite binary quadratic form, and set $D=b^2-4ac$. We define the \emph{Markov constant} of $Q$ by
\[
\mathcal M(Q):=\sup_{(x,y)\in \mathbb Z^2 \setminus \{0\}} \frac{\sqrt{D}}{|Q(x,y)|}\in \mathbb R\cup\{+\infty\},
\]
where this Markov constant is understood to be $+\infty$ if the denominator vanishes for some nonzero lattice point. We define 
\[
\mathcal M := \{\mathcal M(Q) \mid Q\colon\text{real indefinite binary quadratic form}\},
\]
and call it the \emph{Markov spectrum}.
\end{definition}

The same two-sided continued-fraction notation also gives the quadratic forms appearing in the Perron identity for the Markov spectrum. For a bi-infinite positive integer sequence $\mathbf a=(a_n)_{n\in\mathbb Z}$ and $r\in\mathbb Z$, set
\[
\alpha_r(\mathbf a):=[a_r;a_{r+1},a_{r+2},\ldots],
\qquad
\beta_r(\mathbf a):=-[0;a_{r-1},a_{r-2},\ldots],
\]
and define
\[
Q_{\mathbf a,r}(x,y):=(x-\alpha_r(\mathbf a)y)(x-\beta_r(\mathbf a)y).
\]
Then
\[
\ell_r(\mathbf a)=\alpha_r(\mathbf a)-\beta_r(\mathbf a).
\]
For a finite positive integer sequence $T$, set
\[
\alpha_T:=[T^\infty],\qquad
\beta_T:=-[0;(T^\ast)^\infty],\qquad
Q_T(x,y):=(x-\alpha_T y)(x-\beta_T y).
\]
Here $\beta_T$ is the quadratic conjugate of $\alpha_T$, and hence $Q_T$ has rational coefficients.

The following theorem is the Perron identity for the Markov spectrum.
\begin{theorem}[Perron identity for the Markov spectrum \cites{perron1,perron2}]\label{thm:markov-quadratic}
Let $\mathbf a=(a_n)_{n\in\mathbb Z}\in(\mathbb Z_{>0})^\mathbb Z$. For every $r\in\mathbb Z$, the Markov constant of the quadratic form $Q_{\mathbf a,r}$ is
\[
\mathcal M(Q_{\mathbf a,r})=L(\mathbf a)=\sup_{j\in\mathbb Z}\ell_j(\mathbf a),
\]
where this Markov constant may be $+\infty$. In particular, for any finite positive integer sequence $T$, we have
\[
\mathcal M(Q_T) = L({}^\infty T^\infty).
\]
Thus $\mathcal L([T^\infty])=\mathcal M(Q_T)=L({}^\infty T^\infty)$.
\end{theorem}

When $T$ is a strongly admissible sequence $s(t)$, the next theorem identifies a cyclic cut at which this maximum is attained.

\begin{theorem}[\cites{mar1,mar2}]\label{thm:markov-value}For any irreducible fraction $t\in [0,\infty]$, 
\begin{itemize}\setlength{\leftskip}{-10pt} 
\item[(1)] $[s(t)^\infty]\in \mathbb Q[\sqrt{9n_t^2-4}]$, and
\item[(2)] we have
\[
L({}^\infty s(t)^\infty) = \ell(^\infty s(t)|s(t)^\infty) = \ell\left(^\infty s\left(\tfrac{1}{t}\right)\,\middle|\,s\left(\tfrac{1}{t}\right)^\infty\right),
\]
and this value is
\[
\frac{\sqrt{9 n_t^2 - 4}}{n_t}.
\]
\end{itemize}
%Moreover, for any division $^\infty s(t)^\infty=P^\ast Q$ except for $^\infty s(t)s(t)^\infty$ and $^\infty s(\tfrac{1}{t})s(\tfrac{1}{t})^\infty$, we have
%\[L({}^\infty s(t)^\infty)>\ell(P^\ast|Q).\]
\end{theorem}

Here, the value $9n_t^2=(3n_t)^2$ is determined by the trace of $C_t$. Theorem \ref{thm:markov-value} is proved within a more general framework in Section \ref{section:Generalized discrete Markov spectra}. Theorems \ref{thm:lagrange-quadratic}, \ref{thm:markov-quadratic}, and \ref{thm:markov-value} immediately imply the following corollary.
\begin{corollary}\label{cor:markov-spectrum}
For any strongly admissible sequence $s(t)$, we have
\[
\mathcal L([s(t)^\infty]) = \mathcal M(Q_{s(t)})=\frac{\sqrt{9 n_t^2 - 4}}{n_t}.
\]
\end{corollary}

We set
\[
\mathcal M_{d} := \left\{ \frac{\sqrt{9 n^2 - 4}}{n} \ \middle| \ n: \text{Markov number} \right\},
\]
called the \emph{discrete Markov spectrum}.

There is an important characterization of $\mathcal M_d$.
Let
\begin{align*}
\mathcal L_{<3} &:= \{\mathcal L(\alpha) \mid \alpha \in \mathbb R \setminus \mathbb Q \text{ and } \mathcal L(\alpha) < 3\},\\
\mathcal M_{<3} &:=\{\mathcal M(Q) \mid Q: \text{real indefinite quadratic form  and } \mathcal M(Q) < 3\}.
\end{align*}
Then by Corollary \ref{cor:markov-spectrum} and the fact that $\frac{\sqrt{9 n^2 - 4}}{n}<3$, we have $\mathcal M_{d} \subset \mathcal L_{<3}$ and $\mathcal M_{d} \subset \mathcal M_{<3}$. The converse inclusions also hold.
\begin{theorem}[\cites{mar1,mar2}]\label{thm:characterization}
We have $\mathcal M_{d} = \mathcal L_{<3}=\mathcal{M}_{<3}$.
\end{theorem}

\section{Generalized Markov Equation and Generalized Markov Numbers}\label{section:Generalized Markov equation and GM numbers}
From this point on, we work with a generalized Markov equation. More precisely, we consider the \emph{$(k_1,k_2,k_3)$-generalized Markov equation}
\[
x^2 + y^2 + z^2 + k_1 yz + k_2 zx + k_3 xy = (3 + k_1 + k_2 + k_3) xyz,
\]
where $k_1, k_2, k_3 \in \mathbb Z_{\geq 0}$. We abbreviate this equation as the $(k_1,k_2,k_3)$-GM equation. 

Every permutation of a positive integer solution to this equation is called a \emph{$(k_1,k_2,k_3)$-GM triple}, and a positive integer that appears in some $(k_1,k_2,k_3)$-GM triple is called a \emph{$(k_1,k_2,k_3)$-GM number}.

We introduce the \emph{$(k_1,k_2,k_3,\sigma)$-generalized Markov tree} (abbreviated \emph{$(k_1,k_2,k_3,\sigma)$-GM tree}). Let $\mathfrak S_3$ be the symmetric group on 3 elements, acting on $\{1,2,3\}$ from the left. 

Define the $(k_1,k_2,k_3,\sigma)$-GM tree $\mathrm{M}\mathbb T(k_1,k_2,k_3,\sigma)$ for $\sigma \in \mathfrak S_3$ as follows:
\begin{itemize}\setlength{\leftskip}{-10pt}
\item [(1)] The root vertex is
\[
\big((1,\sigma(1)), (k_{\sigma(2)}+2, \sigma(2)), (1,\sigma(3))\big).
\]
\item [(2)] Every vertex $((a,h),(b,i),(c,j))$ has the following two children:
\[
\begin{xy}
(0,0)*+{((a,h),(b,i),(c,j))}="1",
(-40,-15)*+{\left((a,h),\left(\frac{a^2+k_j ab+b^2}{c},j\right),(b,i)\right)}="2",
(40,-15)*+{\left((b,i),\left(\frac{b^2+k_h bc+c^2}{a},h\right),(c,j)\right).}="3", 
\ar@{-}"1";"2"\ar@{-}"1";"3"
\end{xy}
\]
\end{itemize}

\begin{example}\label{ex:CMT(1,2,0,id)}
The first few vertices of $\mathrm{M}\mathbb T(1,2,0,\mathrm{id})$ are as follows:
\begin{align*}
\begin{xy}
(10,0)*+{((1,1),(4,2),(1,3))}="1",
(25,16)*+{((4,2),(21,1),(1,3))}="2",
(25,-16)*+{((1,1),(17,3),(4,2))}="3", 
(65,24)*+{((21,1),(121,2),(1,3))}="4",
(65,8)*+{((4,2),(457,3),(21,1))}="5",
(65,-8)*+{((17,3),(373,1),(4,2))}="6",
(65,-24)*+{((1,1),(81,2),(17,3))}="7", 
(120,28)*+{((121,2),(703,1),(1,3))\cdots}="8",
(120,20)*+{((21,1),(15082,3),(121,2))\cdots}="9",
(120,12)*+{((457,3),(57121,2),(21,1))\cdots}="10",
(120,4)*+{((4,2),(10033,1),(457,3))\cdots}="11",
(120,-4)*+{((373,1),(8185,3),(4,2))\cdots}="12",
(120,-12)*+{((17,3),(38025,2),(373,1))\cdots}="13",
(120,-20)*+{((81,2),(8227,1),(17,3))\cdots}="14",
(120,-28)*+{((1,1),(386,3),(81,2))\cdots}="15",
\ar@{-}"1";"2"\ar@{-}"1";"3"
\ar@{-}"2";"4"\ar@{-}"2";"5"
\ar@{-}"3";"6"\ar@{-}"3";"7"
\ar@{-}"4";"8"\ar@{-}"4";"9"
\ar@{-}"5";"10"\ar@{-}"5";"11"
\ar@{-}"6";"12"\ar@{-}"6";"13"
\ar@{-}"7";"14"\ar@{-}"7";"15"
\end{xy}
\end{align*}        
\end{example}
For any $\tau\in \mathfrak S_3$, set $\tau(x_1,x_2,x_3):=(x_{\tau(1)},x_{\tau(2)},x_{\tau(3)})$. For any $\sigma\in \mathfrak{S}_3$, we denote by $\sigma^{\ast}\in \mathfrak{S}_3$ the unique element distinct from $\sigma$ and satisfying $\sigma(2)=\sigma^\ast(2)$. The following theorem describes all GM triples obtained in this way.

\begin{theorem}\label{thm:all-solution}
Let $(a,b,c)$ be a $(k_1,k_2,k_3)$-GM triple with $b> \max\{a,c\}$, and assume $\tau(a,b,c)$ is a solution to the $(k_1,k_2,k_3)$-GM equation for $\tau\in \mathfrak S_3$. Then there exists a unique pair $\{\sigma,\sigma^{*}\}\subset \mathfrak S_3$ and a unique pair $(v,v^\ast)$ such that $v\in \mathrm{M}\mathbb T(k_1,k_2,k_3,\sigma)$ and $v^\ast\in \mathrm{M}\mathbb T(k_1,k_2,k_3,\sigma^\ast)$, and $v$ and $v^\ast$ are permutations of $((a,\tau^{-1}(1)),(b,\tau^{-1}(2)),(c,\tau^{-1}(3)))$. Furthermore, the position of $v^\ast$ in $\mathrm{M}\mathbb{T}(k_1,k_2,k_3,\sigma^\ast)$ is symmetric, with respect to the center, to the position of $v$ in $\mathrm{M}\mathbb{T}(k_1,k_2,k_3,\sigma)$.
\end{theorem}

\begin{proof}
The descent theorem \cite{gyomatsu}*{Theorem 1} shows that every positive solution other than $(1,1,1)$ occurs in one of the three trees generated by the mutation operations.  The graph isomorphisms in \cite{bana-gyo}*{Proposition 3.8} identify these three trees with the six position-labelled GM trees used here.  The two trees having the same middle position label are mirror images by \cite{bana-gyo}*{Remark 3.9}, while the uniqueness of the labelled vertex with a prescribed strictly largest middle entry is \cite{bana-gyo}*{Corollary 3.10}.  These statements give the asserted pair of vertices and its uniqueness.
\end{proof}

\begingroup\emergencystretch=2em
In each pair appearing in a vertex of $\mathrm{M}\mathbb{T}(k_1, k_2, k_3, \sigma)$, the first component is a $(k_1, k_2, k_3)$-GM number. Moreover, for each triple, the three numbers given by the first components form a solution to the $(k_1, k_2, k_3)$-GM equation when placed in the positions indicated by the second components of the corresponding pairs. These pairs are called \emph{$(k_1,k_2,k_3,\sigma)$-GM-number-position pairs}.
\par\endgroup

\begin{proposition} 
GM trees have the following two properties.
\begin{itemize}\setlength{\leftskip}{-10pt}
    \item [(1)] If $k_1=k_2=k_3$, then the $(k_1,k_2,k_3,\sigma)$-GM tree is independent of $\sigma$ up to the position labels. We often omit $\sigma$ in this case.
    \item [(2)] For any $(k_1,k_2,k_3)$-GM triple $(a,b,c)$, if $b=\max\{a,c\}$, then $(a,b,c)=(1,1,1)$. Consequently, the union
    \[
    \bigcup_{\sigma\in\mathfrak S_3}
    \operatorname{Vert}\!\bigl(\mathrm{M}\mathbb T(k_1,k_2,k_3,\sigma)\bigr)
    \]
    represents every $(k_1,k_2,k_3)$-GM triple other than $(1,1,1)$.
\end{itemize}
\end{proposition}
\begin{proof}
Part~(1) follows immediately from the definition of a GM tree when the three parameters are equal. Part~(2) follows from \cite{gyomatsu}*{Lemma 4} together with Theorem~\ref{thm:all-solution}.
\end{proof}
The following proposition generalizes Proposition \ref{relatively-prime}. 
\begin{proposition}[\cite{gyomatsu}*{Corollary 8}]\label{relatively-prime-generalized}
For any $(k_1,k_2,k_3)$-GM triple $(a,b,c)$, any two of $a,b,c$ are relatively prime.
\end{proposition}

Now, fix $(k_1,k_2,k_3)\in \mathbb Z_{\geq 0}^3$ and $\sigma \in \mathfrak S_3$.

Abusing notation as in the previous section, we denote by $(n_{t},i_t)$ the $(k_1,k_2,k_3,\sigma)$-GM-number-position pair in $\mathrm{M}\mathbb T(k_1,k_2,k_3,\sigma)$ located at the position of the irreducible fraction $t\in [0,\infty]$ in $\mathrm{F}\mathbb{T}$. 

For any $(k_1,k_2,k_3)$-GM triple $(n_{r},n_{t},n_{s})$ in $\mathrm{M}\mathbb T(k_1,k_2,k_3,\sigma)$, consider an integer $u_t$ satisfying:
\[
\begin{cases}
n_{r} u_{t} \equiv n_s \pmod{n_t},\\
0 < u_t < n_t.
\end{cases}
\]
For $t\in(0,\infty)$, one has $n_t>1$: the root value is $n_1=k_{\sigma(2)}+2\geq2$, and every new middle entry in the GM tree is strictly larger than either adjacent entry.  Proposition \ref{relatively-prime-generalized} therefore implies that $n_r$ is invertible modulo $n_t$ and that the residue $n_r^{-1}n_s$ is nonzero.  It has a unique representative in the interval $0<u_t<n_t$.  Moreover, $u_t$ depends only on $t$, since $(r,t,s)$ is determined by $t$ from Proposition \ref{prop:property-farey} (2).

For $t\in (0,\infty)$, we define $u_t$ as above, and set 
\[
u_{\frac{0}{1}}=-k_{\sigma(1)},\quad u_{\frac{1}{0}}=1.
\]
The number $u_t$ is called the \emph{characteristic number} of $t$.

\begin{proposition}\label{rem:dual-remark}
Let $(n^\ast_{\frac{1}{t}},i^*_\frac{1}{t})$ be a $(k_1,k_2,k_3,\sigma^{\ast})$-GM number-position pair associated with $\frac1t$. Then
\[
(n_t,i_t) = (n_\frac{1}{t}^{\ast},i^\ast_{\frac{1}{t}}).
\]
\end{proposition}
\begin{proof}
This follows immediately from Theorem~\ref{thm:all-solution}.
\end{proof}
We abbreviate $k_{i_t}$ as $k_t$. The following lemma generalizes Lemma \ref{lem:t-1/t-relation}.
\begin{lemma}\label{lem:t-1/t-relation-gen}
For any irreducible fraction $t \in (0,\infty)$,
\[
u^{\ast}_{\frac1t} = n_t - u_t-k_t.
\]
\end{lemma}
For $k_1=k_2=k_3$, Lemma \ref{lem:t-1/t-relation-gen} is proved in \cite{gyoda-maruyama-sato}*{Proposition 7.22 (2)}. The general case is proved in the same way. As in the case of Lemma \ref{lem:t-1/t-relation}, replace $v_t^+$ in \cite{gyoda-maruyama-sato} with $u^*_{\frac{1}{t}}$.

\section{Generalized Cohn Matrices and Their Decompositions}\label{section:Generalized Cohn matrices and their decompositions}
This section introduces analogues of Cohn matrices corresponding to GM numbers and derives their continued-fraction matrix decompositions.
\subsection{Generalized Cohn Matrices}
Fix $k_1, k_2, k_3 \in \mathbb Z_{\geq 0}$ and $\sigma \in \mathfrak{S}_3$. We set
\begin{align*}
C_{\frac{0}{1}} &= 
\begin{bmatrix}
3+k_1+k_2+k_3 & -(3+k_1+k_2+k_3)k_{\sigma(1)}-1\\
1 & -k_{\sigma(1)}
\end{bmatrix},\\
C_{\frac{1}{1}} &= 
\begin{bmatrix}
(3+k_1+k_2+k_3)(k_{\sigma(2)}+2)-k_{\sigma(2)}-1 & 2+k_1+k_2+k_3\\
k_{\sigma(2)}+2 & 1
\end{bmatrix},\\
C_{\frac{1}{0}} &= 
\begin{bmatrix}
2+k_1+k_2+k_3-k_{\sigma(3)} & 1+k_1+k_2+k_3-k_{\sigma(3)}\\
1 & 1
\end{bmatrix}.
\end{align*}

For any $(r,t,s)\in \mathrm{F}\mathbb{T}$, we set inductively
\[
C_{r\oplus t} := C_{r}C_t - D_s, \quad C_{t\oplus s} = C_t C_s - D_r,
\]
where
\[
D_r = \begin{bmatrix}
k_r & k_r(3+k_1+k_2+k_3)\\
0 & k_r
\end{bmatrix}.
\]

These matrices coincide with ordinary Cohn matrices when $k_1=k_2=k_3=0$. The matrix $C_t$ is called a \emph{$(k_1,k_2,k_3,\sigma)$-generalized Cohn matrix}. Define the \emph{$(k_1,k_2,k_3,\sigma)$-generalized Cohn tree}:
\[
\mathrm{Co}\mathbb{T}(k_1,k_2,k_3,\sigma) := \mathrm{F}\mathbb{T}\big|_{(r,t,s)\mapsto(C_r,C_t,C_s)}.
\]

\begin{example}
The first few vertices in $\mathrm{Co}\mathbb T(1,2,0,\mathrm{id})$ are:
\begin{align*}
\scalebox{0.95}{\begin{xy}
(0,0)*+{\left(\begin{bmatrix}6 & -7\\1 & -1\end{bmatrix}, \begin{bmatrix}21 & 5\\4 & 1\end{bmatrix}, \begin{bmatrix}5 & 4\\1 & 1\end{bmatrix}\right)}="1",
(40,-16)*+{\left(\begin{bmatrix}6 & -7\\1 & -1\end{bmatrix}, \begin{bmatrix}98 & 23\\17 & 4\end{bmatrix}, \begin{bmatrix}21 & 5\\4 & 1\end{bmatrix}\right)}="2",
(40,16)*+{\left(\begin{bmatrix}21 & 5\\4 & 1\end{bmatrix}, \begin{bmatrix}109 & 83\\21 & 16\end{bmatrix}, \begin{bmatrix}5 & 4\\1 & 1\end{bmatrix}\right)}="3", 
(105,-8)*+{\left(\begin{bmatrix}98 & 23\\17 & 4\end{bmatrix}, \begin{bmatrix}2149 & 507\\373 & 88\end{bmatrix}, \begin{bmatrix}21 & 5\\4 & 1\end{bmatrix}\right)}="4",
(105,-24)*+{\left(\begin{bmatrix}6 & -7\\1 & -1\end{bmatrix}, \begin{bmatrix}467 & 98\\81 & 17\end{bmatrix}, \begin{bmatrix}98 & 23\\17 & 4\end{bmatrix}\right)}="5",
(105,24)*+{\left(\begin{bmatrix}109 & 83\\21 & 16\end{bmatrix}, \begin{bmatrix}626 & 507\\121 & 98\end{bmatrix}, \begin{bmatrix}5 & 4\\1 & 1\end{bmatrix}\right)}="6",
(105,8)*+{\left(\begin{bmatrix}21 & 5\\4 & 1\end{bmatrix}, \begin{bmatrix}2394 & 1823\\457 & 348\end{bmatrix}, \begin{bmatrix}109 & 83\\21 & 16\end{bmatrix}\right)}="7",
\ar@{-}"1";"2"\ar@{-}"1";"3"\ar@{-}"2";"4"\ar@{-}"2";"5"\ar@{-}"3";"6"\ar@{-}"3";"7"
\end{xy}}
\end{align*}
\end{example}

\begin{theorem}\label{thm:Mt-description}
For any irreducible fraction $t\in [0,\infty]$, we have
\[
C_t = \begin{bmatrix}
(3+k_1+k_2+k_3)n_t - k_t - u_t & \dfrac{(3+k_1+k_2+k_3)n_t u_t - k_t u_t - u_t^2 - 1}{n_t} \\
n_t & u_t
\end{bmatrix}.
\] In particular, we have $\det (C_t)=1$ and $\mathrm{tr}(C_t)=(3+k_1+k_2+k_3)n_t-k_t$.
\end{theorem}

The case $k_1=k_2=k_3$ and $t\in[0,1]$ is treated in \cite{gyo-maru}*{Lemma 4.5}; the same argument proves the general case.

\subsection{Continued Fraction Matrix Decompositions of Generalized Cohn Matrices}

As in Section \ref{section:Classical Markov Spectrum}, we use geometric information to decompose $C_t$ into a product of matrices of the form $\begin{bmatrix}a&1\\1&0\end{bmatrix}$. For any irreducible fraction $t\in [0,\infty]$, we construct a \emph{generalized strongly admissible sequence} $s(t)$ 
by the triangle-crossing rule introduced in Section  \ref{section:Classical Markov Spectrum} and the following \emph{edge-crossing rule}.

Fix $(k_1,k_2,k_3)\in\mathbb Z_{\geq0}^3$ and $\sigma\in\mathfrak S_3$. Let $\mathcal E_1,\mathcal E_2,\mathcal E_3$ be the sets of horizontal, diagonal, and vertical unit edges, respectively, of the elementary triangulation, put
\[
\mathcal E:=\mathcal E_1\cup\mathcal E_2\cup\mathcal E_3,
\]
and let $m_e$ denote the midpoint of a unit edge $e$.

\begin{definition}
Let $\gamma$ be an oriented curve segment satisfying the standing finiteness and transversality assumptions of Section~\ref{section:Classical Markov Spectrum}, and assume that
\[
\gamma((0,1))\cap\{m_e\mid e\in\mathcal E\}=\varnothing.
\]
For each edge-crossing occurrence $(e,u)$, let $i$ be the type of $e$ and put $\kappa=k_{\sigma(i)}$. Assign signs as follows:
\begin{itemize}\setlength{\leftskip}{-10pt}
    \item [(i)] If $\kappa>0$ and $m_e$ lies on the left side of the oriented local branch of $\gamma$ at $u$, assign $\kappa$ minus signs ($-$) to this occurrence (see Figure \ref{fig:minus-edge}).
    \begin{figure}[ht]
        \centering
        \begin{tikzpicture}[baseline=0mm]
        \draw (-0.5,0) -- (0.5,0);
        \fill (0,0) circle (1.5pt);
        \draw[red] (-0.5,-0.3) -- (0.5,0.1);
        \node[rotate=205,red] at (0.3,0.01) {$\mathbf <$};
        \end{tikzpicture}
        \hspace{0.5cm}
        \begin{tikzpicture}[baseline=0mm]
        \draw (-0.5,0) -- (0.5,0);
        \fill (0,0) circle (1.5pt);
        \draw[red] (-0.5,-0.1) -- (0.5,0.3);
        \node[rotate=207,red] at (-0.3,-0.01) {$\mathbf >$};
        \end{tikzpicture}
        \hspace{0.5cm}
        \begin{tikzpicture}[baseline=0mm]
        \draw (-0.5,0.5) -- (0.5,-0.5);
        \fill (0,0) circle (1.5pt);
        \draw[red] (-0.5,-0.3) -- (0.5,0);
        \node[rotate=200,red] at (0.29,-0.06) {$\mathbf <$};
        \end{tikzpicture}
        \hspace{0.6cm}
        \begin{tikzpicture}[baseline=0mm]
        \draw (-0.5,0.5) -- (0.5,-0.5);
        \fill (0,0) circle (1.5pt);
        \draw[red] (-0.5,-0) -- (0.5,0.3);
        \node[rotate=200,red] at (-0.29,0.06) {$\mathbf >$};
        \end{tikzpicture}
        \hspace{0.7cm}
        \begin{tikzpicture}[baseline=0mm]
        \draw (0,-0.5) -- (0,0.5);
        \fill (0,0) circle (1.5pt);
        \draw[red] (-0.5,-0.3) -- (0.5,0);
        \node[rotate=200,red] at (0.29,-0.06) {$\mathbf <$};
        \end{tikzpicture}
        \hspace{0.7cm}
        \begin{tikzpicture}[baseline=0mm]
        \draw (0,-0.5) -- (0,0.5);
        \fill (0,0) circle (1.5pt);
        \draw[red] (-0.5,-0) -- (0.5,0.3);
        \node[rotate=200,red] at (-0.29,0.06) {$\mathbf >$};
        \end{tikzpicture}
        \caption{Edge-crossing configurations producing $-$}
        \label{fig:minus-edge}
    \end{figure}
    \item [(ii)] If $\kappa>0$ and $m_e$ lies on the right side of the oriented local branch of $\gamma$ at $u$, assign $\kappa$ plus signs ($+$) to this occurrence (see Figure \ref{fig:plus-edge}).
    \begin{figure}[ht]
        \centering
        \begin{tikzpicture}[baseline=0mm]
        \draw (-0.5,0) -- (0.5,0);
        \fill (0,0) circle (1.5pt);
        \draw[red] (-0.5,-0.3) -- (0.5,0.1);
        \node[rotate=205,red] at (-0.18,-0.18) {$\mathbf >$};
        \end{tikzpicture}
        \hspace{0.5cm}
        \begin{tikzpicture}[baseline=0mm]
        \draw (-0.5,0) -- (0.5,0);
        \fill (0,0) circle (1.5pt);
        \draw[red] (-0.5,-0.1) -- (0.5,0.3);
        \node[rotate=207,red] at (0.23,0.18) {$\mathbf <$};
        \end{tikzpicture}
        \hspace{0.5cm}
        \begin{tikzpicture}[baseline=0mm]
        \draw (-0.5,0.5) -- (0.5,-0.5);
        \fill (0,0) circle (1.5pt);
        \draw[red] (-0.5,-0.3) -- (0.5,0);
        \node[rotate=200,red] at (-0.29,-0.24) {$\mathbf >$};
        \end{tikzpicture}
        \hspace{0.6cm}
        \begin{tikzpicture}[baseline=0mm]
        \draw (-0.5,0.5) -- (0.5,-0.5);
        \fill (0,0) circle (1.5pt);
        \draw[red] (-0.5,-0) -- (0.5,0.3);
        \node[rotate=200,red] at (0.29,0.24) {$\mathbf <$};
        \end{tikzpicture}
        \hspace{0.7cm}
        \begin{tikzpicture}[baseline=0mm]
        \draw (0,-0.5) -- (0,0.5);
        \fill (0,0) circle (1.5pt);
        \draw[red] (-0.5,-0.3) -- (0.5,0);
        \node[rotate=200,red] at (-0.29,-0.24) {$\mathbf >$};
        \end{tikzpicture}
        \hspace{0.7cm}
        \begin{tikzpicture}[baseline=0mm]
        \draw (0,-0.5) -- (0,0.5);
        \fill (0,0) circle (1.5pt);
        \draw[red] (-0.5,-0) -- (0.5,0.3);
        \node[rotate=200,red] at (0.29,0.24) {$\mathbf <$};
        \end{tikzpicture}
        \caption{Edge-crossing configurations producing $+$}
        \label{fig:plus-edge}
    \end{figure}
\end{itemize}
If $\kappa=0$, assign no sign and make no left--right test.
This rule is called the \emph{edge-crossing rule} of $\gamma$ for $(k_1,k_2,k_3,\sigma)$.
\end{definition}

For the sign rules above, take $\overline L_t$ to be any sufficiently small left push-off satisfying the conditions in Section~\ref{section:Classical Markov Spectrum}. Thus every edge crossing is transverse, no crossing occurs at an edge midpoint, the curve leaves its initial lattice point through the lower-left edge, and it terminates in the interior of the upper-right edge at its terminal lattice point. There are only finitely many crossings, and every degenerate crossing of $L_t$ is resolved on the prescribed left side. Hence, if the push-off is sufficiently small, the ordered list of crossed triangles and edges and the side on which the curve passes each crossed-edge midpoint are independent of its choice. Therefore all such push-offs give the same sign word, and the sequence defined below is well defined.

We define the $(k_1,k_2,k_3,\sigma)$-generalized strongly admissible sequence $s(t)$. First, set
\[
s(\tfrac{0}{1})=(1+k_{\sigma(2)}+k_{\sigma(3)},1),\quad s(\tfrac{1}{0})=(1+k_{\sigma(1)}+k_{\sigma(2)},1).
\]
For any irreducible fraction $t\in (0,\infty)$, orient $\overline L_t$ from left to right and read the signs assigned by the triangle-crossing and edge-crossing rules in parameter order, recording every occurrence separately. The maximal consecutive blocks of identical signs in this word are its runs. If their lengths are $a_1,\dots,a_\ell$, define $s(t)=(a_1,\dots,a_\ell)$.

\begin{example}\label{ex:t=2/5}
Let $(k_1,k_2,k_3,\sigma)=(1,2,0,\textrm{id})$ and $t=\frac{2}{5}$. Figure \ref{fig:signextendedpresnake1gen} shows the signs along $\overline L_t$, giving
\[
s(\tfrac{2}{5})=(5,1,3,3,1,5,4,1,3,4).
\]
\begin{figure}[ht]
    \centering
    \includegraphics[scale=0.09]{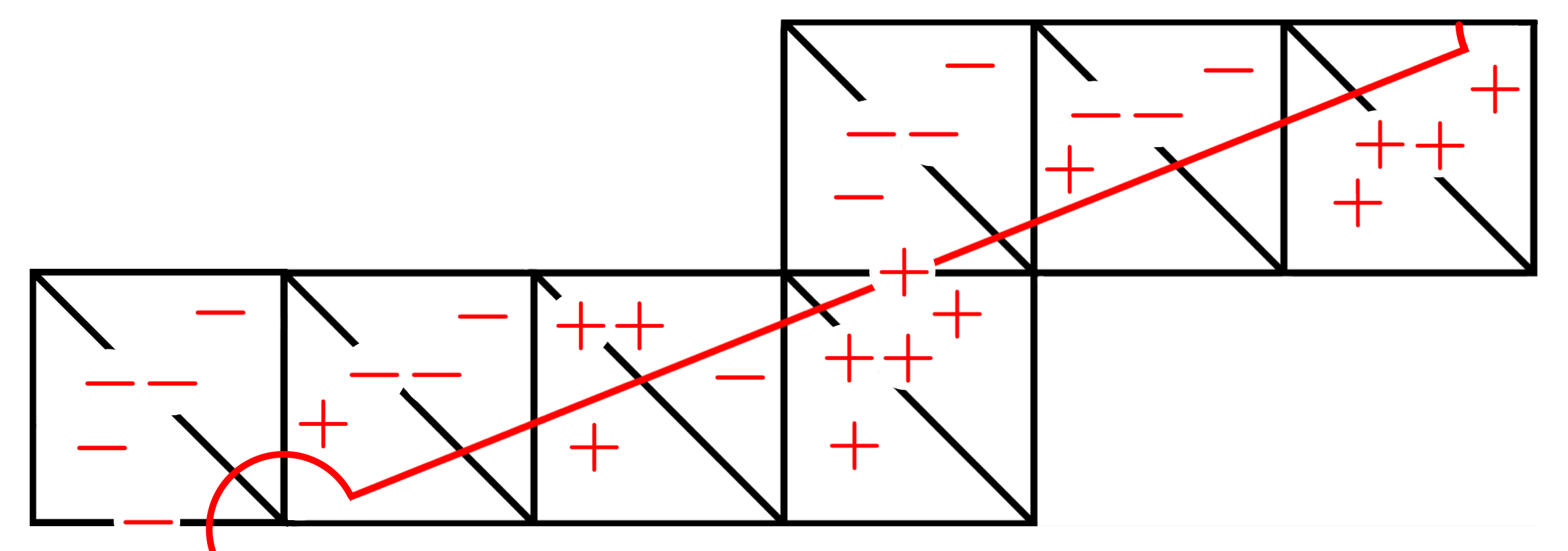}
    \caption{Signs of triangles and edges intersecting $\overline L_t$ for $t=\frac{2}{5}$}
    \label{fig:signextendedpresnake1gen}
\end{figure}
The right endpoint of $\overline{L_t}$ lies within the upper-rightmost edge in $\widetilde {\mathbb R^2}$; therefore no sign is assigned to this edge. Although the unique horizontal edge intersecting $\overline{L_t}$ appears to meet it at its midpoint, in fact, since $\overline{L_t}$ is defined by shifting $L_t$ slightly to the left, the intersection point lies slightly to the left of the midpoint. Thus the corresponding sign is $+$.
\end{example}

\begin{proposition}[Structure of generalized strongly admissible sequences]
\label{rem:difference-(0,1)(1,infty)}
Let $t\in(0,\infty)$ and write $s(t)=(a_1,\dots,a_n)$.  Then the following statements hold.
\begin{itemize}\setlength{\leftskip}{-10pt}
    \item[(1)] The number $n$ of entries is even.  The sequences for $t=\frac01$ and $t=\frac10$ also have an even number of entries.
    \item[(2)] $a_1=2+k_1+k_2+k_3$.
    \item[(3)] If $t=\frac11$, then
    \[
    s(t)=(2+k_1+k_2+k_3,\,2+k_{\sigma(2)}).
    \]
    \item[(4)] If $t\in(0,1)$, then $a_2=1$ and $a_n\neq1$.  If $t\in(1,\infty)$, then $a_2\neq1$ and $a_n=1$.
    \item[(5)] If $t=\frac12$, then $n=4$ and $a_4=a_3+1+k_t$.  If $t=\frac21$, then $n=4$ and $a_2=a_3+1+k_t$.
    \item[(6)] If $t\in(0,\frac12)\cup(\frac12,1)$, then $n\geq6$ and
    \[
    a_3+1=a_n,\qquad
    a_{3+i}=a_{n-i}\quad
    \left(1\leq i\leq\frac n2-3\right),
    \]
    and
    \[
    a_{\frac n2+1}
      =a_{\frac n2+2}+(-1)^{\frac n2+1}k_t.
    \]
    \item[(7)] If $t\in(1,2)\cup(2,\infty)$, then $n\geq6$ and
    \[
    a_2=a_{n-1}+1,\qquad
    a_{2+i}=a_{n-i-1}\quad
    \left(1\leq i\leq\frac n2-3\right),
    \]
    and
    \[
    a_{\frac n2}
      =a_{\frac n2+1}+(-1)^{\frac n2}k_t.
    \]
    \item[(8)] Let $s^\ast(1/t)$ denote the sequence constructed with $\sigma^\ast$.  Then
    \[
    s^\ast(1/t)=(a_1,a_n,a_{n-1},\dots,a_2).
    \]
    Hence ${}^\infty(s^\ast(1/t))^\infty$ is the reverse of ${}^\infty s(t)^\infty$, up to a shift of the cut.  The analogous assertion holds for $t=\frac01$ and $t=\frac10$.
\end{itemize}
\end{proposition}

\begin{proof}
We verify the assertions from the two local sign rules.

Write $t=p/q$ with $p,q>0$ and $\gcd(p,q)=1$, and first retain the individual signs before grouping adjacent equal signs into runs.  We may use a polygonal representative arbitrarily close to $L_t$.  In a small neighborhood of either endpoint, its crossing pattern is determined solely by the pair of incident edges through which the local branch passes.  Inspecting the two incident elementary triangles and the three types of incident unit edges gives the following facts: the first run has two triangle signs together with $k_1+k_2+k_3$ edge signs; the initial and terminal signs are opposite, so the number of runs is even; below slope $1$ the second run has length one, whereas above slope $1$ the last run has length one.  At slope $1$, the same local inspection gives the two runs displayed in (3).  These facts prove (1)--(4), including $a_1=2+k_1+k_2+k_3$.

The remaining relations follow from the half-turn symmetry.  The map
\[
h_t(x)=(q,p)-x
\]
pairs, in reverse order, the triangle and edge crossings of the undeformed segment $L_t$.  It preserves each of the three edge types and their midpoints.  Resolve a paired crossing by the prescribed left push-off.  Apart from the unique crossing at the midpoint of $L_t$, which is not paired by the half-turn, the two resulting local sign blocks have the same length.  The endpoint pair contributes one additional triangle sign on the terminal side when $0<t<1$, and on the initial side when $t>1$.  The unpaired edge crossing is at $(q/2,p/2)$.  Its edge type, and hence the position label of $n_t$, is
\[
i_t=
\begin{cases}
\sigma(1),&p\text{ even and }q\text{ odd},\\
\sigma(2),&p\text{ and }q\text{ odd},\\
\sigma(3),&p\text{ odd and }q\text{ even}.
\end{cases}
\]
Indeed, these are respectively the horizontal, diagonal, and vertical edges containing $(q/2,p/2)$.  The equality with the GM position label holds directly at $2/1$, $1/1$, and $1/2$; when one passes from a Farey-tree vertex to either of its children, the three parity classes rotate in exactly the same way as the position labels in the GM tree.  The equality therefore follows for every $t$ by induction on its depth in the Farey tree.  The middle edge block has length $k_{i_t}=k_t$.

If the middle block meets the endpoint block, which occurs exactly for $t=1/2$ or $t=2$, this comparison gives the two formulas in (5).  Otherwise there is at least one paired run on either side, so $n\geq6$; reading the paired runs from the two ends gives the equalities in (6) and (7).  Whether the middle edge block joins the run immediately before it or the run immediately after it is determined by the alternating signs.  This gives, respectively,
\[
a_{\frac n2+1}-a_{\frac n2+2}=(-1)^{\frac n2+1}k_t
\quad\text{and}\quad
a_{\frac n2}-a_{\frac n2+1}=(-1)^{\frac n2}k_t,
\]
which completes the proof of (5)--(7).

Finally, reflection in the diagonal, $(x,y)\mapsto(y,x)$, sends the crossing pattern of slope $t$ to that of slope $1/t$, reverses the planar left--right convention, interchanges horizontal and vertical edges, and fixes diagonal edges.  The reflected curve lies on the right of the reflected oriented segment.  Apply the half-turn about the midpoint to this curve and orient it again from its initial endpoint to its terminal endpoint.  This produces the prescribed left push-off and reverses the order of the crossings.  The permutation $\sigma^\ast$ interchanges precisely the corresponding first and third edge weights and keeps the second one fixed.  Thus the resulting sign word is obtained by reversing the original sign word and then cyclically shifting it so that the cut lies immediately after its first run.  Its run lengths are $(a_1,a_n,a_{n-1},\dots,a_2)$, proving (8).  For $0/1$ and $1/0$, the same conclusion follows directly from their two-entry definitions.
\end{proof}

\begin{theorem}\label{continued-fraction-theorem2}
For any irreducible fraction $t\in (0,\infty]$, we have $C_t=CF_{s(t)}$.    
\end{theorem}

\begin{figure}[ht]
    \centering
    \includegraphics[scale=0.3]{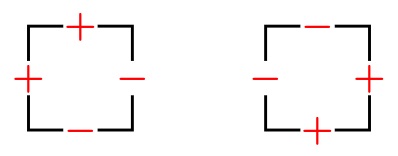}
    \caption{Signed tiles}
    \label{fig:signedtile}
\end{figure}

To prove Theorem \ref{continued-fraction-theorem2}, we introduce the \emph{snake graph}. Assume $a_i\in \mathbb Z_{\geq 1}$. First, recall the snake graph associated with a continued fraction $[a_1,\dots,a_\ell]$ according to \cite{cs18}:

For a continued fraction $[a_1,\dots,a_\ell]$ with $(\ell,a_1)\neq (1,1)$:
\begin{enumerate}\setlength{\leftskip}{-10pt}
    \item Arrange $(a_1 + \cdots + a_\ell)$ signs: the first $a_1$ are $-$, the next $a_2$ are $+$, alternating between $-$ and $+$.
    \item Remove the first and last signs to obtain a tuple $S$ of length $(a_1+\cdots+a_\ell-2)$.
    \item Using $S$, place tiles as in Figure \ref{fig:signedtile}, satisfying:
    \begin{itemize}\setlength{\leftskip}{-10pt}
        \item The first (leftmost) tile is the left one in Figure \ref{fig:signedtile}.
        \item Each new tile is placed to the right or above the preceding tile.
        \item Signs on adjoining edges coincide.
        \item $S$ coincides with the sequence of signs on adjoining parts from left to right.
    \end{itemize}
\end{enumerate}
This graph is the \emph{snake graph}, denoted $\mathcal G[a_1,\dots,a_\ell]$. For the empty continued fraction $[\ ]$, set $\mathcal G[\ ]= \emptyset$. For $[1]$, set $\mathcal G[1]$ as a line segment. $\mathcal G[2]$ is a single tile.

\begin{example}
For $[2,4,2,1]$, the snake graph is shown in Figure \ref{fig:ex-snakegraph}. Indeed, the signs located inside the connected tiles are arranged from the leftmost to the rightmost as follows: there are ($2-1$) consecutive ``$-$" signs, followed by 4 ``$+$" signs, then 2 ``$-$" signs, and finally $(1-1)$ (therefore, no) consecutive ``$+$" signs.
\begin{figure}[ht]
    \centering
    \includegraphics[scale=0.12]{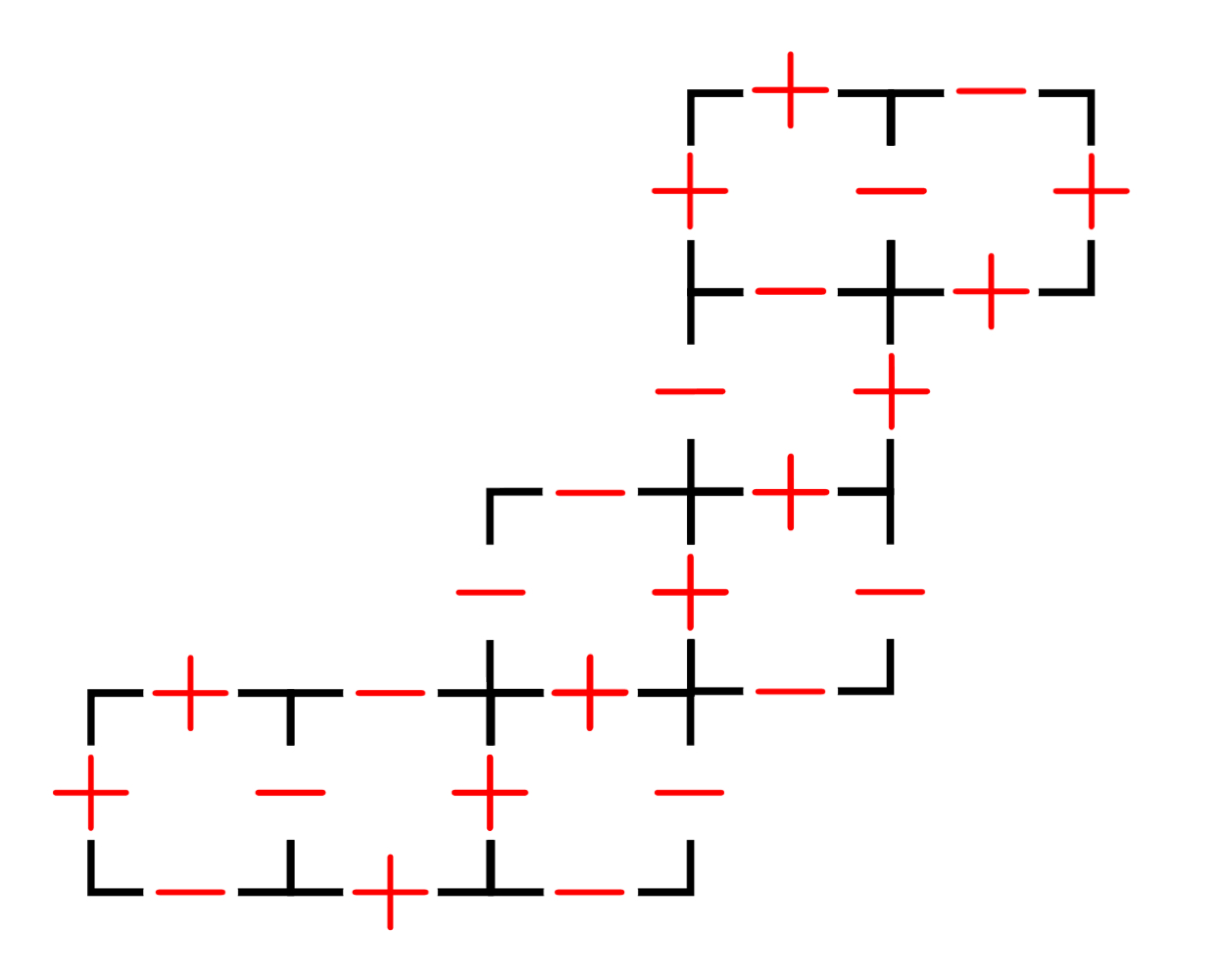}
    \caption{Snake graph associated with $[2,4,2,1]$}
    \label{fig:ex-snakegraph}
\end{figure}
\end{example}

\begin{proposition}
For every continued fraction $[a_1,\dots,a_\ell]$, there is a unique snake graph associated with $[a_1,\dots,a_\ell]$.
\end{proposition}
\begin{proof}
In each tile, the sign on the right edge differs from the sign on the upper edge. Hence the successive attachments of the tiles are uniquely determined.
\end{proof}

Let $G$ be an undirected graph. We recall that a subset $P$ of the edge set of $G$ is called a \emph{perfect matching} of $G$ if each vertex of $G$ is incident to exactly one edge in $P$. We denote by $m(\mathcal G[a_1,\dots,a_\ell])$ the number of perfect matchings of $\mathcal{G} [a_1,\dots,a_\ell]$. We set $m(\mathcal{G}[\ ])=1$.

\begin{example}
By an induction argument, we have $m(\mathcal G [n])=n$. Figure \ref{fig:list-of-perfect-matching} is the list of perfect matchings of $\mathcal{G}[5]$. 
\begin{figure}[ht]
    \centering
\scalebox{3}{\begin{tikzpicture}[baseline=0mm]
\coordinate(1) at (0,0){}; 
\coordinate(2) at (0.2,0){}; 
\coordinate(3) at (0.4,0){}; 
\coordinate(5) at (-0.2,-0.2){}; 
\coordinate(6) at (0,-0.2){}; 
\coordinate(7) at (0.2,-0.2){};
\coordinate(8) at (0.4,-0.2){};
\coordinate(10) at (-0.2,-0.4){}; 
\coordinate(11) at (0,-0.4){}; 
\coordinate(12) at (0.2,-0.4){};
\draw[very thick, red](1) to (2);
\draw(2) to (3);
\draw(5) to (6);
\draw(6) to (7);
\draw(7) to (8);
\draw(10) to (11);
\draw(11) to (12);
\draw[very thick, red](3) to (8);
\draw(2) to (7);
\draw[very thick, red](7) to (12);
\draw(1) to (6);
\draw[very thick, red](6) to (11);
\draw[very thick, red](5) to (10);
\end{tikzpicture}}\hspace{3mm}
\scalebox{3}{\begin{tikzpicture}[baseline=0mm]
\coordinate(1) at (0,0){}; 
\coordinate(2) at (0.2,0){}; 
\coordinate(3) at (0.4,0){}; 
\coordinate(5) at (-0.2,-0.2){}; 
\coordinate(6) at (0,-0.2){}; 
\coordinate(7) at (0.2,-0.2){};
\coordinate(8) at (0.4,-0.2){};
\coordinate(10) at (-0.2,-0.4){}; 
\coordinate(11) at (0,-0.4){}; 
\coordinate(12) at (0.2,-0.4){};
\draw[very thick, red](1) to (2);
\draw(2) to (3);
\draw[very thick, red](5) to (6);
\draw(6) to (7);
\draw(7) to (8);
\draw[very thick, red](10) to (11);
\draw(11) to (12);
\draw[very thick, red](3) to (8);
\draw(2) to (7);
\draw[very thick, red](7) to (12);
\draw(1) to (6);
\draw(6) to (11);
\draw(5) to (10);
\end{tikzpicture}}
\hspace{3mm}
\scalebox{3}{\begin{tikzpicture}[baseline=0mm]
\coordinate(1) at (0,0){}; 
\coordinate(2) at (0.2,0){}; 
\coordinate(3) at (0.4,0){}; 
\coordinate(5) at (-0.2,-0.2){}; 
\coordinate(6) at (0,-0.2){}; 
\coordinate(7) at (0.2,-0.2){};
\coordinate(8) at (0.4,-0.2){};
\coordinate(10) at (-0.2,-0.4){}; 
\coordinate(11) at (0,-0.4){}; 
\coordinate(12) at (0.2,-0.4){};
\draw[very thick, red](1) to (2);
\draw(2) to (3);
\draw(5) to (6);
\draw[very thick, red](6) to (7);
\draw(7) to (8);
\draw(10) to (11);
\draw[very thick, red](11) to (12);
\draw[very thick, red](3) to (8);
\draw(2) to (7);
\draw(7) to (12);
\draw(1) to (6);
\draw(6) to (11);
\draw[very thick, red](5) to (10);
\end{tikzpicture}}
\hspace{3mm}
\scalebox{3}{\begin{tikzpicture}[baseline=0mm]
\coordinate(1) at (0,0){}; 
\coordinate(2) at (0.2,0){}; 
\coordinate(3) at (0.4,0){}; 
\coordinate(5) at (-0.2,-0.2){}; 
\coordinate(6) at (0,-0.2){}; 
\coordinate(7) at (0.2,-0.2){};
\coordinate(8) at (0.4,-0.2){};
\coordinate(10) at (-0.2,-0.4){}; 
\coordinate(11) at (0,-0.4){}; 
\coordinate(12) at (0.2,-0.4){};
\draw(1) to (2);
\draw(2) to (3);
\draw(5) to (6);
\draw(6) to (7);
\draw(7) to (8);
\draw(10) to (11);
\draw[very thick, red](11) to (12);
\draw[very thick, red](3) to (8);
\draw[very thick, red](2) to (7);
\draw(7) to (12);
\draw[very thick, red](1) to (6);
\draw(6) to (11);
\draw[very thick, red](5) to (10);
\end{tikzpicture}}
\hspace{3mm}
\scalebox{3}{\begin{tikzpicture}[baseline=0mm]
\coordinate(1) at (0,0){}; 
\coordinate(2) at (0.2,0){}; 
\coordinate(3) at (0.4,0){}; 
\coordinate(5) at (-0.2,-0.2){}; 
\coordinate(6) at (0,-0.2){}; 
\coordinate(7) at (0.2,-0.2){};
\coordinate(8) at (0.4,-0.2){};
\coordinate(10) at (-0.2,-0.4){}; 
\coordinate(11) at (0,-0.4){}; 
\coordinate(12) at (0.2,-0.4){};
\draw(1) to (2);
\draw[very thick, red](2) to (3);
\draw(5) to (6);
\draw(6) to (7);
\draw[very thick, red](7) to (8);
\draw(10) to (11);
\draw[very thick, red](11) to (12);
\draw(3) to (8);
\draw(2) to (7);
\draw(7) to (12);
\draw[very thick, red](1) to (6);
\draw(6) to (11);
\draw[very thick, red](5) to (10);
\end{tikzpicture}}
\caption{List of perfect matchings of $\mathcal{G}[5]$}\label{fig:list-of-perfect-matching}
\end{figure}
\end{example}

The following proposition follows from the fact that the number of perfect matchings of a snake graph is invariant under congruent transformations of the graph.

\begin{proposition}\label{invariant}
For any positive integer sequence $(a_1,\dots,a_n)$, we have
\[m(\mathcal G[a_1,\dots,a_n])=m(\mathcal G[a_n,\dots,a_1]).\]
\end{proposition}

In \cite{cs18}, the authors establish the following relation between a continued fraction and its associated snake graph:

\begin{theorem}[\cite{cs18}*{Theorem 3.4}]\label{thm:snakegraph-continuedfraction}
The following equality holds:
\[
[a_1,\dots,a_n] = \frac{m(\mathcal G[a_1,\dots,a_n])}{m(\mathcal G[a_2,\dots,a_n])}.
\]
\end{theorem}

Theorem \ref{thm:snakegraph-continuedfraction} includes the case in which the denominator on the right-hand side is $m(\mathcal G[\ ])$.
Throughout this section, blocks are read in the indicated direction of their indices: an ascending block $[a_i,a_{i+1},\dots,a_j]$ is empty when $i>j$, whereas a descending block $[a_i,a_{i-1},\dots,a_j]$ is empty when $i<j$. In either case, the matching number of the empty block is $m(\mathcal G[\ ])=1$. In particular, the descending block $[a_{n-1},\dots,a_4]$ in Theorem \ref{thm:M_t-C_t-combinatorics} is empty when $n=4$, while $[a_n,\dots,a_4]=[a_4]$ in that case.

If $[a_1,\dots,a_n] = \frac{p_n}{q_n}$, then the corresponding product of continued fraction matrices is given by (see \cite{bombieri}*{Section 2.1}):
\[
\begin{bmatrix} a_1 & 1 \\ 1 & 0 \end{bmatrix}
\begin{bmatrix} a_2 & 1 \\ 1 & 0 \end{bmatrix} \cdots
\begin{bmatrix} a_n & 1 \\ 1 & 0 \end{bmatrix}
=
\begin{bmatrix} p_n & p_{n-1} \\ q_n & q_{n-1} \end{bmatrix}.
\]

For $n\geq2$, Theorem \ref{thm:snakegraph-continuedfraction} restates this as
\begin{align}\label{eq:continued-fraction-result}
\begin{bmatrix} a_1 & 1 \\ 1 & 0 \end{bmatrix}
\begin{bmatrix} a_2 & 1 \\ 1 & 0 \end{bmatrix} \cdots
\begin{bmatrix} a_n & 1 \\ 1 & 0 \end{bmatrix}
=
\begin{bmatrix}
m(\mathcal G[a_1,\dots,a_n]) & m(\mathcal G[a_1,\dots,a_{n-1}]) \\
m(\mathcal G[a_2,\dots,a_n]) & m(\mathcal G[a_2,\dots,a_{n-1}])
\end{bmatrix}.
\end{align}

Hence, to prove Theorem \ref{continued-fraction-theorem2}, it suffices to show that for $s(t) = (a_1,\dots,a_n)$,
\begin{align}\label{combinatorics-eq}
C_t =
\begin{bmatrix}
m(\mathcal G[a_1,\dots,a_n]) & m(\mathcal G[a_1,\dots,a_{n-1}]) \\
m(\mathcal G[a_2,\dots,a_n]) & m(\mathcal G[a_2,\dots,a_{n-1}])
\end{bmatrix}
\end{align}
for any $t \in (0,\infty]$. 

On the other hand, \cite{gyoda-maruyama-sato}*{Theorem 7.10} provides an expression for $C_t$ for $t\in(0,1)$. In the notation of the present paper, it gives the following theorem.

\begin{theorem}\label{thm:M_t-C_t-combinatorics}
Let $t\in (0,1)$ be an irreducible fraction. For $s(t) = (a_1,\dots,a_n)$, we have
\begin{align}\label{combinatorics-eq2}
C_t =\scalebox{0.85}{$
\begin{bmatrix}
\begin{gathered}
(3+k_1+k_2+k_3)m(\mathcal G[a_n,\dots,a_3+1])\\
{}-m(\mathcal G[a_n,\dots,a_4])
\end{gathered}
&
\begin{gathered}
(3+k_1+k_2+k_3)m(\mathcal G[a_{n-1},\dots,a_3+1])\\
{}-m(\mathcal G[a_{n-1},\dots,a_4])
\end{gathered}
\\
m(\mathcal{G}[a_n,\dots,a_3+1]) & m(\mathcal{G}[a_{n-1},\dots,a_3+1])
\end{bmatrix}$},
\end{align}
\end{theorem}
\begin{remark}
Note that when comparing Theorem \ref{thm:M_t-C_t-combinatorics} with \cite{gyoda-maruyama-sato}*{Theorem 7.10 (2)}, the positions of the matrix entries differ due to the difference in the convention of Cohn matrices described in Remark \ref{rem:cohn-mannar}.
\end{remark}
Theorem \ref{thm:M_t-C_t-combinatorics} gives the following description of $C_t$ for $t\in (1,\infty)$.
\begin{theorem}\label{thm:M_t-C_t-combinatorics2}
Let $t\in (1,\infty)$ be an irreducible fraction. For $s(t) = (a_1,\dots,a_n)$, we have 
\begin{align}\label{combinatorics-eq3}
C_t =
\begin{bmatrix}
\begin{gathered}
(2+k_1+k_2+k_3)m(\mathcal G[a_n,\dots,a_2])\\
{}+m(\mathcal G[a_n,\dots,a_3])
\end{gathered}
& \ast
\\
m(\mathcal{G}[a_n,\dots,a_2]) & m(\mathcal{G}[a_{n-1},\dots,a_2])
\end{bmatrix}
\end{align}
\end{theorem}
\begin{remark}
The $(1,2)$ entry in \eqref{combinatorics-eq3} is
\[
(2 + k_1 + k_2 + k_3)\, m(\mathcal{G}[a_{n-1},\dots,a_2]) + m(\mathcal{G}[a_{n-1},\dots,a_3]).
\]

\end{remark}
To prove Theorem \ref{thm:M_t-C_t-combinatorics2}, we use two elementary continuant identities.  The first is \cite{banaian-sen}*{Lemma 3}, stated there for numerators of continued fractions; Theorem \ref{thm:snakegraph-continuedfraction} gives the matching-number formulation below.  We prove the second identity directly.

\begin{lemma}[\cite{banaian-sen}*{Lemma 3}]\label{lem:snakegraph-calculation1}
 Let $n\geq 2$. For any positive integer sequence $(a_1,\dots,a_n)$, the following equality holds:
 \[m(\mathcal G[a_1,\dots,a_n])=a_nm(\mathcal G[a_1,\dots,a_{n-1}])+m(\mathcal G[a_1,\dots,a_{n-2}]),\]
 where $m(\mathcal G[\ ])=1$.
\end{lemma}

We call a sequence of the form
$(a_1,\dots,a_{\ell-1},\, a_\ell+k,\, a_\ell,\, a_{\ell-1},\dots,a_1)$
a \emph{semi-palindromic sequence}, where $\ell\geq1$, the integers $a_1,\dots,a_\ell$ and $a_\ell+k$ are positive, and $k\in \mathbb Z$.  When $\ell=1$, the sequence is simply $(a_1+k,a_1)$.  For a sequence $W$ of length at least two, write $W_{\rm L}$ and $W_{\rm R}$ for the sequences obtained by deleting its first and last entry, respectively.

\begin{lemma}\label{lem:snakegraph-calculation2}
Put
\[
S=(a_1,\dots,a_{\ell-1},a_\ell+k,a_\ell,a_{\ell-1},\dots,a_1).
\]
Then
\[
m(\mathcal G[S_{\rm L}])-m(\mathcal G[S_{\rm R}])=(-1)^\ell k.
\]
If the two middle entries are interchanged and
\[
S'=(a_1,\dots,a_{\ell-1},a_\ell,a_\ell+k,a_{\ell-1},\dots,a_1),
\]
then
\[
m(\mathcal G[S'_{\rm L}])-m(\mathcal G[S'_{\rm R}])=(-1)^{\ell+1}k.
\]
For $\ell=1$, these formulas read $m(\mathcal G[a_1])-m(\mathcal G[a_1+k])=-k$ and its negative, respectively.
\end{lemma}
\begin{proof}
Put
\[
A(c)=\begin{bmatrix}c&1\\1&0\end{bmatrix},
\qquad
P=A(a_1)\cdots A(a_{\ell-1}),
\]
where $P$ is the identity matrix if $\ell=1$.  Since every $A(c)$ is symmetric, the continued-fraction matrix of the semi-palindromic sequence is
\[
P A(a_\ell+k)A(a_\ell)P^{\mathsf T}.
\]
Its lower-left entry is the matching number obtained by deleting the first entry of the sequence, and its upper-right entry is the matching number obtained by deleting the last entry; this follows from \eqref{eq:continued-fraction-result}.  Now
\[
A(a_\ell+k)A(a_\ell)
-\bigl(A(a_\ell+k)A(a_\ell)\bigr)^{\mathsf T}
=k\begin{bmatrix}0&1\\-1&0\end{bmatrix}.
\]
For every $2\times2$ matrix $P$,
\[
P\begin{bmatrix}0&1\\-1&0\end{bmatrix}P^{\mathsf T}
=(\det P)\begin{bmatrix}0&1\\-1&0\end{bmatrix},
\]
and here $\det P=(-1)^{\ell-1}$.  Subtracting the upper-right entry from the lower-left entry therefore gives $(-1)^\ell k$, which is the first identity, including the stated interpretation when $\ell=1$.  Interchanging the two central entries replaces $k$ by $-k$ in the skew-symmetric part and gives the second identity.  Equivalently, the latter also follows from Proposition \ref{invariant}.
\end{proof}
\begin{proof}[Proof of Theorem \ref{thm:M_t-C_t-combinatorics2}]
Let $s^\ast(\tfrac{1}{t})=(b_1,\dots,b_\ell)$. Since $b_2=1$ by Proposition \ref{rem:difference-(0,1)(1,infty)} (4), we have
\[
 m(\mathcal{G}[b_\ell,\dots,b_3+1])= m(\mathcal{G}[b_\ell,\dots,b_2]).
\]
By Theorem \ref{thm:Mt-description} and Theorem \ref{thm:M_t-C_t-combinatorics}, it follows that
\[
n^\ast_{\frac1t}=m(\mathcal{G}[b_{\ell},\dots,b_2]), \quad 
u^\ast_{\frac1t}=m(\mathcal{G}[b_{\ell-1},\dots,b_2]).
\]
Moreover, by Proposition \ref{rem:difference-(0,1)(1,infty)} (8), we have $\ell=n$, $b_1=a_1$, and $b_i=a_{n-i+2}$ for $2\leq i\leq n$. Proposition \ref{rem:dual-remark} and Lemma \ref{lem:t-1/t-relation-gen} give
\[n_t=n^\ast_{\frac{1}{t}}, 
\qquad 
u_t=n_t-u^\ast_{\frac{1}{t}}-k_t.
\]
Thus, we have
\begin{align}
n_t&=m(\mathcal{G}[a_{2},\dots,a_n]),\label{eq:n_t}\\
u_t&=m(\mathcal{G}[a_{2},\dots,a_n])-m(\mathcal{G}[a_{3},\dots,a_{n}])-k_t.\label{eq:u_t}
\end{align}
Since $a_n=1$ by Proposition \ref{rem:difference-(0,1)(1,infty)} (4), we obtain
\begin{align}
m(\mathcal{G}[a_{3},\dots,a_{n}])=m(\mathcal{G}[a_{3},\dots,a_{n-1}+1]).\label{eq:a_n=a_n-1+1}
\end{align}
Furthermore, since $(a_2,\dots,a_{n-1}+1)$ is a semi-palindromic sequence by Proposition \ref{rem:difference-(0,1)(1,infty)} (5) and (7), we obtain
\begin{align}\label{eq:m=m+k}
    m(\mathcal{G}[a_{2},\dots,a_{n-2}])=m(\mathcal{G}[a_{3},\dots,a_{n-1}+1])+k_t\stackrel{\eqref{eq:a_n=a_n-1+1}}{=}m(\mathcal{G}[a_{3},\dots,a_{n}])+k_t,
\end{align}
from Lemma \ref{lem:snakegraph-calculation2}. Therefore, substituting \eqref{eq:m=m+k} into \eqref{eq:u_t}, we have
\[
u_t=m(\mathcal{G}[a_{2},\dots,a_n])-m(\mathcal{G}[a_{2},\dots,a_{n-2}]).
\]

Since $a_n=1$ again, Lemma \ref{lem:snakegraph-calculation1} yields
\[
m(\mathcal{G}[a_{2},\dots,a_n])=m(\mathcal{G}[a_2,\dots,a_{n-1}])+m(\mathcal{G}[a_{2},\dots,a_{n-2}]).
\]
Hence, we conclude that
\[
u_t=m(\mathcal{G}[a_{n-1},\dots,a_2]).
\]
By Theorem \ref{thm:Mt-description}, it follows that
\[
C_t=\begin{bmatrix}
    \ast&\ast\\ 
    m(\mathcal{G}[a_n,\dots,a_2]) &m(\mathcal{G}[a_{n-1},\dots,a_2]) 
\end{bmatrix}.
\]
Since $\mathrm{tr}(C_t)=(3+k_1+k_2+k_3)n_t-k_t$, using \eqref{eq:m=m+k}, the $(1,1)$-entry is
\[
(2+k_1+k_2+k_3)\, m(\mathcal{G}[a_n,\dots,a_2])+m(\mathcal{G}[a_n,\dots,a_3]).
\]\end{proof}

The following relation between snake graphs will also be used.

\begin{lemma}\label{lem:m=m-m}
Let $n\geq 3$. For any positive integer sequence $(a_1,1,a_3,\dots,a_n)$, we have
\[
m(\mathcal G[a_1,1,a_3,\dots,a_n]) = (a_1+1) m(\mathcal{G}[a_3+1,\dots,a_n]) - m(\mathcal{G}[a_4,\dots,a_{n}]).
\]
When $n=3$, the last term means $m(\mathcal G[\ ])=1$.
\end{lemma}
\begin{proof}
Write $K(c_1,\dots,c_r)=m(\mathcal G[c_1,\dots,c_r])$ and $K(\,)=1$.  The continuant recurrence, equivalently Lemma \ref{lem:snakegraph-calculation1} together with Proposition \ref{invariant}, gives
\[
K(c_1,\dots,c_r)=c_1K(c_2,\dots,c_r)+K(c_3,\dots,c_r).
\]
It also gives
\[
K(1,a_3,\dots,a_n)=K(a_3+1,a_4,\dots,a_n)
\]
and
\[
K(a_3+1,a_4,\dots,a_n)
=K(a_3,\dots,a_n)+K(a_4,\dots,a_n).
\]
Consequently,
\begin{align*}
K(a_1,1,a_3,\dots,a_n)
&=a_1K(1,a_3,\dots,a_n)+K(a_3,\dots,a_n)\\
&=(a_1+1)K(a_3+1,a_4,\dots,a_n)-K(a_4,\dots,a_n),
\end{align*}
as claimed.  The same calculation applies for $n=3$ with the stated empty-sequence convention.
\end{proof}

\begin{proof}[Proof of Theorem \ref{continued-fraction-theorem2}]
The cases $t = \frac{1}{1}$ and $t=\frac{1}{0}$ can be checked directly. We focus on $t\in (0,1)$. Comparing \eqref{combinatorics-eq} and \eqref{combinatorics-eq2}, it suffices to show
\begin{align}
m(\mathcal G[a_1,\dots,a_n]) &= (3+k_1+k_2+k_3) m(\mathcal{G}[a_n,\dots,a_3+1]) - m(\mathcal{G}[a_{n},\dots,a_4]), \label{eq:1}\\
m(\mathcal G[a_1,\dots,a_{n-1}]) &= (3+k_1+k_2+k_3) m(\mathcal{G}[a_{n-1},\dots,a_3+1]) - m(\mathcal{G}[a_{n-1},\dots,a_4]), \label{eq:2}\\
m(\mathcal G[a_2,\dots,a_n]) &= m(\mathcal{G}[a_n,\dots,a_3+1]), \label{eq:3} \\
m(\mathcal G[a_2,\dots,a_{n-1}]) &= m(\mathcal{G}[a_{n-1},\dots,a_3+1]). \label{eq:4}
\end{align}

Since $a_2=1$, we have 
\[
m(\mathcal{G}[a_n,\dots,a_3+1]) = m(\mathcal{G}[a_n,\dots,a_2]), \quad
m(\mathcal{G}[a_{n-1},\dots,a_3+1]) = m(\mathcal{G}[a_{n-1},\dots,a_2]).
\]
Since $\mathcal{G}[a_n,\dots,a_2]$ and $\mathcal{G}[a_2,\dots,a_n]$ are mirror images of each other as graphs, the values $m(\mathcal{G}[a_n,\dots,a_2])$ and $m(\mathcal{G}[a_2,\dots,a_n])$ are equal.
Therefore, \eqref{eq:3} and \eqref{eq:4} hold. 

By Proposition \ref{rem:difference-(0,1)(1,infty)} (2), Equation \eqref{eq:1} can be rewritten as
\[
m(\mathcal G[2+k_1+k_2+k_3,1,a_3,\dots,a_n]) = (3+k_1+k_2+k_3) m(\mathcal{G}[a_n,\dots,a_3+1]) - m(\mathcal{G}[a_{n},\dots,a_4]).
\]
This is the special case $a_1=2+k_1+k_2+k_3$ of Lemma \ref{lem:m=m-m}. The proof of \eqref{eq:2} is identical; when $t=\frac12$, it uses the $n=3$ case of that lemma for the truncated sequence $(a_1,1,a_3)$.
It remains to prove the statement for $t\in(1,\infty)$. By Proposition \ref{rem:difference-(0,1)(1,infty)} (1), the number of entries of $s(t)$ is even, and hence $\det (CF_{s(t)})=1$. Since $\det C_t=1$ by Theorem \ref{thm:Mt-description}, it suffices to show the equality of the $(1,1)$-, $(2,1)$-, and $(2,2)$-entries. Comparing \eqref{combinatorics-eq} and \eqref{combinatorics-eq3}, it suffices to show
\begin{align}
m(\mathcal G[a_1,\dots,a_n]) &= (2+k_1+k_2+k_3) m(\mathcal{G}[a_n,\dots,a_2]) + m(\mathcal{G}[a_{n},\dots,a_3]), \label{eq:5}\\
m(\mathcal G[a_2,\dots,a_n]) &= m(\mathcal{G}[a_n,\dots,a_2]), \label{eq:6} \\
m(\mathcal G[a_2,\dots,a_{n-1}]) &= m(\mathcal{G}[a_{n-1},\dots,a_2]). \label{eq:7}
\end{align}
Equation \eqref{eq:5} follows from Lemma \ref{lem:snakegraph-calculation1}, while \eqref{eq:6} and \eqref{eq:7} follow from the reversal invariance in Proposition \ref{invariant}.
\end{proof}
\begin{remark}
By definition, $C_{\frac{0}{1}}\neq CF_{s(\frac{0}{1})}$. Nevertheless, Theorem \ref{thm:markov-value-gen} includes $t=\frac{0}{1}$; this case is verified directly.
\end{remark}
\section{Generalized Markov Length and Generalized Markov Distance}\label{section:Generalized Markov length and generalized Markov distance}
This section introduces the key tool for generalizing the discrete Markov spectrum, namely Theorem \ref{thm:distance-theorem}. We retain the fixed parameters $(k_1,k_2,k_3,\sigma)$ and the unit-edge set $\mathcal E$ from the preceding section.

\begin{definition}\label{def:generalized-arc}
Let $\gamma$ be an oriented curve segment in $\widetilde{\mathbb R^2}$ satisfying the standing finiteness and transversality assumptions of Section~\ref{section:Classical Markov Spectrum}. It is called a \emph{generalized arc} if
\[
\gamma(0),\gamma(1)\in\mathbb Z^2,
\qquad
\gamma((0,1))\cap
\left(\mathbb Z^2\cup\{m_e\mid e\in\mathcal E\}\right)=\varnothing,
\]
and if no two consecutive edge-crossing occurrences lie on the same unit edge.
\end{definition}

Self-intersections are permitted. The requirement that $\gamma((0,1))$ avoid every midpoint $m_e$ can always be achieved by an arbitrarily small generic perturbation.

To define the generalized Markov length (GM length for short), we introduce a sign rule for generalized arcs.
\begin{definition}
Let $\gamma$ be an oriented generalized arc. Assign one freely chosen sign $-$ or $+$ to each of its first and last triangle-passage occurrences. If these occurrences coincide, assign only one sign. This is the \emph{endpoint rule} for $\gamma$ (see Figure~\ref{fig:minus-righttriangles-endpoint}).
\begin{figure}[ht]
    \centering
    \begin{tikzpicture}[baseline=0mm]
        \draw (0,0) -- (1,0) -- (0,1) -- cycle;
        \draw[red,thick] (0,0) -- (0.5,0.5);
        \node[rotate=225,red] at (0.33,0.33) {$\mathbf <$};
    \end{tikzpicture}\hspace{-0.3cm}
    \rotatebox{180}{\begin{tikzpicture}[baseline=10mm]
        \draw (0,0) -- (1,0) -- (0,1) -- cycle;
        \draw[red,thick] (0,0) -- (0.5,0.5);
        \node[rotate=225,red] at (0.22,0.22) {$\mathbf >$};
    \end{tikzpicture}}
    \hspace{0.3cm}
    \begin{tikzpicture}[baseline=0mm]
        \draw (0,0) -- (1,0) -- (0,1) -- cycle;
        \draw[red,thick] (0,1) -- (0.5,0);
        \node[rotate=300,red] at (0.38,0.22) {$\mathbf >$};
    \end{tikzpicture}\hspace{-0.3cm}
    \rotatebox{180}{\begin{tikzpicture}[baseline=10mm]
        \draw (0,0) -- (1,0) -- (0,1) -- cycle;
        \draw[red,thick] (0,1) -- (0.5,0);
        \node[rotate=300,red] at (0.22,0.56) {$\mathbf <$};
    \end{tikzpicture}}
    \hspace{0.3cm}
    \begin{tikzpicture}[baseline=0mm]
        \draw (0,0) -- (1,0) -- (0,1) -- cycle;
        \draw[red,thick] (0,0.5) -- (1,0);
        \node[rotate=150,red] at (0.33,0.33) {$\mathbf >$};
    \end{tikzpicture}\hspace{-0.3cm}
    \rotatebox{180}{\begin{tikzpicture}[baseline=10mm]
        \draw (0,0) -- (1,0) -- (0,1) -- cycle;
        \draw[red,thick] (0,0.5) -- (1,0);
        \node[rotate=150,red] at (0.45,0.28) {$\mathbf <$};
    \end{tikzpicture}}
    \hspace{0.3cm}
    \begin{tikzpicture}[baseline=0mm]
        \draw (0,0) -- (1,0) -- (0,1) -- cycle;
        \draw[red,thick] (0,0) -- (0.5,0.5);
        \node[rotate=225,red] at (0.22,0.22) {$\mathbf >$};
    \end{tikzpicture}\hspace{-0.3cm}
    \rotatebox{180}{\begin{tikzpicture}[baseline=10mm]
        \draw (0,0) -- (1,0) -- (0,1) -- cycle;
        \draw[red,thick] (0,0) -- (0.5,0.5);
        \node[rotate=225,red] at (0.33,0.33) {$\mathbf <$};
    \end{tikzpicture}}
    \hspace{0.3cm}
    \begin{tikzpicture}[baseline=0mm]
        \draw (0,0) -- (1,0) -- (0,1) -- cycle;
        \draw[red,thick] (0,1) -- (0.5,0);
        \node[rotate=300,red] at (0.24,0.52) {$\mathbf <$};
    \end{tikzpicture}\hspace{-0.3cm}
    \rotatebox{180}{\begin{tikzpicture}[baseline=10mm]
        \draw (0,0) -- (1,0) -- (0,1) -- cycle;
        \draw[red,thick] (0,1) -- (0.5,0);
        \node[rotate=300,red] at (0.33,0.33) {$\mathbf >$};
    \end{tikzpicture}}
    \hspace{0.3cm}
    \begin{tikzpicture}[baseline=0mm]
        \draw (0,0) -- (1,0) -- (0,1) -- cycle;
        \draw[red,thick] (0,0.5) -- (1,0);
        \node[rotate=150,red] at (0.56,0.22) {$\mathbf <$};
    \end{tikzpicture}\hspace{-0.3cm}
    \rotatebox{180}{\begin{tikzpicture}[baseline=10mm]
        \draw (0,0) -- (1,0) -- (0,1) -- cycle;
        \draw[red,thick] (0,0.5) -- (1,0);
        \node[rotate=150,red] at (0.22,0.38) {$\mathbf >$};
    \end{tikzpicture}}
    \caption{Right-angled triangles with $-$ or $+$ signs assigned according to the endpoint rule.}
    \label{fig:minus-righttriangles-endpoint}
\end{figure}
\end{definition}
\begin{proposition}
The signs assigned by the endpoint rule do not affect the shape of the snake graph constructed from the entire sign sequence. Consequently, either $+$ or $-$ may be chosen in the endpoint rule.
\end{proposition}
\begin{proof}
In the snake-graph construction, the endpoint signs occur only at the two ends of the sign sequence and are not used to determine an attachment between two consecutive tiles. Changing either endpoint sign therefore leaves the shape of the snake graph unchanged.
\end{proof}

Assume that $\gamma$ is a generalized arc. Apply the endpoint rule to the first and last triangle-passage occurrences, the triangle-crossing rule to every other triangle-passage occurrence, and the edge-crossing rule to every edge-crossing occurrence. Read the resulting signs in parameter order along $\gamma$. Every occurrence is recorded separately, even if the same geometric triangle or edge occurs more than once. Denote the resulting finite word by $\varepsilon(\gamma)$. The sequence of lengths of its runs is called the \emph{sign sequence associated with $\gamma$} and is denoted by $s(\gamma)$.

If $s(\gamma)=(a_1,\dots,a_n)$, then the number of perfect matchings $m(\mathcal G[a_1,\dots,a_n])$ is called the \emph{$(k_1,k_2,k_3,\sigma)$-generalized Markov (GM) length of $\gamma$}, denoted by $|\gamma|$.

\begin{proposition}
The GM length of a generalized arc $\gamma$ is independent of its orientation.
\end{proposition}
\begin{proof}
After choosing the endpoint signs compatibly, reversing the orientation reverses the sign sequence. Proposition~\ref{invariant} therefore shows that the two GM lengths are equal.
\end{proof}

\begin{example}
Set $(k_1,k_2,k_3,\sigma)=(1,2,0,\mathrm{id})$. Let $\gamma$ be the red oriented generalized arc shown in Figure \ref{fig:gen-Markov-length}. Then
\[
s(\gamma)=(3,7,1,6,2,1,3,1,2,2,6,5),
\]
and 
\[
|\gamma|=551409.
\]
\begin{figure}[ht]
    \centering
    \includegraphics[scale=0.06]{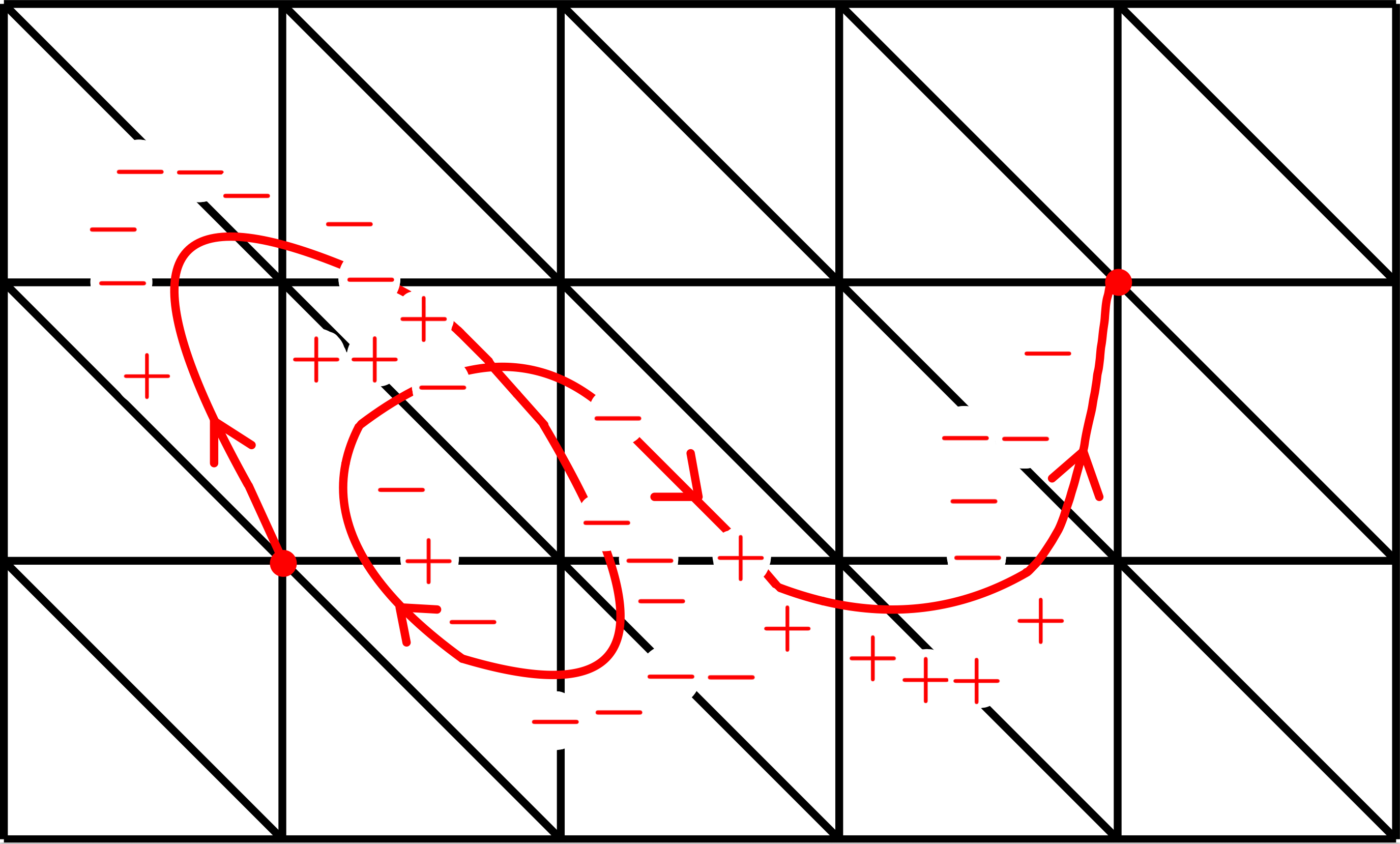} 
    \caption{Generalized arc $\gamma$ and sign assignments.}
    \label{fig:gen-Markov-length}
\end{figure}
\end{example}

For lattice points $A,B$, let $\Gamma(A,B)$ be the set of generalized arcs oriented from $A$ to $B$, that is, with $\gamma(0)=A$ and $\gamma(1)=B$. We define
\[
d(A,B):=\begin{cases}
\min\{|\gamma|\mid\gamma\in\Gamma(A,B)\},&\text{if $A\neq B$},\\
0,&\text{if $A=B$.}
\end{cases}
\]
which is called the \emph{$(k_1,k_2,k_3,\sigma)$-generalized Markov (GM) distance} between $A$ and $B$.
When $A\neq B$, the set $\Gamma(A,B)$ is nonempty by the push-offs defined below, and since its GM lengths are positive integers, the minimum above exists.

Assume $A\neq B$ and orient the segment $AB$ from $A$ to $B$. Keeping the endpoints fixed, take sufficiently small push-offs of $AB$ to its left and right, chosen to be generalized arcs. Denote them by $\gamma^L_{AB}$ and $\gamma^R_{AB}$, respectively. These are the left and right infinitesimal deformations of \cite{banaian}*{Convention~28}.

The following theorem identifies the arcs attaining the minimum that defines $d(A,B)$.

\begin{theorem}\label{thm:distance-theorem}
For any distinct lattice points $A,B$, one has
\[
d(A,B)=|\gamma^R_{AB}|=|\gamma^L_{AB}|.
\]
Consequently, every generalized arc $\gamma$ connecting $A$ and $B$ satisfies $|\gamma|\geq d(A,B)$.
\end{theorem}

\begin{proof}
Our sign rules are the edge-typed generalization of the rules defining the $k$-Markov length in \cite{banaian}. The proof of \cite{banaian}*{Lemma 46 and Theorem 48} is formulated in terms of edge-crossing records and uses the condition that the same unit edge does not occur twice consecutively, as required in Definition~\ref{def:generalized-arc}. By assigning $k_{\sigma(i)}$ to a crossing of an edge of type $i$, that proof applies unchanged to general $(k_1,k_2,k_3)$. Hence every generalized arc from $A$ to $B$ has GM length at least that of either straight push-off. Requiring the interiors of generalized arcs to avoid all edge midpoints does not change this application: every edge crossing is then nondegenerate, and both straight push-offs satisfy this requirement. Consequently,
\[
d(A,B)=|\gamma^R_{AB}|=|\gamma^L_{AB}|.
\]
In the case $k_1=k_2=k_3=0$, see also \cite{llrs}*{Theorem 3.5}.
\end{proof}

\begin{example}
Set $(k_1,k_2,k_3,\sigma)=(1,2,0,\mathrm {id})$. 
By Theorem~\ref{thm:distance-theorem}, the following push-offs attain the GM distance.
\begin{itemize}\setlength{\leftskip}{-10pt}
\item [(1)] 
Let $A=(0,0)$ and $B=(3,2)$. Then $\gamma^L_{AB}$ is the red curve segment shown schematically in Figure \ref{fig:generalized-Markov-distance-2/3}; its infinitesimal leftward displacement at the midpoint is suppressed in the drawing. We have
\[
s(\gamma^L_{AB})=(4,4,5,4)
\]
and
\[
d(A,B)=373.
\]
\begin{figure}[ht]
    \centering
    \includegraphics[scale=0.12]{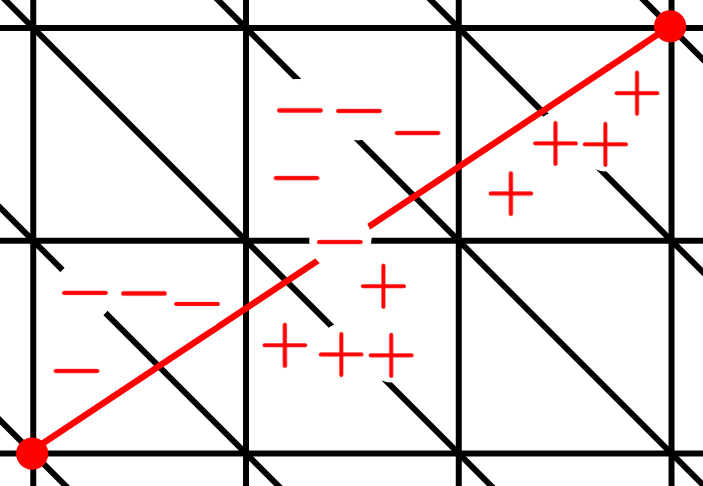} 
    \caption{$(1,2,0,\mathrm {id})$-GM distance between $(0,0)$ and $(3,2)$.}
    \label{fig:generalized-Markov-distance-2/3}
\end{figure}
\item [(2)] 
Let $A=(0,0)$ and $B=(6,4)$. Then $\gamma^L_{AB}$ is the red curve segment shown in Figure \ref{fig:generalized-Markov-distance}. We have
\[
s(\gamma^L_{AB})=(4,5,4,4,5,1,3,5,4,4)
\]
and
\[
d(A,B)=834401.
\]
Note that the intersection point of the horizontal edge in $\widetilde{\mathbb R^2}$ and $\gamma^L_{AB}$ is slightly to the left of the midpoint of the horizontal edge, as in Example~\ref{ex:t=2/5}.
\begin{figure}[ht]
    \centering
    \includegraphics[scale=0.08]{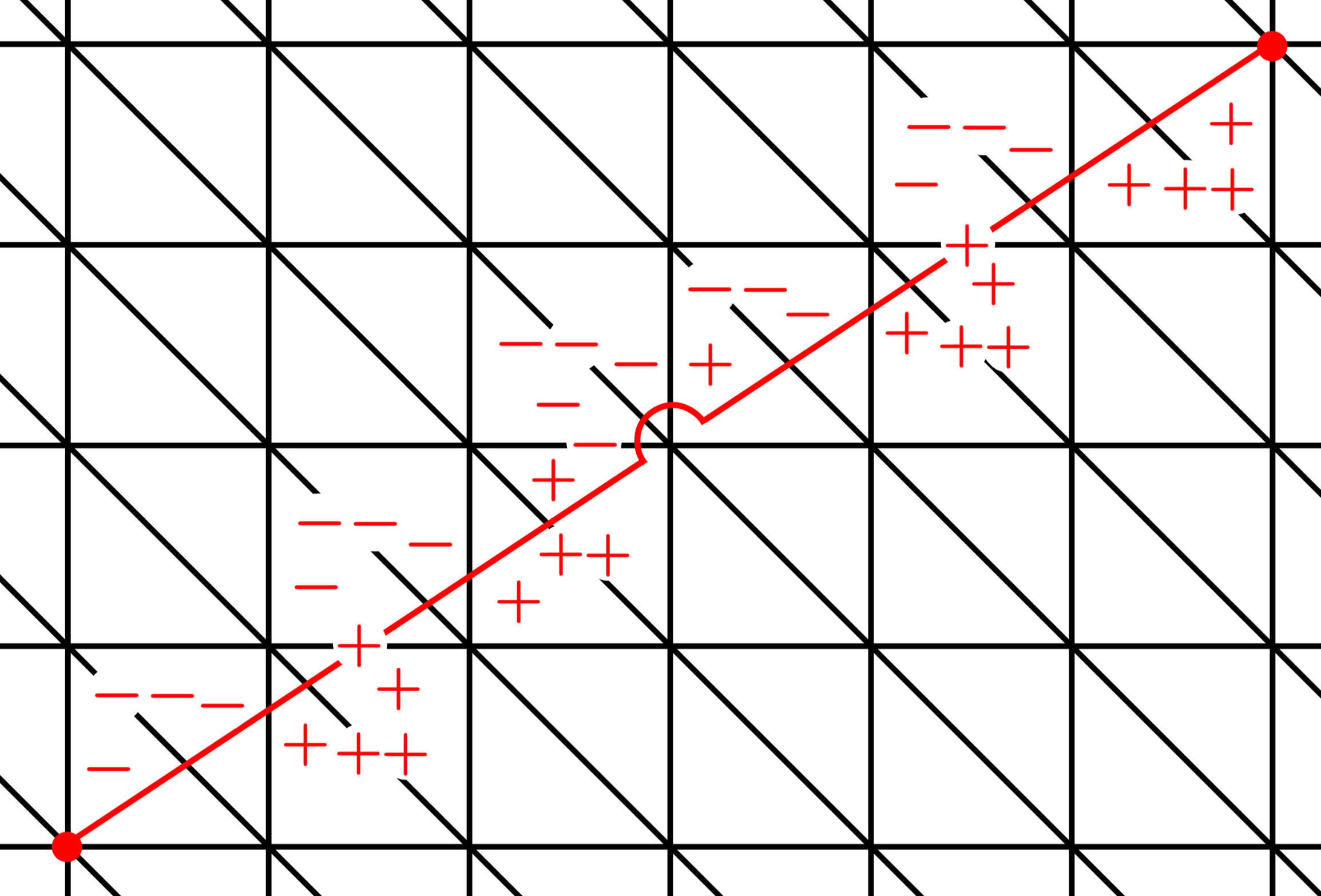} 
    \caption{$(1,2,0,\mathrm {id})$-GM distance between $(0,0)$ and $(6,4)$.}
    \label{fig:generalized-Markov-distance}
\end{figure}
\end{itemize}
\end{example}

\section{Generalized Discrete Markov Spectra}\label{section:Generalized discrete Markov spectra}
This section generalizes Theorem \ref{thm:markov-value} to GM numbers and establishes several properties of the resulting spectra. 
\subsection{Main theorems and their Proofs}
The main theorem of this paper is the following.
\begin{theorem}\label{thm:markov-value-gen}
Fix $(k_1,k_2,k_3)\in \mathbb Z_{\geq 0}^3$ and $\sigma\in \mathfrak S_3$.  
For any irreducible fraction $t\in [0,\infty]$, let $(n_t,i_t)$ denote the corresponding $(k_1,k_2,k_3,\sigma)$-GM number-position pair, and let $s(t)$ be the corresponding generalized strongly admissible sequence. 
Then we have
\[
L({}^\infty s(t)^\infty)=\ell(^\infty s(t)|s(t)^\infty)=\ell(^\infty s^\ast(\tfrac{1}{t})|s^\ast(\tfrac{1}{t})^\infty)=L({}^\infty s^\ast(\tfrac1t)^\infty),
\]
and this value is given by
\[\frac{\sqrt{((3+k_1+k_2+k_3)n_t-k_t)^2-4}}{n_t}.
\]
%Moreover, for any cut $P^\ast|Q$ of $^\infty s(t)^\infty$ except for $^\infty s(t)s(t)^\infty$ (and $^\infty s^\ast(\tfrac{1}{t})s^\ast(\tfrac{1}{t})^\infty$ if $^\infty s(t)^\infty=^\infty s^\ast(\tfrac{1}{t})^\infty$), we have
%\[L({}^\infty s(t)^\infty)>\ell(P^\ast|Q).\]
\end{theorem}

The proof uses the following two lemmas.

\begin{lemma}[See \cite{bombieri2}*{Lemma 17}]\label{lem:trace}
For any finite positive integer sequence $S=(a_1,\dots,a_n)$, we have
\[
\ell(^\infty S|S^\infty)=\frac{\sqrt{(\mathrm{tr}(CF_S))^2-(-1)^n4}}{(CF_{S})_{21}},
\]
where $(CF_S)_{21}$ denotes the $(2,1)$-entry of $CF_S$.
\end{lemma}

For an irreducible fraction $t=p/q\in\mathbb Q_{>0}\setminus\{1\}$, write $s(t)=(a_1,\ldots,a_n)$. By projecting $\widetilde{\mathbb R^2}$ to a triangulation of the once-punctured torus, $\overline{L_t}$ goes to a loop on the torus (see also Remark \ref{rem:one-punctured-torus}). Viewing $\overline{L_t}$ as a loop on this torus, for each $w=(a_{k+1},\dots,a_{k+n-1})$, let $\overline{L_t}(w)$ denote the segment of $\overline{L_t}$ producing the corresponding consecutive runs. Attach endpoint branches as in Table~\ref{table1}, and denote the resulting curve segment by $\widetilde{L_t}(w)$. These endpoint modifications satisfy the following lemma.
\begin{table}[ht]
\centering
\begin{tabular}{|c|c|}
Endpoint of $\overline{L_t}(w)$& $\widetilde{L_t}(w)$-Modification\\
\hline
&\\[-4mm]
\includegraphics[scale=0.04]{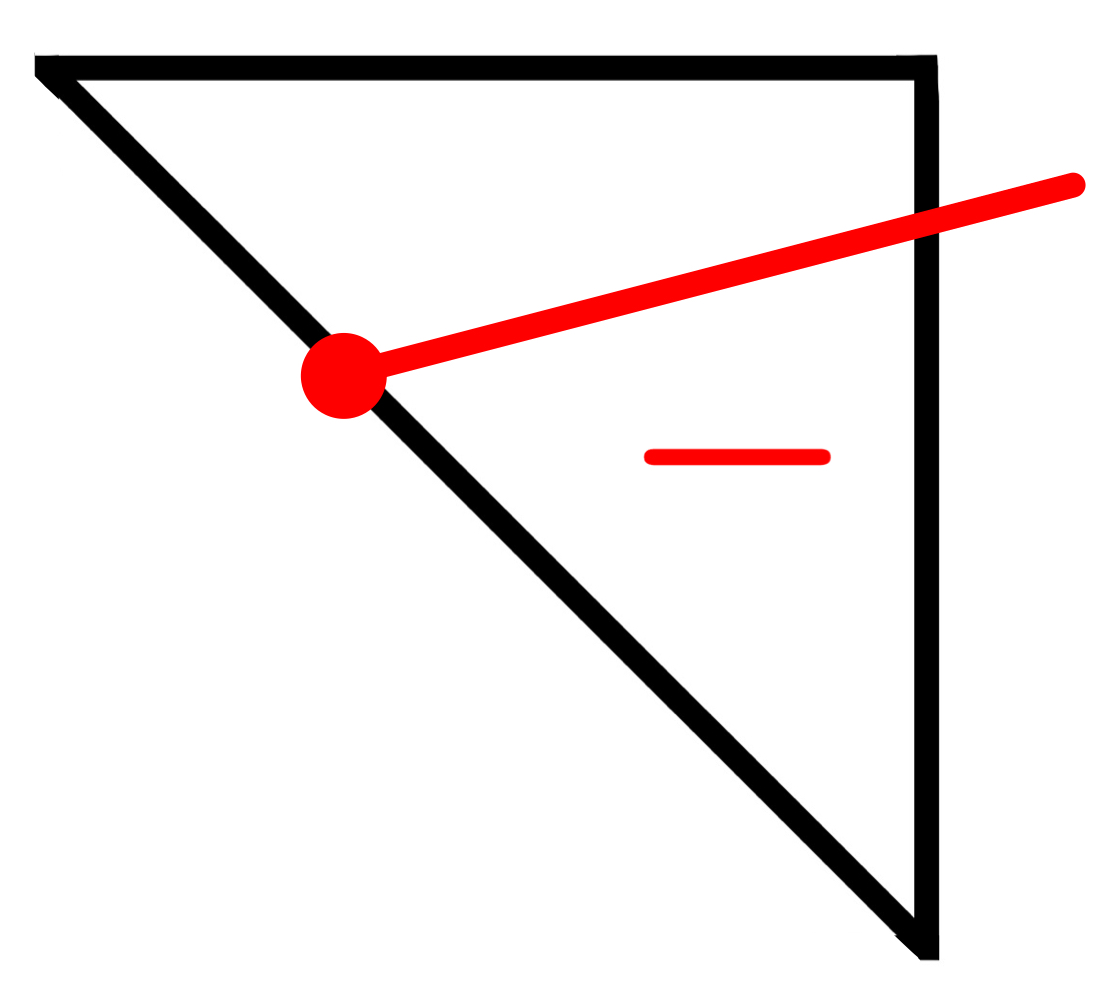} &\includegraphics[scale=0.04]{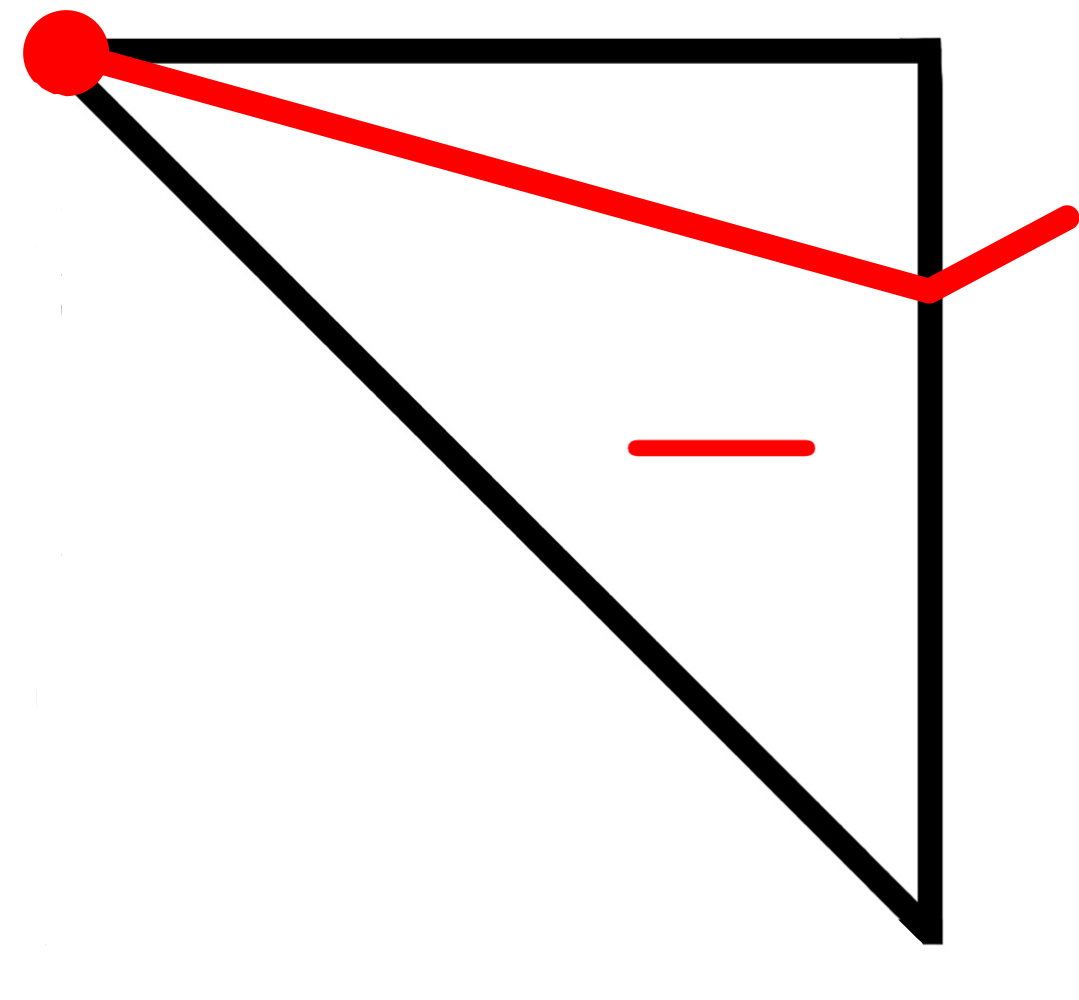}\\
\hline
&\\[-4mm]
\includegraphics[scale=0.04]{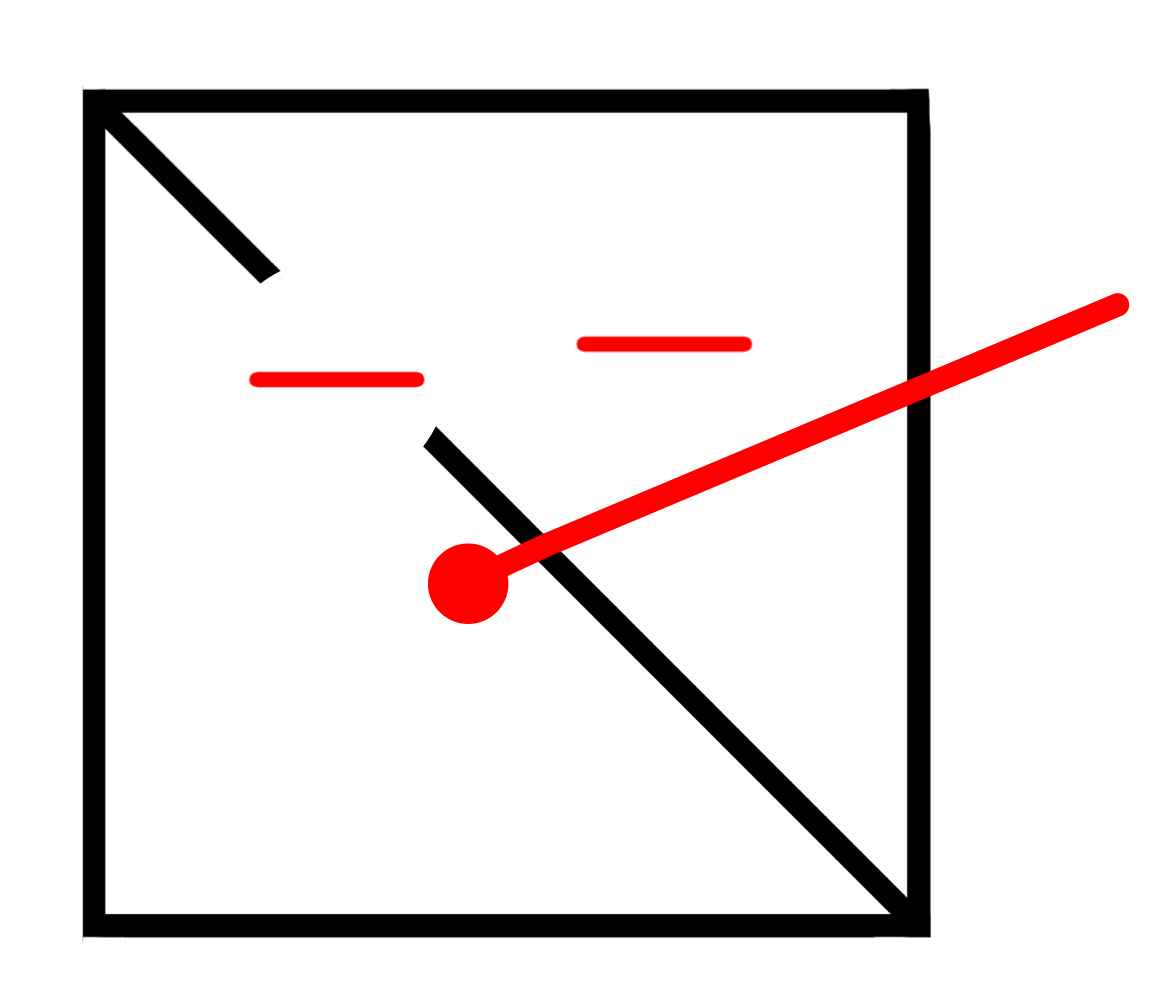} &\includegraphics[scale=0.04]{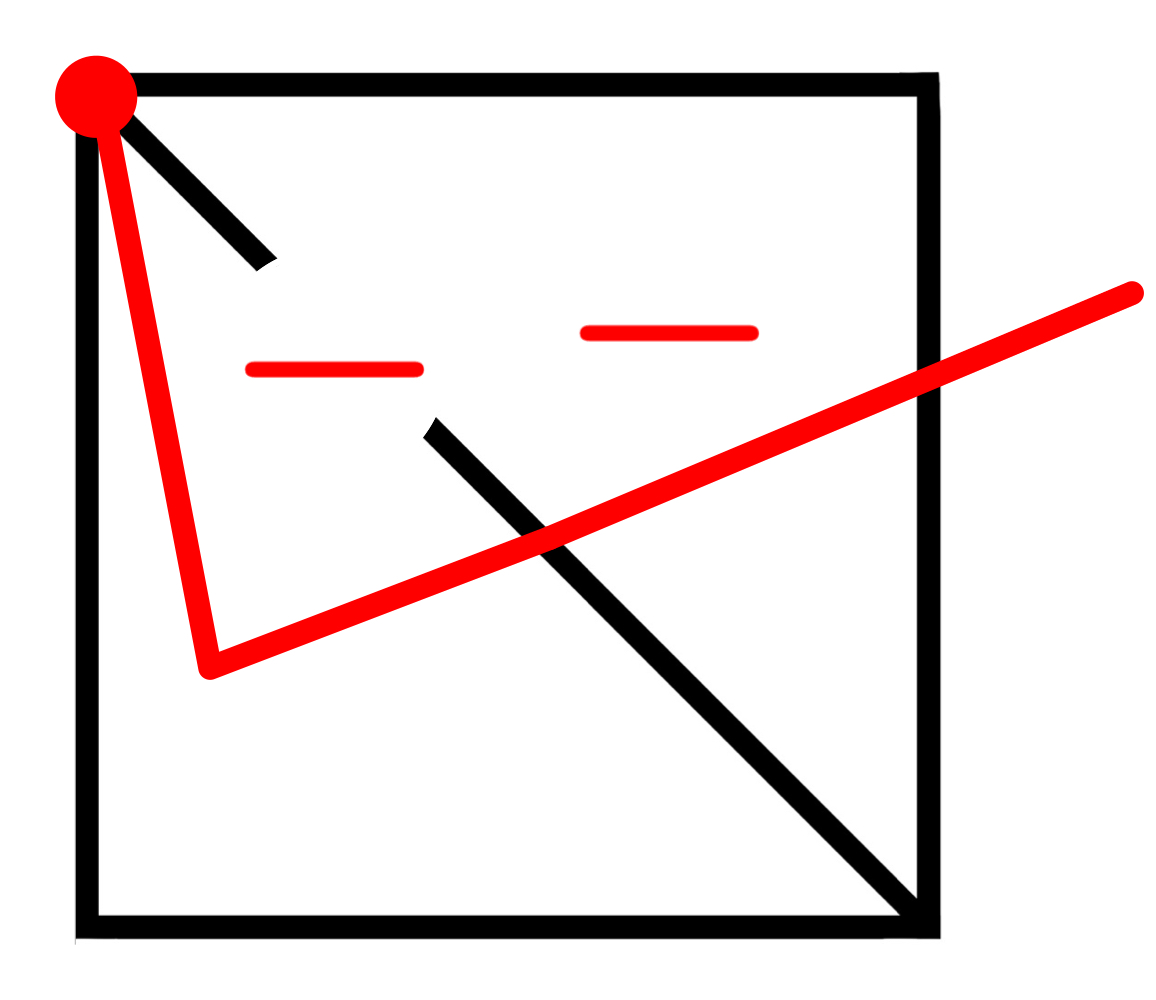}\\
\hline
&\\[-4mm]
\includegraphics[scale=0.04]{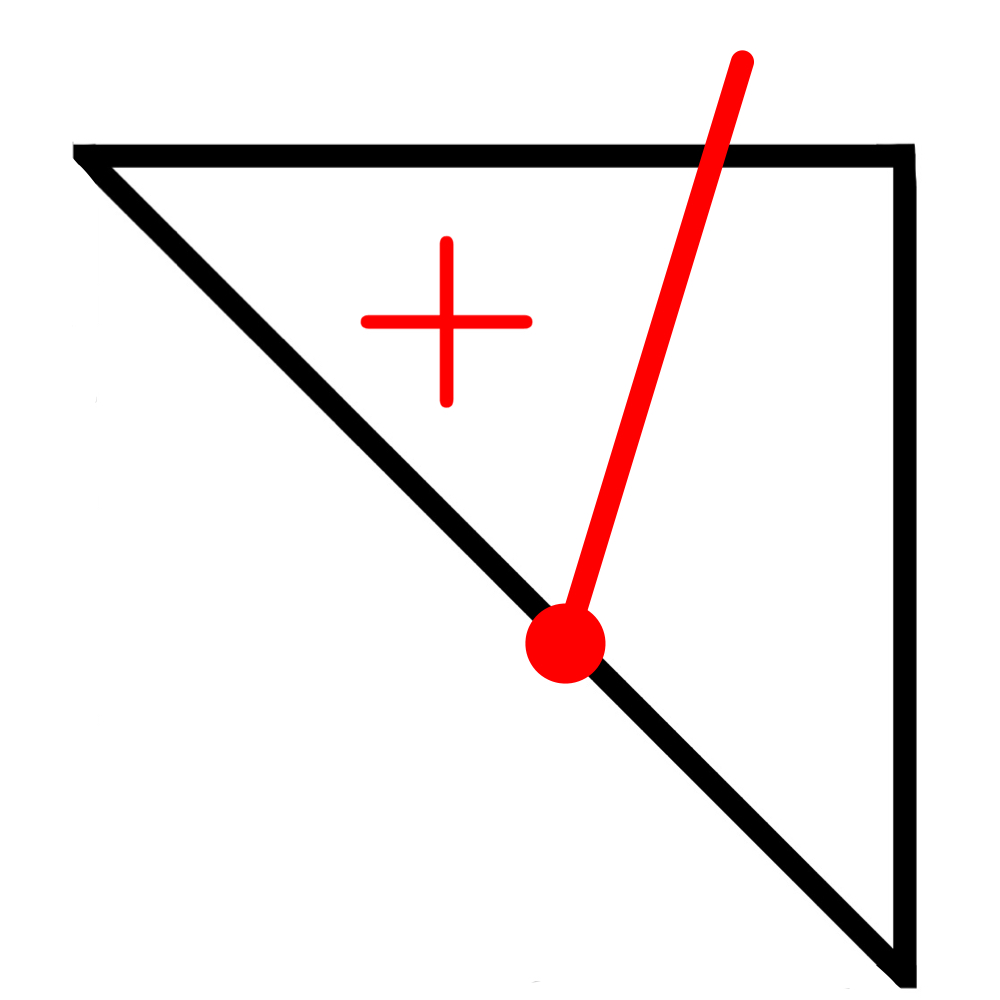} &\includegraphics[scale=0.04]{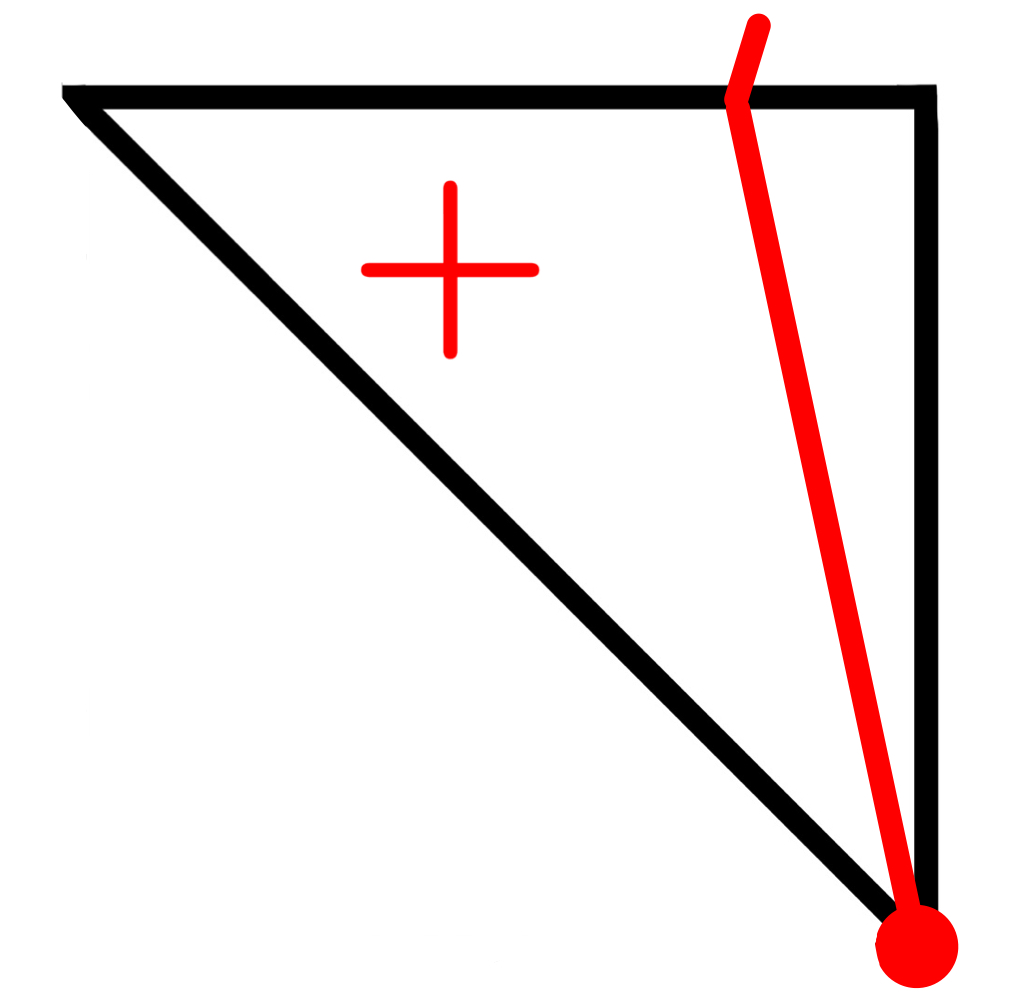}\\
\hline
&\\[-4mm]
\includegraphics[scale=0.04]{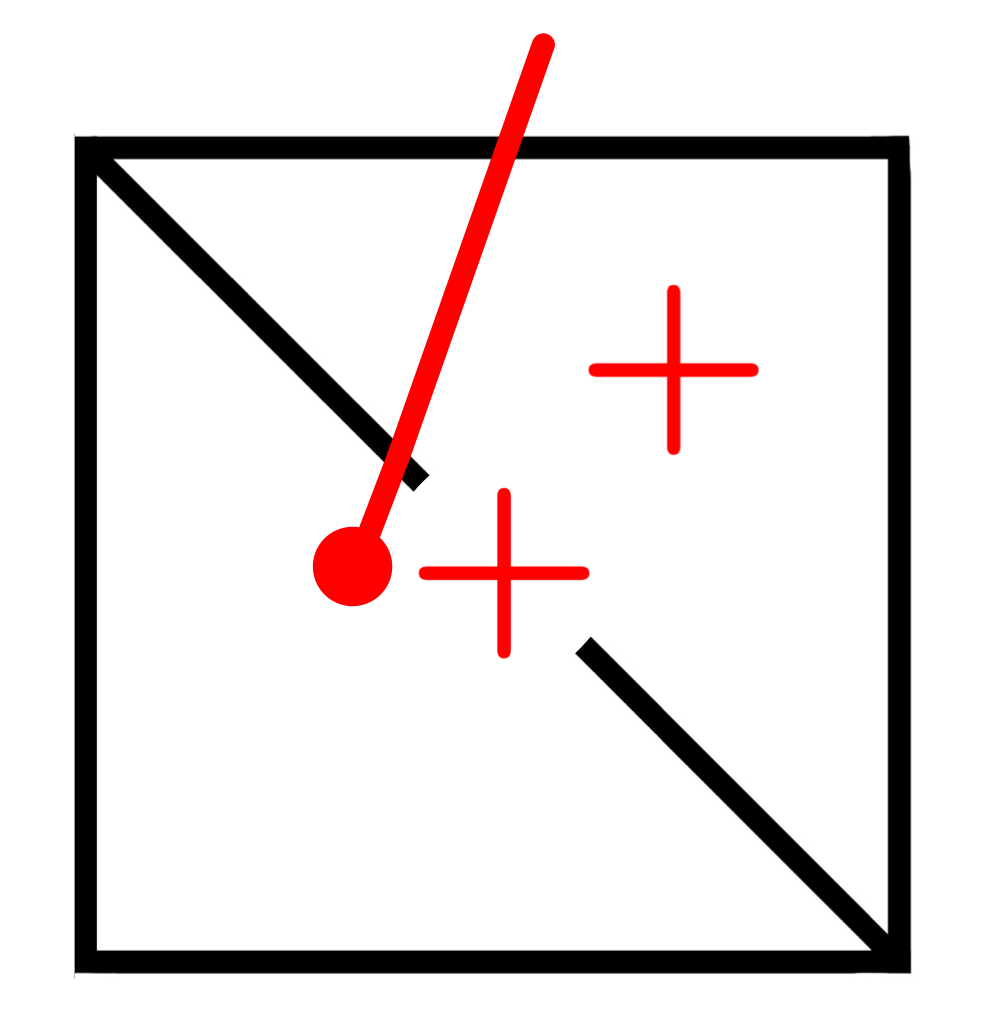} &\includegraphics[scale=0.04]{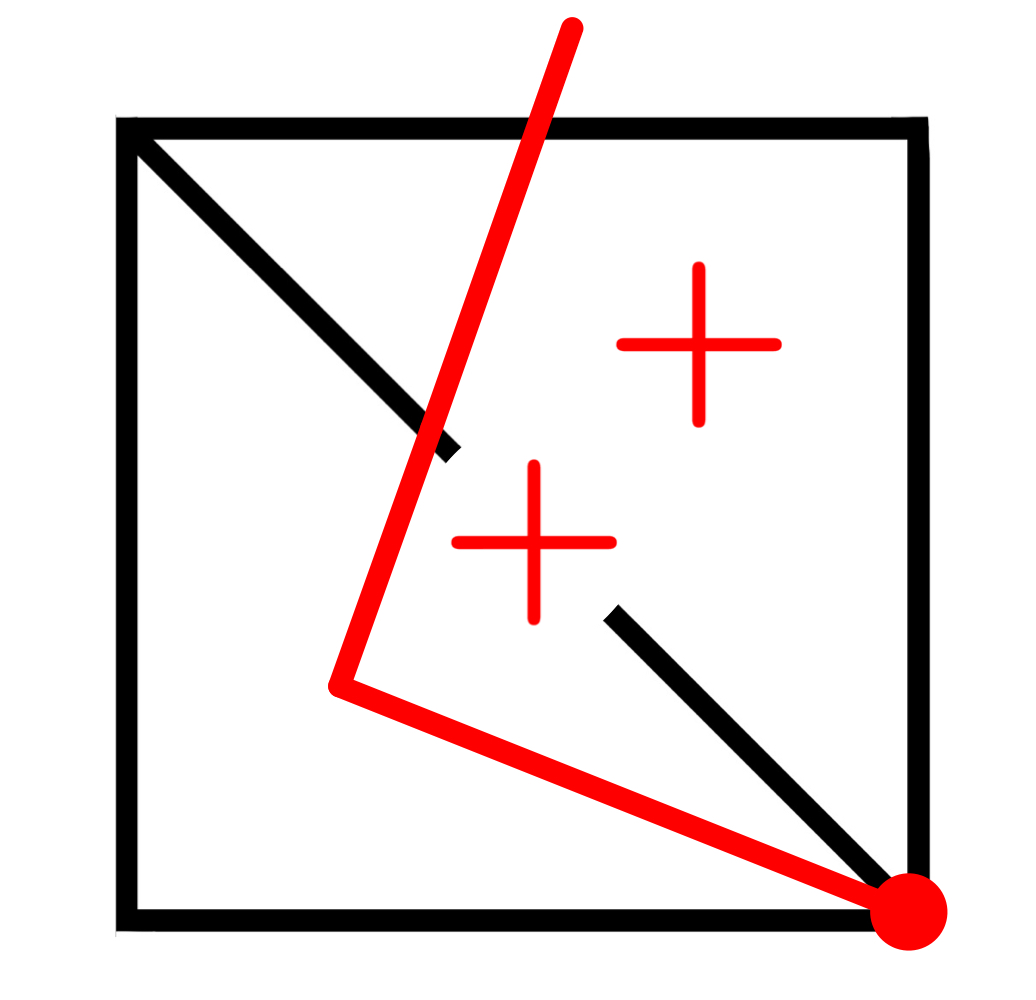}\\
\hline
&\\[-4mm]
\includegraphics[scale=0.04]{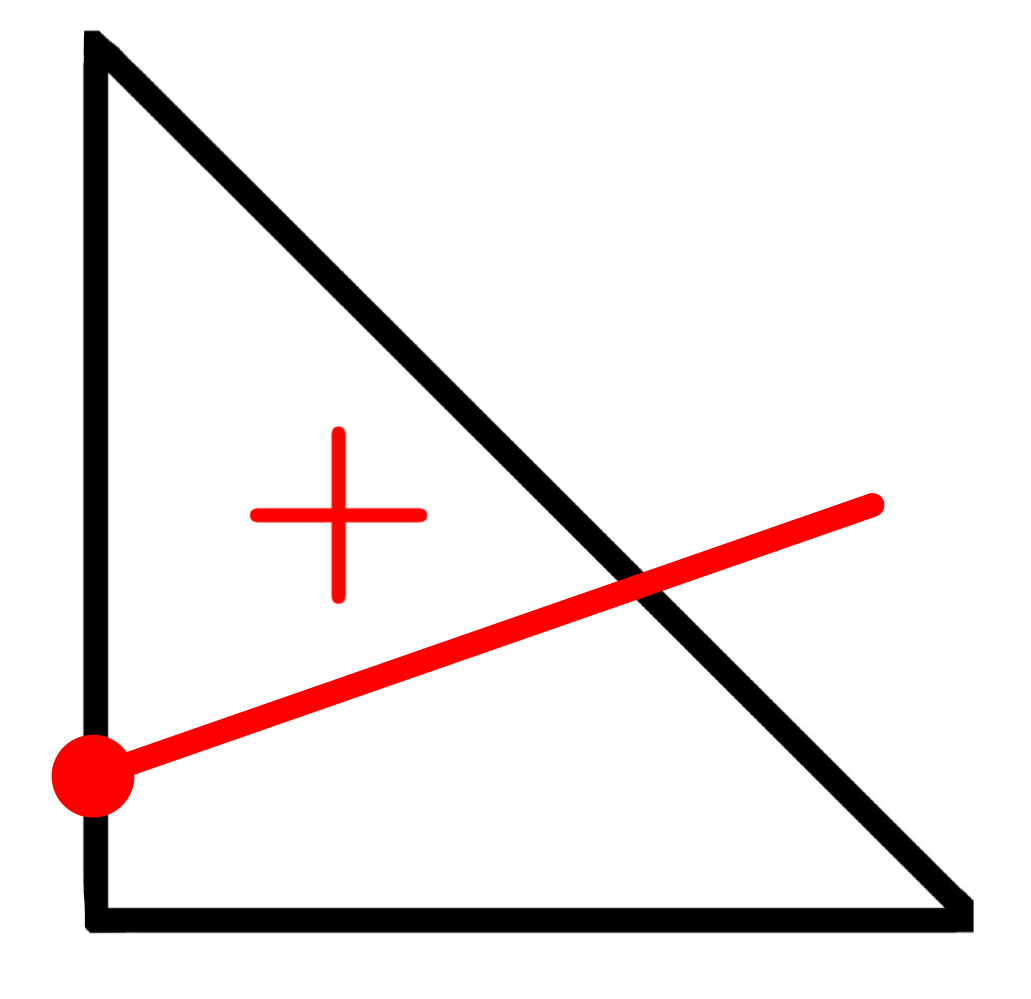} &\includegraphics[scale=0.04]{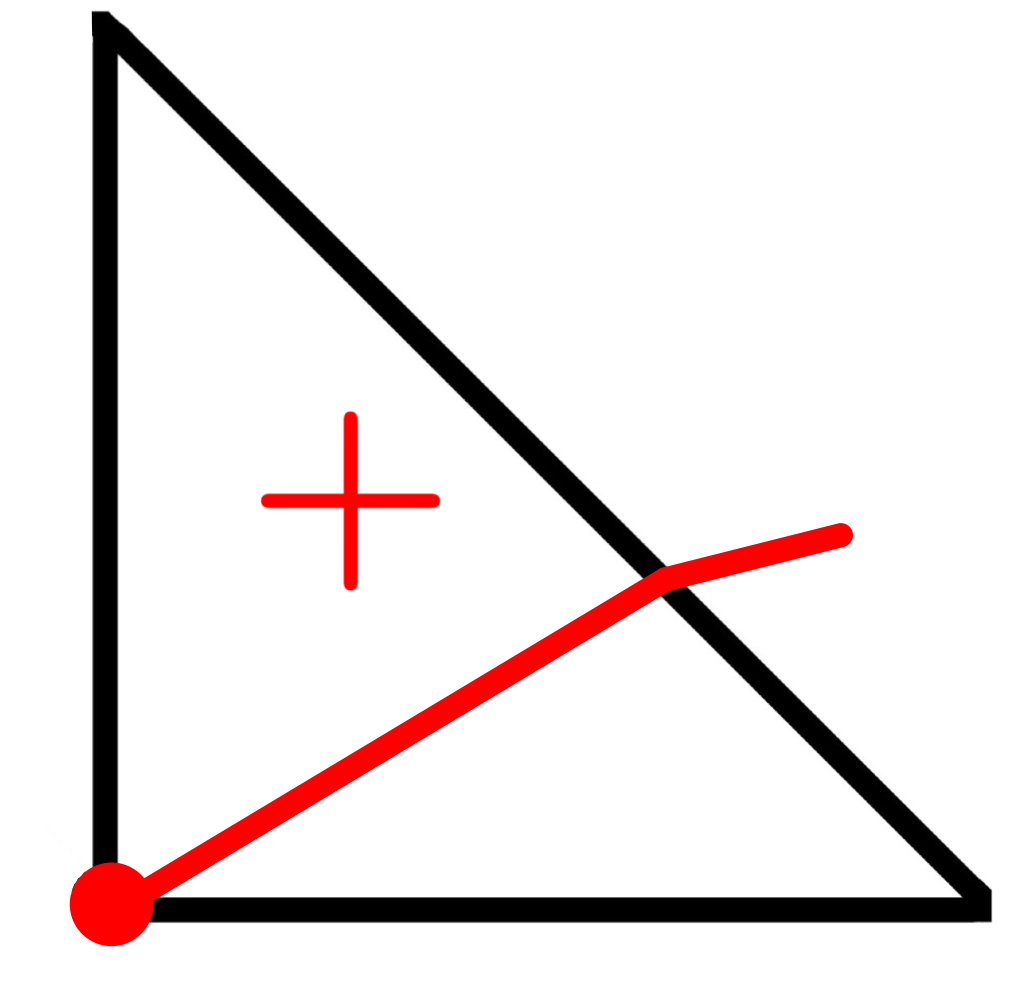}\\
\hline
&\\[-4mm]
\includegraphics[scale=0.04]{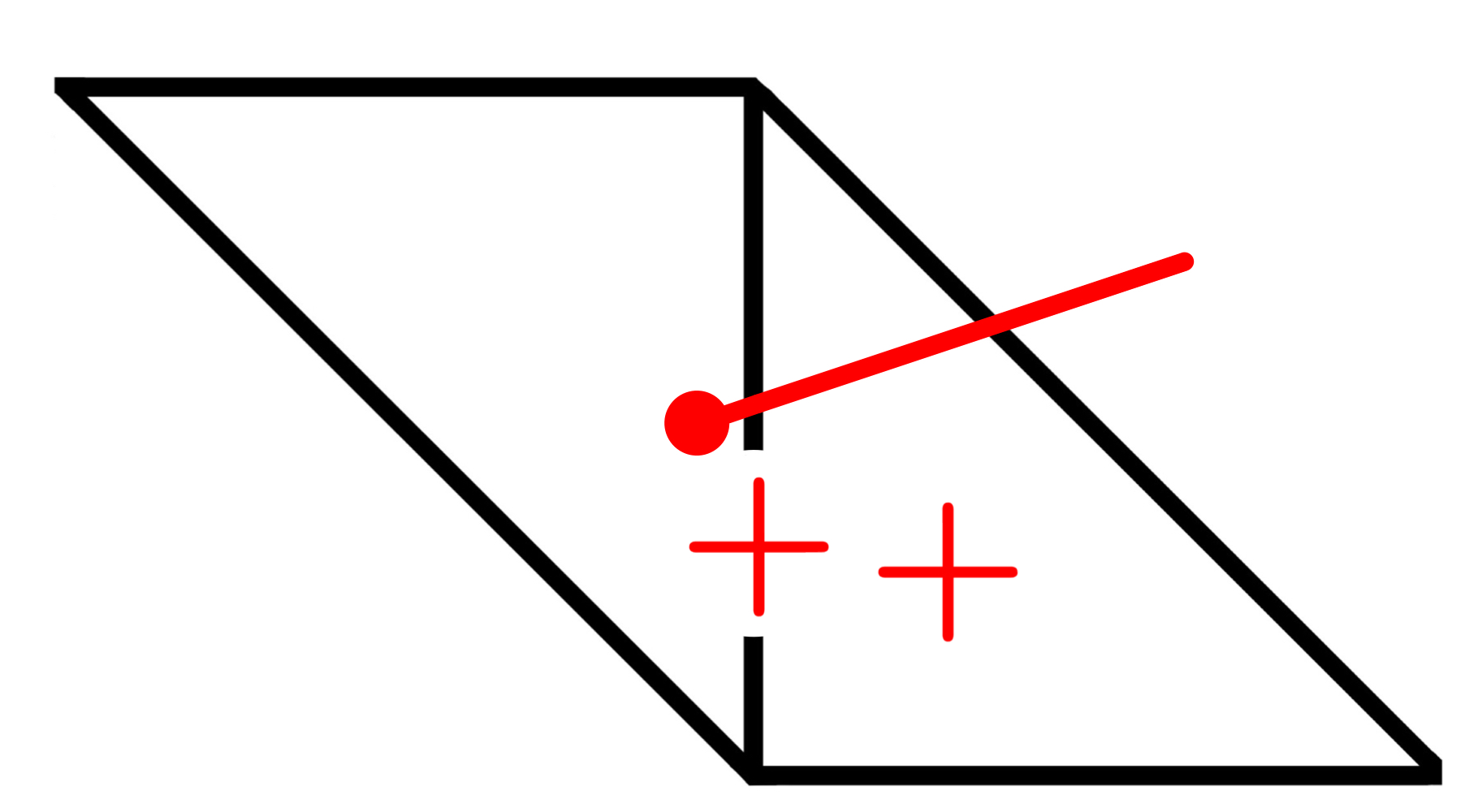} &\includegraphics[scale=0.04]{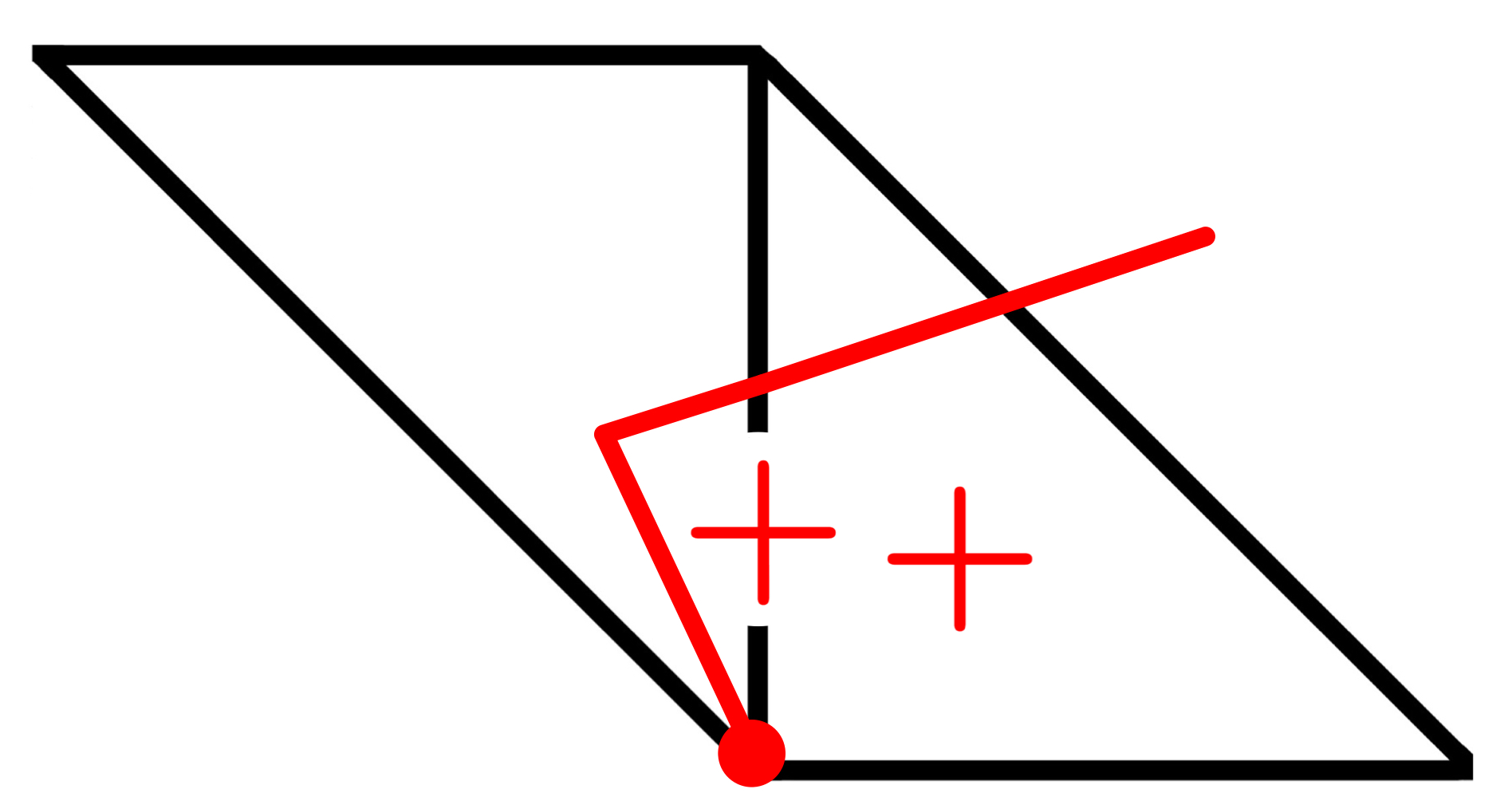}\\
\hline
&\\[-4.7mm]
\includegraphics[scale=0.04]{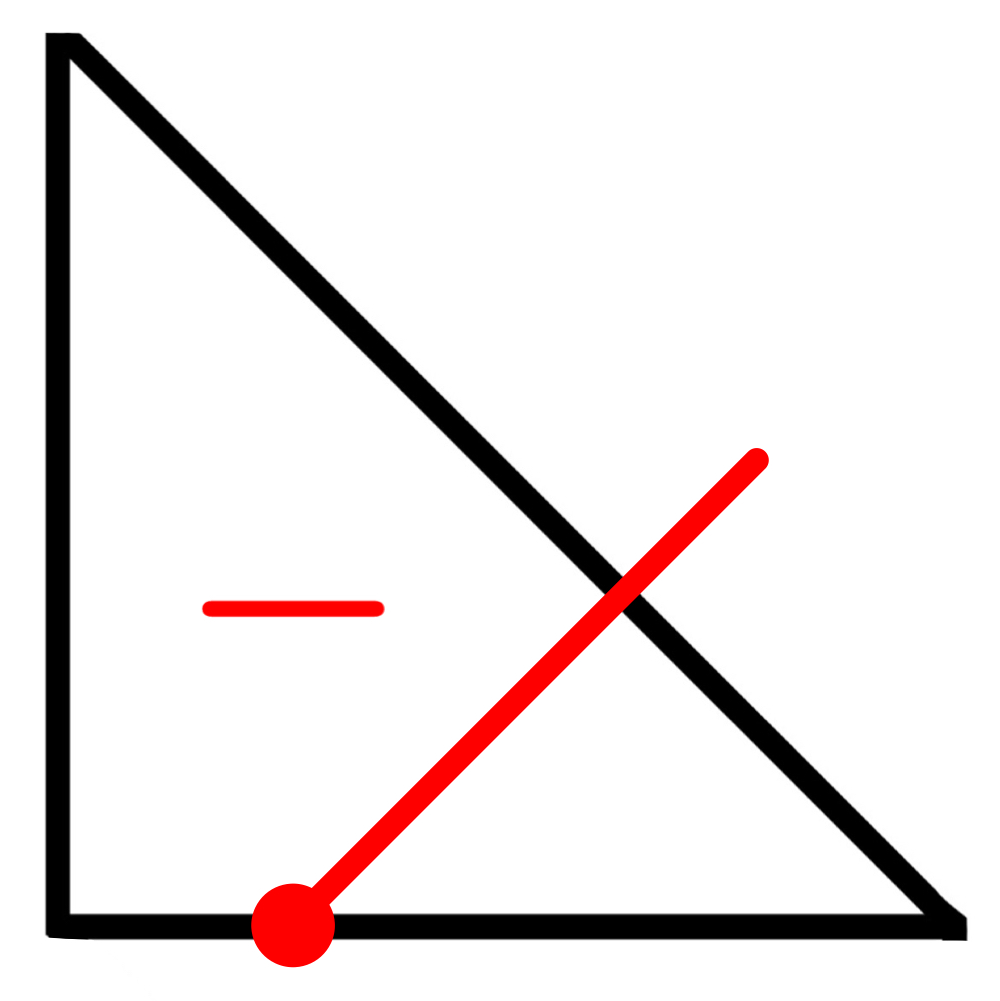} &\includegraphics[scale=0.04]{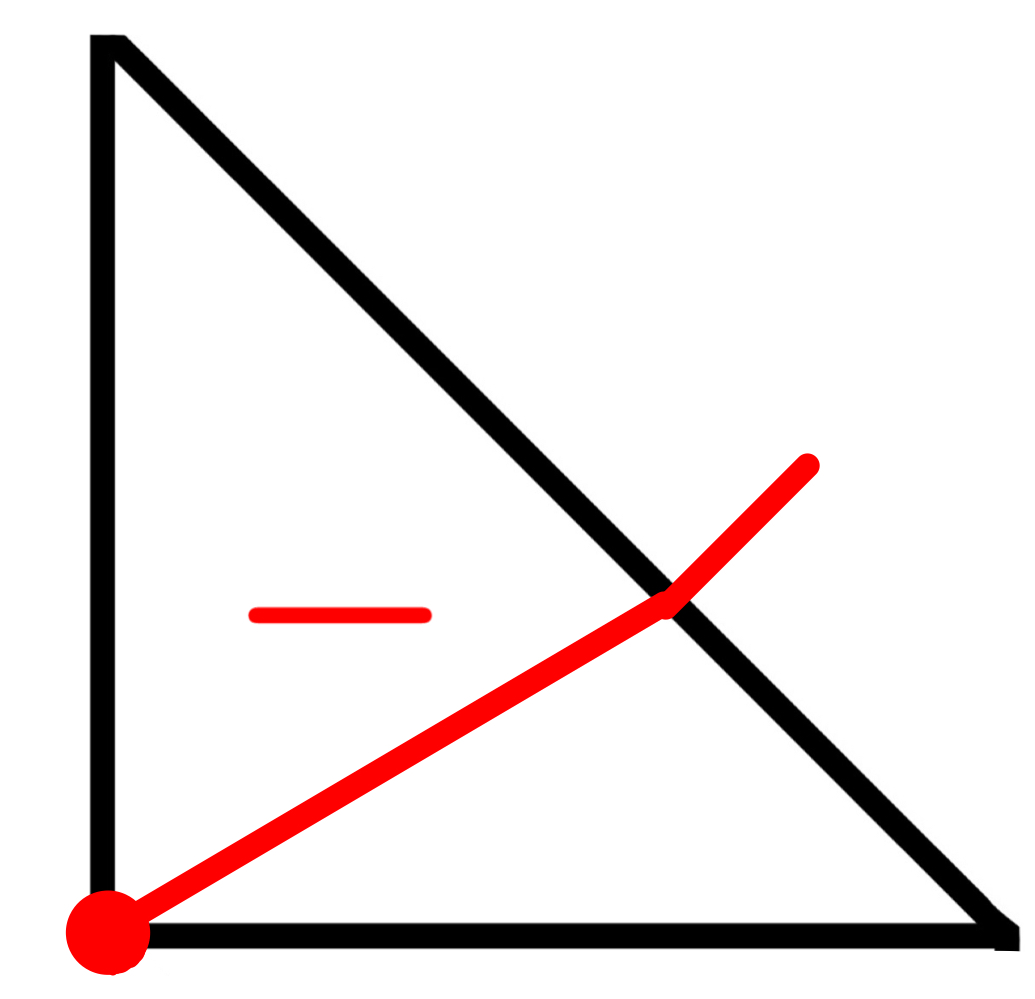}\\[-0.5mm]
\hline
&\\[-4mm]
\includegraphics[scale=0.04]{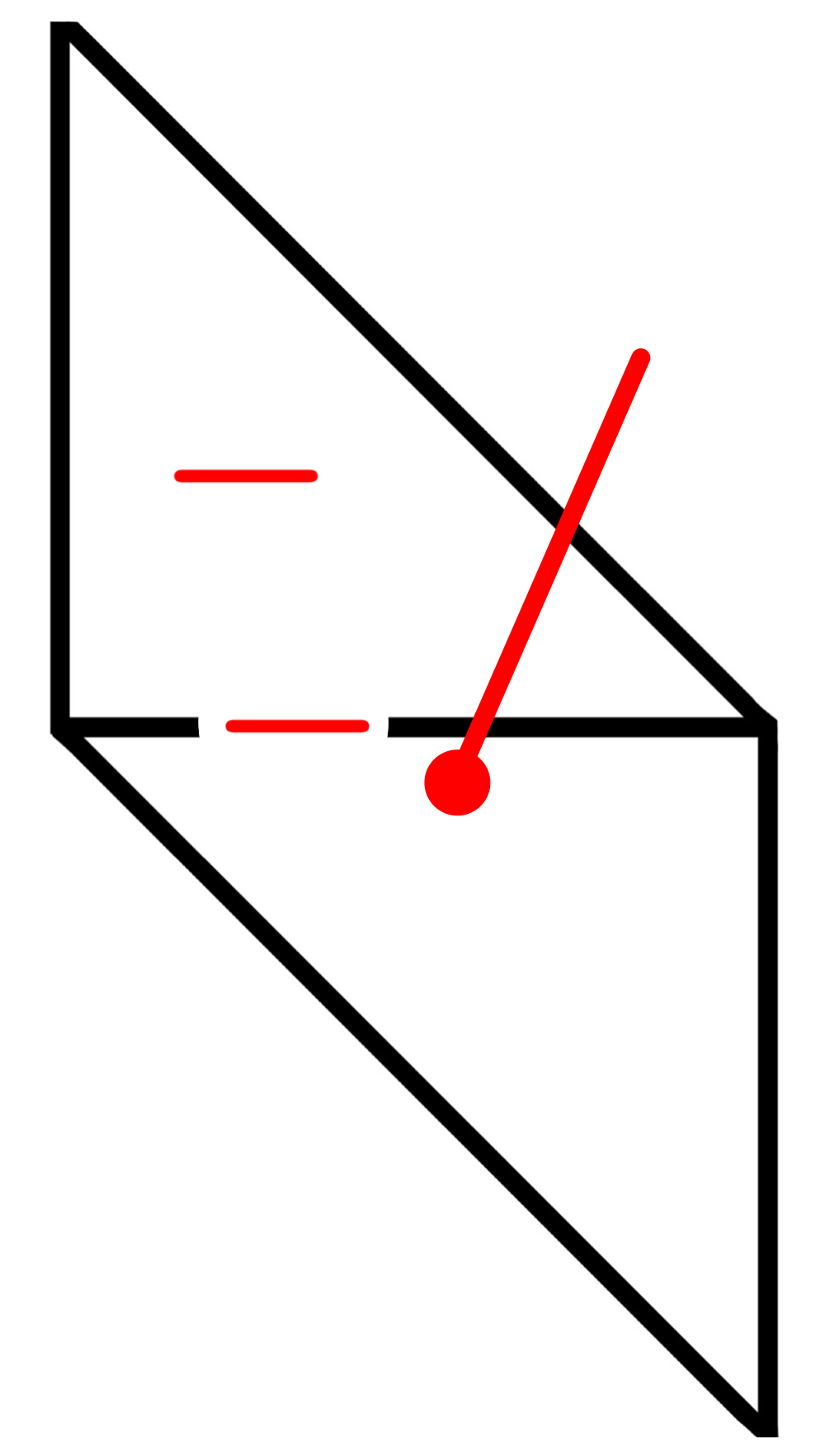} &\includegraphics[scale=0.04]{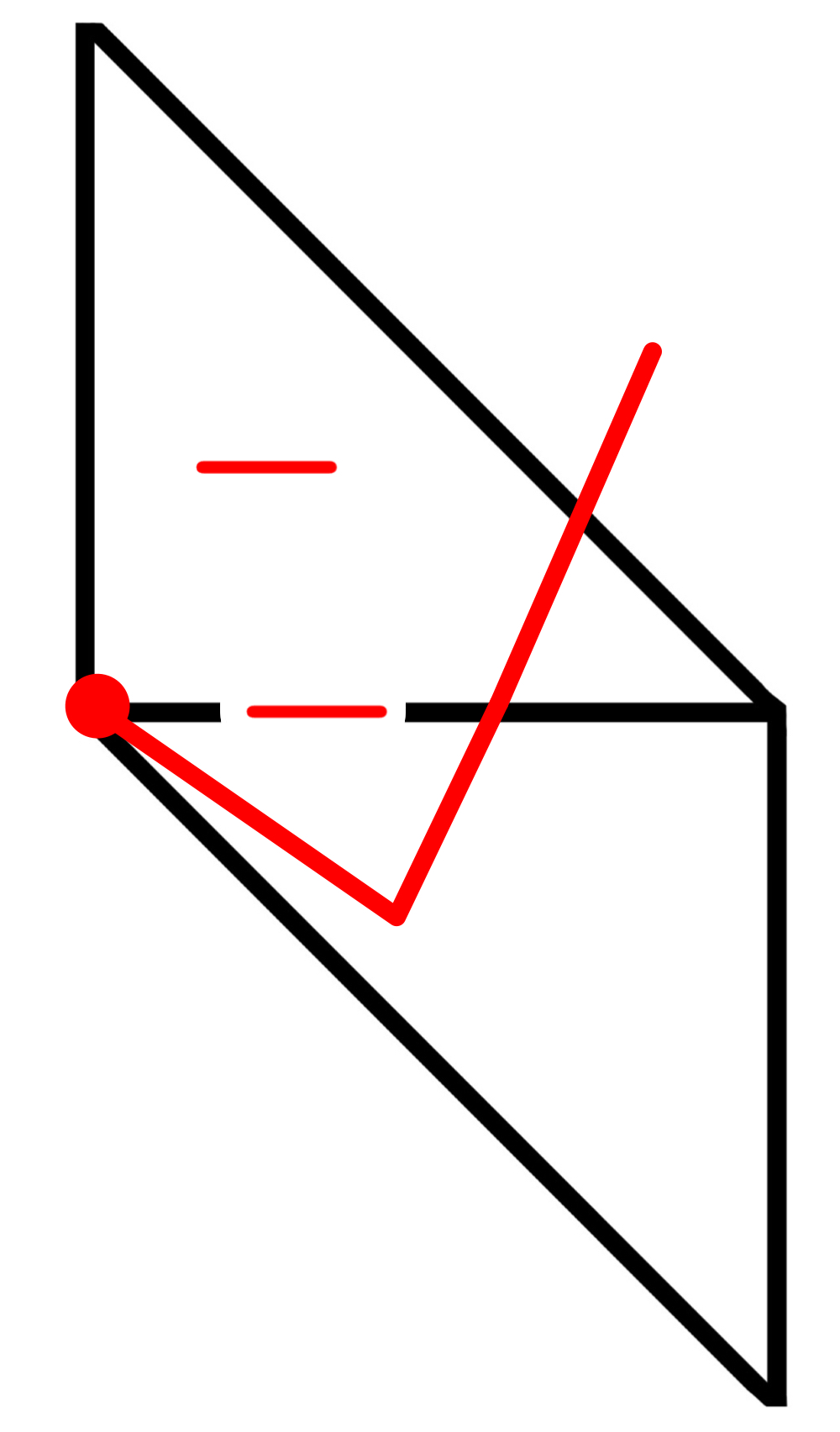}\\
\hline
\end{tabular}
\hspace{3mm}
\begin{tabular}{|c|c|}
Endpoint of $\overline{L_t}(w)$& $\widetilde{L_t}(w)$-Modification\\
\hline
&\\[-4mm]
\includegraphics[scale=0.037]{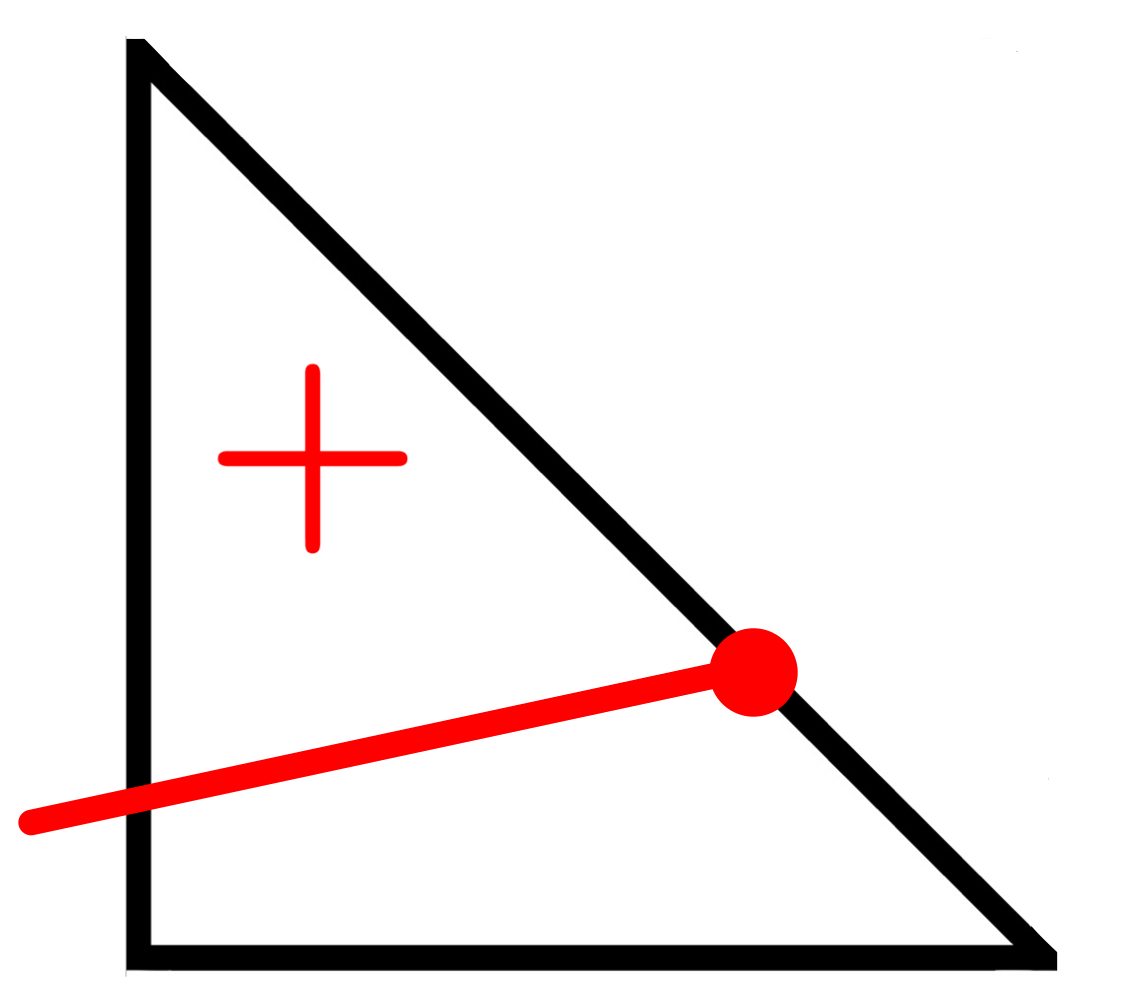} &\includegraphics[scale=0.04]{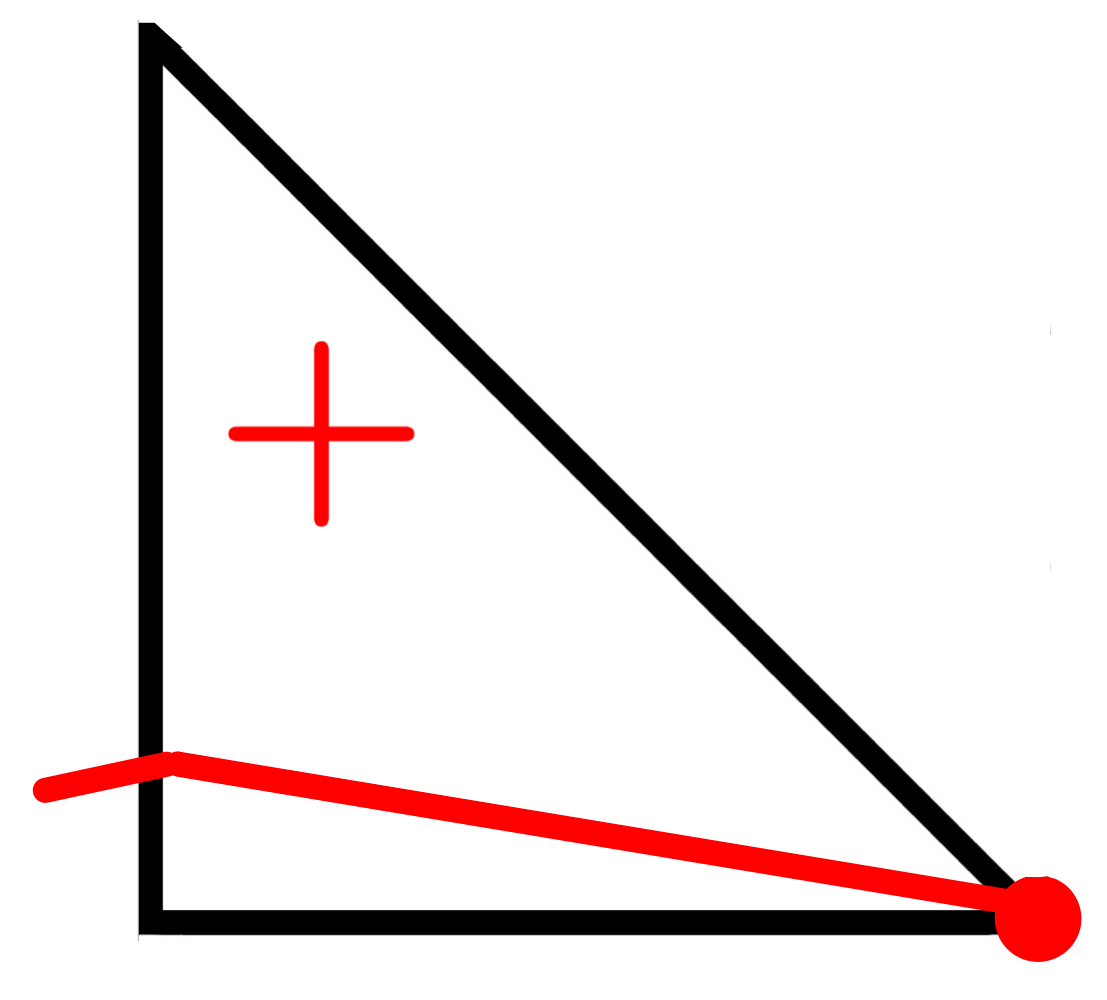}\\
\hline
&\\[-4mm]
\includegraphics[scale=0.036]{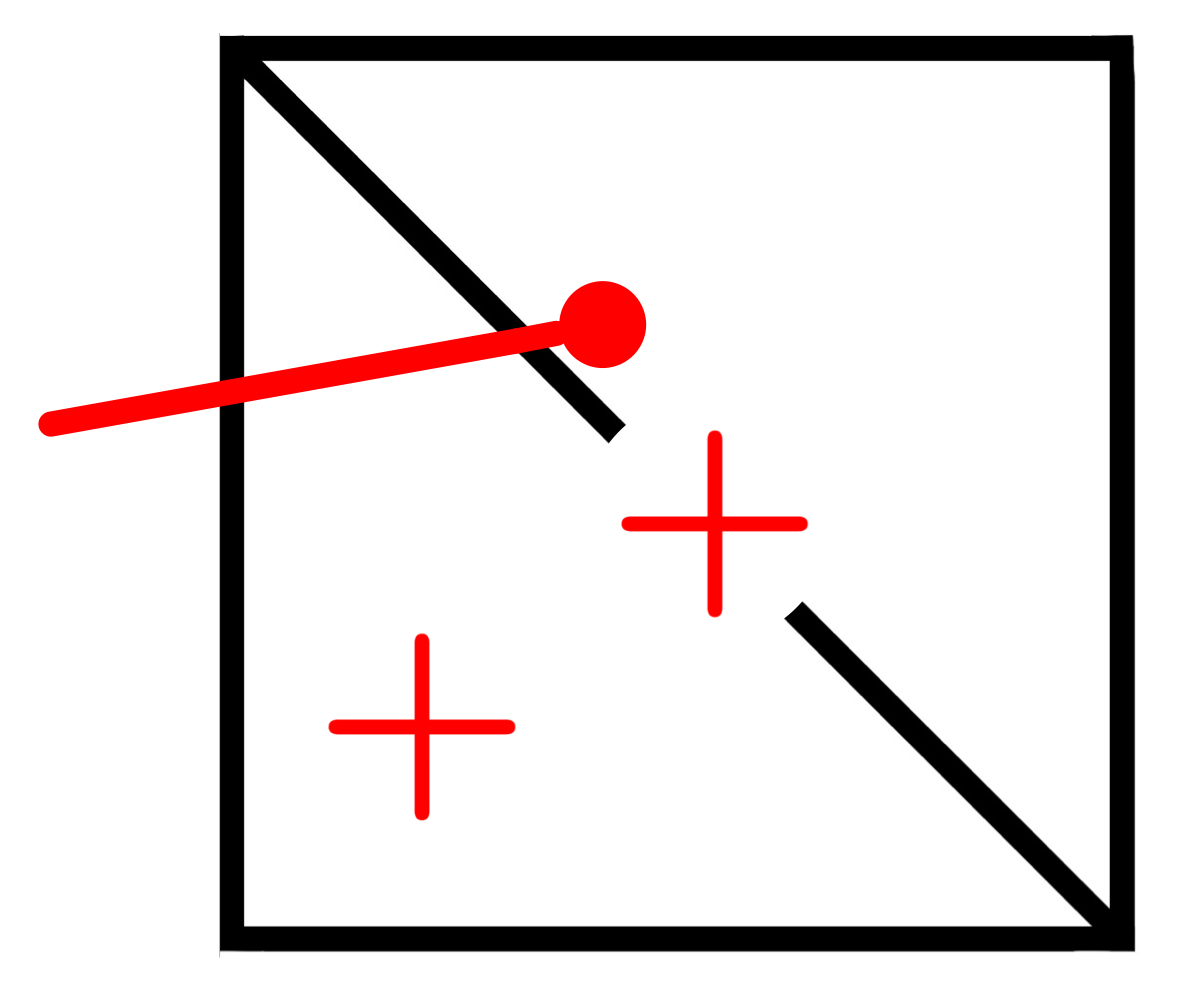} &\includegraphics[scale=0.037]{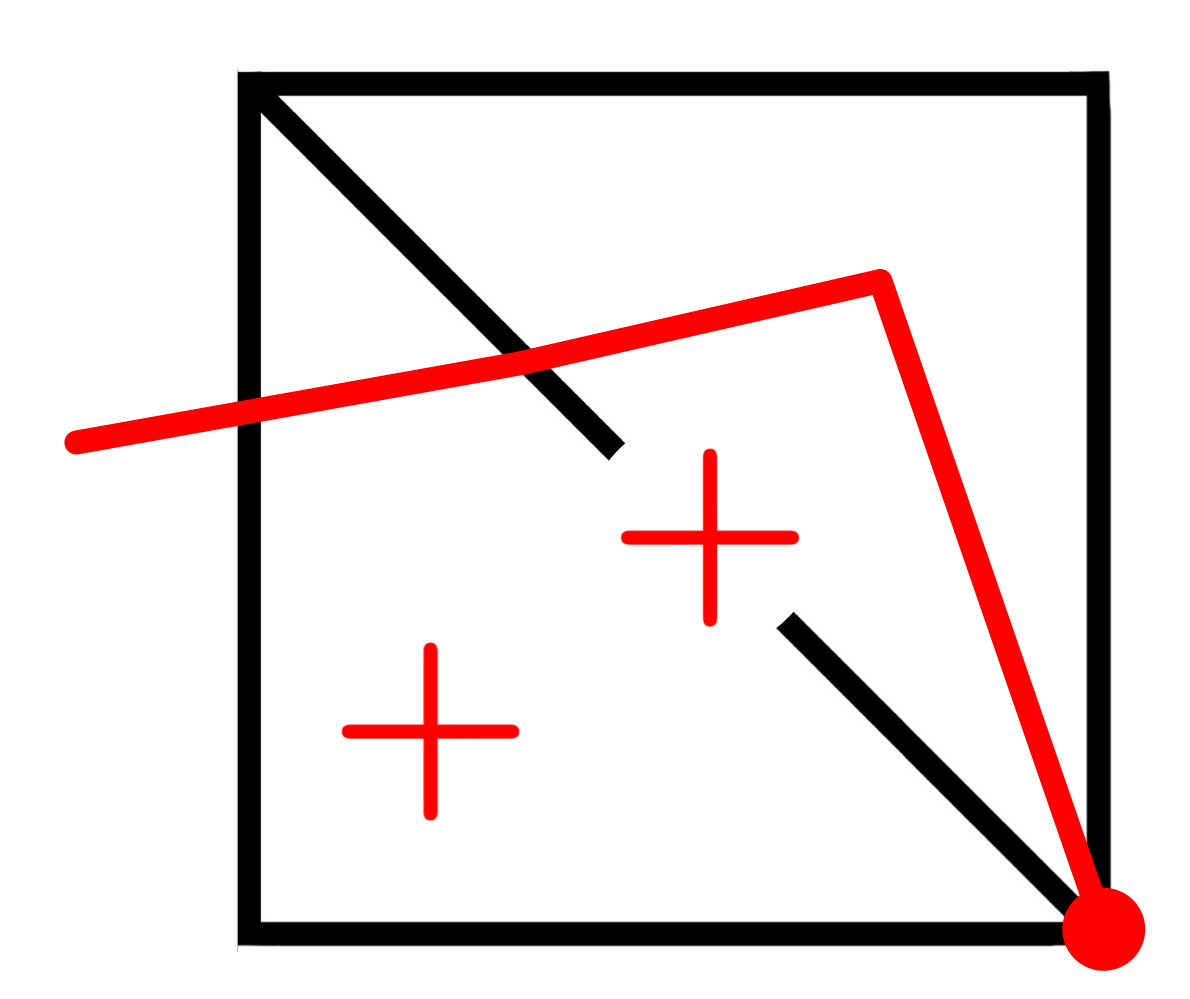}\\
\hline
&\\[-4mm]
\includegraphics[scale=0.038]{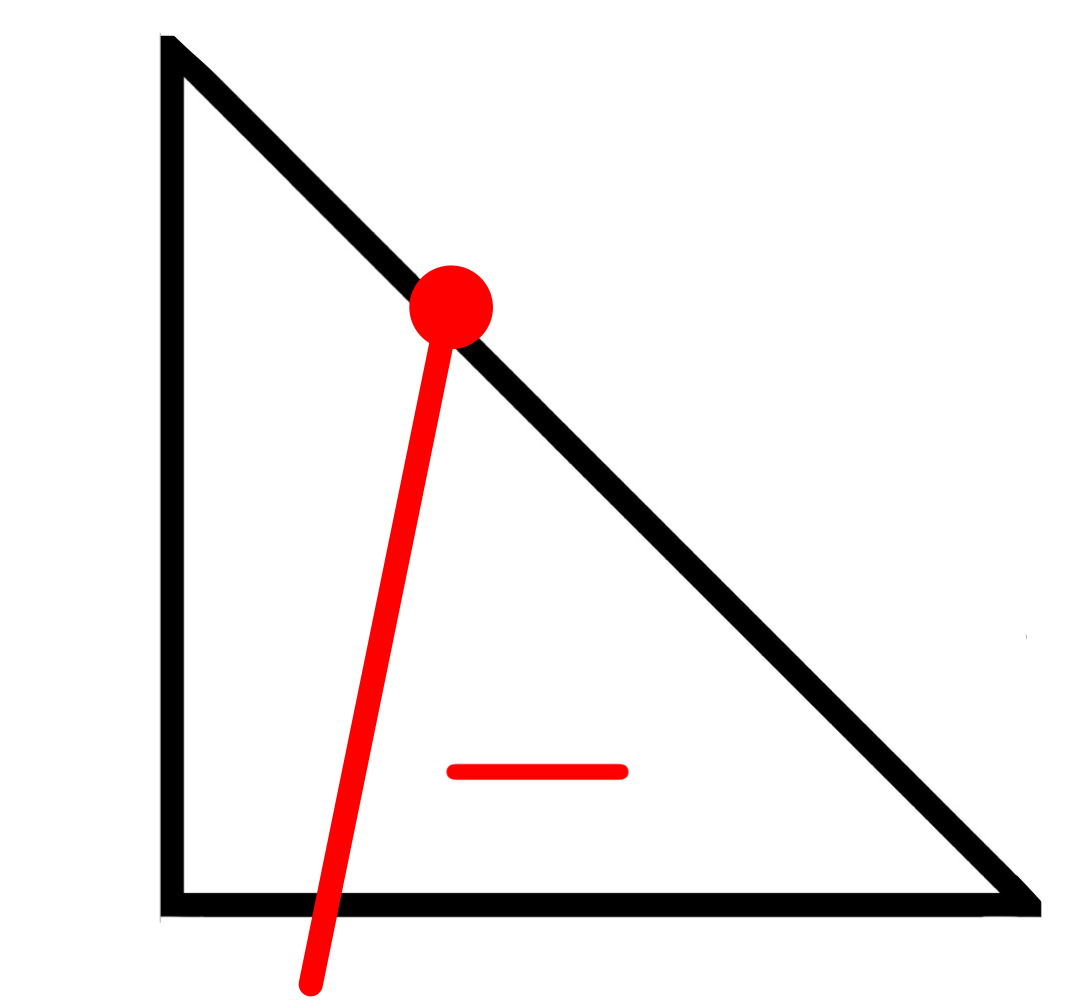} &\includegraphics[scale=0.04]{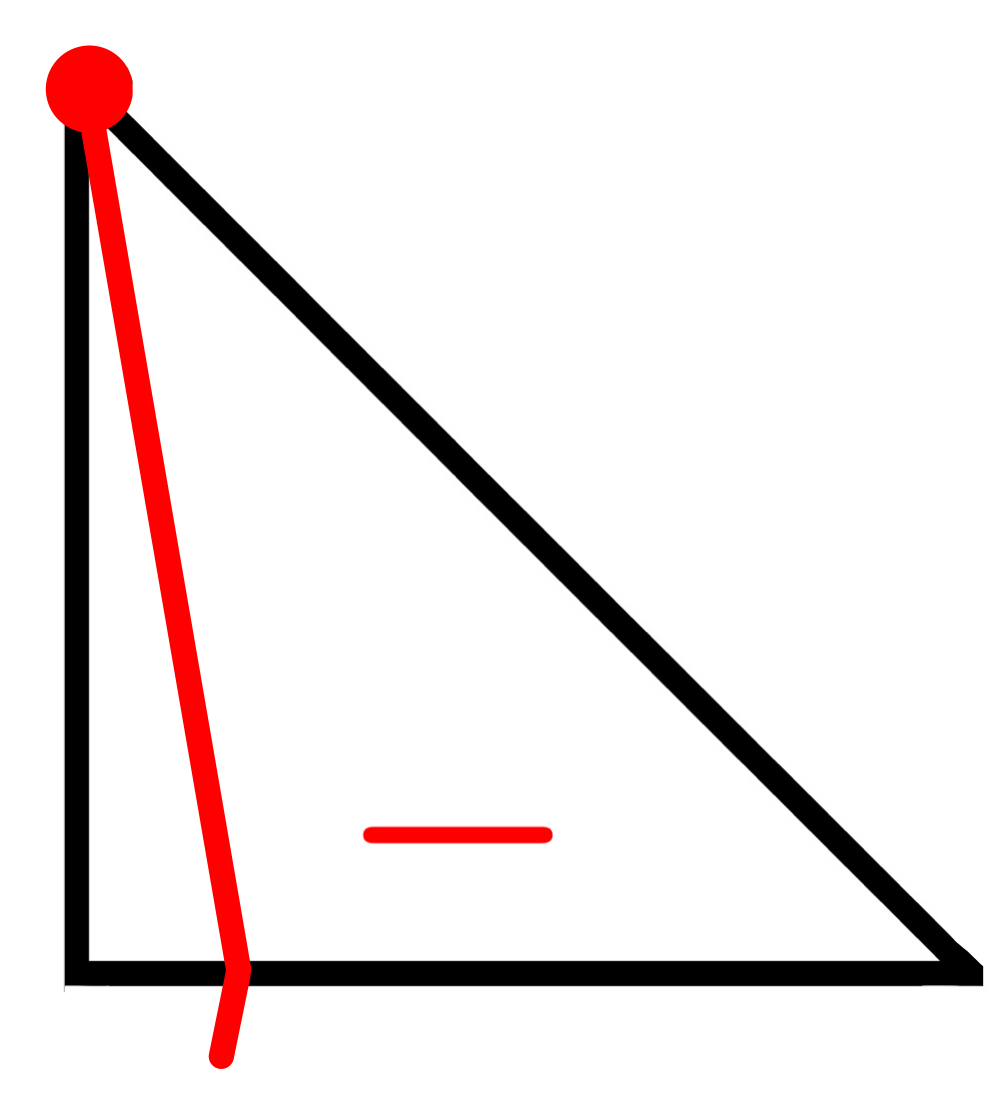}\\[-1mm]
\hline
&\\[-4mm]
\includegraphics[scale=0.04]{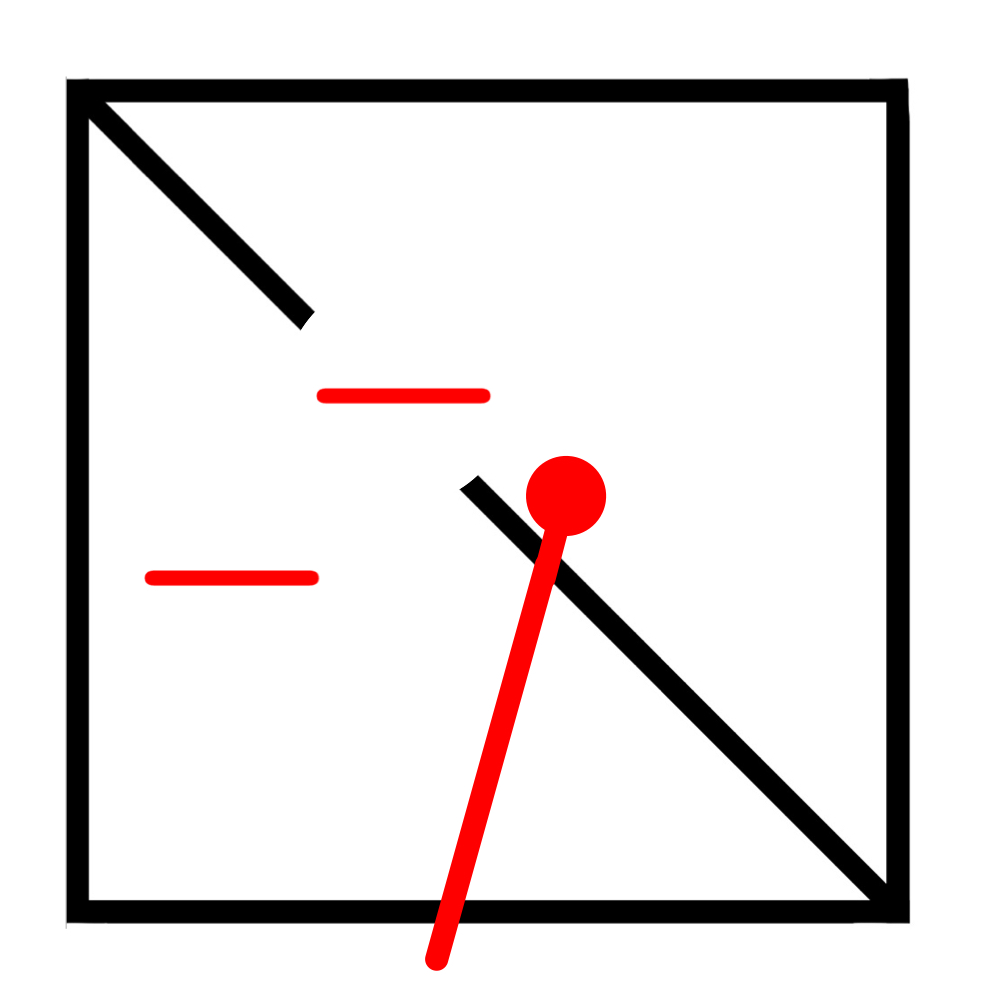} &\includegraphics[scale=0.04]{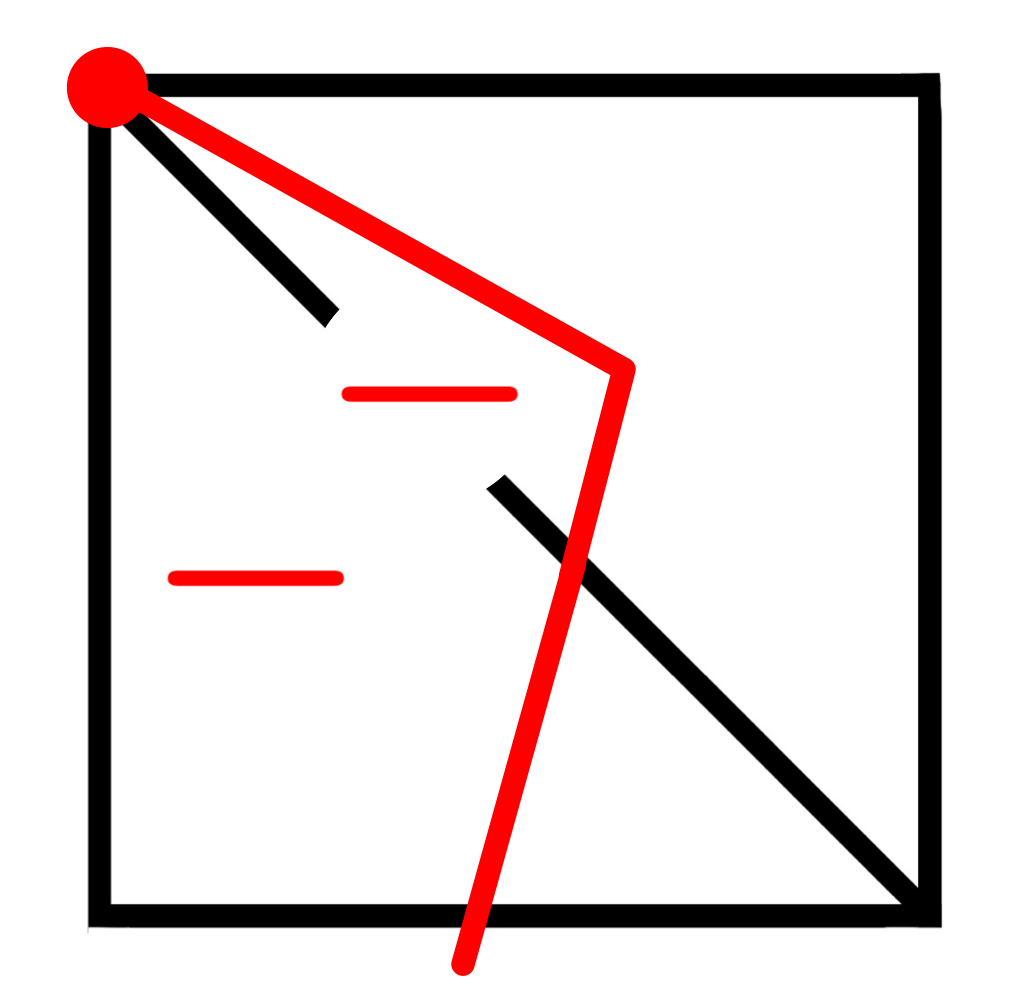}\\
\hline
&\\[-4mm]
\includegraphics[scale=0.04]{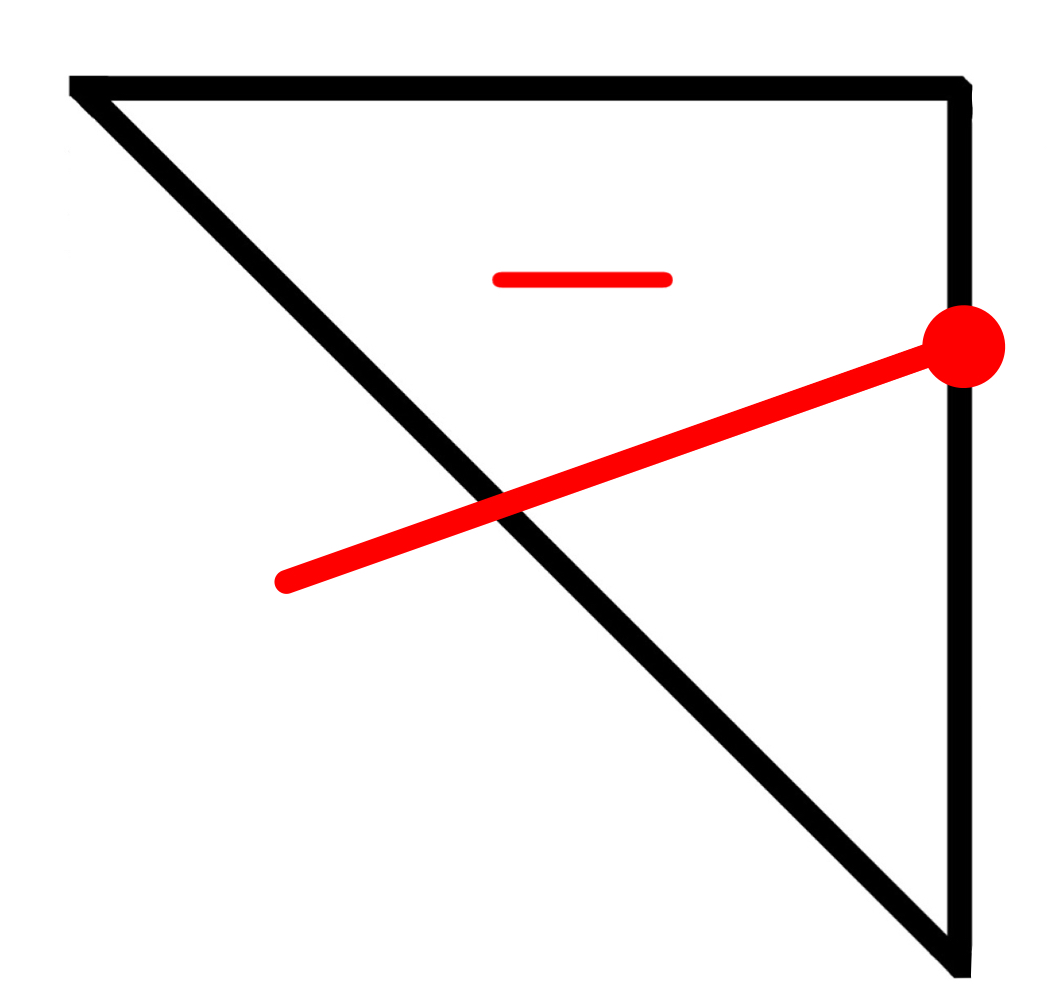} &\includegraphics[scale=0.04]{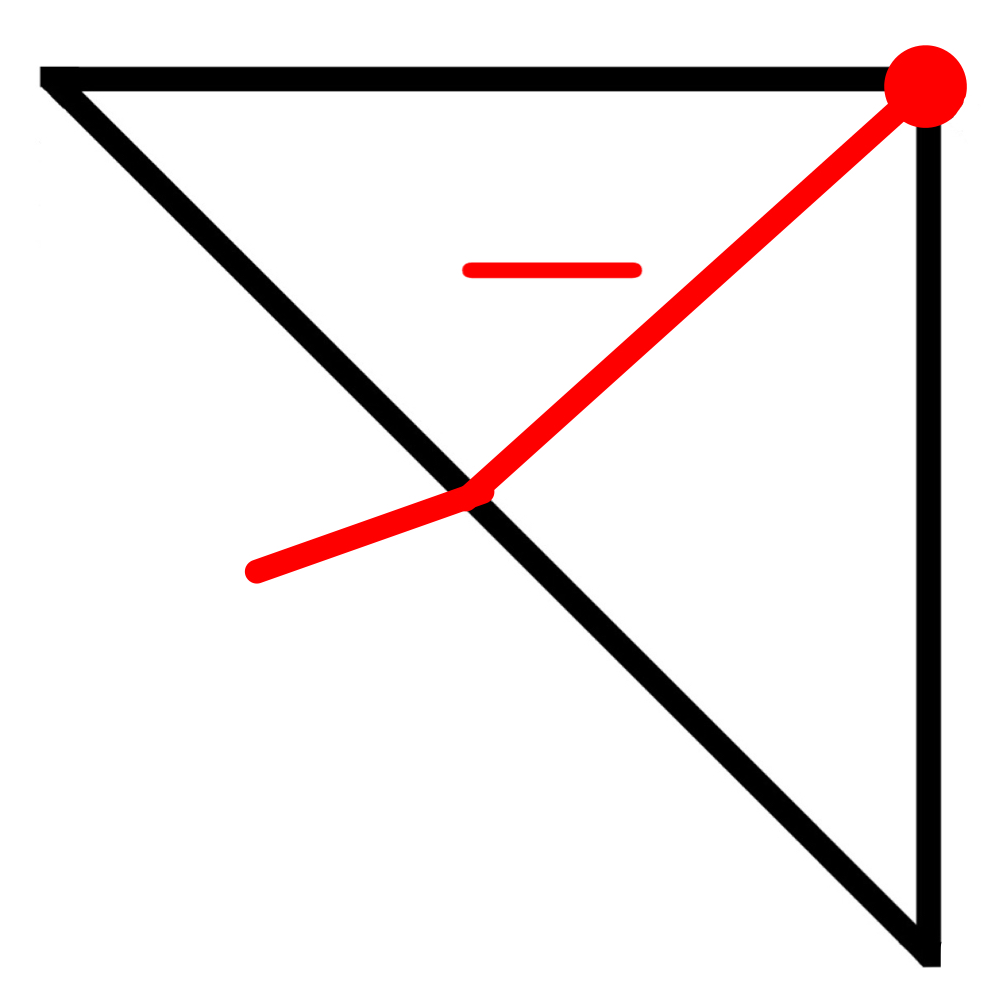}\\
\hline
&\\[-4mm]
\includegraphics[scale=0.04]{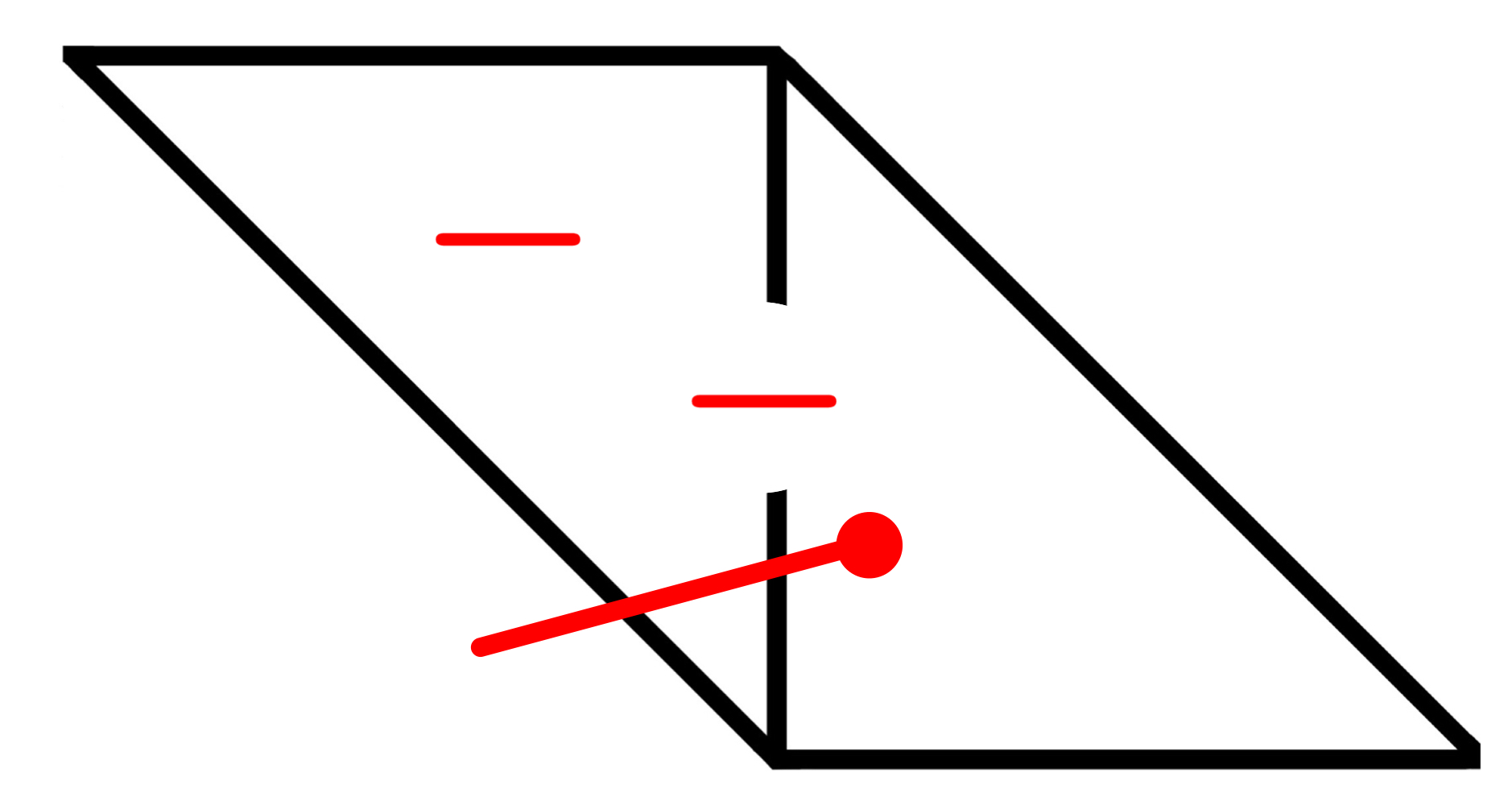} &\includegraphics[scale=0.04]{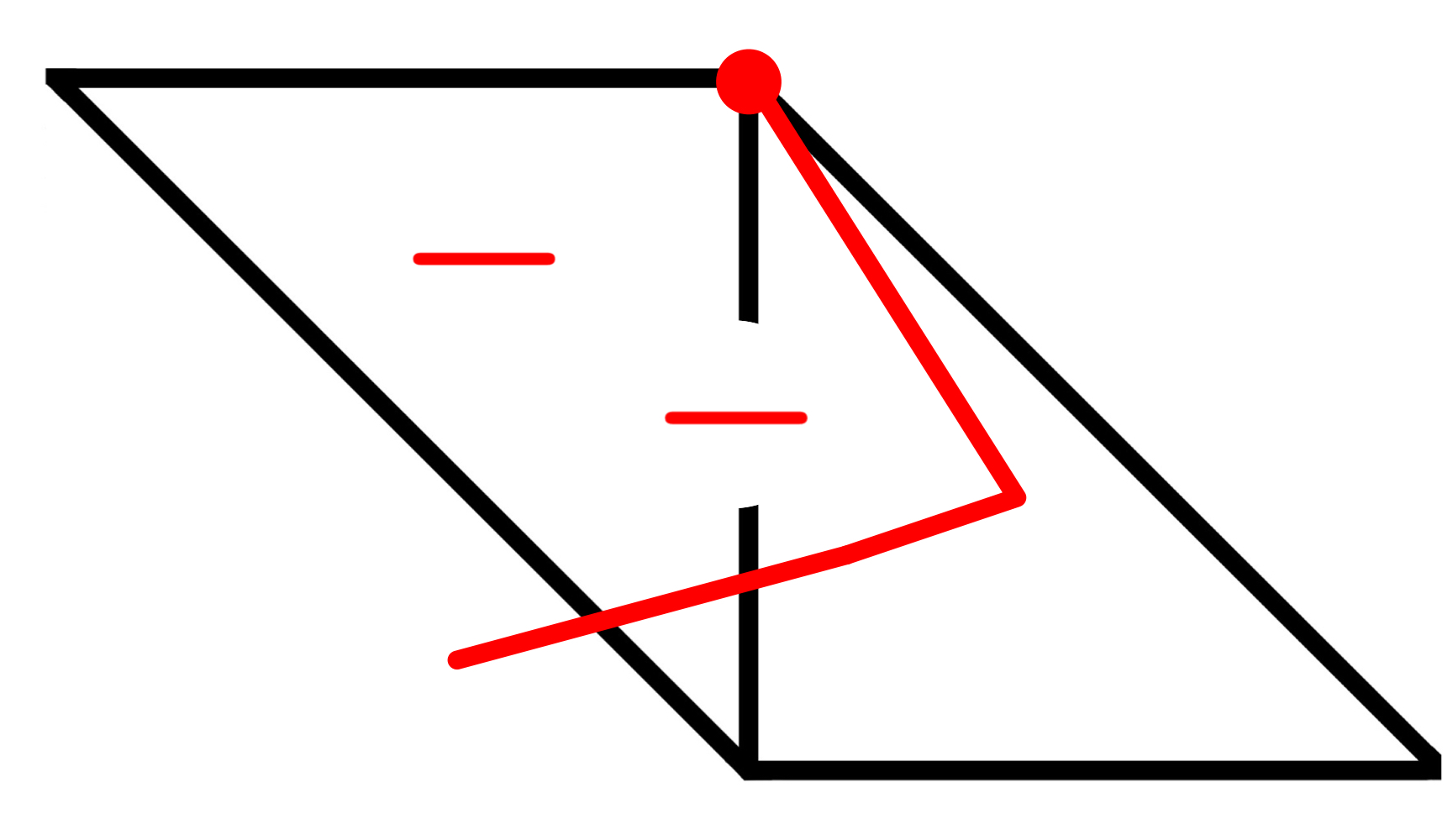}\\
\hline
&\\[-4mm]
\includegraphics[scale=0.04]{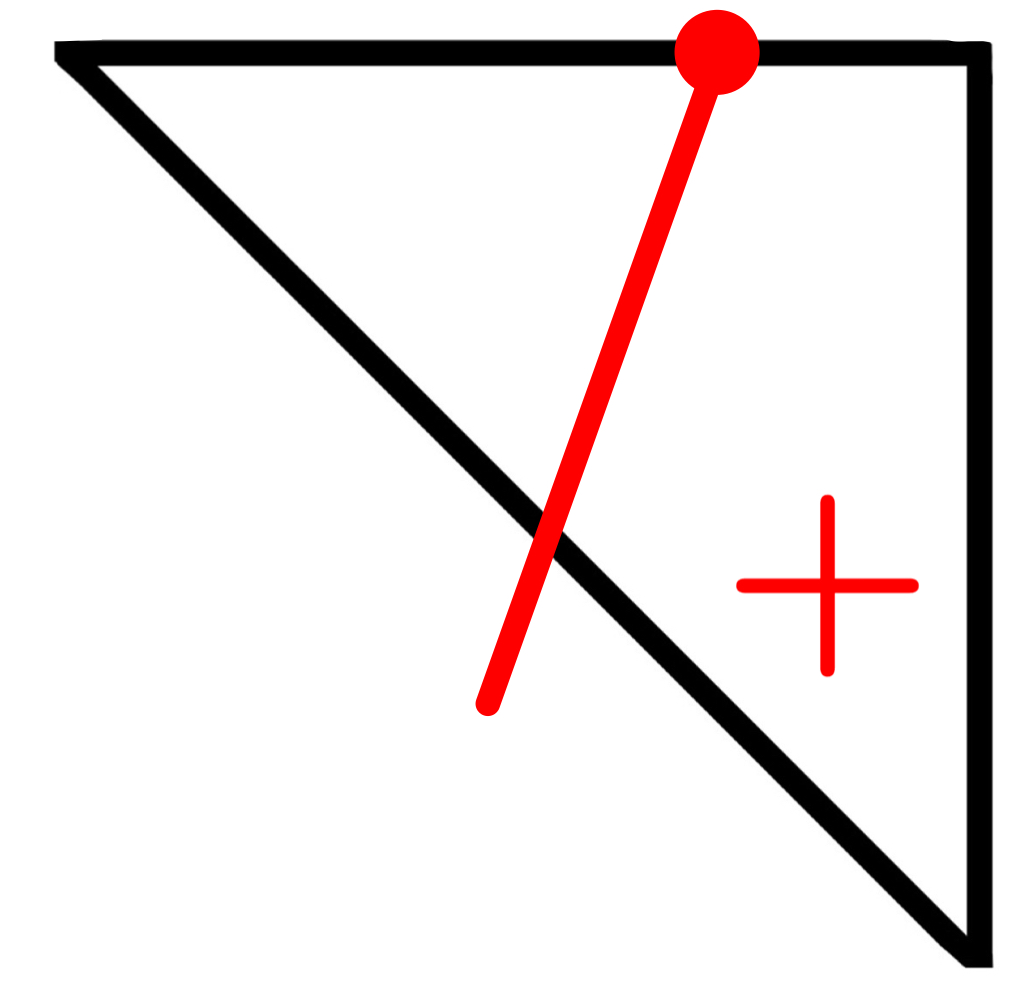} &\includegraphics[scale=0.04]{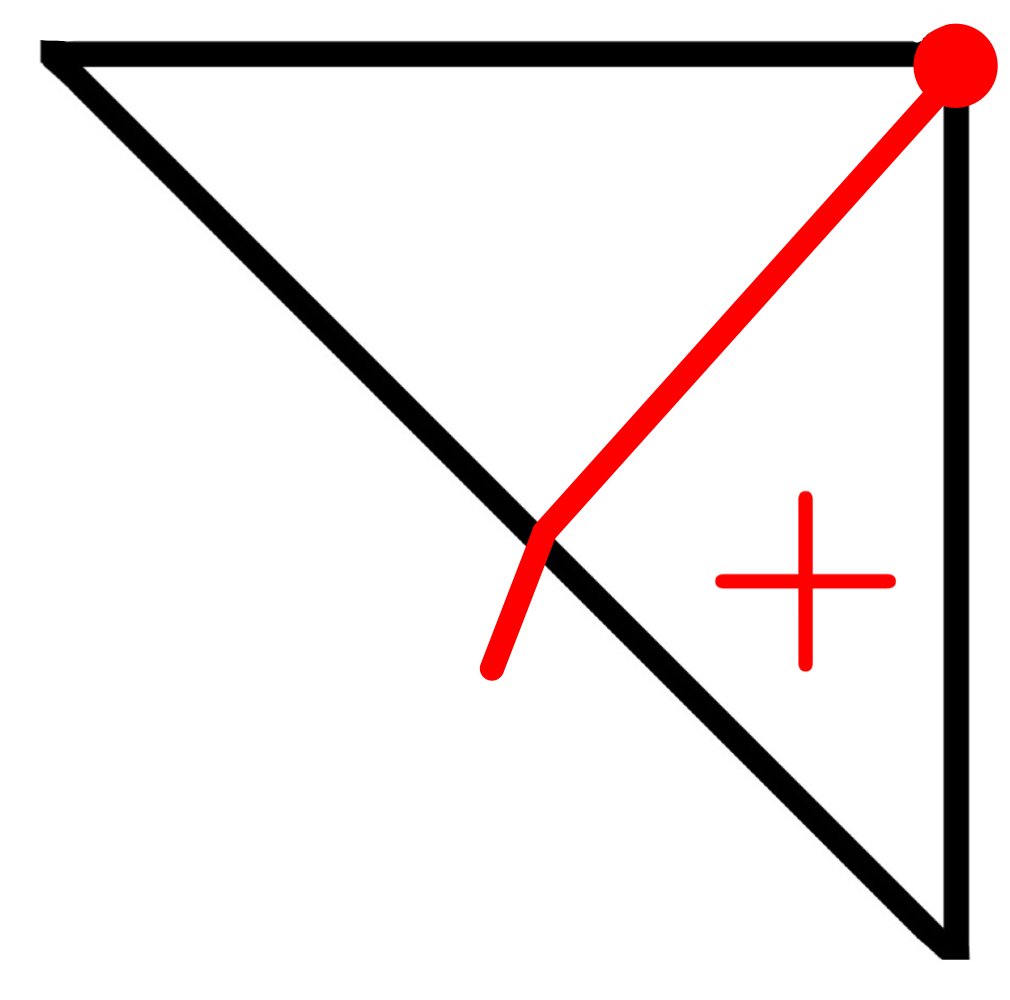}\\
\hline
&\\[-4mm]
\includegraphics[scale=0.04]{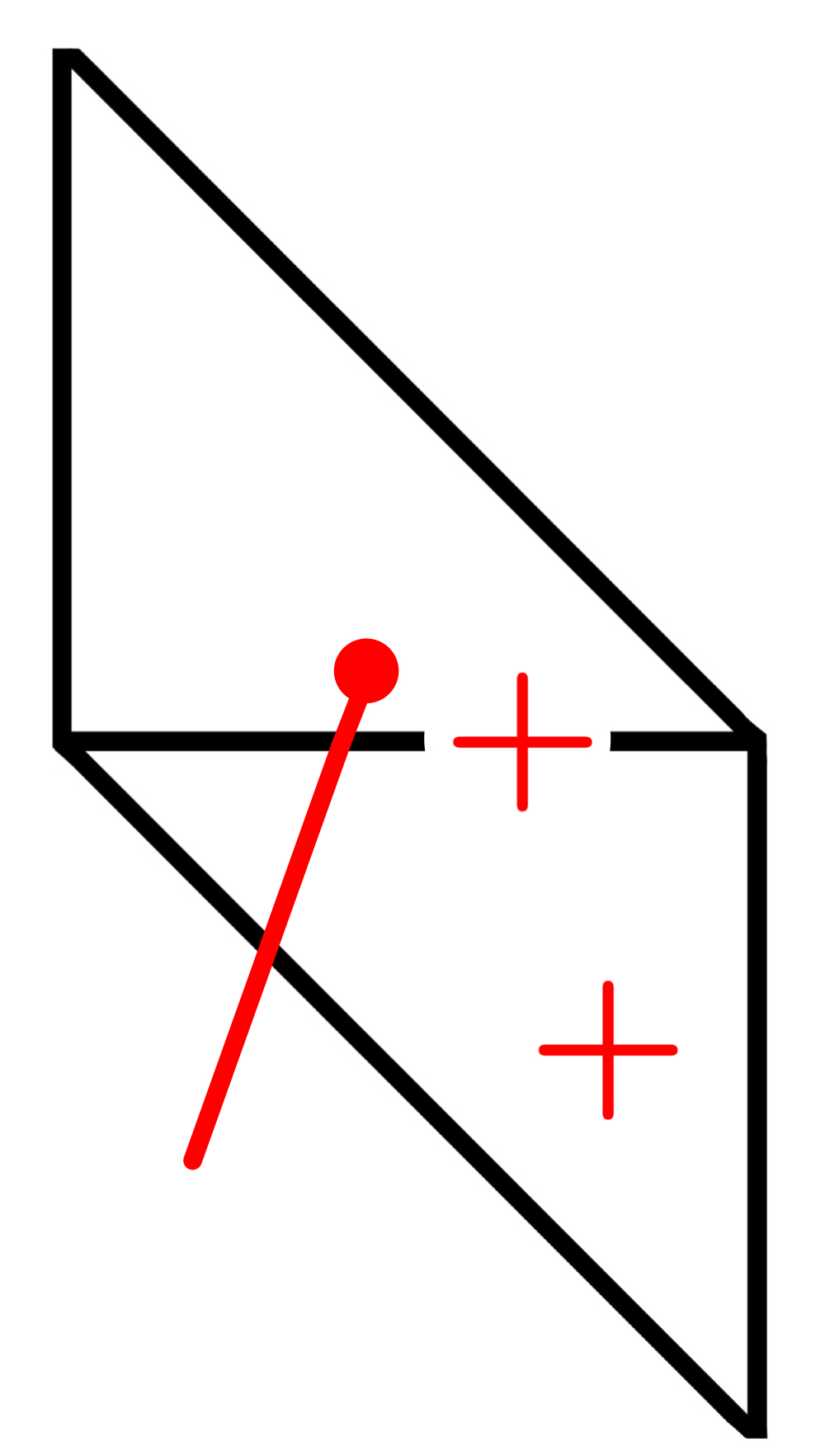} &\includegraphics[scale=0.04]{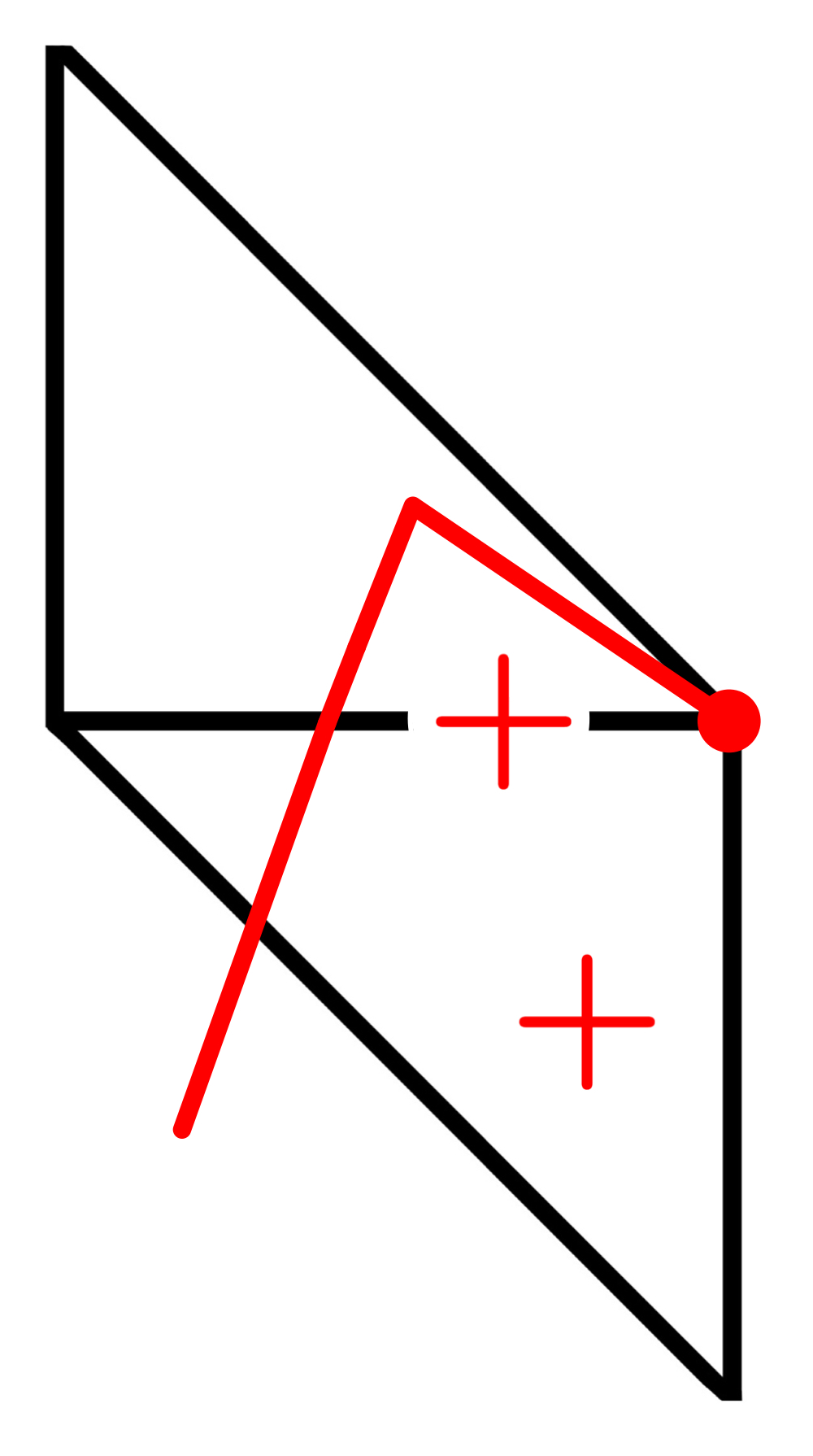}\\
\hline
\end{tabular}
\vspace{2mm}
\caption{Endpoint modification}\label{table1}
\end{table}

\begin{lemma}[Endpoint modification]\label{lem:endpoint-modification}
Let $t=p/q\in\mathbb Q_{>0}\setminus\{1\}$ be irreducible, write $s(t)=(a_1,\ldots,a_n)$, let $A=(0,0)$ and $B=(q,p)$, and let
\[
w=(a_{k+1},a_{k+2},\ldots,a_{k+n-1})
\]
be a cyclic block of length $n-1$, with indices read modulo $n$. Then the construction in Table~\ref{table1} has the following properties.
\begin{enumerate}
    \item Grouping adjacent equal signs of $\widetilde{L_t}(w)$, including its endpoint signs, and recording their block lengths gives $w$. Its lift to $\widetilde{\mathbb R^2}$ beginning at $A$ ends at $B$.
    \item The only possible obstruction to $\widetilde{L_t}(w)$ being a generalized arc is a pair of consecutive crossings of the same unit edge at an endpoint. Deleting each such endpoint passage produces a generalized arc $\widetilde{L_t}(w)_{-}$ from $A$ to $B$.
    \item If a deletion is necessary, then
    \[
    s(\widetilde{L_t}(w)_{-})
    =(b_{k+1},a_{k+2},\ldots,a_{k+n-2},b_{k+n-1}),
    \]
    where $0\leq b_{k+1}\leq a_{k+1}$ and $0\leq b_{k+n-1}\leq a_{k+n-1}$, at least one inequality is strict, and a zero endpoint entry is omitted. Moreover,
    \[
    |\widetilde{L_t}(w)_{-}|<m(\mathcal G[w]).
    \]
    \item For $w_0=(a_2,\ldots,a_n)$ no deletion is required, and $\widetilde{L_t}(w_0)$ is the left push-off $\gamma^L_{AB}$.
\end{enumerate}
\end{lemma}

\begin{proof}
Choose mutually disjoint endpoint disks which meet $\overline{L_t}(w)$ only in its initial or terminal subarc and contain no other lattice point. The endpoint branches in Table~\ref{table1} are drawn inside these disks. They may therefore be chosen transverse to every edge, disjoint from every edge midpoint, and disjoint from the rest of the curve. In particular, they create neither a self-intersection nor an additional winding around the puncture.

The table is exhaustive: a cut lies immediately before or after one of the two possible triangle signs or an edge-sign block, and at the two ends the orientation is reversed. The paired drawings record both choices of the adjacent elementary triangle in each of these local states.

At each end, the branch in the right-hand column of the table continues the unchanged interior branch to the lattice endpoint. The endpoint rule supplies the boundary sign removed when the cyclic sign word was cut, while every edge and triangle sign in the interior is unchanged. Reading from the initial to the terminal endpoint therefore gives exactly the consecutive run lengths
\[
a_{k+1},a_{k+2},\ldots,a_{k+n-1}.
\]
This proves the sign assertion in (1). The modification adds no horizontal or vertical winding, so it preserves the signed crossing numbers of the original torus loop with the horizontal and vertical fundamental cycles, namely $(q,p)$. The lift to $\widetilde{\mathbb R^2}$ starting at $(0,0)$ consequently ends at $(q,p)$.

\begin{figure}[ht]
    \centering
    \includegraphics[scale=0.05]{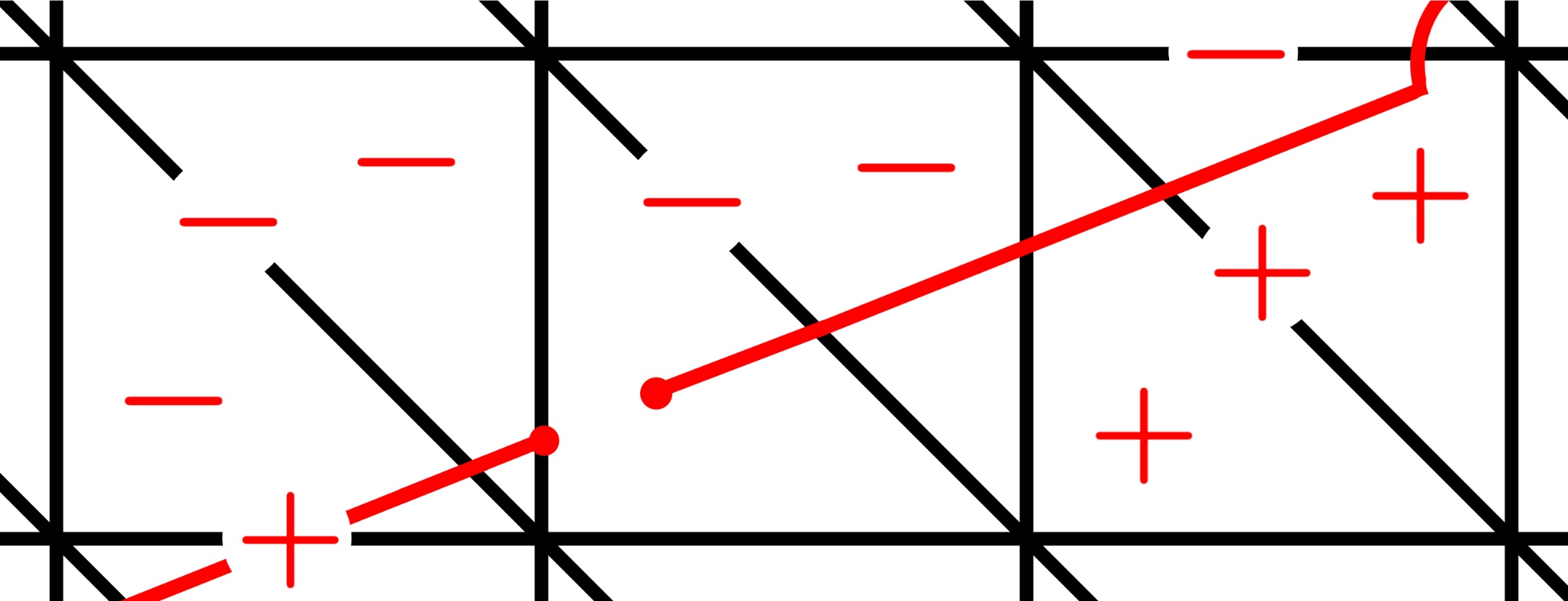}\hspace{3mm}
   \raisebox{1.2em}{$\mapsto$}\hspace{3mm}
    \includegraphics[scale=0.05]{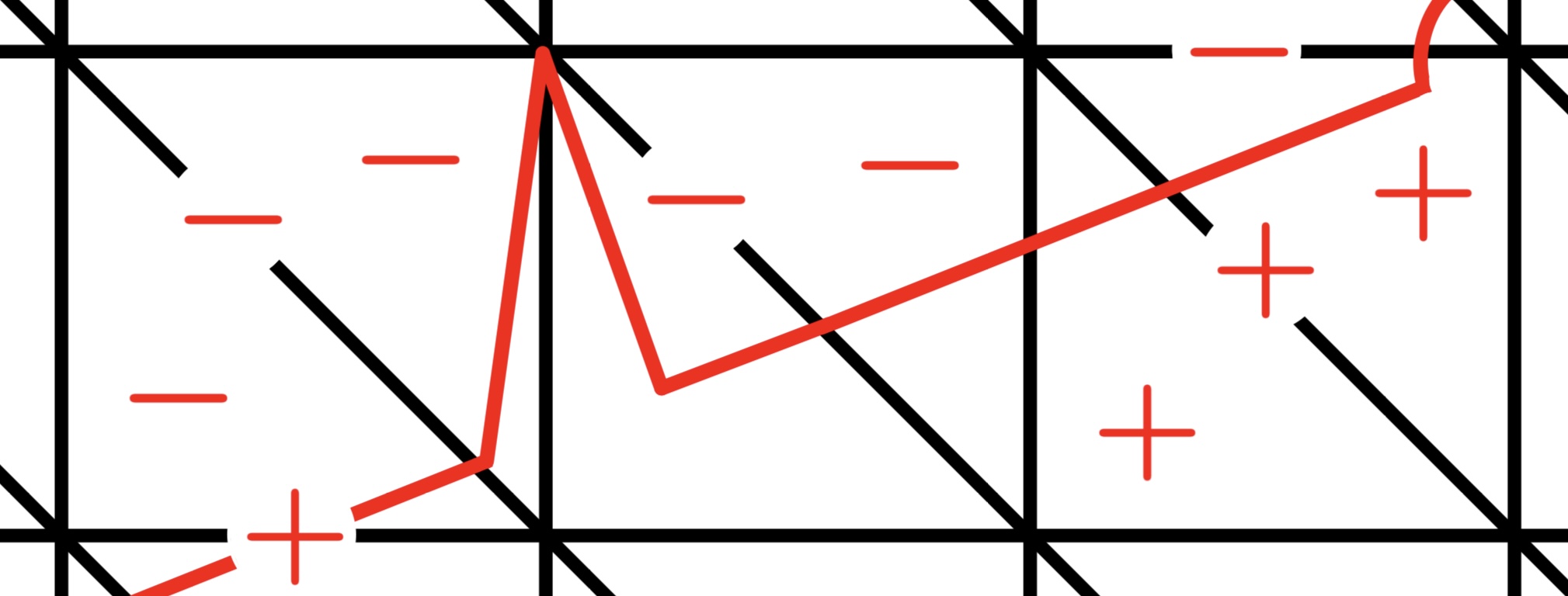}
    \caption{Example of modification from $\overline{L_t}(w)$ to $\widetilde{L_t}(w)$}
    \label{fig:modification}
\end{figure}

The unmodified interior is a subarc of a generic perturbation of the rational line, so it never crosses the same unit edge twice consecutively. Since each endpoint disk meets that subarc in only one connected piece, a repeated pair can occur only when an endpoint branch crosses the same unit edge as the adjacent interior branch. These are exactly the repeated-edge configurations displayed in Table~\ref{table1}. Replacing such a branch by a generic direct endpoint connection on the side displayed in the table removes the repeated crossing, avoids every edge midpoint, and gives the generalized arc in (2).

\begin{figure}[ht]
    \centering
    \includegraphics[scale=0.05]{modification2.jpg}\hspace{3mm}
    \raisebox{1.2em}{$\mapsto$}\hspace{3mm}
    \includegraphics[scale=0.05]{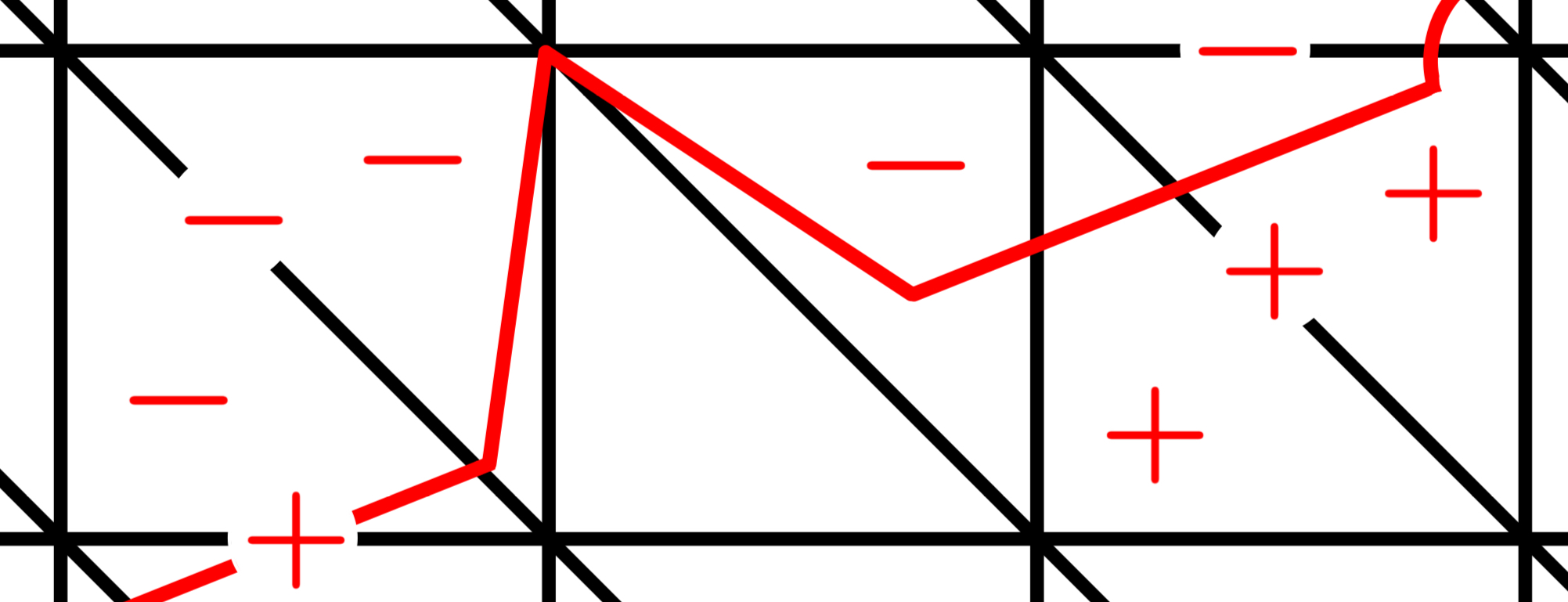}
    \caption{Modification from $\widetilde{L_t}(w)$ to $\widetilde{L_t}(w)_{-}$}
    \label{fig:modification2}
\end{figure}

This deletion affects only the first or last run. It shortens the affected run, and it may remove it completely, which gives the sequence in (3). Since $t\neq1$, Proposition~\ref{rem:difference-(0,1)(1,infty)} (4)--(7) gives $n\geq4$, so the intervening block is nonempty. The continuant recurrence, applied from the left and from the right, shows that the numerator of a positive continued fraction is strictly increasing when either endpoint entry is increased, and that adjoining a positive endpoint entry also strictly increases it. Hence at least one strict endpoint shortening gives
\[
m\!\left(\mathcal G[b_{k+1},a_{k+2},\ldots,a_{k+n-2},b_{k+n-1}]\right)
<m(\mathcal G[w]),
\]
with zero endpoint entries omitted.

Finally, $w_0$ is obtained by cutting the cyclic word immediately after the distinguished first run of length $a_1$. The endpoint branches in the corresponding row of Table~\ref{table1} join the two ends on the left of the oriented segment and introduce no repeated edge. They are precisely the endpoint pieces of $\gamma^L_{AB}$. This proves (4).
\end{proof}

\begin{proof}[Proof of Theorem \ref{thm:markov-value-gen}]
When $t=\frac{0}{1}$ or $\frac{1}{0}$, the conclusion follows directly from Lemma \ref{lem:trace}. We assume that $t\in (0,\infty)$. By Proposition \ref{rem:difference-(0,1)(1,infty)} (1), the number of entries of $s(t)$ is even. Hence Theorems \ref{thm:Mt-description} and \ref{continued-fraction-theorem2}, together with Lemma \ref{lem:trace}, give
\[
\ell(^\infty s(t)|s(t)^\infty)=\frac{\sqrt{((3+k_1+k_2+k_3)n_t-k_t)^2-4}}{n_t}.
\]
Combining this fact with Proposition \ref{rem:dual-remark} and Theorem \ref{thm:Mt-description}, we deduce that
\[
\ell(^\infty s(t)|s(t)^\infty)=\ell(^\infty s^\ast(\tfrac{1}{t})|s^\ast(\tfrac{1}{t})^\infty).
\]
By Proposition~\ref{rem:difference-(0,1)(1,infty)} (8), the two periodic sequences are reverses of one another. Reversal exchanges the two one-sided continued fractions at the corresponding cut and therefore preserves their sum; in particular, it preserves $L$. Consequently, once the designated cut of ${}^\infty s(t)^\infty$ is shown to be maximal, the same is automatically true for the displayed cut of the dual periodic sequence.
Thus it remains to prove that
\[
L({}^\infty s(t)^\infty)=\ell(^\infty s(t)|s(t)^\infty),
\]
that is, $\ell(^\infty s(t)|s(t)^\infty)$ is maximal among all values $\ell(P^\ast|Q)$ for cuts $P^\ast|Q$ of $^\infty s(t)^\infty$.
We set $s(t)=(a_1,a_2,\dots,a_n)$. 
If $t=1$, then Proposition~\ref{rem:difference-(0,1)(1,infty)} (3) gives
\[
s(1)=(a_1,a_2)=(2+k_1+k_2+k_3,\,2+k_{\sigma(2)}),
\qquad a_2\leq a_1.
\]
There are only two cuts. Their continued-fraction matrices have the same trace, while their $(2,1)$-entries are respectively $a_2$ and $a_1$. Lemma~\ref{lem:trace} therefore shows that the displayed cut, whose denominator is $a_2$, is maximal. This proves the theorem for $t=1$. Henceforth assume $t\neq1$; then $n\geq4$.
The possible cuts of $^{\infty}s(t)^{\infty}$ are given by
\begin{align*}
\dots a_1,a_2\dots a_{n-1},a_n&\;|\;a_1,a_2,\dots,a_{n-1},a_n\dots,\\
\dots a_2,a_3\dots a_n,a_1&\;|\;a_2,a_3,\dots,a_n,a_1\dots,\\
&\vdots\\
\dots a_n,a_1\dots a_{n-2},a_{n-1}&\;|\;a_n,a_1,\dots,a_{n-2},a_{n-1}\dots,
\end{align*}
therefore there are $n$ such cuts. Since $n$ is even by Proposition \ref{rem:difference-(0,1)(1,infty)} (1), Lemma \ref{lem:trace} shows that it suffices to prove that $\ell(^\infty s(t)|s(t)^\infty)$ is the maximum value in
\[F_t=\left\{\frac{\sqrt{(\mathrm{tr}(CF_S))^2-4}}{(CF_{S})_{21}}\ \middle|\ S=(a_{k}\dots,a_{k+n-1}) \text{ for some $k$ }\right\},\]
where indices of $a_i$ are taken modulo $n$. The matrices $CF_{(a_i,a_{i+1},\dots,a_{i-1})}$ are related by conjugation, hence all have the same trace. Therefore, the maximum value in $F_t$ has the minimal denominator $(CF_{(a_i,a_{i+1},\dots,a_{i-1})})_{21}$. By \eqref{eq:continued-fraction-result}, we have
\[(CF_{(a_i,a_{i+1},\dots,a_{i+n-1})})_{21}=m(\mathcal G[a_{i+1},\dots,a_{i+n-1}]).\]
Therefore, to prove the claim, it suffices to show that
$m(\mathcal G[a_2,\dots,a_n])$ is the minimal element in
\[
N_t=\{m(\mathcal G[w]) \mid w=(a_{k+1},a_{k+2},\dots,a_{k+n-1}) \ \text{for some $k$}\ \}.
\]
Write $t=p/q$ in lowest terms and set $A=(0,0)$ and $B=(q,p)$. By Lemma~\ref{lem:endpoint-modification} (4) and Theorem~\ref{thm:distance-theorem},
\[
m(\mathcal G[a_2,\dots,a_n])=|\gamma^L_{AB}|=d(A,B).
\]
For an arbitrary cyclic choice $w=(a_{k+1},\dots,a_{k+n-1})$, if $\widetilde{L_t}(w)$ is already a generalized arc, then
\[
m(\mathcal G[a_2,\dots,a_n])=d(A,B)
\leq |\widetilde{L_t}(w)|=m(\mathcal G[w]).
\]
Here the last equality is Lemma~\ref{lem:endpoint-modification} (1). Otherwise, Lemma~\ref{lem:endpoint-modification} (2)--(3) and the distance theorem give
\[
m(\mathcal G[a_2,\dots,a_n])
\leq m\!\left(\mathcal G[b_{k+1},a_{k+2},\dots,a_{k+n-2},b_{k+n-1}]\right)
<m(\mathcal G[w]).
\]
Thus $m(\mathcal G[a_2,\dots,a_n])$ attains the minimal value in $N_t$, which completes the proof.
\end{proof}
This gives the following corollary, which identifies the corresponding Lagrange constant and Markov constant.
\begin{corollary}\label{cor:markov-lagrange}
Fix $(k_1,k_2,k_3)\in \mathbb Z_{\geq 0}^3$ and $\sigma\in \mathfrak S_3$.  
For any irreducible fraction $t\in [0,\infty]$, let $(n_t,i_t)$ denote the corresponding $(k_1,k_2,k_3,\sigma)$-GM number-position pair, and let $s(t)$ be the corresponding generalized strongly admissible sequence. Then we have \[[s(t)^\infty]\in \mathbb Q[\sqrt{((3+k_1+k_2+k_3)n_t-k_t)^2-4}]\] and
\[
\mathcal L([s(t)^\infty])=\mathcal M(Q_{s(t)})=\frac{\sqrt{((3+k_1+k_2+k_3)n_t-k_t)^2-4}}{n_t}.
\]
\end{corollary}
\begin{proof}
We first show that
\[
[s(t)^\infty]\in \mathbb Q[\sqrt{((3+k_1+k_2+k_3)n_t-k_t)^2-4}].
\]
Let $\alpha=[s(t)^\infty]$. Then
\[\alpha=[s(t),\alpha].\]
By the general theory of continued-fraction matrices (see \cite{bombieri2}*{Appendix A}), we have
\[\alpha=CF_{s(t)}(\alpha):=\dfrac{a\alpha+b}{c\alpha+d},\]
where $CF_{s(t)}=\begin{bmatrix}a&b\\c&d
\end{bmatrix}$, and thus
\[c\alpha^2-(a-d)\alpha-b=0.\]Therefore, we have
\[\alpha=\dfrac{a-d\pm\sqrt{(a-d)^2+4bc}}{2c}=\dfrac{a-d\pm\sqrt{(a+d)^2-4(ad-bc)}}{2c}.\]
Since $\alpha>0$, we have 
\[\alpha=\dfrac{a-d+\sqrt{(a+d)^2-4(ad-bc)}}{2c}.\]
When $t=\frac{0}{1}$, the conclusion can be checked directly. We assume that $t\in (0,\infty].$ 
By Theorems \ref{thm:Mt-description} and \ref{continued-fraction-theorem2}, we have
\[\alpha=\dfrac{(3+k_1+k_2+k_3)n_t - k_t - 2u_t+\sqrt{((3+k_1+k_2+k_3)n_t - k_t)^2-4}}{2n_t}.\] Therefore, we have
\[\alpha\in \mathbb Q[\sqrt{((3+k_1+k_2+k_3)n_t-k_t)^2-4}].\]
The asserted equality of the Lagrange constant and the Markov constant follows from Theorems \ref{thm:markov-value-gen}, \ref{thm:lagrange-quadratic} and \ref{thm:markov-quadratic}.
\end{proof}
\begin{example}
Set $(k_1,k_2,k_3,\sigma)=(1,2,0,\textrm{id})$ and $t=\tfrac{2}{5}$. Then $i_\frac{2}{5}=1$ and
\[
s\!\left(\tfrac{2}{5}\right)=(5,1,3,3,1,5,4,1,3,4).
\]
Thus the sequences $w$ that contribute to $m(\mathcal G[w])\in N_{\frac25}$ are
\begin{align*}
&(1,3,3,1,5,4,1,3,4),\ (3,3,1,5,4,1,3,4,5),\ (3,1,5,4,1,3,4,5,1),
\ (1,5,4,1,3,4,5,1,3),\\ 
&(5,4,1,3,4,5,1,3,3),\ (4,1,3,4,5,1,3,3,1),\ (1,3,4,5,1,3,3,1,5),\ (3,4,5,1,3,3,1,5,4),\\ 
&(4,5,1,3,3,1,5,4,1),\ (5,1,3,3,1,5,4,1,3).
\end{align*}
The corresponding values $m(\mathcal G[w])$ are
\[
8227,\ 32957,\ 12039,\ 12041,\ 32937,\ 8261,\ 9997,\ 31881,\ 12199,\ 11127,
\]
and the minimal number in $N_{\frac25}$ is $8227$. The corresponding arcs $\overline{L_t}(w)$ and $\widetilde{L_t}(w)$ are as in Table \ref{table2}. Since
\[CF_{s(\frac{2}{5})}=\begin{bmatrix}
    m(\mathcal G[5,1,3,3,1,5,4,1,3,4])&m(\mathcal G[5,1,3,3,1,5,4,1,3])\\
    m(\mathcal G[1,3,3,1,5,4,1,3,4])&m(\mathcal G[1,3,3,1,5,4,1,3])
\end{bmatrix}=\begin{bmatrix}
    47431&11127\\8227&1930
\end{bmatrix},\]
we have
\[[s(\tfrac{2}{5})^\infty]=\dfrac{\sqrt{2436508317}+45501}{16454},\quad L({}^\infty s(\tfrac{2}{5})^\infty)=\dfrac{\sqrt{2436508317}}{8227}.\]
Moreover, we have
\[Q_{s(\frac{2}{5})}=x^2 - \frac{45501}{8227}xy - \frac{11127}{8227}y^2,\]
and thus
\[\mathcal L\left(\frac{\sqrt{2436508317}+45501}{16454}\right)=\mathcal M\left(x^2 - \frac{45501}{8227}xy - \frac{11127}{8227}y^2\right)=\dfrac{\sqrt{2436508317}}{8227}.\]
\begin{table}[ht]
\centering
\begin{tabular}{|>{\centering\arraybackslash}m{4cm}|>{\centering\arraybackslash}m{4cm}|>{\centering\arraybackslash}m{4cm}|>{\centering\arraybackslash}m{1.5cm}|}
$w$&$\overline{L_t}(w)$& $\widetilde{L_t}(w)$&$m(\mathcal G[w])$\\
\hline
&&&\\[-4mm]
$(1,3,3,1,5,4,1,3,4)$&\includegraphics[scale=0.044]{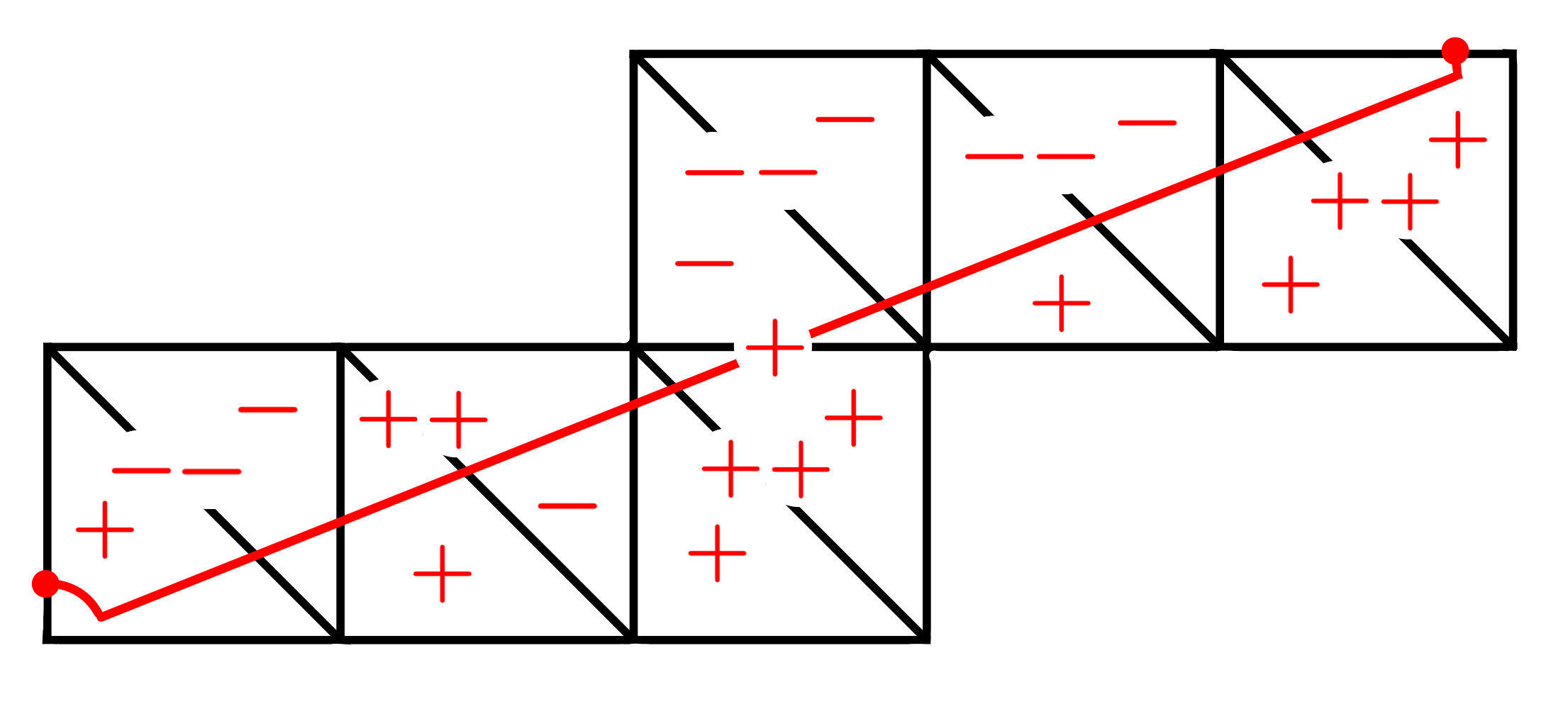} &\includegraphics[scale=0.044]{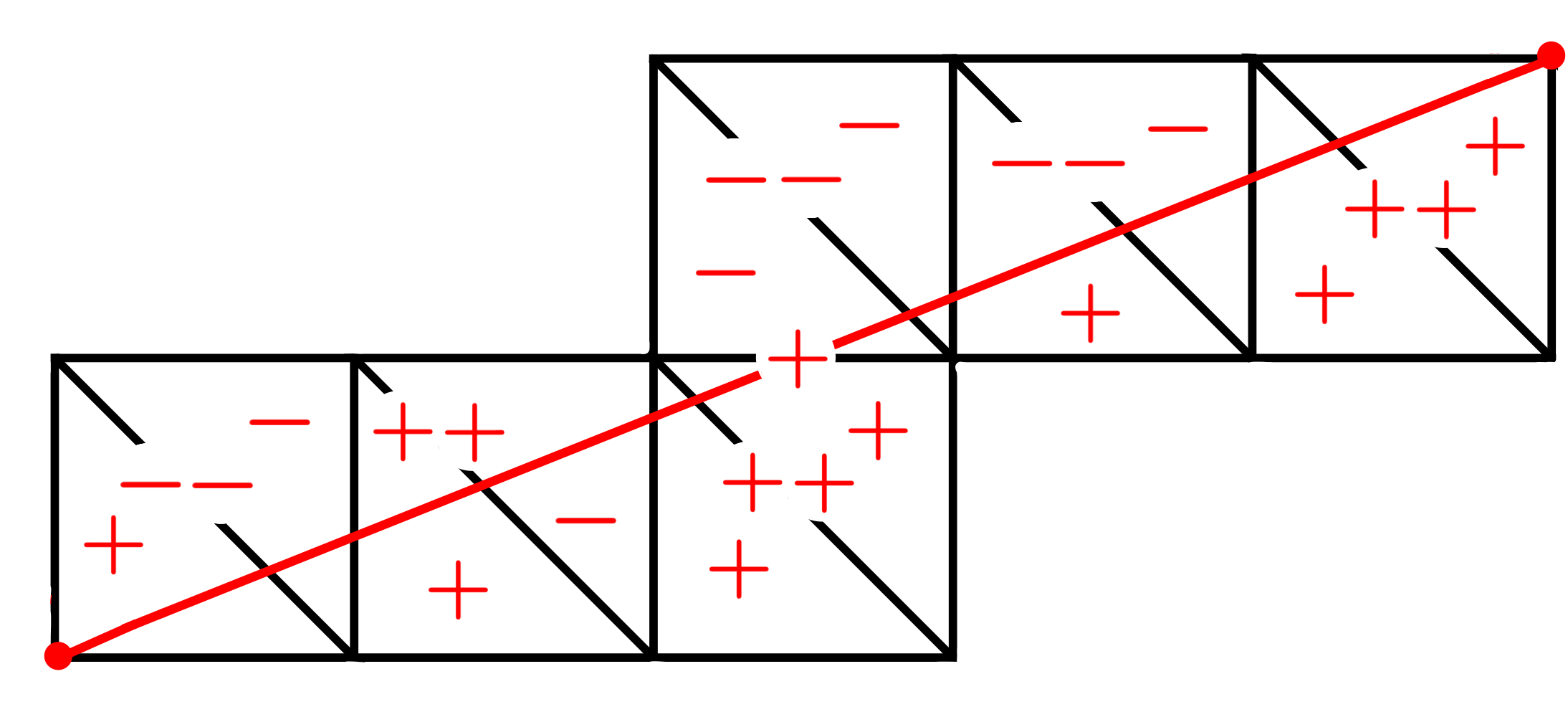}&$8227$\\
\hline
&&&\\[-4mm]
$(3,3,1,5,4,1,3,4,5)$&\includegraphics[scale=0.06]{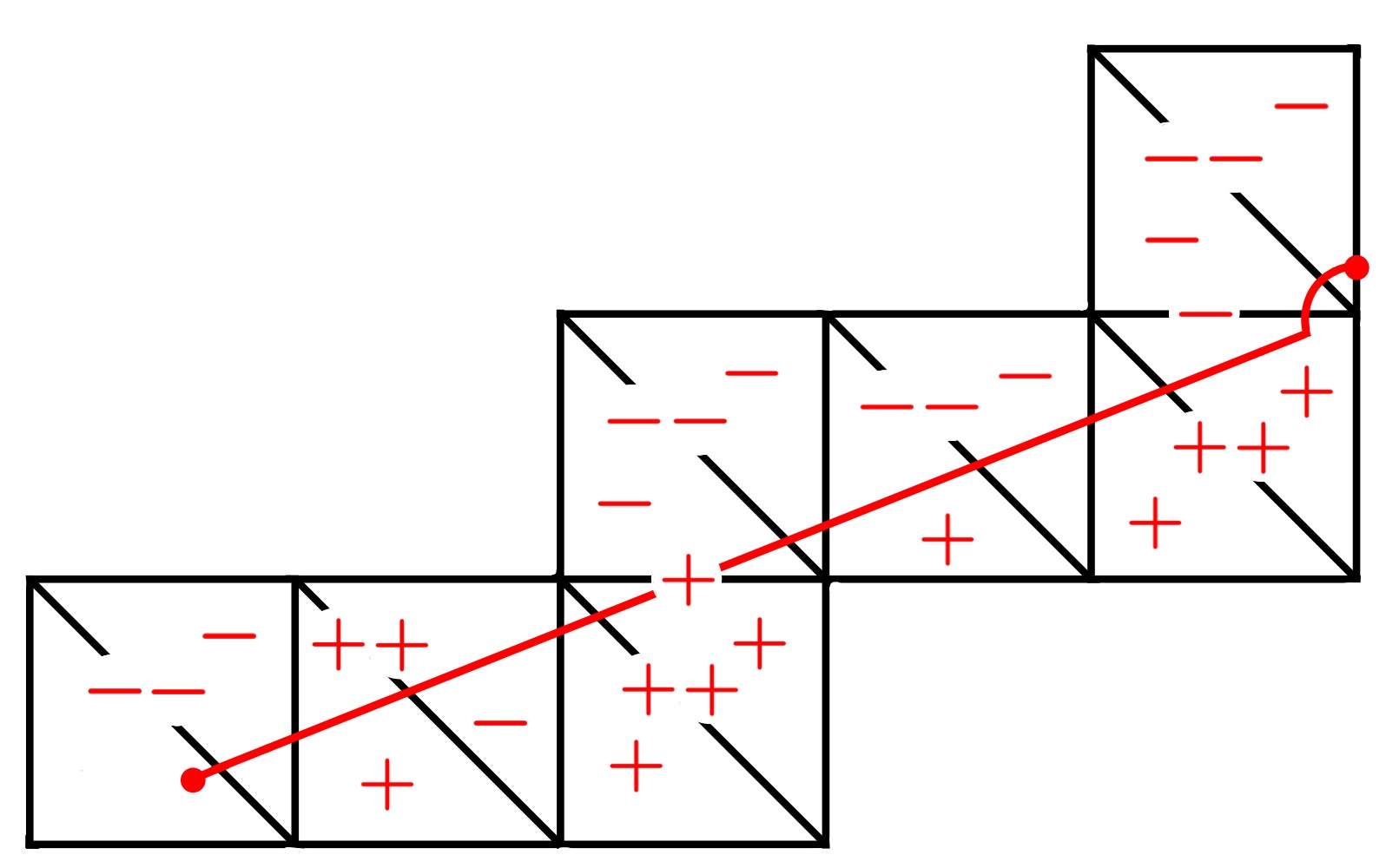} &\includegraphics[scale=0.06]{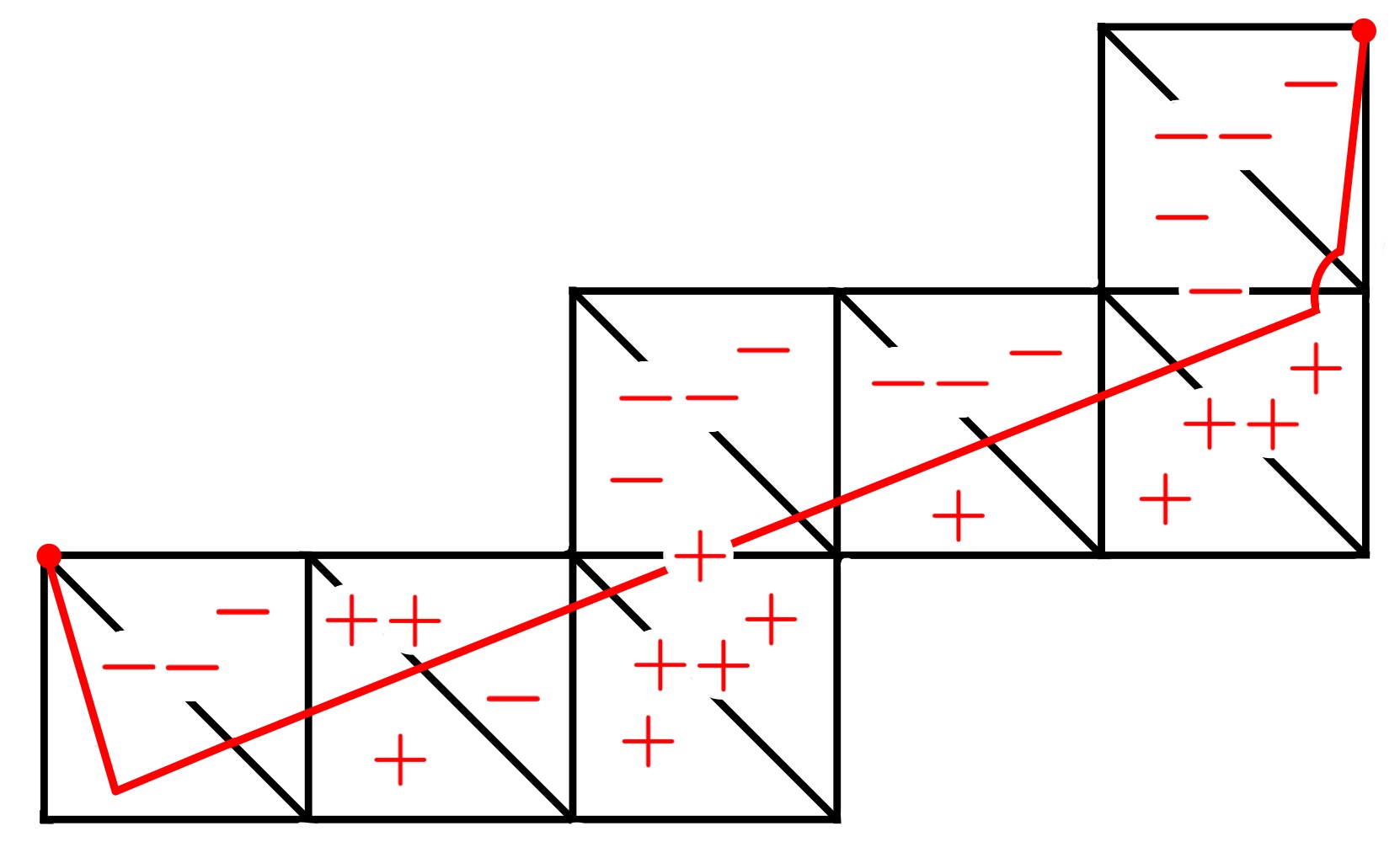}&$32957$\\
\hline
&&&\\[-4mm]
$(3,1,5,4,1,3,4,5,1)$&\includegraphics[scale=0.06]{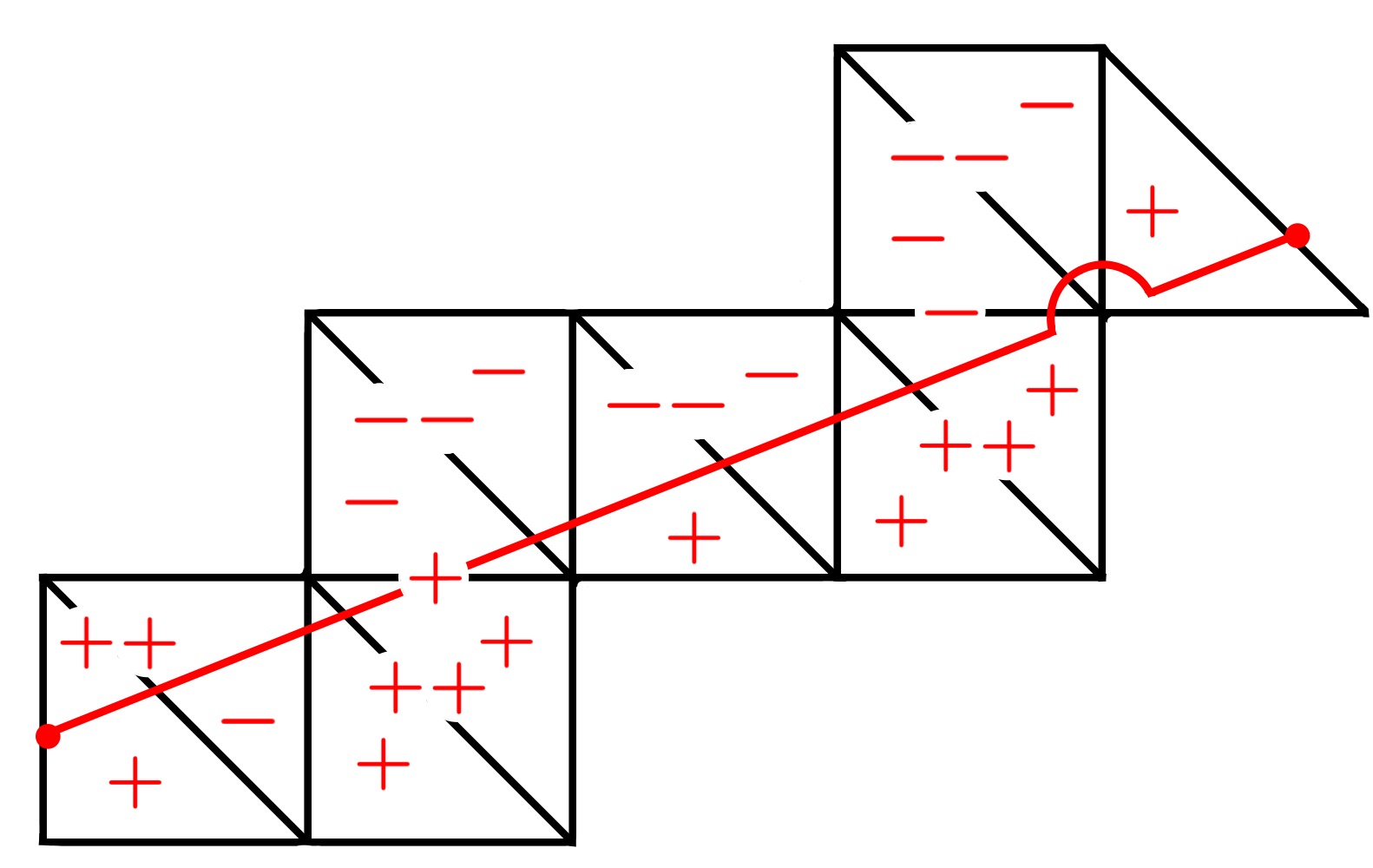} &\includegraphics[scale=0.06]{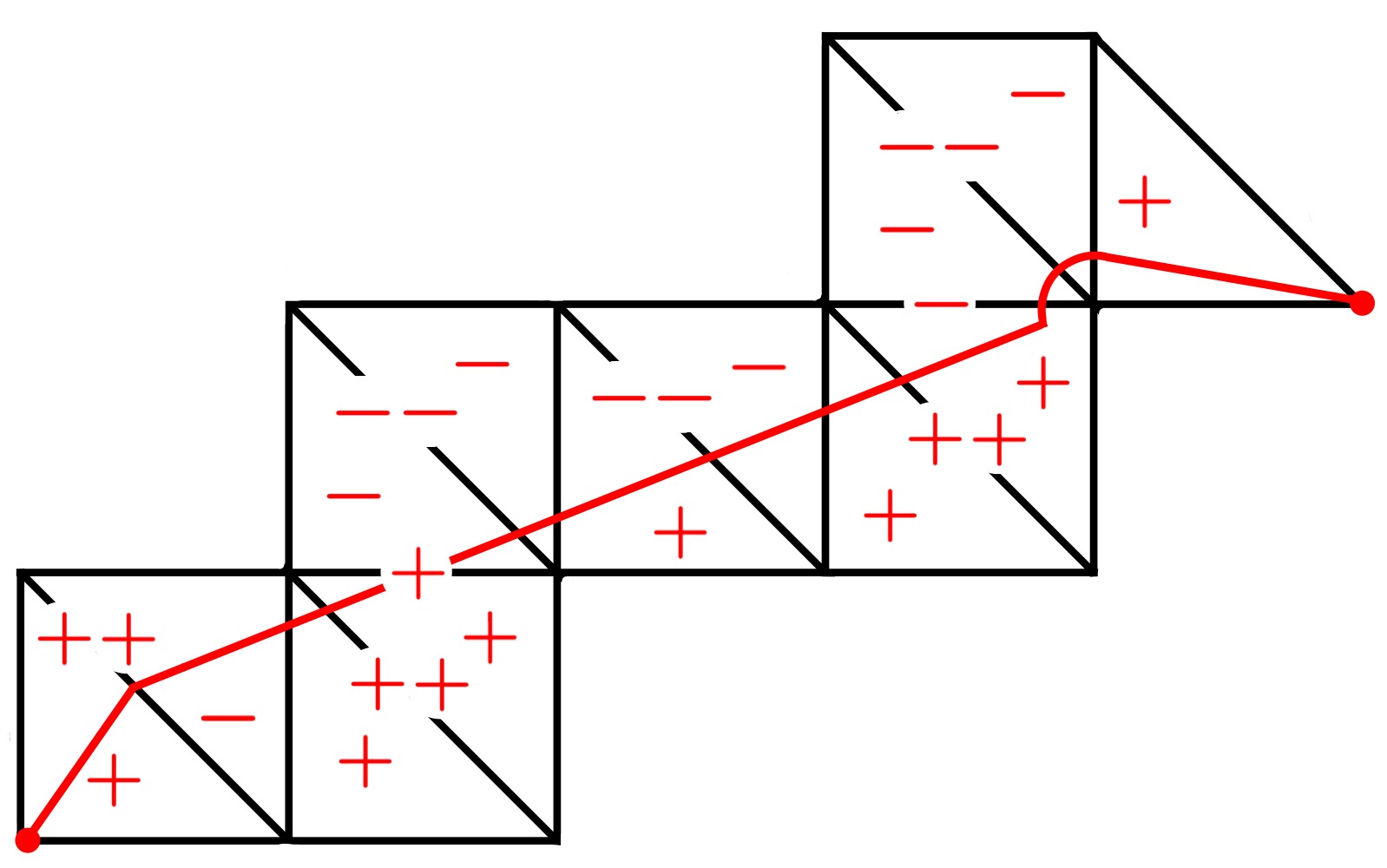}&$12039$\\
\hline
&&&\\[-4mm]
$(1,5,4,1,3,4,5,1,3)$&\includegraphics[scale=0.06]{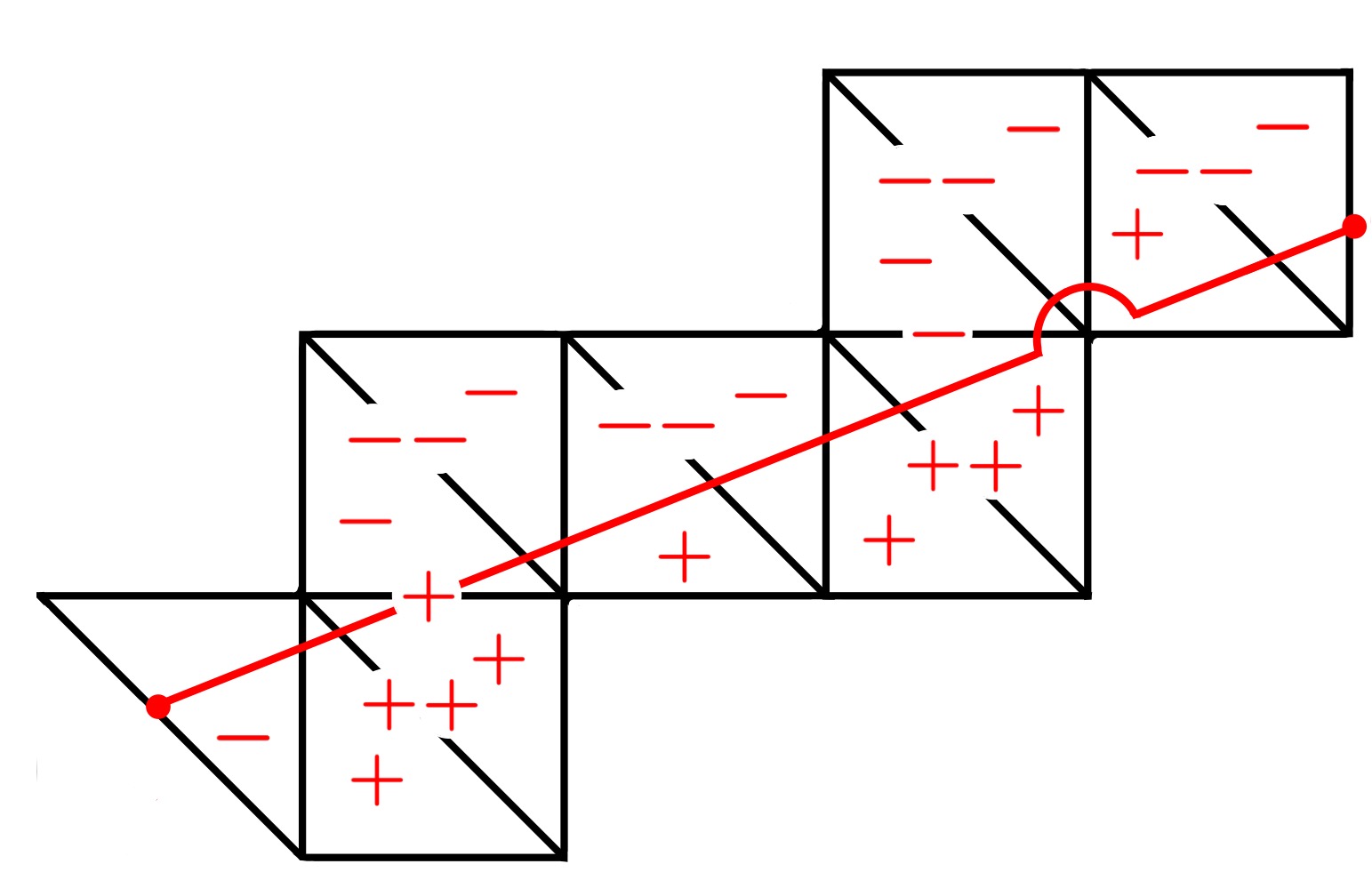} &\includegraphics[scale=0.06]{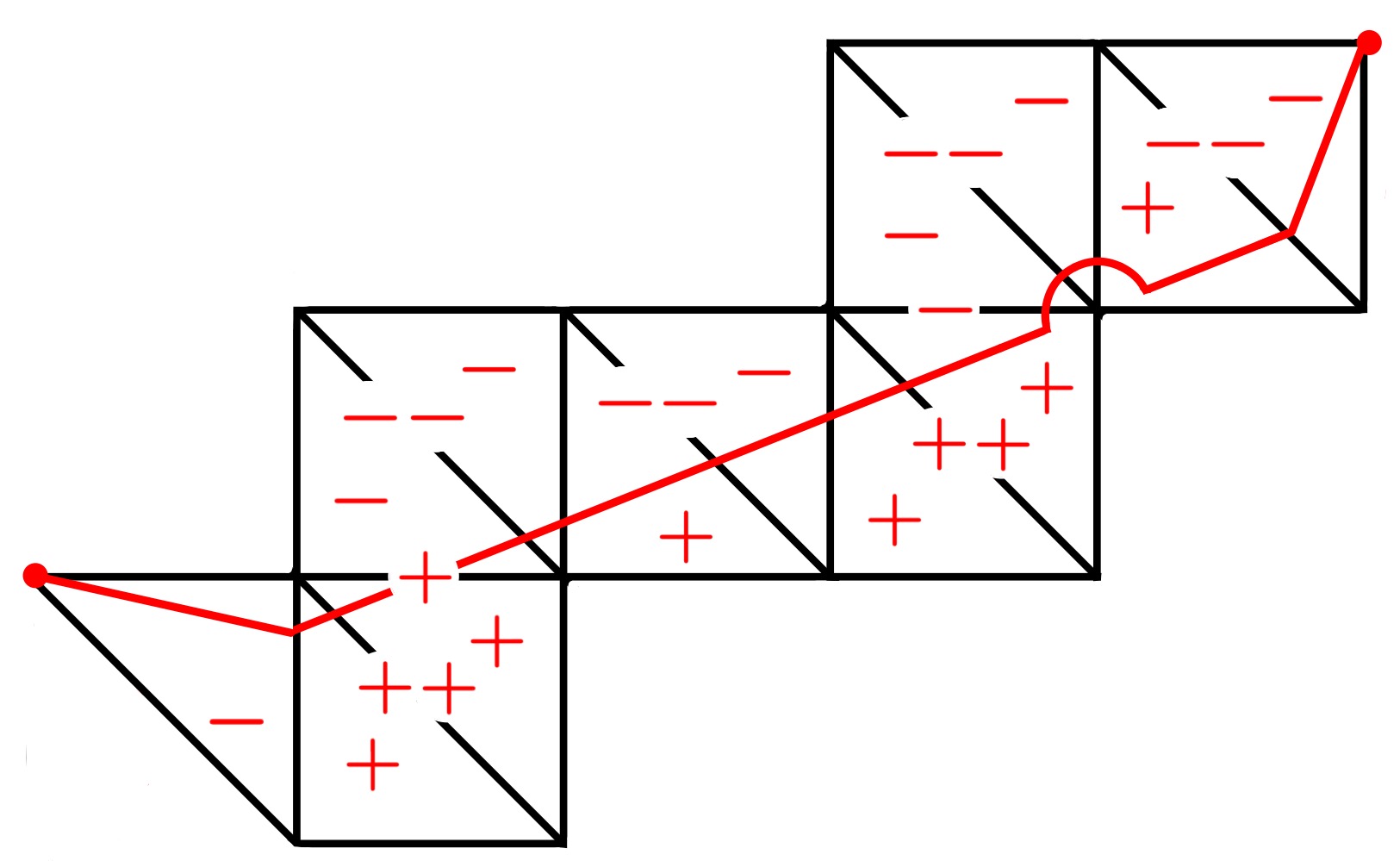}&$12041$\\
\hline
&&&\\[-4mm]
$(5,4,1,3,4,5,1,3,3)$&\includegraphics[scale=0.06]{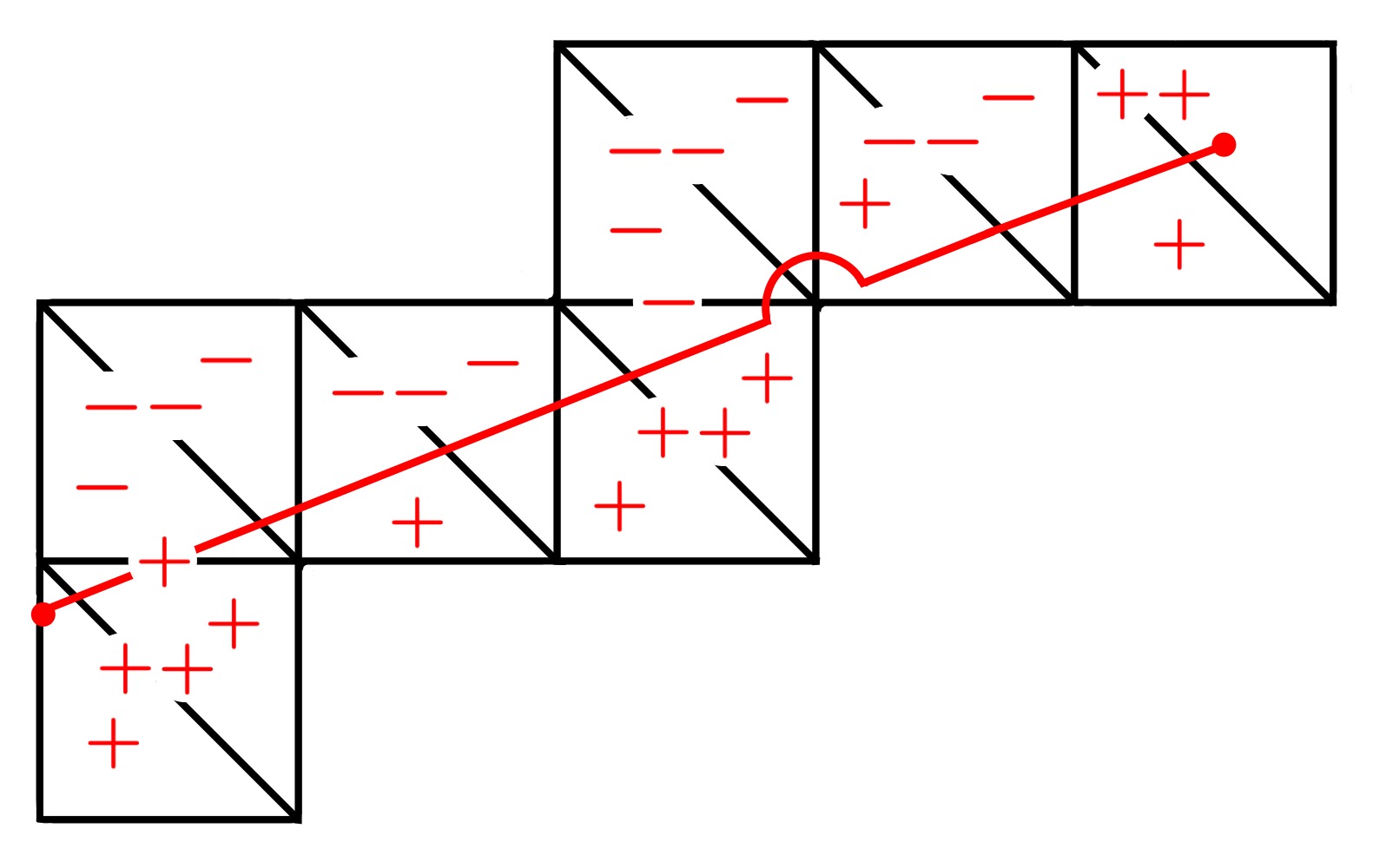} &\includegraphics[scale=0.06]{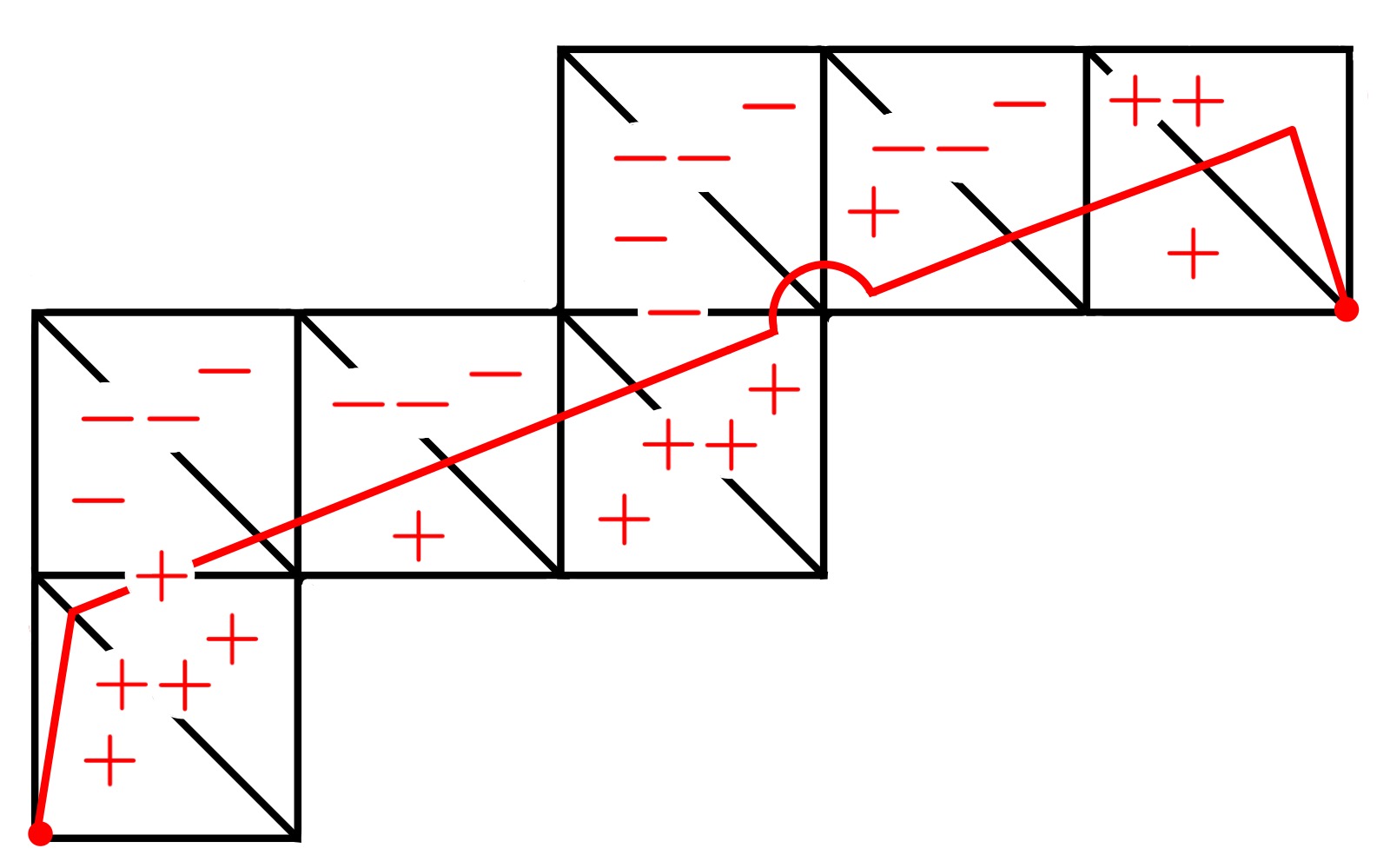}&$32937$\\
\hline
&&&\\[-4mm]
$(4,1,3,4,5,1,3,3,1)$&\includegraphics[scale=0.044]{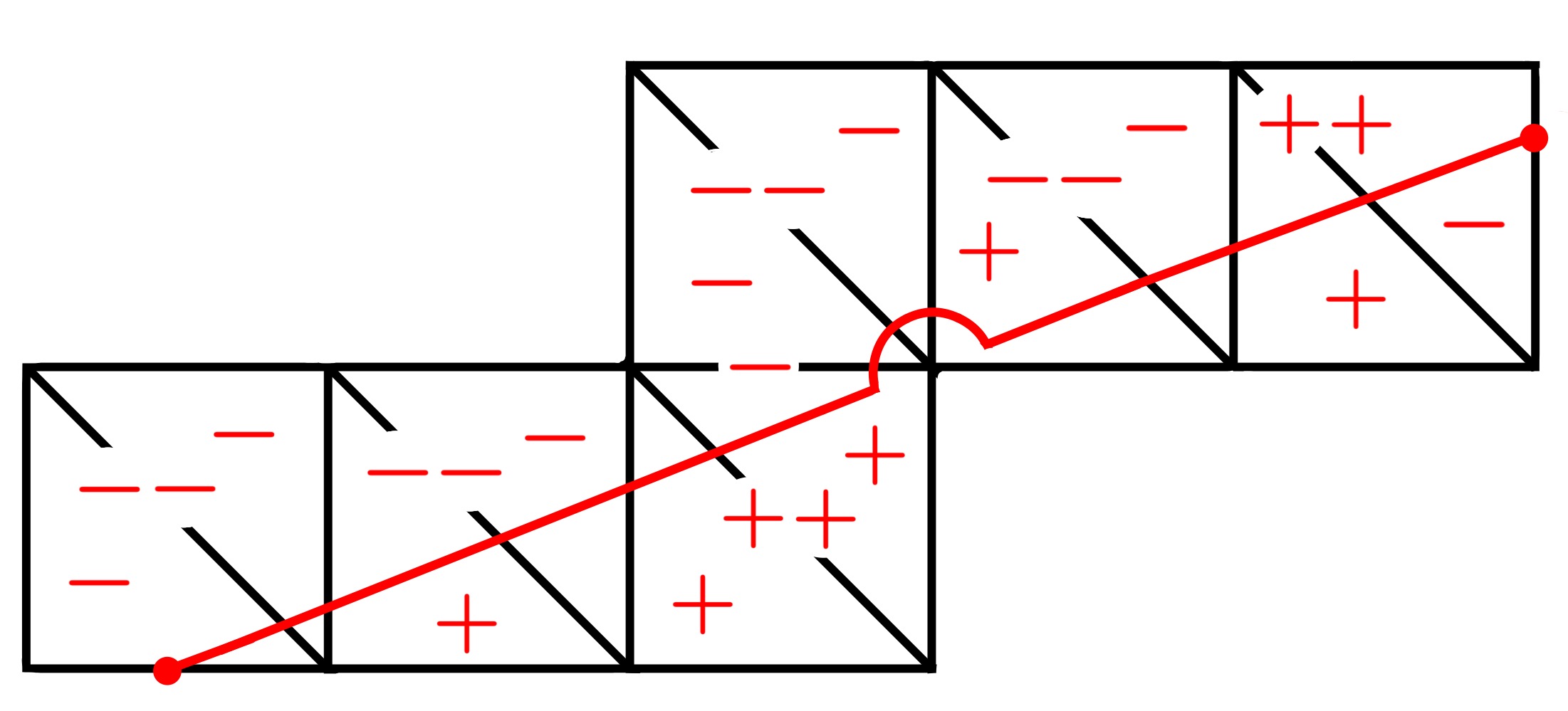} &\includegraphics[scale=0.044]{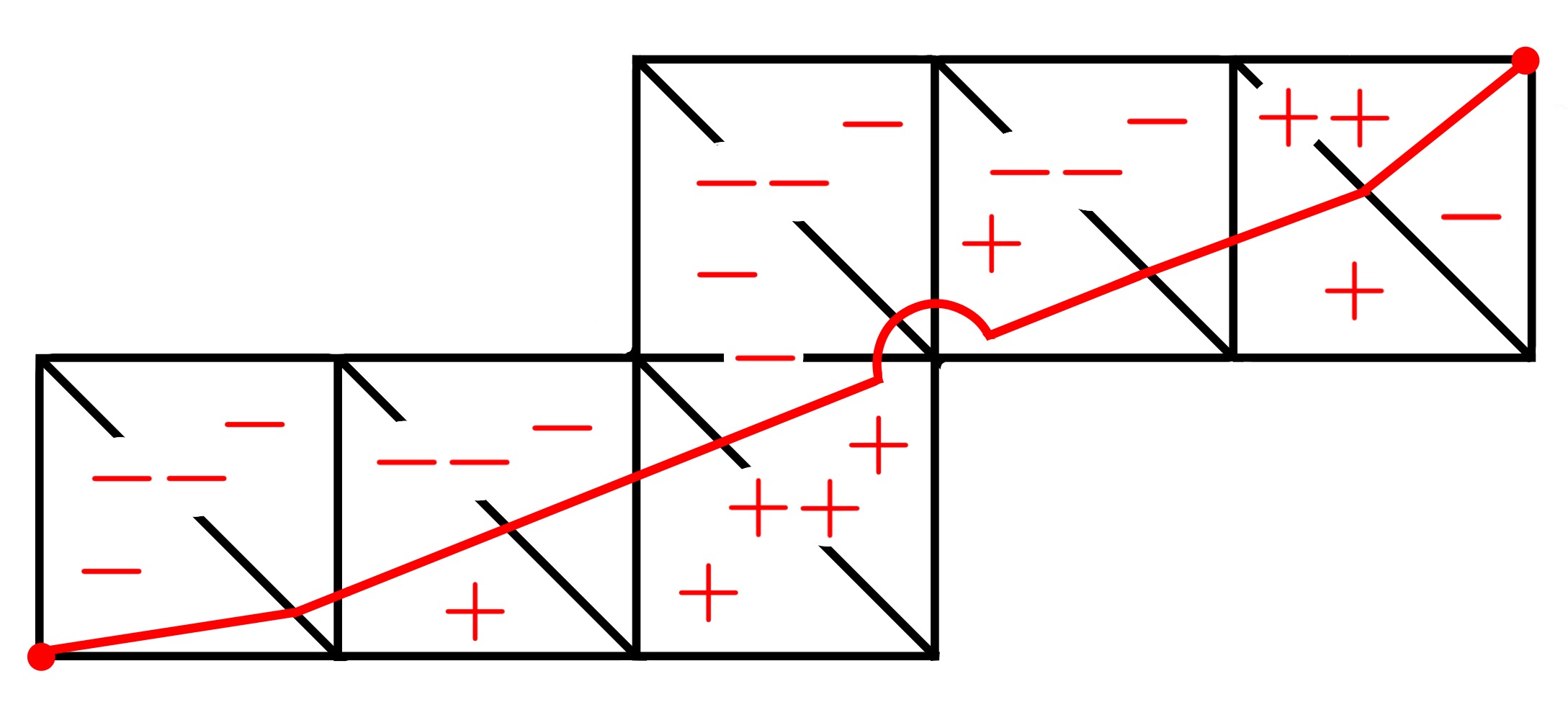}&8261\\
\hline
&&&\\[-4mm]
$(1,3,4,5,1,3,3,1,5)$&\includegraphics[scale=0.06]{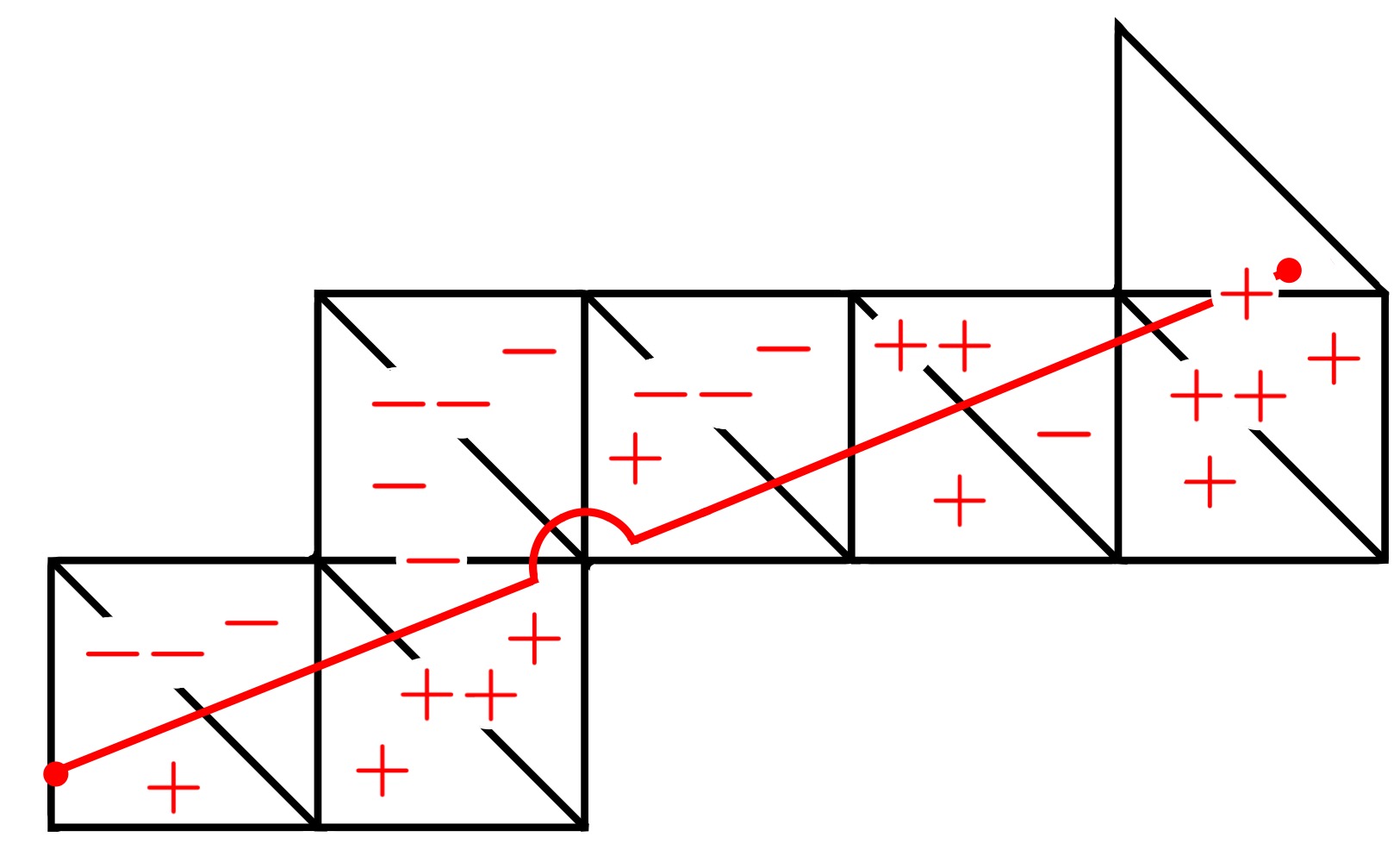} &\includegraphics[scale=0.06]{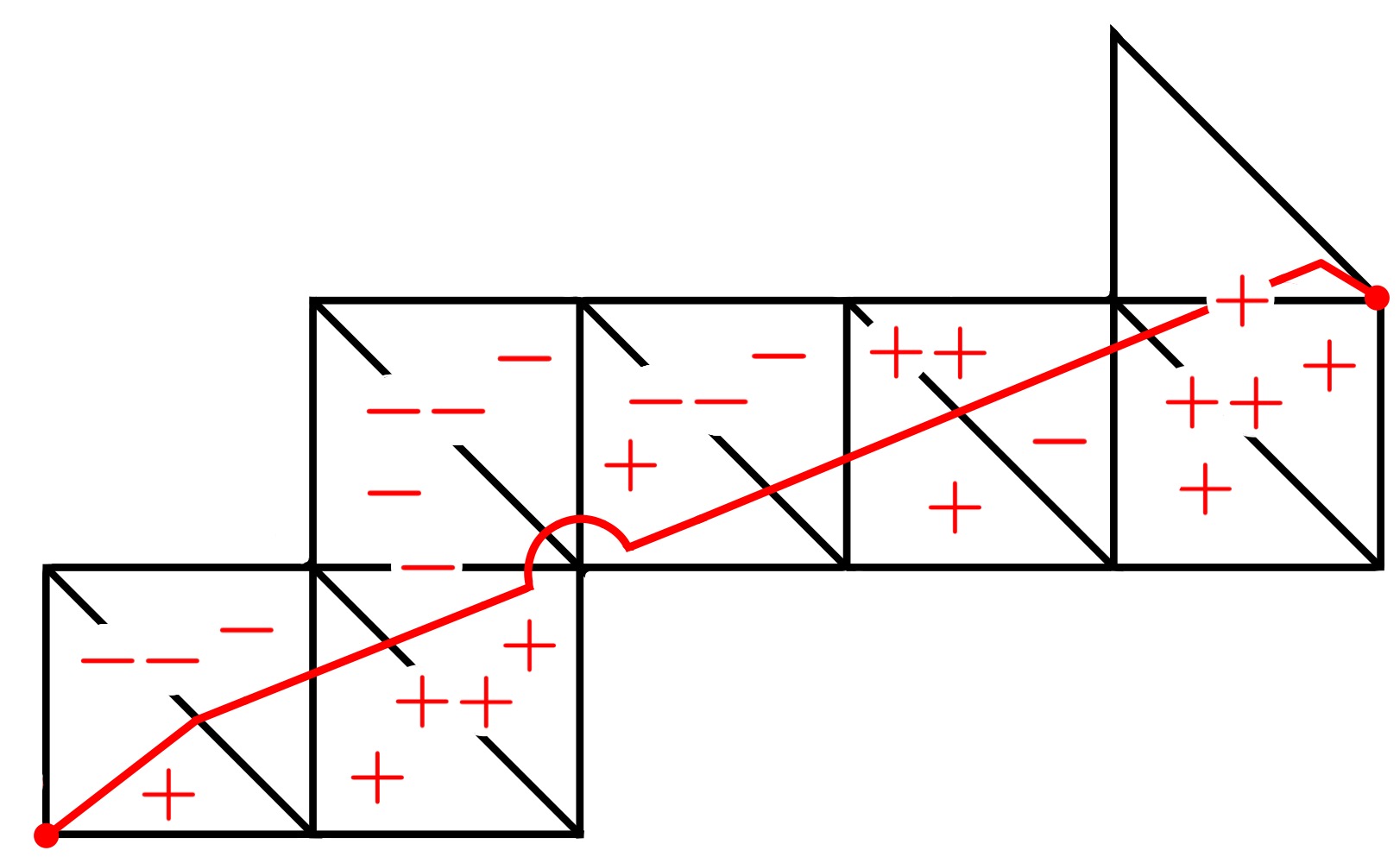}&$9997$\\
\hline
&&&\\[-4mm]
$(3,4,5,1,3,3,1,5,4)$&\includegraphics[scale=0.06]{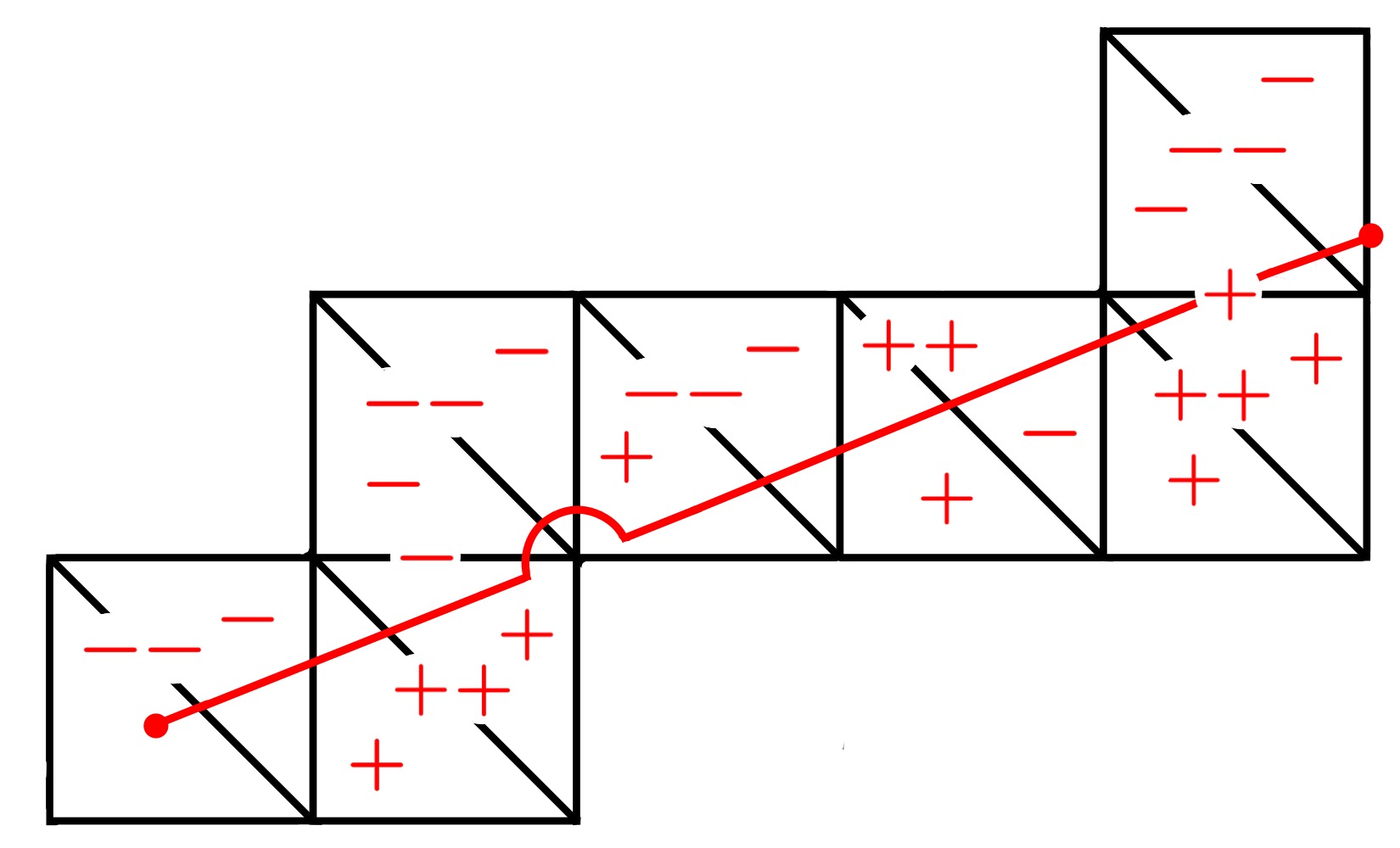} &\includegraphics[scale=0.06]{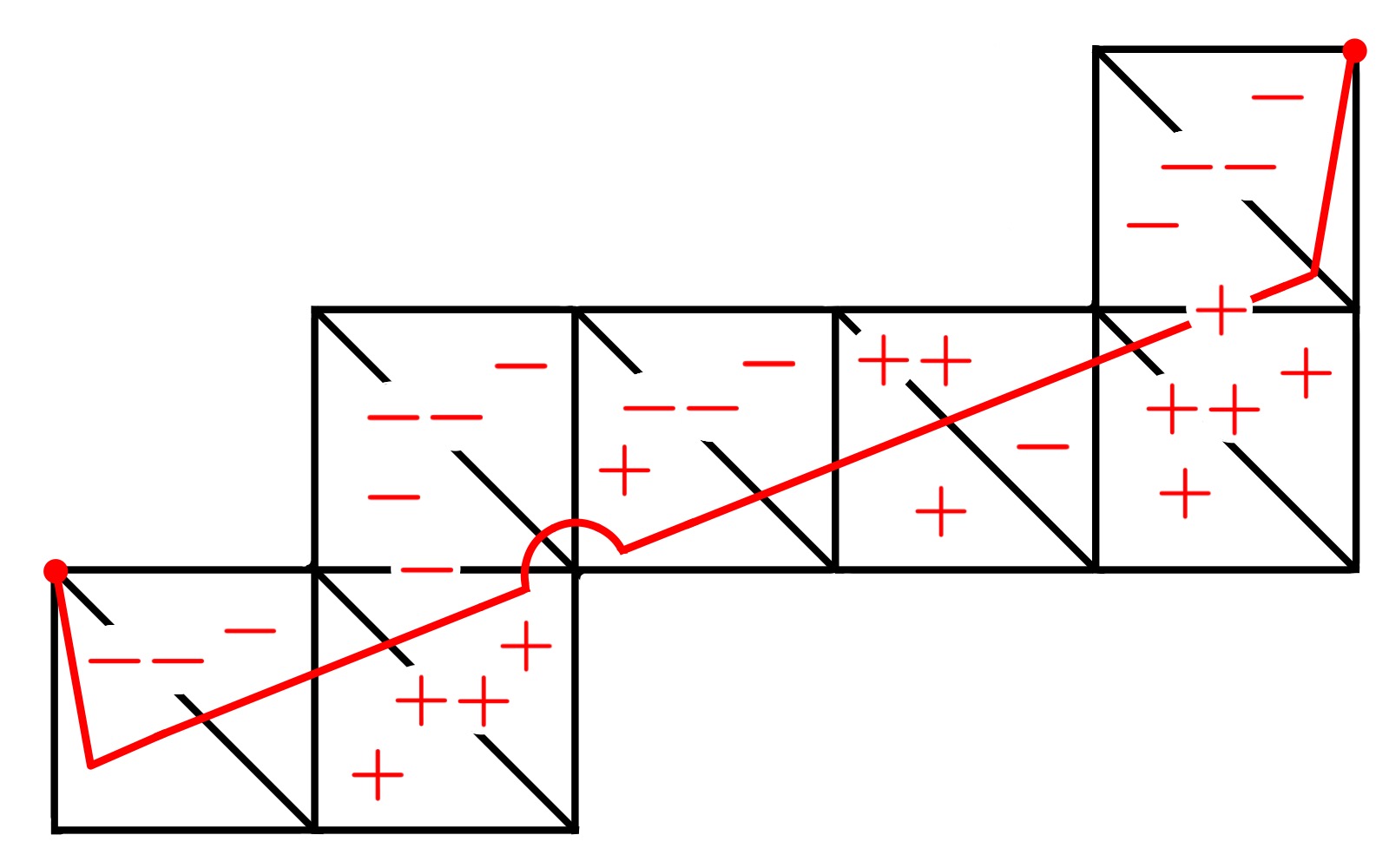}&$31881$\\
\hline
&&&\\[-4mm]
$(4,5,1,3,3,1,5,4,1)$&\includegraphics[scale=0.06]{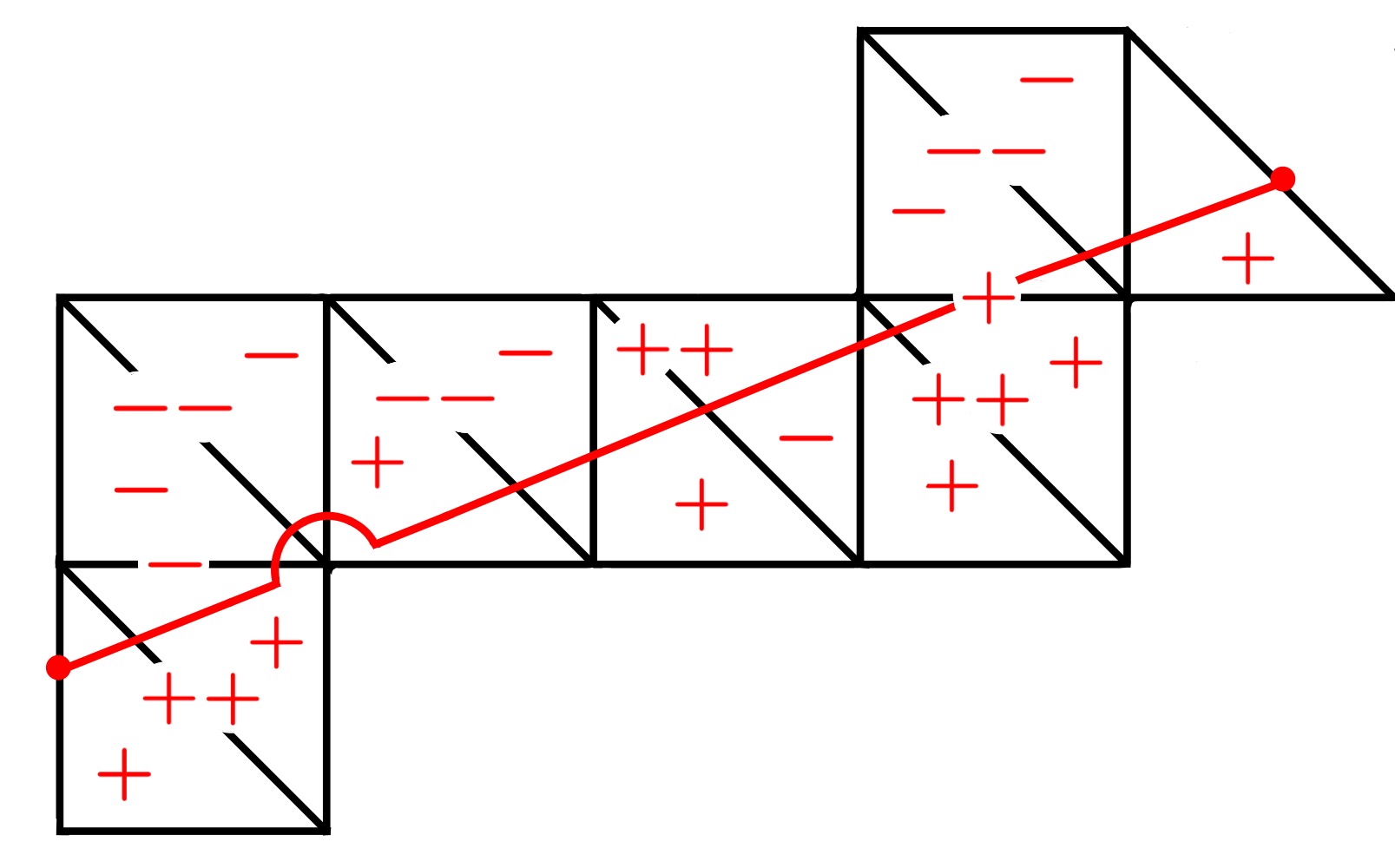} &\includegraphics[scale=0.06]{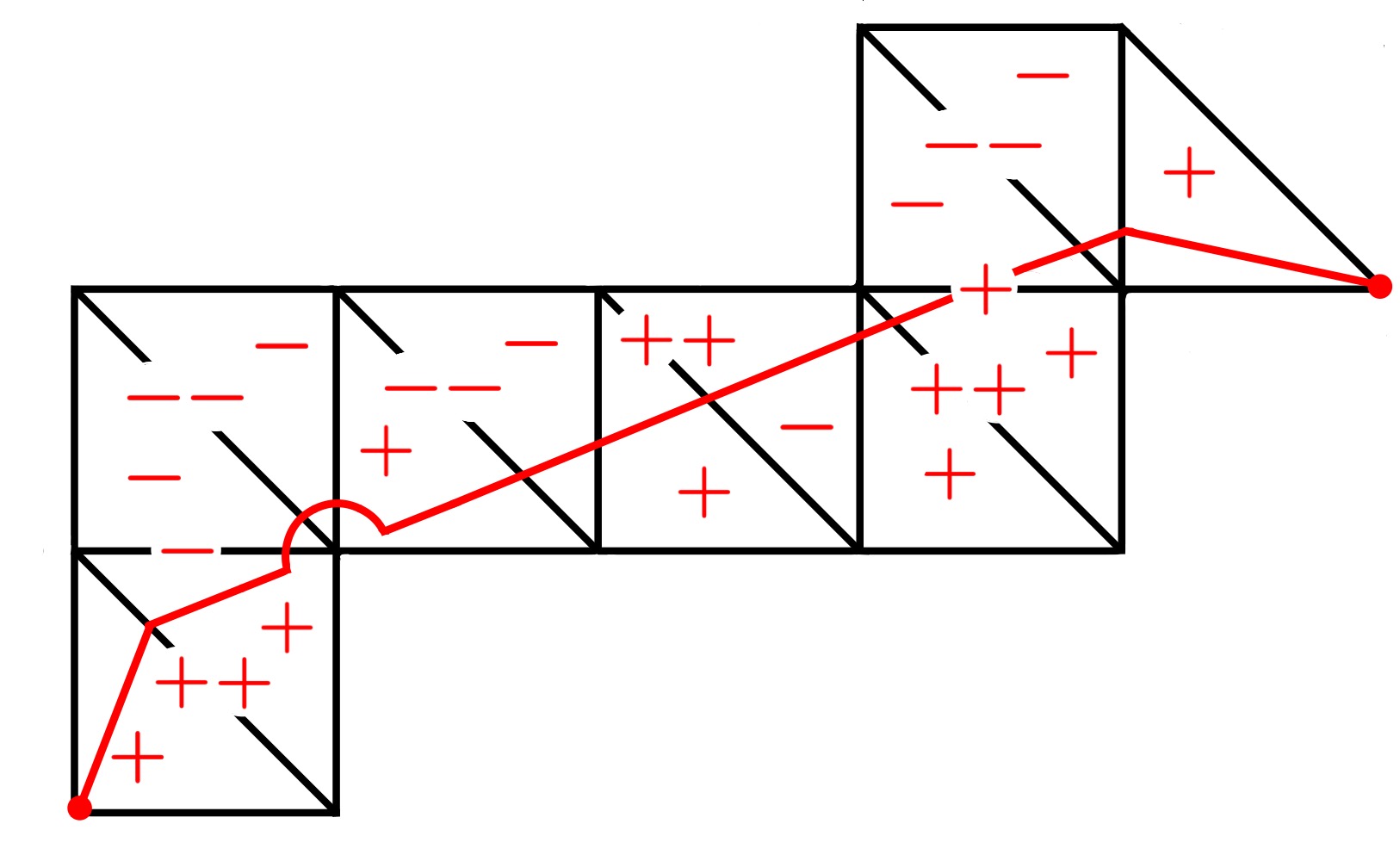}&$12199$\\
\hline
&&&\\[-4mm]
$(5,1,3,3,1,5,4,1,3)$&\includegraphics[scale=0.06]{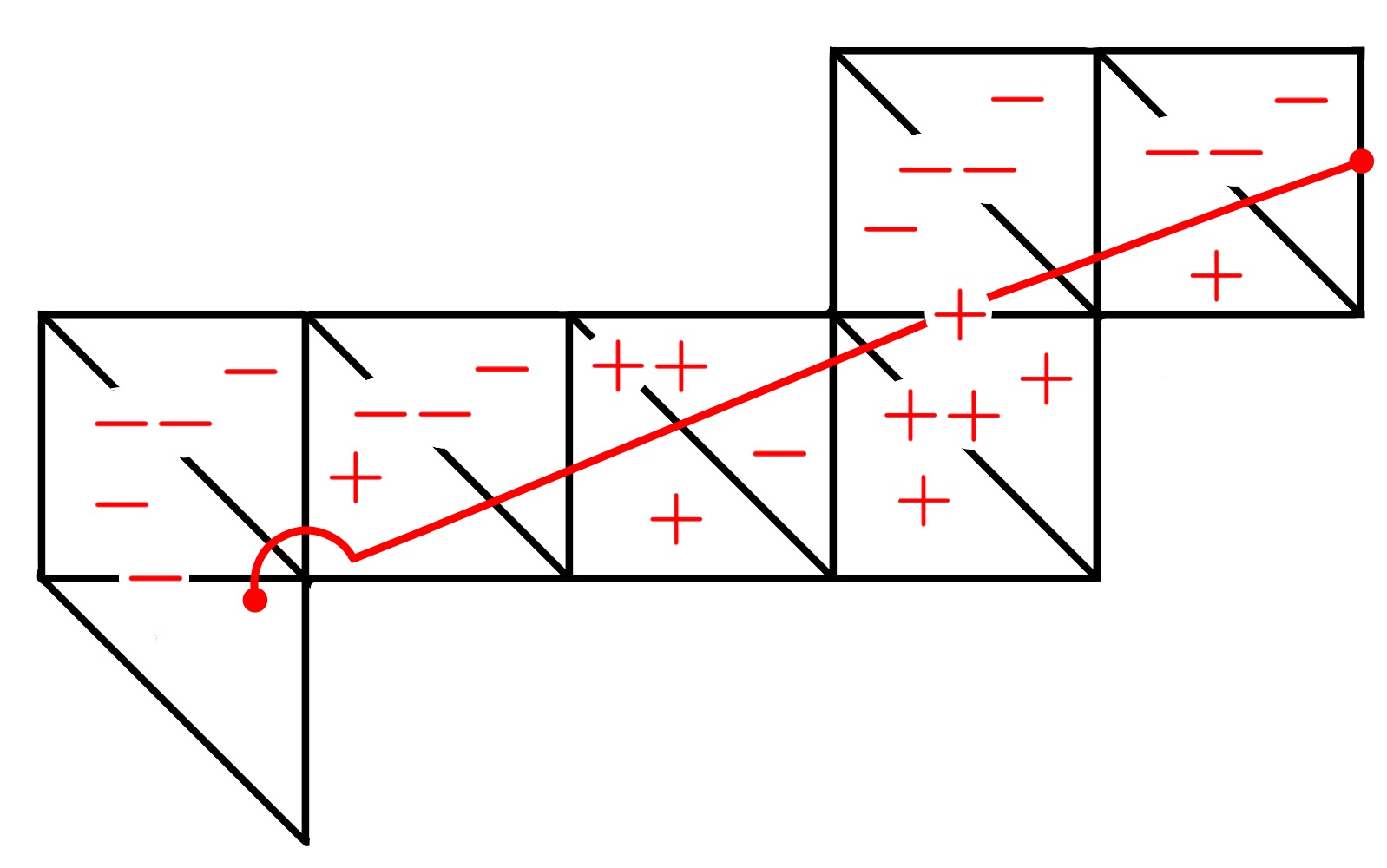} &\includegraphics[scale=0.06]{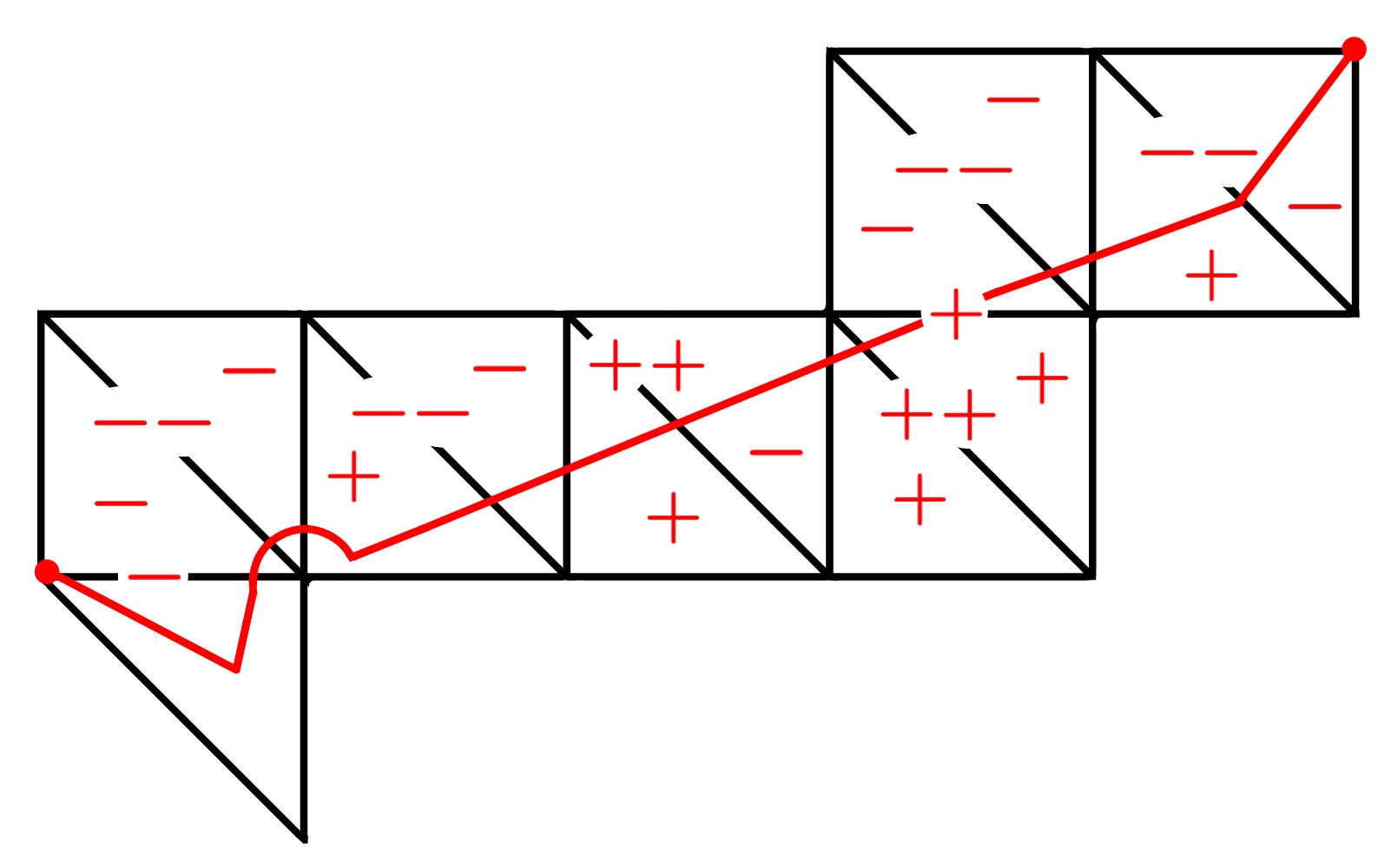}&$11127$\\
\hline
\end{tabular}
\vspace{2mm}
\caption{$\overline{L_t}(w)$ and $\widetilde{L_t}(w)$}\label{table2}
\end{table}
\end{example}
\begin{remark}
The following observations are useful for computing concrete examples.
\begin{itemize}\setlength{\leftskip}{-10pt}
    \item [(1)]
For given $(k_1,k_2,k_3)\in \mathbb Z_{\geq 0}^3$ and $\sigma\in \mathfrak S_3$, the position $i_t$ corresponds to the edge of $\widetilde{\mathbb{R}^2}$ through which $L_t$ passes exactly at its midpoint. Namely, if this edge is horizontal, then $i_t =\sigma(1)$; if it is diagonal, then $i_t = \sigma(2)$; and if it is vertical, then $i_t = \sigma(3)$. Moreover, $k_t$ coincides with the number of signs associated with this edge. 
\item[(2)]
For a positive integer sequence $w=(b_1,\dots,b_n)$, one can compute $m(\mathcal G[w])$ by evaluating the continued fraction $[b_1,\dots,b_n]$ and taking its numerator in lowest terms, by Theorem~\ref{thm:snakegraph-continuedfraction}. Equivalently, one can take the denominator of $[b_0,b_1,\dots,b_n]$ in lowest terms for any positive integer $b_0$.
\end{itemize}
\end{remark}

\subsection{Properties of Generalized Discrete Markov Spectra}
Fix $(k_1,k_2,k_3)\in\mathbb Z_{\geq 0}^3$. For any $n\in\mathbb Z_{\geq 1}$ and $i\in \{1,2,3\}$, we define
\[\Delta(n,i):=((3+k_1+k_2+k_3)n-k_i)^2-4.\] 
We denote 
\[
\mathcal M_{k_1,k_2,k_3,\sigma}:=\left\{\frac{\sqrt{\Delta(n_t,i_t)}}{n_t}\ \middle|\ \begin{aligned} (n_t,i_t)\text{ is a $(k_1,k_2,k_3,\sigma)$-GM number-position pair}\end{aligned}\right\}
\]
and consider 
\[\mathcal M_{k_1,k_2,k_3}:=\bigcup_{\sigma\in \mathfrak S_3}\mathcal M_{k_1,k_2,k_3,\sigma}.\]
We call the set $\mathcal M_{k_1,k_2,k_3}$ the \emph{$(k_1,k_2,k_3)$-generalized discrete Markov spectrum}. For any $\sigma\in \mathfrak S_3$, one has $\mathcal M_{k_1,k_2,k_3}=\mathcal M_{k_{\sigma(1)},k_{\sigma(2)},k_{\sigma(3)}}.$ Moreover, if $k_1=k_2=k_3=0$, then $\mathcal M_{0,0,0}$ coincides with the ordinary discrete Markov spectrum $\mathcal M_{d}$.

Corollary \ref{cor:markov-lagrange} immediately gives the following consequence.
\begin{corollary}\label{cor:inL}
For any $k_1,k_2,k_3\in \mathbb Z_{\geq 0}$, $\mathcal{M}_{k_1,k_2,k_3}\subset \mathcal L$ and $\mathcal{M}_{k_1,k_2,k_3}\subset \mathcal M$ hold.
\end{corollary}

\begin{example}
We list some quadratic irrational numbers $\alpha$ and the corresponding Lagrange constants $\mathcal L(\alpha)$ obtained from Corollary \ref{cor:markov-lagrange}. Table \ref{table000} shows the case $(k_1,k_2,k_3)=(0,0,0)$, Tables \ref{table001-1}, \ref{table001-2}, and \ref{table001-3} show the case $(k_1,k_2,k_3)=(0,0,1)$, Tables \ref{table011-1}, \ref{table011-2}, and \ref{table011-3} show the case $(k_1,k_2,k_3)=(0,1,1)$. Table \ref{table111} shows the case $(k_1,k_2,k_3)=(1,1,1)$, and Table \ref{table222} shows the case $(k_1,k_2,k_3)=(2,2,2)$. Table \ref{table120-1} shows the case $(k_1,k_2,k_3,\sigma)=(1,2,0, \mathrm{id})$. For each list, the elements are arranged in ascending order of the corresponding GM numbers.
\begin{table}[ht]
\centering
\begin{tabular}{|>{\centering\arraybackslash}m{1cm}|>{\centering\arraybackslash}m{5cm}|>{\centering\arraybackslash}m{4cm}|>{\centering\arraybackslash}m{1cm}|>{\centering\arraybackslash}m{3cm}|}
$t$&$s(t)$& $\alpha=[s(t)^\infty]$ &$n_t$&$\mathcal L(\alpha)$\\
\hline
&&&&\\[-3mm]
$\dfrac{0}{1}$&$(1,1)$&$\dfrac{\sqrt{5}+1}{2}$&$1$&$\sqrt{5}$\\[4.8mm]
\hline
&&&&\\[-3mm]
$\dfrac{1}{1}$&$(2,2)$&${\sqrt{2}+1}$&$2$&$2\sqrt{2}$\\[3mm]
\hline
&&&&\\[-3mm]
$\dfrac{1}{2}$&$(2,1,1,2)$&$\dfrac{\sqrt{221}+11}{10}$&$5$&$\dfrac{\sqrt{221}}{5}$\\[4.8mm]
\hline
&&&&\\[-3mm]
$\dfrac{1}{3}$&$(2,1,1,1,1,2)$&$\dfrac{\sqrt{1517}+29}{26}$&$13$&$\dfrac{\sqrt{1517}}{13}$\\[4.8mm]
\hline
&&&&\\[-3mm]
$\dfrac{2}{3}$&$(2,1,1,2,2,2)$&$\dfrac{\sqrt{7565}+63}{58}$&$29$&$\dfrac{\sqrt{7565}}{29}$\\[4.8mm]
\hline
&&&&\\[-3mm]
$\dfrac{1}{4}$&$(2,1,1,1,1,1,1,2)$&$\dfrac{5\sqrt{26}+19}{17}$&$34$&$\dfrac{10\sqrt{26}}{17}$\\[4.8mm]
\hline
&&&&\\[-3mm]
$\dfrac{1}{5}$&$(2,1,1,1,1,1,1,1,1,2)$&$\dfrac{\sqrt{71285}+199}{178}$&$89$&$\dfrac{\sqrt{71285}}{89}$\\[4.8mm]
\hline
&&&&\\[-3mm]
$\dfrac{3}{4}$&$(2,1,1,2,2,2,2,2)$&$\dfrac{\sqrt{257045}+367}{338}$&$169$&$\dfrac{\sqrt{257045}}{169}$\\[4.8mm]
\hline
\end{tabular}
\vspace{2mm}
\caption{Case $(k_1,k_2,k_3)=(0,0,0)$}\label{table000}
\end{table}
\begin{table}[ht]
\centering
\begin{tabular}{|>{\centering\arraybackslash}m{1cm}|>{\centering\arraybackslash}m{5cm}|>{\centering\arraybackslash}m{4cm}|>{\centering\arraybackslash}m{1cm}|>{\centering\arraybackslash}m{3cm}|}
$t$&$s(t)$& $\alpha=[s(t)^\infty]$ &$n_t$&$\mathcal L(\alpha)$\\
\hline
&&&&\\[-3mm]
$\dfrac{0}{1}$&$(2,1)$&$\sqrt{3}+1$&$1$&$2\sqrt{3}$\\[3mm]
\hline
&&&&\\[-3mm]
$\dfrac{1}{1}$&$(3,2)$&$\dfrac{\sqrt{15}+3}{2}$&$2$&$\sqrt{15}$\\[4.8mm]
\hline
&&&&\\[-3mm]
$\dfrac{1}{2}$&$(3,1,1,3)$&$\dfrac{5\sqrt{29}+23}{14}$&$7$&$\dfrac{5\sqrt{29}}{7}$\\[4.8mm]
\hline
&&&&\\[-3mm]
$\dfrac{1}{3}$&$(3,1,2,1,1,3)$&$\dfrac{7\sqrt{51}+43}{25}$&$25$&$\dfrac{14\sqrt{51}}{25}$\\[4.8mm]
\hline
&&&&\\[-3mm]
$\dfrac{2}{3}$&$(3,1,1,3,3,2)$&$\dfrac{\sqrt{11235}+83}{53}$&$53$&$\dfrac{2\sqrt{11235}}{53}$\\[4.8mm]
\hline
&&&&\\[-3mm]
$\dfrac{1}{4}$&$(3,1,2,1,1,2,1,3)$&$\dfrac{\sqrt{15293}+107}{62}$&$93$&$\dfrac{\sqrt{15293}}{31}$\\[4.8mm]
\hline
&&&&\\[-3mm]
$\dfrac{1}{5}$&$(3,1,2,1,2,1,1,2,1,3)$&$\dfrac{3\sqrt{53207}+599}{346}$&$346$&$\dfrac{3\sqrt{53207}}{173}$\\[4.8mm]
\hline
&&&&\\[-3mm]
$\dfrac{3}{4}$&$(3,1,1,3,2,3,3,2)$&$\dfrac{\sqrt{308765}+435}{278}$&$417$&$\dfrac{2\sqrt{308765}}{278}$\\[4.8mm]
\hline
\end{tabular}
\vspace{1mm}
\caption{Case $(k_1,k_2,k_3,\sigma)=(0,0,1,\mathrm{id})$}\label{table001-1}
\end{table}
\begin{table}
\centering
\begin{tabular}{|>{\centering\arraybackslash}m{1cm}|>{\centering\arraybackslash}m{5cm}|>{\centering\arraybackslash}m{4cm}|>{\centering\arraybackslash}m{1cm}|>{\centering\arraybackslash}m{3cm}|}
$t$&$s(t)$& $\alpha=[s(t)^\infty]$ &$n_t$&$\mathcal L(\alpha)$\\
\hline
&&&&\\[-3mm]
$\dfrac{0}{1}$&$(2,1)$&$\sqrt{3}+1$&$1$&$2\sqrt{3}$\\[3mm]
\hline
&&&&\\[-3mm]
$\dfrac{1}{1}$&$(3,3)$&$\dfrac{\sqrt{13}+3}{2}$&3&$\sqrt{13}$\\[4.8mm]
\hline
&&&&\\[-3mm]
$\dfrac{1}{2}$&$(3,1,2,3)$&$\dfrac{\sqrt{399}+17}{10}$&$10$&$\dfrac{\sqrt{399}}{5}$\\[4.8mm]
\hline
&&&&\\[-3mm]
$\dfrac{1}{3}$&$(3,1,2,2,1,3)$&$\dfrac{\sqrt{21605}+127}{74}$&$37$&$\dfrac{\sqrt{21605}}{37}$\\[4.8mm]
\hline
&&&&\\[-3mm]
$\dfrac{2}{3}$&$(3,1,2,3,3,3)$&$\dfrac{\sqrt{47523}+185}{109}$&$109$&$\dfrac{2\sqrt{47523}}{109}$\\[4.8mm]
\hline
&&&&\\[-3mm]
$\dfrac{1}{4}$&$(3,1,2,1,2,2,1,3)$&$\dfrac{5\sqrt{3003}+237}{137}$&$137$&$\dfrac{10\sqrt{3003}}{137}$\\[4.8mm]
\hline
&&&&\\[-3mm]
$\dfrac{1}{5}$&$(3,1,2,1,2,2,1,2,1,3)$&$\dfrac{\sqrt{4173845}+1769}{1022}$&$511$&$\dfrac{\sqrt{4173845}}{511}$\\[4.8mm]
\hline
&&&&\\[-3mm]
$\dfrac{3}{4}$&$(3,1,2,3,3,3,3,3)$&$\dfrac{\sqrt{5654883}+2018}{1189}$&$1189$&$\dfrac{2\sqrt{5654883}}{1189}$\\[4.8mm]
\hline
\end{tabular}
\vspace{1mm}
\caption{Case $(k_1,k_2,k_3,\sigma)=(0,0,1,(1\ 2 \ 3))$}\label{table001-2}
\end{table}
\begin{table}
\centering
\begin{tabular}{|>{\centering\arraybackslash}m{1cm}|>{\centering\arraybackslash}m{5cm}|>{\centering\arraybackslash}m{4cm}|>{\centering\arraybackslash}m{1cm}|>{\centering\arraybackslash}m{3cm}|}
$t$&$s(t)$& $\alpha=[s(t)^\infty]$ &$n_t$&$\mathcal L(\alpha)$\\
\hline
&&&&\\[-3mm]
$\dfrac{0}{1}$&$(1,1)$&$\dfrac{\sqrt{5}+1}{2}$&$1$&$\sqrt{5}$\\[4.8mm]
\hline
&&&&\\[-3mm]
$\dfrac{1}{1}$&$(3,2)$&$\dfrac{\sqrt{15}+3}{2}$&$2$&$\sqrt{15}$\\[4.8mm]
\hline
&&&&\\[-3mm]
$\dfrac{1}{2}$&$(3,1,1,2)$&$\dfrac{3\sqrt{11}+8}{5}$&$5$&$\dfrac{6\sqrt{11}}{5}$\\[4.8mm]
\hline
&&&&\\[-3mm]
$\dfrac{1}{3}$&$(3,1,1,1,1,2)$&$\dfrac{15\sqrt{3}+21}{13}$&$13$&$\dfrac{30\sqrt{3}}{13}$\\[4.8mm]
\hline
&&&&\\[-3mm]
$\dfrac{1}{4}$&$(3,1,1,1,1,1,1,2)$&$\dfrac{\sqrt{4623}+55}{34}$&$34$&$\dfrac{\sqrt{4623}}{17}$\\[4.8mm]
\hline
&&&&\\[-3mm]$\dfrac{2}{3}$&$(3,1,1,3,2,2)$&$\dfrac{\sqrt{2669}+41}{26}$&$39$&$\dfrac{\sqrt{2669}}{13}$\\[4.8mm]
\hline
&&&&\\[-3mm]
$\dfrac{1}{5}$&$(3,1,1,1,1,1,1,1,1,2)$&$\dfrac{\sqrt{31683}+144}{89}$&$89$&$\dfrac{2\sqrt{31683}}{89}$\\[4.8mm]
\hline
&&&&\\[-3mm]
$\dfrac{1}{6}$&$(3,1,1,1,1,1,1,1,1,1,1,2)$&$\dfrac{\sqrt{217155}+377}{233}$&$233$&$\dfrac{2\sqrt{217155}}{233}$\\[4.8mm]
\hline
\end{tabular}
\vspace{1mm}
\caption{Case $(k_1,k_2,k_3,\sigma)=(0,0,1,(1\ 3 \ 2))$}\label{table001-3}
\end{table}
\begin{table}[ht]
\centering
\begin{tabular}{|>{\centering\arraybackslash}m{1cm}|>{\centering\arraybackslash}m{5cm}|>{\centering\arraybackslash}m{4cm}|>{\centering\arraybackslash}m{1cm}|>{\centering\arraybackslash}m{3cm}|}
$t$&$s(t)$& $\alpha=[s(t)^\infty]$ &$n_t$&$\mathcal L(\alpha)$\\
\hline
&&&&\\[-3mm]
$\dfrac{0}{1}$&$(3,1)$&$\dfrac{\sqrt{21}+3}{2}$&$1$&$\sqrt{21}$\\[4.8mm]
\hline
&&&&\\[-3mm]
$\dfrac{1}{1}$&$(4,3)$&$\dfrac{4\sqrt{3}+6}{3}$&$3$&$\dfrac{8\sqrt{3}}{3}$\\[4.8mm]
\hline
&&&&\\[-3mm]
$\dfrac{1}{2}$&$(4,1,2,4)$&$\dfrac{\sqrt{1023}+29}{13}$&$13$&$\dfrac{2\sqrt{1023}}{13}$\\[4.8mm]
\hline
&&&&\\[-3mm]
$\dfrac{1}{3}$&$(4,1,3,2,1,4)$&$\dfrac{3\sqrt{2567}+139}{61}$&$61$&$\dfrac{6\sqrt{2567}}{61}$\\[4.8mm]
\hline
&&&&\\[-3mm]
$\dfrac{2}{3}$&$(4,1,2,4,4,3)$&$\dfrac{\sqrt{49506}+195}{89}$&$178$&$\dfrac{2\sqrt{49506}}{89}$\\[4.8mm]
\hline
&&&&\\[-3mm]
$\dfrac{1}{4}$&$(4,1,3,1,2,3,1,4)$&$\dfrac{44\sqrt{273}+666}{291}$&$291$&$\dfrac{88\sqrt{273}}{291}$\\[4.8mm]
\hline
&&&&\\[-3mm]
$\dfrac{1}{5}$&$(4,1,3,1,3,2,1,3,1,4)$&$\dfrac{531\sqrt{43}+3191}{1393}$&$1393$&$\dfrac{1062\sqrt{43}}{1393}$\\[4.8mm]
\hline
&&&&\\[-3mm]
$\dfrac{3}{4}$&$(4,1,2,4,3,4,4,3)$&$\dfrac{2\sqrt{9600702}+5431}{2479}$&$2479$&$\dfrac{4\sqrt{9600702}}{2479}$\\[4.8mm]
\hline
\end{tabular}
\vspace{1mm}
\caption{Case $(k_1,k_2,k_3,\sigma)=(0,1,1,\mathrm{id})$}\label{table011-1}
\end{table}
\begin{table}
\centering
\begin{tabular}{|>{\centering\arraybackslash}m{1cm}|>{\centering\arraybackslash}m{5cm}|>{\centering\arraybackslash}m{4cm}|>{\centering\arraybackslash}m{1cm}|>{\centering\arraybackslash}m{3cm}|}
$t$&$s(t)$& $\alpha=[s(t)^\infty]$ &$n_t$&$\mathcal L(\alpha)$\\
\hline
&&&&\\[-3mm]
$\dfrac{0}{1}$&$(2,1)$&$\sqrt{3}+1$&$1$&$2\sqrt{3}$\\[3mm]
\hline
&&&&\\[-3mm]
$\dfrac{1}{1}$&$(4,3)$&$\dfrac{4\sqrt{3}+6}{3}$&$3$&$\dfrac{8\sqrt{3}}{3}$\\[4.8mm]
\hline
&&&&\\[-3mm]
$\dfrac{1}{2}$&$(4,1,2,3)$&$\dfrac{2\sqrt{39}+11}{5}$&$10$&$\dfrac{4\sqrt{39}}{5}$\\[4.8mm]
\hline
&&&&\\[-3mm]
$\dfrac{1}{3}$&$(4,1,2,2,1,3)$&$\dfrac{\sqrt{8463}+82}{37}$&$37$&$\dfrac{2\sqrt{8463}}{37}$\\[4.8mm]
\hline
&&&&\\[-3mm]
$\dfrac{1}{4}$&$(4,1,2,1,2,2,1,3)$&$\dfrac{\sqrt{469221}+611}{274}$&$137$&$\dfrac{\sqrt{469221}}{137}$\\[4.8mm]
\hline
&&&&\\[-3mm]
$\dfrac{2}{3}$&$(4,1,2,4,3,3)$&$\dfrac{2\sqrt{30102}+305}{139}$&$139$&$\dfrac{4\sqrt{30102}}{139}$\\[4.8mm]
\hline
&&&&\\[-3mm]
$\dfrac{1}{5}$&$(4,1,2,1,2,2,1,2,1,3)$&$\dfrac{6\sqrt{45298}+1140}{511}$&$511$&$\dfrac{12\sqrt{45298}}{511}$\\[4.8mm]
\hline
&&&&\\[-3mm]
$\dfrac{2}{5}$&$(4,1,2,2,1,4,3,1,2,3)$&$\dfrac{22\sqrt{43662}+4050}{1839}$&$1839$&$\dfrac{44\sqrt{43662}}{1839}$\\[4.8mm]
\hline
\end{tabular}
\vspace{1mm}
\caption{Case $(k_1,k_2,k_3,\sigma)=(0,1,1,(1\ 2 \ 3))$}\label{table011-2}
\end{table}
\begin{table}
\centering
\begin{tabular}{|>{\centering\arraybackslash}m{1cm}|>{\centering\arraybackslash}m{5cm}|>{\centering\arraybackslash}m{4cm}|>{\centering\arraybackslash}m{1cm}|>{\centering\arraybackslash}m{3cm}|}
$t$&$s(t)$& $\alpha=[s(t)^\infty]$ &$n_t$&$\mathcal L(\alpha)$\\
\hline
&&&&\\[-3mm]
$\dfrac{0}{1}$&$(2,1)$&$\sqrt{3}+1$&$1$&$2\sqrt{3}$\\[3mm]
\hline
&&&&\\[-3mm]
$\dfrac{1}{1}$&$(4,2)$&$\sqrt6+2$&$2$&$2\sqrt{6}$\\[3mm]
\hline
&&&&\\[-3mm]
$\dfrac{1}{2}$&$(4,1,1,3)$&$\dfrac{12\sqrt{2}+15}{7}$&$7$&$\dfrac{24\sqrt{2}}{7}$\\[4.8mm]
\hline
&&&&\\[-3mm]
$\dfrac{1}{3}$&$(4,1,2,1,1,3)$&$\dfrac{\sqrt{15621}+111}{50}$&$25$&$\dfrac{\sqrt{15621}}{25}$\\[4.8mm]
\hline
&&&&\\[-3mm]
$\dfrac{2}{3}$&$(4,1,1,4,3,2)$&$\dfrac{4\sqrt{1743}+138}{67}$&$67$&$\dfrac{8\sqrt{1743}}{67}$\\[4.8mm] 
\hline
&&&&\\[-3mm]
$\dfrac{1}{4}$&$(4,1,2,1,1,2,1,3)$&$\dfrac{\sqrt{53823}+207}{93}$&$93$&$\dfrac{2\sqrt{53823}}{93}$\\[4.8mm]
\hline
&&&&\\[-3mm]
$\dfrac{1}{5}$&$(4,1,2,1,2,1,1,2,1,3)$&$\dfrac{12\sqrt{1299}+386}{173}$&$346$&$\dfrac{24\sqrt{1299}}{173}$\\[4.8mm]
\hline
&&&&\\[-3mm]
$\dfrac{3}{4}$&$(4,1,1,4,2,3,4,2)$&$\dfrac{\sqrt{2729103}+1356}{661}$&$661$&$\dfrac{2\sqrt{2729103}}{661}$\\[4.8mm]
\hline
\end{tabular}
\vspace{1mm}
\caption{Case $(k_1,k_2,k_3,\sigma)=(0,1,1,(1\ 3 \ 2))$}\label{table011-3}
\end{table}
\begin{table}[ht]
\centering
\begin{tabular}{|>{\centering\arraybackslash}m{1cm}|>{\centering\arraybackslash}m{5cm}|>{\centering\arraybackslash}m{4cm}|>{\centering\arraybackslash}m{1cm}|>{\centering\arraybackslash}m{3cm}|}
$t$&$s(t)$& $\alpha=[s(t)^\infty]$ &$n_t$&$\mathcal L(\alpha)$\\
\hline
&&&&\\[-3mm]
$\dfrac{0}{1}$&$(3,1)$&$\dfrac{\sqrt{21}+3}{2}$&$1$&$\sqrt{21}$\\[4.8mm]
\hline
&&&&\\[-3mm]
$\dfrac{1}{1}$&$(5,3)$&$\dfrac{\sqrt{285}+15}{6}$&$3$&$\dfrac{\sqrt{285}}{3}$\\[4.8mm]
\hline
&&&&\\[-3mm]
$\dfrac{1}{2}$&$(5,1,2,4)$&$\dfrac{5\sqrt{237}+71}{26}$&$13$&$\dfrac{5\sqrt{237}}{13}$\\[4.8mm]
\hline
&&&&\\[-3mm]
$\dfrac{1}{3}$&$(5,1,3,2,1,4)$&$\dfrac{11\sqrt{1101}+339}{122}$&$61$&$\dfrac{11\sqrt{1101}}{61}$\\[4.8mm]
\hline
&&&&\\[-3mm]
$\dfrac{2}{3}$&$(5,1,2,5,4,3)$&$\dfrac{\sqrt{1692597}+1167}{434}$&$217$&$\dfrac{\sqrt{1692597}}{217}$\\[4.8mm]
\hline
&&&&\\[-3mm]
$\dfrac{1}{4}$&$(5,1,3,1,2,3,1,4)$&$\dfrac{\sqrt{3045021}+1623}{582}$&$291$&$\dfrac{\sqrt{3045021}}{291}$\\[4.8mm]
\hline
&&&&\\[-3mm]
$\dfrac{1}{5}$&$(5,1,3,1,3,2,1,3,1,4)$&$\dfrac{\sqrt{69839445}+7775}{2786}$&$1393$&$\dfrac{\sqrt{69839445}}{1393}$\\[4.8mm]
\hline
&&&&\\[-3mm]
$\dfrac{3}{4}$&$(5,1,2,5,3,4,5,3)$&$\dfrac{\sqrt{485629365}+19735}{7346}$&$3673$&$\dfrac{\sqrt{485629365}}{3673}$\\[4.8mm]
\hline
\end{tabular}
\vspace{1mm}
\caption{Case $(k_1,k_2,k_3)=(1,1,1)$}\label{table111}
\end{table}
\begin{table}[ht]
\centering
\begin{tabular}{|>{\centering\arraybackslash}m{1cm}|>{\centering\arraybackslash}m{5cm}|>{\centering\arraybackslash}m{4cm}|>{\centering\arraybackslash}m{1cm}|>{\centering\arraybackslash}m{3cm}|}
$t$&$s(t)$& $\alpha=[s(t)^\infty]$ &$n_t$&$\mathcal L(\alpha)$\\
\hline
&&&&\\[-3mm]
$\dfrac{0}{1}$&$(5,1)$&$\dfrac{3\sqrt{5}+5}{2}$&$1$&$3\sqrt{5}$\\[4.8mm]
\hline
&&&&\\[-3mm]
$\dfrac{1}{1}$&$(8,4)$&$3\sqrt{2}+4$&$4$&$6\sqrt{2}$\\[4mm]
\hline
&&&&\\[-3mm]
$\dfrac{1}{2}$&$(8,1,3,6)$&$\dfrac{3\sqrt{221}+43}{10}$&$25$&$\dfrac{3\sqrt{221}}{5}$\\[4.8mm]
\hline
&&&&\\[-3mm]
$\dfrac{1}{3}$&$(8,1,5,3,1,6)$&$\dfrac{3\sqrt{1517}+113}{26}$&$169$&$\dfrac{3\sqrt{1517}}{13}$\\[4.8mm]
\hline
&&&&\\[-3mm]
$\dfrac{2}{3}$&$(8,1,3,8,6,4)$&$\dfrac{3\sqrt{7565}+247}{58}$&$841$&$\dfrac{3\sqrt{7565}}{29}$\\[4.8mm]
\hline
&&&&\\[-3mm]
$\dfrac{1}{4}$&$(8,1,5,1,3,5,1,6)$&$\dfrac{15\sqrt{26}+74}{17}$&$1156$&$\dfrac{30\sqrt{26}}{17}$\\[4.8mm]
\hline
&&&&\\[-3mm]
$\dfrac{1}{5}$&$(8,1,5,1,5,3,1,5,1,6)$&$\dfrac{3\sqrt{71285}+775}{178}$&$7921$&$\dfrac{3\sqrt{71285}}{89}$\\[4.8mm]
\hline
&&&&\\[-3mm]
$\dfrac{3}{4}$&$(8,1,3,8,4,6,8,4)$&$\dfrac{3\sqrt{257045}+1439}{338}$&$28561$&$\dfrac{3\sqrt{257045}}{169}$\\[4.8mm]
\hline
\end{tabular}
\vspace{1mm}
\caption{Case $(k_1,k_2,k_3)=(2,2,2)$}\label{table222}
\end{table}
\begin{table}[ht]
\centering
\begin{tabular}{|>{\centering\arraybackslash}m{1cm}|>{\centering\arraybackslash}m{5cm}|>{\centering\arraybackslash}m{4cm}|>{\centering\arraybackslash}m{1cm}|>{\centering\arraybackslash}m{3cm}|}
$t$&$s(t)$& $\alpha=[s(t)^\infty]$ &$n_t$&$\mathcal L(\alpha)$\\
\hline
&&&&\\[-3mm]
$\dfrac{0}{1}$&$(3,1)$&$\dfrac{\sqrt{21}+3}{2}$&$1$&$\sqrt{21}$\\[4.8mm]
\hline
&&&&\\[-3mm]
$\dfrac{1}{1}$&$(5,4)$&$\dfrac{{\sqrt{30}+5}}{2}$&$4$&$\sqrt{30}$\\[4mm]
\hline
&&&&\\[-3mm]
$\dfrac{1}{2}$&$(5,1,3,4)$&$\dfrac{10\sqrt{26}+47}{17}$&$17$&$\dfrac{20\sqrt{26}}{17}$\\[4.8mm]
\hline
&&&&\\[-3mm]
$\dfrac{1}{3}$&$(5,1,3,3,1,4)$&$\dfrac{\sqrt{723}+25}{9}$&$81$&$\dfrac{2\sqrt{723}}{9}$\\[4.8mm]
\hline
&&&&\\[-3mm]
$\dfrac{2}{3}$&$(5,1,3,5,4,4)$&$\dfrac{\sqrt{5004165}+2061}{746}$&$373$&$\dfrac{\sqrt{5004165}}{373}$\\[4.8mm]
\hline
&&&&\\[-3mm]
$\dfrac{1}{4}$&$(5,1,3,1,3,3,1,4)$&$\dfrac{\sqrt{1340963}+1077}{386}$&$386$&$\dfrac{\sqrt{1340963}}{193}$\\[4.8mm]
\hline
&&&&\\[-3mm]
$\dfrac{1}{5}$&$(5,1,3,1,3,3,1,3,1,4)$&$\dfrac{\sqrt{16635}+120}{43}$&$1849$&$\dfrac{2\sqrt{16635}}{43}$\\[4.8mm]
\hline
&&&&\\[-3mm]
$\dfrac{3}{4}$&$(5,1,3,5,4,4,5,4)$&\scalebox{0.98}{$\dfrac{2\sqrt{150737006}+22602}{8185}$}&$8185$&$\dfrac{4\sqrt{150737006}}{8185}$\\[4.8mm]
\hline
\end{tabular}
\vspace{2mm}
\caption{Case $(k_1,k_2,k_3,\sigma)=(1,2,0,\mathrm{id})$}\label{table120-1}
\end{table}
\end{example}
The following proposition relates $(0,0,0)$-GM triples and $(2,2,2)$-GM triples by squaring their entries.
\begin{proposition}[\cite{gyomatsu}*{Theorem 11}]\label{square-markov}
If $(a,b,c)$ is a $(0,0,0)$-GM triple, then $(a^2,b^2,c^2)$ is a $(2,2,2)$-GM triple. Conversely, if $(A,B,C)$ is a $(2,2,2)$-GM triple, then $(\sqrt{A},\sqrt{B},\sqrt{C})$ is a $(0,0,0)$-GM triple.
\end{proposition}

Proposition \ref{square-markov} gives the following relation between the two spectra.

\begin{theorem}\label{thm:0-2-relation}
If $r\in \mathcal M_{0,0,0}$, then $3r\in \mathcal M_{2,2,2}$. Conversely, if $R\in \mathcal M_{2,2,2}$, then $\frac R3\in \mathcal M_{0,0,0}$.
\end{theorem}

\begin{proof}
We first prove the forward implication. Since $r\in \mathcal M_{0,0,0}$, there exists a $(0,0,0)$-GM number $n$ such that
\[
r=\frac{\sqrt{9n^2-4}}{n}.
\]
Thus
\[
3r=\frac{3\sqrt{9n^2-4}}{n}=\frac{\sqrt{81n^4-36n^2}}{n^2}=\frac{\sqrt{(9n^2-2)^2-4}}{n^2}.
\]
By Proposition \ref{square-markov}, $n^2$ is a $(2,2,2)$-GM number, and hence $3r\in \mathcal M_{2,2,2}$.  
The reverse implication follows by reversing the above computation.
\end{proof}
The duality $(n_t,i_t)=(n^\ast_{\frac1t},i^\ast_{\frac1t})$ gives the following reductions in the definition of $\mathcal M_{k_1,k_2,k_3}$. One may restrict either $\sigma$ to the alternating group $\mathfrak A_3$ or the fraction labels $t$ to $[0,1]$ without changing the spectrum. If two of $k_1,k_2,k_3$ are equal, both restrictions may be imposed simultaneously.

Put $K=3+k_1+k_2+k_3$ and, for $\kappa\in\{k_1,k_2,k_3\}$, set
\[
\lambda_\kappa(x):=\frac{\sqrt{(Kx-\kappa)^2-4}}{x}\qquad (x\geq1).
\]
Then $\lambda_\kappa$ is strictly increasing. Indeed,
\[
\frac{d}{dx}\lambda_\kappa(x)^2
=\frac{2(K\kappa x-\kappa^2+4)}{x^3}>0,
\]
where the last inequality follows from $K\geq \kappa+3$. Moreover, $\lambda_\kappa(x)<K$ for every $x\geq1$.
Indeed,
\[
K^2-\lambda_\kappa(x)^2=\frac{2K\kappa x-\kappa^2+4}{x^2}>0.
\]

GM trees have the following growth property.
\begin{lemma}[\cite{gyomatsu}*{Corollary 6 and Remark 7}]\label{lem:GM-tree-growth}
If
\[
((a,h),(b,i),(c,j))
\]
is a vertex of a GM tree, then $b>\max\{a,c\}$, and each of the two new middle entries in its children is strictly larger than $b$.
\end{lemma}

For any $(k_1,k_2,k_3)\in \mathbb Z_{\geq 0}^3$, the minimal $(k_1,k_2,k_3)$-GM number is $1$. If $k_i,k_j$ are the two smallest elements of $\{k_1,k_2,k_3\}$ and $k_\ell$ is the largest one, the monotonicity above gives
\[
\lambda_{k_t}(n_t)\geq \lambda_{k_t}(1)
\geq K-k_t-1\geq K-k_\ell-1=2+k_i+k_j.
\]
Here the middle inequality follows from $K-k_t\geq3$ and
$(K-k_t)^2-4\geq(K-k_t-1)^2$.
Together with $\lambda_{k_t}(n_t)<K$, this proves
\[
\mathcal{M}_{k_1,k_2,k_3}\subset [2+k_i+k_j,3+k_1+k_2+k_3],
\]
where $k_i,k_j$ are the two smallest elements among $\{k_1,k_2,k_3\}$.

We next characterize generalized discrete Markov spectra in the transition interval.
\begin{theorem}\label{thm:freiman-interval}
Set
\[\mathcal M'=\bigcup_{(k_1,k_2,k_3)\in \mathbb Z_{\geq 0}^3}\mathcal M_{k_1,k_2,k_3}.\]
Then we have
\[\mathcal M'\cap [3,c_F)=(\mathcal M_{0,0,1}\setminus\{\sqrt 5\})\cup\{2\sqrt 5\},\]
where $c_F$ is the Freiman constant, i.e.,
\[
c_F:=\dfrac{2221564096 + 283748 \sqrt{462}}{491993569} \approx 4.5278295661\dots.
\]
\end{theorem}
\begin{proof}
Without loss of generality, assume $k_1\leq k_2\leq k_3$. We first prove $\mathcal M'\cap [3,c_F)\supset(\mathcal M_{0,0,1}\setminus\{\sqrt 5\})\cup\{2\sqrt 5\}$. The inclusion $\mathcal M_{0,0,1}\setminus\{\sqrt 5\}\subset \mathcal M'$ is immediate. Moreover, we have $2\sqrt 5\in \mathcal M_{0,0,2}\subset \mathcal M'$. Indeed, $(4,3)$ is a $(0,0,2,\mathrm{(1\ 2\ 3)})$-GM number-position pair, and the corresponding element of $\mathcal M_{0,0,2}$ is
\[\frac{\sqrt{(5\cdot 4-2)^2-4}}{4}=\frac{\sqrt{320}}{4}=2\sqrt{5}.\]

We prove $(\mathcal M_{0,0,1}\setminus\{\sqrt 5\})\cup\{2\sqrt5\}\subset [3,c_F)$. The element $2\sqrt{5}$ lies in $[3,c_F)$. The functions
\[f(x)=\frac{\sqrt{(4x)^2-4}}{x},\quad g(x)=\frac{\sqrt{(4x-1)^2-4}}{x}\] are both strictly increasing on $x\geq 1$, and
\[\lim_{x\to \infty}f(x)=\lim_{x\to\infty}g(x)=4.\]
Every element of $\mathcal M_{0,0,1}$ is of the form $f(n_t)$ or $g(n_t)$ according as $k_t=0$ or $1$, so all these elements are smaller than $4<c_F$. The value $g(1)=\sqrt5$ occurs when the entry $1$ is one of the two end entries of a root vertex, whereas $f(1)=2\sqrt3$. By Lemma~\ref{lem:GM-tree-growth}, an entry equal to $2$ can only originate as the middle entry of a root; in the present case $k_t=0$. Thus, apart from the occurrence $g(1)=\sqrt5$, the smallest possible arguments are $1$ for $f$ and $3$ for $g$. Since $f(1)=2\sqrt3<g(3)=\sqrt{13}$, the minimal element of $\mathcal M_{0,0,1}\setminus\{\sqrt 5\}$ is $2\sqrt3$. Therefore $\mathcal M_{0,0,1}\setminus\{\sqrt 5\}\subset [3,c_F)$.

We next prove $\mathcal M'\cap [3,c_F)\subset(\mathcal M_{0,0,1}\setminus\{\sqrt 5\})\cup\{2\sqrt{5}\}$. First,
\[
\mathcal M_{0,0,0}\cap[3,c_F)=\emptyset,
\qquad
\mathcal M_{0,0,1}\cap[3,c_F)=\mathcal{M}_{0,0,1}\setminus\{\sqrt{5}\}.
\]
Thus the cases $(k_1,k_2,k_3)=(0,0,0)$ and $(0,0,1)$ have already been treated. We now consider the remaining cases in turn. First, assume $k_1=k_2=0$ and $k_3=k\in \mathbb Z_{\geq2}$. Suppose a number-position pair with new position $j$ is created below a vertex
\[
((a,h),(b,i),(c,j)),\qquad b>\max\{a,c\}.
\]
The new number in the left child satisfies
\[
\frac{a^2+k_jab+b^2}{c}
>k_ja+b\geq k_j+2.
\]
Indeed, $b/c>1$, $a\geq1$, and $b\geq2$. The identical estimate for the right child is
\[
\frac{b^2+k_hbc+c^2}{a}>k_hc+b\geq k_h+2.
\]
Consequently, every number-position pair $(n_t,i_t)$ created below a root satisfies $n_t>k_t+2$.

Since $k_1=k_2$, duality permits both restrictions simultaneously: it is enough to consider the three trees with $\sigma\in\mathfrak A_3$ and labels in $[0,1]$. At their roots, the pairs satisfying $k_t=k$ have numbers $1$ or $k+2$, and both values occur. It follows that
\[
\min\{n_t:k_t=k\}=1,\qquad
\min\{n_t>1:k_t=k\}=k+2.
\]
The minimum of $n_t$ subject to $k_t=0$ is also $1$. Since $\lambda_0$ is increasing, every value with $k_t=0$ is at least
\[
\lambda_0(1)=\sqrt{(k+3)^2-4}\geq\sqrt{21}>c_F.
\]
The value with $k_t=k$ and $n_t=1$ is $\sqrt5<3$, while every such value with $n_t>1$ is at least
\[
\lambda_k(k+2)=\sqrt{k^2+4k+8}.
\]
This equals $2\sqrt5$ when $k=2$, and it is at least $\sqrt{29}>c_F$ when $k\geq3$. For $k=2$, any entry $n_t>4$ is an integer at least $5$, so monotonicity gives
\[
\lambda_2(n_t)\geq\lambda_2(5)
=\frac{\sqrt{(5\cdot5-2)^2-4}}{5}
=\sqrt{21}>c_F.
\]
Therefore,
\[\mathcal{M}_{0,0,2}\cap[3,c_F)=\{2\sqrt 5\},\]
and for any $k\geq 3$,
\[\mathcal{M}_{0,0,k}\cap[3,c_F)=\emptyset.\]

Next, assume that $k_1=0$ and $k_2=k_3=1$. If $n_t=1$, the corresponding value is $2\sqrt3$ when $k_t=1$, and it is $\sqrt{21}>c_F$ when $k_t=0$. By Lemma~\ref{lem:GM-tree-growth}, an entry $n_t=2$ can only be a root middle entry; for these parameters such an entry has $k_t=0$, and its spectral value is $\sqrt{24}>c_F$. The estimate $n_t>k_t+2$ for every pair created below a root shows that every non-root entry with $k_t=1$ is greater than $3$. Hence the smallest entry greater than $1$ with $k_t=1$ is the root middle entry $n_t=3$, whose spectral value is $8\sqrt3/3>c_F$. The monotonicity of $\lambda_0,\lambda_1$ and Lemma~\ref{lem:GM-tree-growth} then show that every larger entry gives a value outside $[3,c_F)$. Therefore
\[\mathcal{M}_{0,1,1}\cap[3,c_F)=\{2\sqrt{3}\},\]
and $2\sqrt{3}\in \mathcal {M}_{0,0,1}$. 

It remains to consider the parameter triples that are neither of the form $(0,0,k)$ nor equal to $(0,1,1)$. Set $K:=3+k_1+k_2+k_3$. For these triples, $K\geq6$, and at least two of $k_1,k_2,k_3$ are positive. Hence $\sum_{r\ne i}k_r\geq1$, or equivalently $k_i\leq K-4$, for every $i$. If
\[\dfrac{\sqrt{(Kn_t-k_t)^2-4}}{n_t}<5,\]
then by assumption we have
\[\dfrac{\sqrt{(Kn_t-K+4)^2-4}}{n_t}\leq \dfrac{\sqrt{(Kn_t-k_t)^2-4}}{n_t}<5.\]
Squaring the inequality between the leftmost and rightmost terms and rearranging gives
\[(K^2-25)n_t^2+(-2K^2+8K)n_t+(K^2-8K+12)<0.\]
Since $K\geq 6$, we have $K^2-25>0$ and thus
\[\dfrac{K(K-4)-\sqrt{29K^2-200K+300}}{K^2-25}<n_t<\dfrac{K(K-4)+\sqrt{29K^2-200K+300}}{K^2-25}.\]
Now, set $D(x)=29x^2-200x+300$. For $x\geq6$, the quantity $13x^2+44x-600$ is positive, and
\begin{align*}
&(13x^2+44x-600)^2-121D(x)\\
&\qquad=(x-6)(169x^3+2158x^2-4225x-53950)\geq0.
\end{align*}
The last cubic is positive and increasing on $[6,\infty)$. Hence
$11\sqrt{D(x)}\leq13x^2+44x-600$, which is equivalent to
\[
F(x):=\dfrac{x(x-4)+\sqrt{D(x)}}{x^2-25}\leq F(6)=\frac{24}{11}.
\]
Therefore, we have
\[n_t<\dfrac{6(6-4)+\sqrt{29\cdot 36-200\cdot6+300}}{36-25}=\frac{24}{11}.\]
Therefore, if 
\[\dfrac{\sqrt{((3+k_1+k_2+k_3)n_t-k_t)^2-4}}{n_t}\in [3,c_F),\]
then $n_t$ must be $1$ or $2$. First, we assume $n_t=1$. If $\sqrt{(3+k_1+k_2+k_3-k_t)^2-4}\in [3,c_F)$, then $k_1+k_2+k_3-k_t=1$.  We have thus $k_1=0$, $k_2=1$, $k_t=k_3$ and
\[\sqrt{(3+k_1+k_2+k_3-k_t)^2-4}=2\sqrt{3}\in \mathcal M_{0,0,1}.\]
Second, we assume that $n_t=2$. By Lemma~\ref{lem:GM-tree-growth}, an entry equal to $2$ cannot be created below a root, because every newly created middle entry is strictly larger than the preceding middle entry, which is at least $2$. Hence it originates as a root middle entry $k_{\sigma(2)}+2$, and $k_t=0$. Thus $k_1=0$. Then we have
\[\dfrac{\sqrt{((3+k_1+k_2+k_3)n_t-k_t)^2-4}}{n_t}=\dfrac{\sqrt{(2\cdot(3+k_2+k_3))^2-4}}{2}.\]
Since this is neither a case of the form $(0,0,k)$ nor the case $(0,1,1)$, we have $k_2+k_3\geq3$. Hence
\[\dfrac{\sqrt{(2\cdot(3+k_2+k_3))^2-4}}{2}\notin [3,c_F).\]
Conversely, let $k\geq2$ and $(k_1,k_2,k_3)=(0,1,k)$. Choose $\sigma$ with $\sigma(1)=3$. The number-position pair $(1,3)$ then occurs at the left endpoint of the root. For this pair, $n_t=1$ and $k_t=k$, and the corresponding value is
\[
\sqrt{(4+k-k)^2-4}=2\sqrt3.
\]
The preceding argument therefore shows that, if $k\geq2$, then
\[\mathcal{M}_{0,1,k} \cap[3,c_F)=\{2\sqrt3\}\subset\mathcal M_{0,0,1},\]
and if $(k_1,k_2,k_3)$ does not have the form $(0,0,0),(0,0,k)$ or $(0,1,k)$, then we have
\[\mathcal{M}_{k_1,k_2,k_3} \cap[3,c_F)=\emptyset.\]
This finishes the proof.
\end{proof}
\begin{proposition}\label{rem:non-contain}
The set $\mathcal M'\cap[3,c_F)$ does not cover all Lagrange constants contained in the transition interval. More precisely,
\[
\mathcal M'\cap [3,c_F)\subsetneq \mathcal L\cap [3,c_F)
\]
holds.
\end{proposition}
\begin{proof}
Set $s=(1,1,1,2,2,2)$ and
\[
\alpha=[(1,1,1,2,2,2)^\infty]=\frac{2\sqrt{210}+17}{29}.
\]
Using Theorem~\ref{thm:lagrange-quadratic} and Lemma~\ref{lem:trace}, we have
\[
\mathcal L(\alpha) =\frac{4\sqrt{210}}{19}
\approx3.05081615709\dots,
\]
and $\mathcal L(\alpha)\in\mathcal L\cap[3,c_F)$. On the other hand, the minimal element of $\mathcal M'\cap[3,c_F)$ is
\[2\sqrt{3}\approx3.46410161514\dots,\]
and $\mathcal L(\alpha)\notin \mathcal M'\cap[3,c_F)$.
\end{proof}

Theorem \ref{thm:characterization} states that the discrete Markov spectrum coincides with the part of the Markov spectrum below $3$. A corresponding theorem for GM numbers is not yet known.

\begin{question}\label{ques:characterization}
Is there a characterization of $\mathcal M_{k_1,k_2,k_3}$, or of their union, that generalizes Markov’s theorem?    
\end{question}

\subsection{Boundary Values from Irrational Slopes}\label{subsec:irrational-slope-boundary}
\begingroup\emergencystretch=2em
In this subsection, we keep $(k_1,k_2,k_3)\in\mathbb Z_{\geq 0}^3$ and $\sigma\in\mathfrak S_3$ fixed, and put
\[
K:=3+k_1+k_2+k_3.
\]
We consider only positive slopes, in accordance with the definition of the generalized strongly admissible sequence $s(t)$ for $t\in[0,\infty]\cap\mathbb Q$. A line of positive slope is oriented in the direction in which the $x$-coordinate increases.
\par\endgroup

Let $\mathcal P$ denote the set consisting of all lattice points and all midpoints of edges in $\widetilde{\mathbb R^2}$. We call a line \emph{regular} if it does not pass through any point of $\mathcal P$.

\begin{definition}
Let $l$ be an oriented regular line of positive slope. Apply the triangle-crossing and edge-crossing rules to $l$, and arrange the resulting signs in the order in which $l$ passes through them. Decompose the resulting bi-infinite sign word into runs, that is, maximal consecutive blocks of identical signs. By recording their lengths, we obtain a two-sided infinite sequence of positive integers, denoted by
\[
\mathbf b(l)=(b_n)_{n\in\mathbb Z}.
\]
The origin of the index is chosen arbitrarily, so $\mathbf b(l)$ is determined only up to shift.
\end{definition}

For an irrational slope, the preceding sequence is genuinely indexed by all of $\mathbb Z$. Indeed, the line meets the locally finite triangulation in a discrete sequence which is unbounded in both directions. Modulo $\mathbb Z^2$, an irrational-slope line is dense in the torus. Choose a short line segment in the torus that the projected line crosses transversely.  The intersection points corresponding to the two triangle-passage signs form two nonempty open subintervals of this segment.  Density therefore implies that each sign occurs infinitely often in both directions. Edge crossings insert blocks of at most $\max\{k_1,k_2,k_3\}$ signs and cannot eliminate these sign changes. Every run is therefore finite, and the changes of sign are unbounded in both directions.

When a finite integer block
\[
W=(b_r,b_{r+1},\ldots,b_s)
\]
of $\mathbf b(l)$ is treated geometrically, we retain not only the finite sign word consisting of the runs represented by $W$, but also the sign immediately preceding its first run and the sign immediately following its last run. We call the resulting finite sign word the \emph{extended sign block} associated with $W$ and denote it by $\widehat W$. The two additional signs are opposite to the signs in the adjacent runs. Thus, if $\widehat W$ is preserved under a perturbation of the line, the first and last runs cannot merge with neighboring runs, and the integer block $W$ itself is preserved exactly.

We use the functions $\ell_r(\mathbf a)$ and $L(\mathbf a)$ for bi-infinite positive integer sequences as defined in Section~\ref{section:Classical Markov Spectrum}. Both shifting the index origin and reversing the orientation leave $L(\mathbf b(l))$ unchanged.

The goal of this subsection is the following theorem.

\begin{theorem}\label{thm:irrational-slope-boundary-value}
Let $(k_1,k_2,k_3)\in\mathbb Z_{\geq 0}^3$ and $\sigma\in\mathfrak S_3$ be fixed, and put $K=3+k_1+k_2+k_3$. If $l$ is a regular line of positive irrational slope in $\widetilde{\mathbb R^2}$, then
\[
L(\mathbf b(l))=K.
\]
\end{theorem}

We prepare several elementary lemmas.

\begin{lemma}\label{lem:rational-line-period}
For every positive irreducible fraction $t=p/q$, there is a regular affine line $l_t^0$ of slope $t$, arbitrarily close to the oriented line through $L_t$, such that
\[
\mathbf b(l_t^0)={}^\infty s(t)^\infty
\]
up to a shift of the indices.
\end{lemma}

\begin{proof}
Choose a sufficiently small generic translate of the supporting line of $L_t$ to its left. The exceptional translates which meet a lattice point or an edge midpoint form a discrete set, so the translate may be chosen regular. Its sign pattern is invariant under translation by $(q,p)$. Take a half-open segment between a point and its translate by $(q,p)$, choosing the cut so that it corresponds to the initial--terminal cut in the definition of $\overline L_t$. By that definition, for a sufficiently small left translate, grouping adjacent equal signs on this fundamental segment and recording their block lengths gives $s(t)$.

It remains only to check the seam between consecutive periods. Proposition~\ref{rem:difference-(0,1)(1,infty)} (1) says that $s(t)$ has an even number of entries, equivalently, that one period of the sign word has an even number of runs. Since successive runs have opposite signs, the last sign of one period is opposite to the first sign of the next; hence no two runs merge at the seam. The sequence of run lengths for the full affine line is therefore exactly ${}^\infty s(t)^\infty$, up to the choice of index origin.
\end{proof}

\begin{lemma}\label{lem:finite-block-continuity-L}
For every $\varepsilon>0$, there exists $N\geq1$ such that, if two bi-infinite positive integer sequences $\mathbf a=(a_n)_{n\in\mathbb Z}$ and $\mathbf c=(c_n)_{n\in\mathbb Z}$ satisfy
\[
a_i=c_i\qquad (-N\leq i\leq N),
\]
then
\[
|\ell_0(\mathbf a)-\ell_0(\mathbf c)|<\varepsilon.
\]
\end{lemma}

\begin{proof}
The value $\ell_0(\mathbf a)$ is defined as
\[
[a_0;a_1,a_2,\ldots]+[0;a_{-1},a_{-2},\ldots].
\]
For the right infinite continued fraction $[a_0;a_1,a_2,\ldots]$, consider two values whose first $N+1$ partial quotients coincide. If $q_N$ and $q_{N-1}$ are the denominators of the convergents determined by this common initial block, then the standard formula for convergents shows that the difference between the two values is at most $1/(q_N(q_N+q_{N-1}))$. Since all partial quotients are at least $1$, the denominator $q_N$ is bounded from below by the denominator obtained when the first $N+1$ partial quotients are all $1$. Hence this difference is bounded above by a quantity independent of the sequence and tending to $0$ as $N\to\infty$. The same argument applies to the left infinite continued fraction $[0;a_{-1},a_{-2},\ldots]$. Taking $N$ sufficiently large makes the sum of the two errors smaller than $\varepsilon$.
\end{proof}

\begin{lemma}\label{lem:irrational-block-rational-approx}
Let $l$ be a regular line of positive irrational slope. Then every finite block appearing in $\mathbf b(l)$ also appears in the periodic sequence
\[
{}^\infty s(t)^\infty=(\ldots,s(t),s(t),s(t),\ldots)
\]
for some irreducible fraction $t\in(0,\infty)\cap\mathbb Q$ whose denominator may be chosen arbitrarily large. More precisely, there is a regular line $l_t$ of slope $t$ such that the block appears in $\mathbf b(l_t)$ and $\mathbf b(l_t)$ is a shift of ${}^\infty s(t)^\infty$.
\end{lemma}

\begin{proof}
Fix a finite block $W$ of $\mathbf b(l)$ and its extended sign block $\widehat W$. The block $\widehat W$ is determined by a finite segment of $l$, namely by the order in which that segment meets triangles and edges of $\widetilde{\mathbb R^2}$, the local configurations of these incidences, and the signs assigned to them. The set $\mathcal P$ is locally finite, and $l$ is regular. Hence the finite segment has positive distance from the finite set of relevant exceptional points in $\mathcal P$. Consequently, its incidence data and the extended sign block $\widehat W$ are unchanged if the slope and intercept of the line are changed sufficiently slightly.

Write the slope and intercept of $l$ as $\tau$ and $\theta$, respectively. Choose an irreducible rational number $t=p/q$ sufficiently close to $\tau$, with $q$ arbitrarily large, and choose the regular affine line $l_t^0$ supplied by Lemma~\ref{lem:rational-line-period}. Thus $\mathbf b(l_t^0)$ is a shift of ${}^\infty s(t)^\infty$.

Translating $l_t^0$ by an integer vector $(m,n)$ changes its intercept by $n-tm$. Since $t=p/q$ with $(p,q)=1$,
\[
\{n-tm\mid m,n\in\mathbb Z\}=q^{-1}\mathbb Z.
\]
Thus, by taking $q$ sufficiently large, an integer translate of $l_t^0$ can be chosen whose intercept is as close to $\theta$ as required. Denote this translate by $l_t$. Integer translations preserve $\mathcal P$ and both sign rules, so $l_t$ is regular and $\mathbf b(l_t)$ is still a shift of ${}^\infty s(t)^\infty$. Taking the slope and translated intercept sufficiently close to $(\tau,\theta)$, as in the first paragraph of the proof, makes $l_t$ produce the same extended sign block $\widehat W$. The two boundary signs in $\widehat W$ ensure that the integer block $W$ is preserved exactly. Therefore $W$ appears in ${}^\infty s(t)^\infty$.
\end{proof}

\begin{lemma}\label{lem:same-irrational-slope-same-language}
Let $l$ and $l'$ be regular lines with the same positive irrational slope. Then the set of finite blocks appearing in $\mathbf b(l)$ coincides with the set of finite blocks appearing in $\mathbf b(l')$.
\end{lemma}

\begin{proof}
Let the common slope be $\tau\notin\mathbb Q$. Translating $l$ by an integer vector $(m,n)$ changes its intercept by $n-\tau m$. By the density of irrational rotations, these changes are dense modulo $1$.

Fix a finite block $W$ appearing in $\mathbf b(l)$ and retain its extended sign block $\widehat W$. Because $l$ is regular, the finite segment producing $\widehat W$ is stable under sufficiently small parallel translations. By the density above, some integer translate of $l$ can be made sufficiently close to $l'$ in the relevant finite region. The translated segment then gives the same extended sign block $\widehat W$. Since integer translations preserve the sign rules, and since the two boundary signs prevent the first and last runs from merging with neighboring runs, the integer block $W$ appears in $\mathbf b(l')$. The reverse inclusion is proved in the same way.
\end{proof}

\begin{lemma}\label{lem:rational-values-approach-boundary}
For every positive irreducible fraction $t\in(0,\infty)\cap\mathbb Q$, one has
\[
L({}^\infty s(t)^\infty)<K.
\]
Moreover, if $(t_j)_{j\geq0}$ is a sequence of distinct positive irreducible fractions converging to an irrational number $\tau$, then
\[
\lim_{j\to\infty}L({}^\infty s(t_j)^\infty)=K.
\]
\end{lemma}

\begin{proof}
By Theorem~\ref{thm:markov-value-gen}, if $(n_t,i_t)$ is the $(k_1,k_2,k_3,\sigma)$-GM number-position pair corresponding to $t$, then
\[
L({}^\infty s(t)^\infty)
=
\frac{\sqrt{(K n_t-k_t)^2-4}}{n_t}.
\]
The right-hand side is strictly smaller than $K$.

Now suppose that $t_j\to\tau\notin\mathbb Q$ and that the $t_j$ are distinct. Write $t_j=p_j/q_j$ in lowest terms. Then $p_j+q_j\to\infty$, since otherwise a subsequence would take values in a finite set of irreducible fractions, contradicting either distinctness or convergence to an irrational number.

Write
\[
s(t_j)=(a^{(j)}_1,\ldots,a^{(j)}_{m_j}).
\]
One period of $L_{t_j}$ has at least $p_j+q_j-1$ triangle-crossing occurrences. Each such occurrence contributes one sign before adjacent equal signs are grouped, and the lengths of the resulting runs form $s(t_j)$. Consequently,
\[
\sum_{r=1}^{m_j}a^{(j)}_r\geq p_j+q_j-1\longrightarrow\infty.
\]
On the other hand, for every $r$, with indices read cyclically,
\[
a^{(j)}_r<\ell_r({}^\infty s(t_j)^\infty)
\leq L({}^\infty s(t_j)^\infty)<K.
\]
Since the $a^{(j)}_r$ are positive integers, they are all at most $K-1$. It follows from the preceding divergent sum that $m_j\to\infty$.

By Theorem~\ref{continued-fraction-theorem2}, $n_{t_j}$ is the $(2,1)$-entry of $CF_{s(t_j)}$, equivalently the denominator of the finite continued fraction $[a^{(j)}_1;\ldots,a^{(j)}_{m_j}]$. The recurrence for convergent denominators and the positivity of all partial quotients give a lower bound by the corresponding Fibonacci denominator. Hence $m_j\to\infty$ implies $n_{t_j}\to\infty$. As $k_{t_j}$ belongs to the fixed finite set $\{k_1,k_2,k_3\}$, it follows that
\[
\frac{\sqrt{(K n_{t_j}-k_{t_j})^2-4}}{n_{t_j}}
=
\sqrt{\left(K-\frac{k_{t_j}}{n_{t_j}}\right)^2-\frac{4}{n_{t_j}^2}}
\longrightarrow K.
\]
\end{proof}

\begin{proof}[Proof of Theorem~\ref{thm:irrational-slope-boundary-value}]
We first prove $L(\mathbf b(l))\leq K$. Fix $r\in\mathbb Z$. It suffices to prove $\ell_r(\mathbf b(l))\leq K$.

Let $\varepsilon>0$. By Lemma~\ref{lem:finite-block-continuity-L}, if $N$ is sufficiently large, then the value of $\ell_r$ is determined up to an error less than $\varepsilon$ by the finite block on the index interval $[r-N,r+N]\cap\mathbb Z$,
\[
(b_{r-N},\ldots,b_r,\ldots,b_{r+N}).
\]
By Lemma~\ref{lem:irrational-block-rational-approx}, this finite block also appears in ${}^\infty s(t)^\infty$ for some rational slope $t$. Hence, for some position $j$,
\[
\ell_r(\mathbf b(l))\leq \ell_j({}^\infty s(t)^\infty)+\varepsilon.
\]
By Lemma~\ref{lem:rational-values-approach-boundary}, we obtain
\[
\ell_r(\mathbf b(l))
\leq L({}^\infty s(t)^\infty)+\varepsilon
<K+\varepsilon.
\]
Since $\varepsilon>0$ was arbitrary, $\ell_r(\mathbf b(l))\leq K$. Since $r$ was arbitrary, $L(\mathbf b(l))\leq K$.

We next prove $L(\mathbf b(l))\geq K$. Let $\tau$ be the slope of $l$, and choose a sequence of distinct positive irreducible fractions $t_j\in(0,\infty)\cap\mathbb Q$ converging to $\tau$. By Lemma~\ref{lem:rational-values-approach-boundary},
\[
L({}^\infty s(t_j)^\infty)\longrightarrow K.
\]
\begingroup\emergencystretch=2em
For each $j$, choose a position in the periodic sequence ${}^\infty s(t_j)^\infty$ at which the value $L({}^\infty s(t_j)^\infty)$ is attained, and shift the sequence so that this position becomes $0$. Denote the shifted sequence by $\mathbf c^{(j)}$. Thus
\[
\ell_0(\mathbf c^{(j)})=L({}^\infty s(t_j)^\infty).
\]
\par\endgroup

We first observe that the entries of all $\mathbf c^{(j)}$ are uniformly bounded. For every $n\in\mathbb Z$,
\[
c^{(j)}_n<\ell_n(\mathbf c^{(j)})\leq L(\mathbf c^{(j)})<K.
\]
Since $K$ is an integer and $c^{(j)}_n$ is a positive integer, we have
\[
c^{(j)}_n\in\{1,2,\ldots,K-1\}
\]
for all $j$ and $n$.

Order the integers as $0,1,-1,2,-2,\ldots$. By repeatedly passing to subsequences and then taking the diagonal subsequence, we may assume that, for every fixed $n\in\mathbb Z$, the value $c^{(j)}_n$ is eventually constant. Define
\[
c_n:=\lim_{j\to\infty}c^{(j)}_n,
\qquad
\mathbf c=(c_n)_{n\in\mathbb Z}.
\]
We claim that $\ell_0(\mathbf c^{(j)})\to\ell_0(\mathbf c)$. Indeed, given $\varepsilon>0$, choose $N$ by Lemma~\ref{lem:finite-block-continuity-L}. For all sufficiently large $j$, the two sequences $\mathbf c^{(j)}$ and $\mathbf c$ agree on $[-N,N]$, and hence $|\ell_0(\mathbf c^{(j)})-\ell_0(\mathbf c)|<\varepsilon$. Therefore $\ell_0(\mathbf c^{(j)})\to\ell_0(\mathbf c)$. Since the left-hand side also tends to $K$, we get
\[
\ell_0(\mathbf c)=K.
\]

We now show that every finite block of $\mathbf c$ also appears in $\mathbf b(l)$. It is enough to consider blocks indexed by symmetric intervals about $0$, since every finite block of $\mathbf c$ is a subblock of one of them. Fix $N\geq0$ and set
\[
W=(c_{-N},c_{-N+1},\ldots,c_N).
\]
By construction of the diagonal subsequence, for all sufficiently large $j$ this block appears as
\[
(c^{(j)}_{-N},c^{(j)}_{-N+1},\ldots,c^{(j)}_N)=W.
\]
Thus $W$ appears, for infinitely many $j$, in the periodic sequence obtained from a regular rational-slope line of slope $t_j$. For each such occurrence, retain the corresponding extended sign block. Since $W$ fixes the lengths of all its consecutive sign runs, the entire extended sign block is determined by the sign of its first run: the run signs alternate, and the two added boundary signs are opposite to the adjacent run signs. Hence there are at most two possible extended sign blocks. After passing to a subsequence, we may assume that all of them are the same; denote the common block by $\widehat W$. Moreover, every triangle passage contributes one sign, so the number of triangle passages, and hence the number of intervening edge crossings, is bounded by the length of $\widehat W$. Up to integer translation, only finitely many ordered local crossing configurations can therefore occur. Passing to a further subsequence, we may assume that this complete local crossing data, including the crossed edge types and the side of each relevant edge midpoint, is the same for all occurrences.

Let $\gamma_j$ be a finite line segment of slope $t_j$ that reads $\widehat W$. Since the triangular decomposition and the sign rules are periodic under integer translations, we may translate $\gamma_j$ by an integer vector without changing the extended sign block it reads. We do this so that the midpoint of $\gamma_j$ lies in the fundamental square $[0,1]^2$.

By compactness of $[0,1]^2$, after passing to a subsequence we may assume that these midpoints converge to some point $P\in[0,1]^2$. The block $\widehat W$ has a fixed finite number of signs, and $t_j\to\tau\in(0,\infty)$. Hence the lengths of the segments $\gamma_j$ may be chosen uniformly bounded. Passing to a further subsequence, we may therefore assume that the $\gamma_j$ converge to a finite line segment $\gamma$ of slope $\tau$ passing through $P$. Let $l_{\theta_0}$ be the full line containing $\gamma$, written as
\[
l_{\theta_0}\colon y=\tau x+\theta_0.
\]

Since $\tau$ is irrational, $l_{\theta_0}$ passes through at most one point of $\mathcal P$: two distinct such points have coordinate difference in $(\frac12\mathbb Z)^2$, which would force the slope between them to be rational. We distinguish whether this possible exceptional point lies on the finite segment relevant to $\widehat W$.

If it does not, then the incidence data and the signs defining $\widehat W$ are locally constant under sufficiently small parallel translations. Thus there is an open interval $I$ containing $\theta_0$ such that every line
\[
l_\theta\colon y=\tau x+\theta,\qquad \theta\in I,
\]
reads the same extended sign block $\widehat W$ in the relevant finite region.

If the exceptional point does lie on the relevant segment, then none of the approximating segments $\gamma_j$ passes through it. After taking a further subsequence, all $\gamma_j$ pass on the same side of that point. Translating $\gamma$ slightly toward this side preserves the other incidences and reproduces the local configuration of the $\gamma_j$ at the exceptional point. Hence, in this case, there is a one-sided open interval
\[
I=(\theta_0,\theta_0+\delta)
\quad\text{or}\quad
I=(\theta_0-\delta,\theta_0)
\]
such that every $l_\theta$ with $\theta\in I$ reads $\widehat W$.

In either case, the set of intercepts for which a line of slope $\tau$ meets a point of $\mathcal P$ is
\[
E_\tau=\{b-\tau a\mid (a,b)\in\mathcal P\},
\]
which is countable. Choose $\theta\in I\setminus E_\tau$. Then $l_\theta$ is a regular line of slope $\tau$ and reads the extended sign block $\widehat W$, hence it contains the integer block $W$ exactly.

By Lemma~\ref{lem:same-irrational-slope-same-language}, the set of finite blocks produced by regular lines of slope $\tau$ is independent of the intercept. Hence every finite block of $\mathbf c$ appears in the original sequence $\mathbf b(l)$.

Finally let $\varepsilon>0$. Choose $N$ by Lemma~\ref{lem:finite-block-continuity-L}, and set
\[
W=(c_{-N},c_{-N+1},\ldots,c_N).
\]
By the preceding argument, $W$ appears in $\mathbf b(l)$, so for some position $q$,
\[
(b_{q-N},b_{q-N+1},\ldots,b_{q+N})=W.
\]
The lemma then gives
\[
|\ell_q(\mathbf b(l))-\ell_0(\mathbf c)|<\varepsilon,
\]
and hence $\ell_q(\mathbf b(l))>K-\varepsilon$.
Therefore $L(\mathbf b(l))\geq K-\varepsilon$. Since $\varepsilon>0$ was arbitrary, $L(\mathbf b(l))\geq K$. Combining the two inequalities gives $L(\mathbf b(l))=K$.
\end{proof}

We now show that the same boundary value is realized as a Lagrange constant.

\begin{corollary}\label{cor:irrational-slope-lagrange-value}
Let $(k_1,k_2,k_3)\in\mathbb Z_{\geq 0}^3$ and $\sigma\in\mathfrak S_3$ be fixed, and put $K=3+k_1+k_2+k_3$. Let $l$ be a regular line of positive irrational slope, and let
\[
\mathbf b(l)=(b_n)_{n\in\mathbb Z}
\]
be the bi-infinite sequence of positive integers obtained from $l$ by the triangle-crossing and edge-crossing rules. For any $r\in\mathbb Z$, put
\[
    \alpha_r:=[b_r;b_{r+1},b_{r+2},\ldots].
\]
Then
\[
    \mathcal L(\alpha_r)=K.
\]
In particular, $K\in\mathcal L$.
\end{corollary}

\begin{proof}
It is enough to prove the case $r=0$, since the other cases are obtained by shifting the indices of $\mathbf b(l)$. Put
\[
    \alpha:=[b_0;b_1,b_2,\ldots].
\]
By the Perron identity (Theorem~\ref{thm:lagrange-quadratic}), we have
\[
    \mathcal L(\alpha)
    =
    \limsup_{n\to\infty}
    \left(
        [b_{n+1};b_{n+2},b_{n+3},\ldots]
        +
        [0;b_n,b_{n-1},\ldots,b_1]
    \right).
\]
On the other hand, for the bi-infinite sequence $\mathbf b(l)$ we have
\[
    \ell_{n+1}(\mathbf b(l))
    =
    [b_{n+1};b_{n+2},b_{n+3},\ldots]
    +
    [0;b_n,b_{n-1},\ldots].
\]
The two continued fractions
\[
    [0;b_n,b_{n-1},\ldots,b_1]
    \quad\text{and}\quad
    [0;b_n,b_{n-1},\ldots]
\]
have the same first $n$ partial quotients. Since every partial quotient is at least $1$, the convergent-denominator estimate for regular continued fractions gives
\[
    \left|
    [0;b_n,b_{n-1},\ldots,b_1]
    -
    [0;b_n,b_{n-1},\ldots]
    \right|
    \longrightarrow 0
    \qquad (n\to\infty).
\]
Therefore
\[
    \mathcal L(\alpha)
    =
    \limsup_{n\to\infty}\ell_{n+1}(\mathbf b(l)).
\]
By Theorem~\ref{thm:irrational-slope-boundary-value}, $L(\mathbf b(l))=K$. Hence $\ell_n(\mathbf b(l))\leq K$ for all $n\in\mathbb Z$, and so
\[
    \mathcal L(\alpha)
    =
    \limsup_{n\to\infty}\ell_{n+1}(\mathbf b(l))
    \leq K.
\]

We prove the reverse inequality. Let $\varepsilon>0$. Since $L(\mathbf b(l))=K$, there exists $r_0\in\mathbb Z$ such that
\[
    \ell_{r_0}(\mathbf b(l))>K-\varepsilon.
\]
By Lemma~\ref{lem:finite-block-continuity-L}, choose $N\geq1$ sufficiently large and set
\[
    W:=(b_{r_0-N},b_{r_0-N+1},\ldots,b_{r_0+N})
\].
Then every bi-infinite positive integer sequence $\mathbf d=(d_n)_{n\in\mathbb Z}$ satisfying
\[
    (d_{-N},d_{-N+1},\ldots,d_N)=W
\]
obeys
\[
    |\ell_0(\mathbf d)-\ell_{r_0}(\mathbf b(l))|<\varepsilon.
\]

We now show that the finite block $W$ appears in the positive direction in $\mathbf b(l)$ infinitely many times. Retain the extended sign block $\widehat W$ associated with $W$, and write
\[
l\colon y=\tau x+\theta.
\]
Choose a finite segment of $l$ which reads $\widehat W$, and enlarge it slightly so that its $x$-coordinate range is contained in
\[
[m_0-M,m_0+M]
\]
for some $m_0\in\mathbb Z$ and $M\geq1$. Since $l$ is regular, all incidences relevant to this segment have positive distance from the exceptional set $\mathcal P$. It follows that there is an open neighborhood $I$ of
\[
\xi_0:=\tau m_0+\theta\pmod 1
\]
in $\mathbb R/\mathbb Z$ with the following property: whenever $m\in\mathbb Z$ satisfies
\[
\tau m+\theta\pmod 1\in I,
\]
the part of $l$ with $x\in[m-M,m+M]$ has, up to an integer translation, the same crossing order, local configurations, and signs as the reference segment. It therefore reads the same extended sign block $\widehat W$, and hence contains the integer block $W$ exactly.

Since $\tau$ is irrational, the positive rotation orbit
\[
\{\tau m+\theta\pmod 1\mid m\in\mathbb Z_{\geq0}\}
\]
is dense in $\mathbb R/\mathbb Z$. Every tail of this orbit is a translate of the same dense orbit and is therefore dense as well. We may therefore choose integers
\[
m_1<m_2<m_3<\cdots,
\qquad m_{i+1}>m_i+2M,
\]
such that $\tau m_i+\theta\pmod 1\in I$ for every $i$. The portions of $l$ with $x$-coordinate ranges $[m_i-M,m_i+M]$ are pairwise disjoint and ordered in the positive $x$-direction, so they give distinct positive-direction occurrences of $W$. Hence there are positions
\[
    q_1<q_2<q_3<\cdots
\]
such that
\[
    (b_{q_i-N},b_{q_i-N+1},\ldots,b_{q_i+N})=W
\]
for every $i$. The choice of $N$ then gives
\[
    \ell_{q_i}(\mathbf b(l))>K-2\varepsilon.
\]
Consequently,
\[
    \limsup_{n\to\infty}\ell_{n+1}(\mathbf b(l))\geq K-2\varepsilon.
\]
Since $\varepsilon>0$ was arbitrary, we obtain
\[
    \limsup_{n\to\infty}\ell_{n+1}(\mathbf b(l))\geq K.
\]
Together with the opposite inequality, this proves
\[
    \mathcal L(\alpha)
    =
    \limsup_{n\to\infty}\ell_{n+1}(\mathbf b(l))
    =K.
\]
\end{proof}

The same boundary value is also realized as a Markov constant of explicit real indefinite quadratic forms.

\begin{corollary}\label{cor:irrational-slope-markov-value}
Let $(k_1,k_2,k_3)\in\mathbb Z_{\geq 0}^3$ and $\sigma\in\mathfrak S_3$ be fixed, and put $K=3+k_1+k_2+k_3$. Let $l$ be a regular line of positive irrational slope, and let
\[
\mathbf b(l)=(b_n)_{n\in\mathbb Z}
\]
be the bi-infinite sequence of positive integers obtained from $l$ by the triangle-crossing and edge-crossing rules. For any $r\in\mathbb Z$, set
\[
\alpha_r=[b_r;b_{r+1},b_{r+2},\ldots],
\qquad
\beta_r=-[0;b_{r-1},b_{r-2},\ldots],
\]
and
\[
Q_r(x,y)=(x-\alpha_r y)(x-\beta_r y).
\]
Then
\[
\mathcal M(Q_r)=K.
\]
In particular, $K\in\mathcal M$.
\end{corollary}

\begin{proof}
By Theorem~\ref{thm:markov-quadratic}, the Markov constant of $Q_r$ is $L(\mathbf b(l))$. Theorem~\ref{thm:irrational-slope-boundary-value} gives $L(\mathbf b(l))=K$.
\end{proof}

\section{Generalized Uniqueness Questions}\label{section:Generalized Uniqueness Conjecture}
In the classical theory, Frobenius's uniqueness conjecture is equivalent both to a formulation in terms of Lagrange constants and to an injectivity statement for the Markov tree. After generalization, these become three distinct questions, none of which holds in full generality. The injectivity conjecture restricted to $k_1=k_2=k_3$, however, remains open. This section states the three questions and their counterexamples, and then formulates the surviving equal-parameter conjecture. We begin with Frobenius's uniqueness conjecture.

\begin{conjecture}[\cite{frobenius}]\label{uniqueness-conjecture}
 For any Markov number $c$, there exists a unique Markov triple $(a,b,c)$ such that $a\leq b\leq c$.    
\end{conjecture}

This conjecture is equivalent to the following (see \cite{aig}*{Corollary 9.30}):

\begin{conjecture}\label{uniqueness-conjecture2}
 For any element $L\in \mathcal M_{0,0,0}$, if $L=\mathcal L(\alpha)=\mathcal L(\beta)$, $\alpha$ and $\beta$ are $GL(2,\mathbb Z)$-equivalent. That is, there exists 
 \[
 \begin{bmatrix}
     a&b\\c&d
 \end{bmatrix}\in GL(2,\mathbb Z)
 \quad\text{such that}\quad
 \alpha=\frac{a\beta+b}{c\beta+d}.
 \]
\end{conjecture}

As a generalization of Conjecture \ref{uniqueness-conjecture}, the author and Matsushita \cite{gyomatsu} posed the following question:

\begin{question}\label{uniqueness-question}
 For any $(k_1,k_2,k_3)$-GM number $c$, is there a unique $(k_1,k_2,k_3)$-GM triple $(a,b,c)$ such that $a\leq b\leq c$?   
\end{question}

One may also consider the following question as a generalization of Conjecture \ref{uniqueness-conjecture2}.

\begin{question}\label{uniqueness-question2}
 For any element $L\in \mathcal M_{k_1,k_2,k_3}$, if $L=\mathcal L(\alpha)=\mathcal L(\beta)$, are $\alpha$ and $\beta$ $GL(2,\mathbb Z)$-equivalent?
\end{question}
After generalization, Questions \ref{uniqueness-question} and \ref{uniqueness-question2} are no longer equivalent. We give counterexamples to both.

A counterexample to Question \ref{uniqueness-question} arises when $k_1,k_2,k_3$ are distinct (\cite{gyomatsu}*{Remark 10}). In fact, the two sorted triples $(1,17,81)$ and $(2,7,81)$ are distinct $(1,2,0)$-GM triples with the same largest entry. They are represented by the ordered positive solutions $(1,81,17)$ and $(7,81,2)$, respectively, of the $(1,2,0)$-GM equation. Thus Question \ref{uniqueness-question} has a negative answer.

Question \ref{uniqueness-question2} also has a negative answer. We first give a sufficient criterion for obtaining counterexamples from reversed periods. Suppose that $k_1,k_2,k_3$ are not all zero. For a generalized strongly admissible sequence $s(t)$, set $\alpha=[s(t)^\infty]$ and 
$\beta=[s^*(\frac1t)^\infty]$. By Theorems \ref{thm:lagrange-quadratic} and \ref{thm:markov-value-gen}, we have $\mathcal{L}(\alpha)=\mathcal{L}(\beta)$. 
By Proposition \ref{rem:difference-(0,1)(1,infty)} (8), the periodic parts of $s(t)^\infty$ and $s^*(\frac1t)^\infty$ 
are reverse-ordered with respect to each other. By Serret's theorem (see, for example, \cite{aig}), two irrational numbers are $GL(2,\mathbb Z)$-equivalent if and only if their regular continued fraction expansions have a common tail. For purely periodic quadratic irrationals, this is equivalent to saying that their primitive period words---that is, their shortest period words, which are not repetitions of shorter words---agree up to a cyclic permutation. Thus, whenever the primitive period words of $s(t)^\infty$ and $s^*(\frac1t)^\infty$ are not cyclic permutations of one another, $\alpha$ and $\beta$ are not $GL(2,\mathbb Z)$-equivalent.
For an explicit counterexample, consider the two GM triples used above to answer Question \ref{uniqueness-question}. Set
\[
L=\frac{\sqrt{((3+0+1+2)81-2)^2-4}}{81}=\frac{2\sqrt{723}}{9}.
\]
Here, there exist at least two distinct quadratic irrationals $\alpha,\beta$ such that $L=\mathcal L(\alpha)=\mathcal L(\beta)$: one arising from the $(1,2,0)$-GM triple $(1,81,17)$, and another from the triple $(7,81,2)$. The triple $(1,81,17)$ corresponds to a vertex in $\mathrm{M}\mathbb T(1,2,0,\mathrm{id})$, where the fraction label of $81$ is $\tfrac{1}{3}$. The associated quadratic irrational is
\[
\alpha := [(5,1,3,3,1,4)^\infty] = \frac{\sqrt{723}+25}{9}.
\]
On the other hand, the triple $(7,81,2)$ corresponds to a vertex in $\mathrm{M}\mathbb T(1,2,0,(1\ 2\ 3))$, where the fraction label of $81$ is $\tfrac{2}{3}$. The associated quadratic irrational is
\[
\beta := [(5,1,1,5,3,2)^\infty] = \frac{\sqrt{723}+23}{9}.
\]
The two displayed period words are primitive and are not cyclic permutations of one another: even their multisets of entries, namely
\[
\{5,4,3,3,1,1\}\qquad\text{and}\qquad
\{5,5,3,2,1,1\},
\]
are different. Serret's theorem therefore gives that $\alpha$ and $\beta$ are not $GL(2,\mathbb Z)$-equivalent. The two period words are also not cyclic permutations of one another's reversals, so this example is distinct from the reversal construction based on $s(t)$ and $s^*(1/t)$. In particular, Question \ref{uniqueness-question2} has a negative answer.

Conjecture \ref{uniqueness-conjecture} also has an equivalent formulation as injectivity within a fixed Markov tree. The original injectivity conjecture is as follows.

\begin{conjecture}[\cite{aig}*{Uniqueness conjecture IV}]\label{conj:injectivity-conjecture}
Fix $\sigma\in \mathfrak S_3$ and consider $\mathrm{M}\mathbb T(0,0,0,\sigma)$. For any irreducible fractions $t,s\in[0,1]$, let $n_t$ and $n_s$ be Markov numbers with fraction labelings $t$ and $s$. If $n_t=n_s$, then $t=s$. 
\end{conjecture}

This conjecture is equivalent to Conjecture \ref{uniqueness-conjecture}. The two triples $(1,81,17)$ and $(7,81,2)$ arise from different GM trees, so they do not disprove injectivity within a fixed GM tree. We therefore consider the following generalization of that injectivity statement.

\begin{question}\label{uniqueness-question-gen}
Fix $\sigma\in \mathfrak S_3$ and consider $\mathrm{M}\mathbb T(k_1,k_2,k_3,\sigma)$. For any irreducible fractions $t,s\in[0,1]$, let $n_t$ and $n_s$ be $(k_1,k_2,k_3)$-GM numbers with fraction labelings $t$ and $s$. If $n_t=n_s$, then $t=s$. 
\end{question}
When $k_1=k_2=k_3$, the tree $\mathrm{M}\mathbb{T}(k_1,k_2,k_3,\sigma)$ is independent of $\sigma$ up to the position labels. The equivalence with Question~\ref{uniqueness-question} follows from Theorem~\ref{thm:all-solution} and Lemma~\ref{lem:GM-tree-growth}: every non-root solution is represented, up to the mirror symmetry of the tree, by a unique labelled vertex whose middle entry is the largest entry of that solution. Thus two fraction labels with the same middle entry yield the same obstruction to uniqueness, and conversely; the root cases are immediate. Hence a positive answer to Question \ref{uniqueness-question-gen} is equivalent to a positive answer to Question \ref{uniqueness-question} for equal parameters, and this case remains open. When the parameters are not required to be equal, however, the following infinite family gives counterexamples to Question \ref{uniqueness-question-gen}.

The following proposition is due to Shoma Nakabayashi.

\begin{proposition}\label{prop:injectivity-counterexamples}
Let the Fibonacci sequence be defined by $F_0=0$, $F_1=1$, and $F_{n+1}=F_n+F_{n-1}$ for $n\geq1$. For every $r\in\mathbb Z_{\geq0}$, set
\[
k_r:=\frac{F_{60r+11}-29}{10}.
\]
Then $k_r$ is a nonnegative integer, and for $(k_1,k_2,k_3,\sigma)=(0,0,k_r,(1\ 3\ 2))$ one has
\[
n_{\frac1{30r+5}}=n_{\frac23}=F_{60r+11}.
\]
In particular, Question \ref{uniqueness-question-gen} has a negative answer, and there are infinitely many counterexamples of this form.
\end{proposition}

\begin{proof}
Fix $(k_1,k_2,k_3,\sigma)=(0,0,k,(1\ 3\ 2))$. The root of the GM tree is
\[
((1,3),(2,1),(1,2)).
\]
Set $x_m:=n_{1/m}$. Along the left branch, the GM-tree recurrence gives
\[
x_1=2,\qquad x_2=5,\qquad
x_{m+1}=\frac{1+x_m^2}{x_{m-1}}\quad(m\geq2).
\]
Induction using the identity
\[
F_{2m-1}F_{2m+3}=F_{2m+1}^2+1
\]
therefore gives
\[
n_{1/m}=x_m=F_{2m+1}\qquad(m\geq1).
\]
On the other hand, the left child of the root is $((1,3),(5,2),(2,1))$, so the new middle entry in its right child is
\[
n_{2/3}=\frac{5^2+k\cdot5\cdot2+2^2}{1}=10k+29.
\]

Computing the recurrence modulo $10$ gives $(F_{60},F_{61})\equiv(0,1)\pmod {10}$. Hence $F_{n+60}\equiv F_n\pmod {10}$, and for every $r\geq0$,
\[
F_{60r+11}\equiv F_{11}=89\equiv9\pmod {10}.
\]
Thus $k_r$ is an integer, and $F_{60r+11}\geq F_{11}=89$ gives $k_r\geq6$. Taking $m=30r+5$ and $k=k_r$, we obtain
\[
n_{1/(30r+5)}=F_{60r+11}=10k_r+29=n_{2/3}.
\]
Moreover, $1/(30r+5)\ne2/3$. Since $F_{60r+11}$ is strictly increasing with $r$, the integers $k_r$ are pairwise distinct, giving infinitely many counterexamples. For $r=0$, one has $k_0=6$ and $n_{1/5}=n_{2/3}=89$.
\end{proof}

The counterexamples in Proposition \ref{prop:injectivity-counterexamples} have two parameters equal to $0$ and the third positive, and therefore do not concern the equal-parameter case. Thus the following conjecture remains open.

\begin{conjecture}\label{conj:equal-parameter-injectivity}
Fix $k\in\mathbb Z_{\geq0}$ and $\sigma\in\mathfrak S_3$, and consider $\mathrm{M}\mathbb T(k,k,k,\sigma)$. For any irreducible fractions $t,s\in[0,1]$, let $n_t$ and $n_s$ be the $(k,k,k)$-GM numbers with fraction labels $t$ and $s$, respectively. Then $n_t=n_s$ implies $t=s$.
\end{conjecture}

Conjecture \ref{conj:equal-parameter-injectivity} remains open, but the weaker monotonicity statement known as the \emph{(generalized) Aigner conjecture} holds.
\begin{theorem}[\cites{mcshane,llrs,banaian,banaian-sen}]\label{gen-aigner}
Fix $\sigma\in \mathfrak S_3$ and $k\in \mathbb Z_{\geq 0}$, and consider $\mathrm{M}\mathbb T(k,k,k,\sigma)$. For any irreducible fractions $t,s\in[0,1]$, let $n_t$ and $n_s$ be $(k,k,k)$-GM numbers with fraction labelings $t$ and $s$. Then the following statements hold.
\begin{itemize}\setlength{\leftskip}{-10pt}   
    \item[(1)] if $t$ and $s$ have the same numerator and $t<s$, then we have $n_t>n_s$, 
    \item[(2)] if $t$ and $s$ have the same denominator and $t<s$, then we have $n_t<n_s$, 
    \item [(3)] if $t$ and $s$ have the same sum of numerator and denominator and $t<s$, we have $n_t>n_s$.
    \end{itemize}
\end{theorem}
Part (3) does not extend to arbitrary unequal parameters. For example, take
\[
(k_1,k_2,k_3,\sigma)=(0,0,1,(1\ 3\ 2)).
\]
Then Table \ref{table001-3} gives $n_{\frac14}=34<39=n_{\frac23}$, although $\frac14<\frac23$ and both fractions have numerator--denominator sum $5$.
\bibliography{myrefs}
\end{document}